\documentclass{amsart}

\makeatletter

\newtheorem{theorem}{Theorem}[section]
\newtheorem{proposition}{Proposition}[section]

\newtheorem{lemma}[theorem]{Lemma}
\newtheorem{corollary}{Corollary}[theorem]
\theoremstyle{definition}
\newtheorem{definition}[theorem]{Definition}

\usepackage[margin=1.5in]{geometry} 
\usepackage{stackrel}

\usepackage{hyperref}
\usepackage{xcolor,colortbl}
\usepackage{hyperref}
\hypersetup{
	colorlinks=true,
	linkcolor=blue,
	filecolor=magenta,      
	urlcolor=cyan,
}

\definecolor{Gray}{gray}{0.5}
\definecolor{LightCyan}{rgb}{0.88,1,1}

\usepackage{amsmath}
\usepackage{accents}

\newcommand{\expect}{{\rm I\!E}}
\newcommand{\real}{{\rm I\!R}}

\newcommand{\integer}{{\rm I\!N}}
\newcommand{\prob}{{\rm I\!P}}

\newcommand{\ball}{{\rm I\!B}}

\usepackage{amsmath,amsthm,amssymb,bm}
\usepackage{algpseudocode}
\usepackage{mathtools}
\usepackage{listings}
\usepackage{color} 
\usepackage{graphicx}
\usepackage{easybmat}
\usepackage{bm}
\usepackage{subfigure}
\usepackage{hyperref}
\usepackage{algorithm}
\usepackage{dsfont}
\usepackage{enumitem}
\usepackage{amsmath}
\usepackage{algpseudocode}
\usepackage{bbm}
\usepackage{dsfont}
\usepackage{epstopdf}
\usepackage{graphicx,subfigure}
\usepackage[labelfont=bf]{caption}
\usepackage[export]{adjustbox}
\usepackage{amsbsy}
\usepackage{romannum}
\usepackage{tikz}
\usepackage{lettrine}

\let\exp\relax

\DeclareMathOperator{\exp}{exp}

\algnewcommand{\IfThenElse}[3]{
	\State \algorithmicif\ #1\ \algorithmicthen\ #2\ \algorithmicelse\ #3}

\let \expect \relax
\newcommand{\expect}{{\rm I\!E}}

\usepackage{filecontents}
\usepackage{mathrsfs}

\newcommand{\df}{\stackrel{\text{def}}{=}}

\theoremstyle{remark}

\usepackage{accents}

\DeclareMathOperator{\relint}{relint}

\newcommand{\vertiii}[1]{{\left\vert\kern-0.25ex\left\vert #1 
		\right\vert\kern-0.25ex\right\vert\kern}}

\usepackage{hyperref}
\hypersetup{
	colorlinks=true,
	linkcolor=[RGB]{110,25,25},
	citecolor=[RGB]{25,25,110},
	filecolor=magenta,      
	urlcolor=black,
}
\usepackage{lipsum}

\numberwithin{equation}{section}

\begin{document}
	
	\title[A Mean-Field Theory for Learning the Sch\"{o}nberg Measure of Radial Basis Functions] {A Mean-Field Theory for Learning the Sch\"{o}nberg Measure of Radial Basis Functions}

	

	\author[Masoud Badiei Khuzani, et al.]{Masoud Badiei Khuzani$^{\dagger,\ddagger}$, Yinyu Ye$^{\dagger}$, Sandy Napel$^{\ddagger,\ast}$,  Lei Xing$^{\ddagger,\ast}$
}

	\address{Stanford University, 450 Serra Mall, Stanford, CA 94305\\
		\\	
	}
	\curraddr{}
	\email{mbadieik,yyye,snapel,lei@stanford.edu}

	{\let\thefootnote\relax\footnotetext{
			\hspace{-5mm}$^{\dagger}$Department of Management Science, 	Stanford University, CA 94043, USA.\\
			$^{\ast}$Department of Electrical Engineering, Stanford University, CA 94043, USA.\\
			$^{\ddagger}$Department of Radiation Oncology, Stanford University, CA 94043, USA.
			\vspace{2mm}
		\vspace{-.2mm}
	}}

\begin{abstract}
We develop and analyze a projected particle Langevin optimization method to learn the distribution in the Sch\"{o}nberg integral representation of the radial basis functions from training samples.  More specifically, we characterize a distributionally robust optimization method with respect to the Wasserstein distance to optimize the distribution in the Sch\"{o}nberg integral representation. To provide theoretical performance guarantees, we analyze the scaling limits of a \textit{projected} particle online (stochastic) optimization method in the mean-field regime. In particular, we prove that in the scaling limits, the empirical measure of the Langevin particles converges to the law of a \textit{reflected} It\^{o} diffusion-drift process. The distinguishing feature of the derived process is that its drift component is also a function of the law of the underlying process. Using It\^{o} lemma for semi-martingales and Grisanov's change of measure for the Wiener processes, we then derive a Mckean-Vlasov type partial differential equation (PDE) with Robin boundary conditions that describes the evolution of the empirical measure of the projected Langevin particles in the mean-field regime. In addition, we establish the existence and uniqueness of the steady-state solutions of the derived PDE in the weak sense. We apply our learning approach to train radial kernels in the kernel locally sensitive hash (LSH) functions, where the training data-set is generated via a $k$-mean clustering method on a small subset of data-base. We subsequently apply our kernel LSH with a trained kernel for image retrieval task on MNIST data-set, and demonstrate the efficacy of our kernel learning approach. We also apply our kernel learning approach in conjunction with the kernel support vector machines (SVMs) for classification of benchmark data-sets.
\end{abstract}
	\maketitle
\pagestyle{plain}

\section{Introduction}
\lettrine{\textbf{R}}{adial} basis functions (RBFs) are key tools in statistical machine learning to approximate multivariate functions by linear combinations of terms based on a single univariate function. Such non-parametric regression schemes are typically used to approximate multivariate functions or high-dimensional data \cite{powell2001radial},\cite{cheney1966introduction} which are only known at a finite number of points or too difficult to evaluate otherwise, so that then evaluations of the approximating function can take place often and efficiently. The accuracy and performance of such techniques, however, relies to a large extent on a good choice of RBF kernel that captures the structure of data. Standard regression methods based on RBF kernels requires the input of a user-defined kernel---a drawback if a good representation of underlying data is unknown a priori. While statistical model selection techniques such as cross-validation or jackknife \cite{efron1981jackknife} can conceptually resolve those statistical model selection issues, they typically slow down the training process on complex data-sets as they repeatedly refit the model. It is thus imperative to devise efficient learning algorithms to facilitate such model selection problems.

In this paper we put forth a novel optimization framework to learn a good radial kernel from training data. Our kernel learning approach is based on a distributionally robust optimization problem  \cite{delage2010distributionally,gao2016distributionally} to learn a good distribution for the Sch\"{o}nberg integral representation of the radial basis functions from training samples  \cite[Thm. 1]{schoenberg1938metric}. Since optimization with respect to the distribution of RBF kernels is intractable, we consider a Monte Carlo sample average approximation to obtain a solvable finite dimensional optimization problem with respect to the samples of the distribution. We then use a projected particle Langevin optimization method to solve the approximated finite dimensional optimization problem. We provide a theoretical guarantee for the consistency of the finite sample-average approximations. Based on a mean-field analysis, we also show the consistency of the proposed projected particle Langevin optimization method in the sense that when the number of particles tends to infinity, the empirical distribution of the Langevin particles follows the path of the gradient descent flow of the distributionally robust optimization problem on the Wasserstein manifold.

\subsection{Related works}
The proposed kernel learning approach in this paper is closely related to the previous work of the authors in \cite{khuzani2019mean}. Therein, we proposed a particle stochastic gradient descent method in conjunction with the random feature model of Rahimi and Recht  \cite{rahimi2008random,rahimi2009weighted} to optimize the distribution of the random features in generative and discriminative machine learning models using training samples. We also showed numerically that compared to the importance sampling techniques for kernel learning (\textit{e.g.}, \cite{sinha2016learning}), the particle SGD in conjunction with the kernel SVMs yields a lower training and test errors on benchmark data-sets. Nevertheless,  we observed that the particle SGD scales rather poorly with data dimension and the number of random feature samples, rendering it inapplicable for discriminative analysis of high dimensional data-sets.  In this work, we address the scalability issue by optimizating the kernel over the sub-class of the radial kernels which includes important special cases such as the Gaussian, inverse multiquadrics, and Mat\'{e}rn kernels. While we characterize a distributionally robust optimization framework similar to \cite{khuzani2019mean}, in this work we optimize distributions defined on the  real line instead of high dimensional distributions in \cite{khuzani2019mean}.  

This work is also closely related to the copious literature on the kernel model selection problem, see, \textit{e.g.}, \cite{lanckriet2004learning,cortes2010two,bach2004multiple,gonen2011multiple}. For classification problems using kernel SVMs, Cortes, \textit{et al}. studied a kernel learning procedure from the class of mixture of base kernels.  They have also studied the generalization bounds of the proposed methods. The same authors have also studied a two-stage kernel learning in \cite{cortes2010two} based on a notion of the kernel alignment. The first stage of this technique consists of learning a kernel that is a convex combination of a set of base kernels. The second stage consists of using the learned kernel with a standard kernel-based learning algorithm such as SVMs to select a prediction hypothesis.  In \cite{lanckriet2004learning}, the authors have proposed a semi-definite programming for the kernel-target alignment problem. However, semi-definite programs based on the interior point method scale rather poorly to a large number of base kernels.  

Our proposed kernel learning framework is related to the work Sinha and Duchi \cite{sinha2016learning} for learning shift invariant kernels with random features. Therein, the authors have proposed a distributionally robust optimization for the importance sampling of random features using $f$-divergences. In \cite{khuzani2019mean}, we proposed a particle stochastic gradient descent to directly optimize the samples of the random features in the distributionally robust optimization framework, instead of optimizing the weight (importance) of the samples. However, for the high dimensional data-sets the particles are high dimensional, and the particle SGD scales poorly. In this paper, we address the scalibility issue by restricting the kernel class to the radial kernels.

The mean-field description of SGD dynamics has been studied in several prior works for different information processing tasks. Wang \textit{et al.} \cite{wang2017scaling} consider the problem of online learning for the principal component analysis (PCA), and analyze the scaling limits of different online learning algorithms based on the notion of \textit{finite exchangeability}.  In their seminal papers, Montanari and co-authors  \cite{mei2018mean,javanmard2019analysis, mei2019mean} consider the scaling limits of SGD for training a two-layer neural network, and characterize the related Mckean-Vlasov PDE for the limiting distribution of the empirical measure associated with the weights of the input layer. They also establish the uniqueness and existence of the solution for the PDE using the connection  between Mckean-Vlasov type PDEs and the gradient flows on the Wasserstein manifolds established by Otto \cite{otto2001geometry}, and Jordan, Kinderlehrer, and Otto \cite{jordan1998variational}. Similar mean-field type results for two-layer neural networks are also studied recently in \cite{rotskoff2018neural,sirignano2018mean}

\subsection{Paper Outline} The rest of this paper is organized as follows:

\begin{itemize}[leftmargin=*]
	\item \textit{Empirical Risk Minimization in Reproducing Kernel Hilbert Spaces}: In Section \ref{Section:Preliminaries of Kernel Methods in Classification and Regression}, we review some preliminaries regarding the empirical risk minimization in reproducing kernel Hilbert spaces. We also provide the notion of the kernel-target alignment for optimizing the kernel in support vector machines (SVMs). We then characterize a distributionally robust optimization problem for multiple kernel learning. 
	
	\item \textit{Theoretical Results}: In Section \ref{Section:Main_Results}, we provide the theoretical guarantees for the performance of our kernel learning algorithm. In particular, we establish the non-asymptotic consistency of the finite sample approximations. We also analyze the scaling limits of the kernel learning algorithm.
	
	\item \textit{Empirical Evaluation on Synthetic and Benchmark Data-Sets}: In Section \ref{Section:Performance Evaluation on Synthetic and Benchmark Data-Sets}, we evaluate the performance of  our proposed kernel learning model on synthetic and benchmark data-sets. In particular, we analyze the performance of our kernel learning approach for the hypothesis testing problem. We also apply our proposed kernel learning approach to develop hash codes for image query task from large data-bases. 
\end{itemize}

\section{Preliminaries and the Optimization Problem for Kernel Learning
} 
\label{Section:Preliminaries of Kernel Methods in Classification and Regression}

In this section, we review preliminaries of kernel methods in classification and regression problems. 

\subsection{Reproducing Kernel Hilbert Spaces and Kernel Alignment Optimization}

Let $\mathcal{X}$ be a metric space. A \textit{Mercer kernel} on $\mathcal{X}$ is a continuous and symmetric function $K:\mathcal{X}\times \mathcal{X}\rightarrow \real$ such that for any finite set of points $\{\bm{x}_{1},\cdots,\bm{x}_{N}\}\subset \mathcal{X}$, the kernel matrix $(K(\bm{x}_{i},\bm{x}_{j}))_{1\leq i,j\leq N}$ is positive semi-definite.  

The \textit{reproducing kernel Hilbert space} (RKHS) $\mathcal{H}_{K}$ associated with the kernel $K$ is  the completion of the linear span of the set of functions $\{K_{\bm{x}}\df K(\bm{x},\cdot), \bm{x}\in \mathcal{X}\}$ with the inner product structure $\langle\cdot,\cdot \rangle_{\mathcal{H}_{K}}$ defined by $\langle K_{\bm{x}_{0}},K_{\bm{x}_{1}} \rangle_{\mathcal{H}_{K}}=K(\bm{x}_{0},\bm{x}_{1})$. That is 
\begin{align}
\left\langle \sum_{i}\alpha_{i}K_{\bm{x}_{i}},\sum_{j}\beta_{j}K_{\bm{x}_{j}} \right\rangle_{\mathcal{H}_{K}}=\sum_{i,j}\alpha_{i}\beta_{j}K(\bm{x}_{i},\bm{x}_{j}).
\end{align}
The \textit{reproducing property} takes the following form
\begin{align}
\langle K_{\bm{x}},f\rangle_{\mathcal{H}_{K}}=f(\bm{x}), \quad \forall \bm{x}\in \mathcal{X}, f\in \mathcal{H}_{K}.
\end{align} 
In the classical supervised learning models, we are given $n$ feature vectors and their corresponding uni-variate class labels  $(\bm{x}_{1},y_{1}),\cdots,(\bm{x}_{n},y_{n})\sim_{\text{i.i.d.}} P_{\bm{X},Y}$,  $(\bm{x}_{i},y_{i})\in \mathcal{X}\times \mathcal{Y}\subset \real^{d}\times \real$. For the binary classification and regression tasks, the target spaces is given by $\mathcal{Y}=\{-1,1\}$ and $\mathcal{Y}=\real$, respectively. Given a loss function $\ell:\mathcal{Y}\times \real\rightarrow \real$, a classifier $f$ is learned from the function class $\mathcal{F}$ via the minimization of the regularized empirical risk 
\begin{align}
\label{Eq:Empirical_Loss_Minimization}
\inf_{f\in \mathcal{F}}\widehat{R}[f]\df \dfrac{1}{n}\sum_{i=1}^{n} \ell(y_{i},f(\bm{x}_{i}))+\dfrac{\lambda}{2}\|f\|_{\mathcal{F}}^{2},
\end{align}
where $\|\cdot\|_{\mathcal{F}}$ is a function norm, and $\lambda>0$ is the parameter of the regularization. Consider a Reproducing Kernel Hilbert Space (RKHS) $\mathcal{H}_{K}$ with the kernel function $K:\mathcal{X}\times \mathcal{X}\rightarrow \real$, and suppose $\mathcal{F}=\mathcal{H}_{K}\oplus 1$. Then, using the expansion $f(\bm{x})=\omega_{0}+\sum_{i=1}^{n-1}\omega_{i}K(\bm{x},\bm{x}_{i})$, and optimization over the kernel class $\mathcal{K}$ yields the following primal and dual optimization problems
\begin{subequations}
	\begin{align}
	&\text{Primal:}\min_{\bm{\omega}\in \real^{n}}\max_{K\in \mathcal{K}} \dfrac{1}{n}\sum_{i=1}^{n} \ell\left(y_{i},\omega_{0}+\sum_{i=1}^{n-1}\omega_{i}K(\bm{x},\bm{x}_{i}) \right)+\dfrac{\lambda}{2}\|\bm{\omega}\|_{2}^{2},\\
	&\text{Dual:}\max_{\bm{\alpha}\in \real^{n}}\min_{K\in \mathcal{K}}-\sum_{i=1}^{n}\ell^{\ast}(\gamma_{i},y_{i})-\dfrac{1}{2\lambda}\bm{\alpha}^{T}\bm{K}\bm{\alpha},
	\end{align}
\end{subequations}
respectively, where $\ell^{\ast}(\beta,y)=\sup_{z\in \real}\{\beta z-\ell(\beta,y)\}$ is the Fenchel's conjugate, and $\bm{K}\df (K(\bm{x}_{i},\bm{x}_{j}))_{1\leq i,j\leq n}$ is the kernel Gram matrix.

In the particular case of the soft margin SVMs classifier $\ell(y,z)=[1-yz]_{+}\df \max\{0,1-yz\}$, the primal and dual optimizations take the following forms
\begin{subequations}
	\label{Eq:Suggest}
	\begin{align}
	&\text{Primal:}\min_{\bm{\omega}\in \real^{n}}\max_{K\in \mathcal{K}}\dfrac{1}{n}\sum_{i=1}^{n}\left[1-\omega_{0}y_{i}-\sum_{i=1}^{n-1}\omega_{i}y_{i}K(\bm{x},\bm{x}_{i})\right]_{+}\hspace{-3mm}+\dfrac{\lambda}{2}\|\bm{\omega}\|_{2}^{2},  \\ 
	&\text{Dual:}\max_{\bm{\beta}\in \real^{n}:\langle \bm{\beta},\bm{y} \rangle=0,\bm{0}\preceq \bm{\beta} \preceq C\bm{1}}\min_{K\in \mathcal{K}}  \langle \bm{\beta},\bm{1}\rangle-\dfrac{1}{2}\mathrm{Tr}(\bm{K}(\bm{\beta}\odot \bm{y})(\bm{\beta}\odot \bm{y})^{T}),
	\end{align}
\end{subequations}
where $\odot$ is the Hadamard (element-wise) product of vectors. The form of the dual optimization in Eq. \eqref{Eq:Suggest} suggests that for a fixed dual vector $\bm{\beta}\in \real^{n}$, the optimal kernel can be computed by optimizing the following $U$-statistics known as the \textit{kernel-target alignment}, \textit{i.e.}
\begin{align}
\label{Eq:kernel_target}
&\max_{K\in \mathcal{K}} \dfrac{2}{n(n-1)} \sum_{1\leq i<j\leq n}y_{i}y_{j}K(\bm{x}_{i},\bm{x}_{j}).
\end{align}
In this paper, we focus on the class of radial kernels \footnote{Compare the kernel class in Eq. \eqref{Eq:Focus} with that of $\mathcal{K}\df \{K:\mathcal{X}\times \mathcal{X}\rightarrow \real:K(\bm{x}_{i},\bm{x}_{j})=\psi(\bm{x}_{i}-\bm{x}_{j}), \psi\in C^{1}(\real)\}$ we considered in \cite{khuzani2019mean}, where $\mathcal{K}$ is the set of all translation invariant kernels. } 
\begin{align}
\label{Eq:Focus}
\mathcal{K}\df \{K:\mathcal{X}\times \mathcal{X}\rightarrow \real: K(\bm{x}_{i},\bm{x}_{j})=\psi(\|\bm{x}_{i}-\bm{x}_{j}\|_{2}), \psi \in C^{1}(\real) \},
\end{align}
which includes the following important special cases:
\begin{itemize}[leftmargin=.2in]
	\item \textit{Gaussian}:$K(\bm{x},\tilde{\bm{x}})=\exp\left(-\dfrac{\|\bm{x}-\tilde{\bm{x}}\|_{2}^{2}}{2\sigma}\right)$, for $\sigma>0$,	
	\item \textit{Inverse multiquadrics}: $K(\bm{x},\tilde{\bm{x}})=\left(c^{2}+\|\bm{x}-\tilde{\bm{x}} \|_{2}^{2}\right)^{-\gamma}$ for $c,\gamma>0$,
	
	\item \textit{Mat\`{e}rn}: $K(\bm{x},\tilde{\bm{x}})= {c^{2\tau -d}\over \Gamma\big(\tau-{d\over 2}\big)2^{\tau-1-{d\over 2}}}\left({\|\bm{x}-\tilde{\bm{x}}\|_{2}^{2}\over c} \right)^{\tau-{d\over 2}}\mathcal{K}_{{d\over 2}-\tau}(c\|\bm{x}-\tilde{\bm{x}}\|_{2})$, where $\mathcal{K}_{\alpha}$ is the modified Bessel function of the third kind, and $\Gamma$ is the Euler Gamma function.
\end{itemize}

The radial kernel $K(\bm{x}_{i},\bm{x}_{j})=\psi(\|\bm{x}_{i}-\bm{x}_{j}\|_{2})$ is a positive semi-definite kernel (a.k.a. Mercer kernel) if the uni-variate function $\psi:\real_{+}\rightarrow \real$ is positive semi-definite. A uni-variate function which is positive semi-definite or positive definite on every $\real^{d}$ admits an integral representation due to Sch\"{o}enberg \cite[Thm. 1]{schoenberg1938metric}:
\begin{theorem}\textsc{(I. J. Sch\"{o}enberg \cite[Thm. 1]{schoenberg1938metric})}
 A continuous function $\psi:\real_{+}\rightarrow \real$  is positive semi-definite  if and only if it admits the following integral representation
\begin{align}
\psi(r)=\int_{0}^{\infty} e^{-tr^{2}}\mathrm{d}\mu(t),
\end{align}
for a finite positive Borel measure $\mu$ on $\real_{+}$. Moreover, if $\mathrm{supp}(\mu)\not= \{0\}$ then $\psi$ is positive definite, where $\mathrm{supp}(\mu)$ denotes the support of the measure $\mu\in \mathcal{M}_{+}(\real_{+})$.
\end{theorem}

From the Sch\"{o}nberg's representation theorem \cite{schoenberg1938metric}, the following integral representation for the radial kernels follows
\begin{align}
\label{Eq:Sch}
K(\bm{x},\widehat{\bm{x}})=\int_{0}^{\infty}e^{-\xi\|\bm{x}-\widehat{\bm{x}}\|_{2}^{2}}\mu(\mathrm{d}\xi), \quad \forall\bm{x},\widehat{\bm{x}}\in \real^{d}, \mu\in \mathcal{M}_{+}(\real_{+}),
\end{align}
where $\mathcal{M}_{+}(\real_{+})$ is the set of all finite  non-negative
Borel measures on $\real_{+}$. Due to the integral representation of Equation \eqref{Eq:Sch}, a kernel function $K$ is completely characterized in terms of the probability measure $\mu$.  Therefore, we can reformulate the kernel-target alignment as an optimization with respect to the distribution of the kernel
\begin{align}
\label{Eq:Kernel_Target_Alignment}
\sup_{\mu\in \mathcal{P}} \widehat{E}_{0}(\mu)\df \dfrac{2}{n(n-1)}\sum_{1\leq i<j\leq n}y_{i}y_{j}\int_{0}^{\infty} e^{-\xi \|\bm{x}_{i}-\bm{x}_{j}\|_{2}^{2}}\mu(\mathrm{d}\xi),
\end{align}
where $\mathcal{P}\subset \mathcal{M}_{+}(\real_{+})$ is a distribution sub-set. In the sequel, we consider a distribution ball $\mathcal{P}\df \{\mu\in \mathcal{M}(\Xi): W_{2}(\mu,\nu) \leq R\}$ with the radius $R>0$ and the (user-defined) center $\mu_{0}\in \mathcal{M}_{+}(\Xi)$, where $\Xi\subseteq \real_{+}$ is the support of the distributions. Moreover, given the metric space $(\mathcal{X},d)$, for the measures $\mu,\nu\in \mathcal{M}_{+}(\mathcal{X})$, $W_{p}(\mu,\nu)$ is the $p$-Wasserstein metric defined below
\begin{align}
\label{Eq:Wasserstein}
W_{p}(\mu,\nu)\df  \left(\inf_{\pi \in \mathcal{C}(\mu,\nu)} \int_{\mathcal{X}\times \mathcal{X}}d^{p}(\bm{\xi}_{1},\bm{\xi}_{2}) \mathrm{d}\pi(\bm{\xi}_{1},\bm{\xi}_{2})\right)^{1\over p},
\end{align}
where the  infimum is taken with respect to all couplings $\pi$ of the measures $\mu,\nu\in \mathcal{M}_{+}(\real^{d})$, and $\Pi(\mu,\nu)$ is the set of all measures for which $\mu$ and $\nu$ are marginals, \textit{i.e.},
\begin{align}
\mathcal{C}(\mu,\nu)\df \Big\{\pi\in \mathcal{M}_{+}(\real^{d}\times \real^{d}):T^{1}_{\#}\pi=\mu, T^{2}_{\#}\pi=\nu\Big\},
\end{align}
for all the maps $T^{1}(\bm{\xi}_{1},\bm{\xi}_{2})=\bm{\xi}_{2}$ and $T^{2}(\bm{\xi}_{1},\bm{\xi}_{2})=\bm{\xi}_{1}$, and $T^{1}_{\#}\pi$ and $T^{2}_{\#}\pi$ are the push-forwards of $\pi$.

Alternatively, we can recast the distributional optimization problem as a risk minimization aiming to match the output of the target kernel with the ideal kernel $\bm{K}_{\ast}=\bm{y}\bm{y}^{T}$ on the training data-set
\begin{align}
\label{Eq:Tends_to}
\inf_{\mu\in \mathcal{P}} \widehat{E}_{\gamma}(\mu)\df \dfrac{2}{n(n-1)\gamma}\sum_{1\leq i<j\leq n}\left( \gamma y_{i}y_{j}-\int_{0}^{\infty} e^{-\xi \|\bm{x}_{i}-\bm{x}_{j}\|_{2}^{2}}\mu(\mathrm{d}\xi)\right)^{2}.
\end{align}
As $\gamma\rightarrow +\infty$, the optimization problem in Eq. \eqref{Eq:Tends_to} tends to that of Eq. \eqref{Eq:Kernel_Target_Alignment}. The distributional optimization in Eq. \eqref{Eq:Tends_to} is infinite dimensional.  To characterize a finite dimensional optimization problem, we instead optimize the samples (particles) of the target distribution. In particular, we consider the independent identically distributed samples $ \xi^{1},\cdots,\xi^{N}\sim_{\text{i.i.d.}} \mu$, and let $\bm{\xi}\df (\xi^{1},\cdots,\xi^{N})$. Let $\xi_{0}^{1},\cdots,\xi_{0}^{N}\sim_{\text{i.i.d.}} \mu_{0}$. Then, we consider the following empirical risk function
\begin{align}
\label{Eq:Tends_to}
\inf_{\widehat{\mu}^{N}\in \mathcal{P}^{N}} \widehat{E}_{\gamma}\big(\widehat{\mu}^{N}\big)&\df \dfrac{2}{n(n-1)\gamma}\sum_{1\leq i<j\leq n}\left( \gamma y_{i}y_{j}-{1\over N}\sum_{k=1}^{N}e^{-\xi^{k} \|\bm{x}_{i}-\bm{x}_{j}\|_{2}^{2}}\right)^{2}.
\end{align}
where $\widehat{\mu}^{N}(\xi)\df {1\over N}\sum_{k=1}^{N}\delta_{\xi^{k}}(\xi)$ and $\widehat{\mu}_{0}^{N}(\zeta)\df {1\over N}\sum_{k=1}^{N}\delta_{\xi_{0}^{k}}(\zeta)$, and
\begin{align}
\mathcal{P}^{N}\df \Big\{\widehat{\mu}^{N}\in \mathcal{M}_{+}(\Xi):W_{2}(\widehat{\mu}^{N},\widehat{\mu}_{0}^{N})\leq R \Big\}.
\end{align}

\subsection{Surrogate Loss Function}
Using $h\geq 0$ for the constraint $W_{2}(\widehat{\mu}^{N},\widehat{\nu}^{N})\leq R $, a partial Lagrangian is
\begin{align}
\inf_{\mu\in \mathcal{M}_{+}(\real_{+})}\sup_{h\in \real_{+}} \widehat{J}_{h}(\mu)\df  \widehat{E}_{\gamma}(\widehat{\mu}^{N})+{h\over 2}\big(W^{2}_{2}(\widehat{\mu}^{N},\widehat{\mu}_{0}^{N})-R^{2}\big),
\end{align}
where $h>0$ is the Lagrange multiplier. The Wasserstein distance in Eq. \eqref{Eq:Tends_to}, the minimization in Eq. \eqref{Eq:Tends_to} is computationally prohibitive. Therefore, we instead minimize a surrogate loss function. To obtain the surrogate loss function,  we define the empirical measure of the samples (particles) $\xi_{m}^{1},\cdots,\xi_{m}^{N}$ at each iteration $m=0,1,2,\cdots,T$ as follows
\begin{align}
\widehat{\mu}^{N}_{m}(\xi)\df \dfrac{1}{N}\sum_{k=1}^{N}\delta(\xi-\xi_{m}^{k}).
\end{align}
Then, the surrogate loss for the Wasserstein distance in Eq. \eqref{Eq:Tends_to} can be obtained as follows. We denote the joint distribution of the particle-pair $(\xi_{m}^{i},\zeta^{j})$ by $\pi_{ij}$. Then, we obtain that
	\begin{align}
	 \label{Eq:Lagrangina}
	W_{2}^{2}\big(\widehat{\mu}_{m}^{N},\widehat{\mu}_{0}^{N}\big)&=\arg\min_{\bm{\pi}\in \mathcal{C}_{+}} \sum_{i,j=1}^{N}\pi_{ij}|\xi_{m}^{i}-\xi_{0}^{j}|^{2},
	\end{align}
where the set of couplings are given by 
\begin{align}
\label{Eq:Couplings}
\mathcal{C}_{+}\df \Bigg\{\bm{\pi}\in \real_{+}^{N\times N}: &\sum_{i=1}^{N}\pi_{ij}=\dfrac{1}{N},\sum_{j=1}^{N}\pi_{ij}=\dfrac{1}{N}, \quad i,j=1,2,\cdots,N\Bigg\}.
\end{align}
is a linear program that is challenging to solve in practice for a large number of particles $N$. To improve the computational efficiency, Cuturi \cite{cuturi2013sinkhorn} have proposed to add an entropic regularization to the Wasserstein distance. The resulting \textit{Sinkhorn divergence} solves a regularized version of Equation \eqref{Eq:Lagrangina}:
\begin{align}
	\bm{\pi}^{\ast}&\df \arg\min_{\bm{\pi}\in \mathcal{C}} \sum_{i,j=1}^{N}\pi_{ij}|\xi_{m}^{i}-\xi_{0}^{j}|^{2}-\varepsilon H(\bm{\pi}),
\end{align}
where $H(\bm{\pi})\df -\sum_{i,j=1}^{N}\pi_{ij}\log(\pi_{ij})$ is the entropic barrier function enforcing the non-negativity constraint on the entries $\pi_{ij}$'s, with the regularization parameter $\varepsilon>0$. Moreover, $\mathcal{C}$ has a similar definition as in Eq. \eqref{Eq:Couplings}, except that $(\pi_{ij})_{i,j}\in \real^{N\times N}$, and thus $\mathcal{C}_{+}\subseteq \mathcal{C}$.The Sikhorn divergence is now computed as follows
\begin{align}
\label{Eq:Without_Entropic_Regularization}
	W_{\varepsilon,2}^{2}\big(\widehat{\mu}_{m}^{N},\widehat{\mu}_{0}^{N}\big)&=\sum_{i,j=1}^{N}\pi_{ij}^{\ast}|\xi_{m}^{i}-\xi_{0}^{j}|^{2},
\end{align}
Notice that the entropic regularization term $H(\bm{\pi})$ is absent from Eq. \eqref{Eq:Without_Entropic_Regularization}.\footnote{The divergence in Eq. \eqref{Eq:Without_Entropic_Regularization} is sometimes referred to as the \textit{sharp Sinkhorn divergence} to differentiate it from its regularized counterpart $\widetilde{W}^{2}_{\varepsilon,2}(\widehat{\mu}_{m}^{N},\widehat{\mu}_{0}^{N})\df \sum_{i,j=1}^{N}\pi_{ij}^{\ast}|\xi_{m}^{i}-\zeta^{j}|-\varepsilon H(\bm{\pi}^{\ast})$; see, \textit{e.g.}, \cite{luise2018differential}.}

We thus consider the following surrogate loss function
\begin{align}
\label{Eq:Surrogate_Loss}
\inf_{\bm{\xi}\in \real_{+}^{N}}\widehat{J}^{N}_{h,\varepsilon}(\bm{\xi})\df \widehat{E}_{\gamma}(\bm{\xi})+{h\over 2}W_{\varepsilon,2}^{2}(\widehat{\mu}_{m}^{N},\widehat{\mu}_{0}^{N}).
\end{align}
The loss function $\widehat{J}_{h,\varepsilon}^{N}(\bm{\xi})$ acts as a proxy for the empirical loss $\widehat{J}_{h}^{N}(\bm{\xi})$ that we wish to optimize. 

After introducing the Lagrange multipliers $\bm{\lambda}\df (\lambda_{i})_{1\leq i\leq N}$ and $\tilde{\bm{\lambda}}\df (\tilde{\lambda}_{i})_{1\leq i\leq N}$ for the constraints on the marginals of $\pi_{ij}$, we obtain that
\begin{align}
\nonumber
\mathcal{L}(\bm{\pi};\bm{\lambda},\bm{\gamma})\df & \sum_{i,j=1}^{N}\varepsilon\pi_{ij}\log(\pi_{ij})+\sum_{i,j=1}^{N}\pi_{ij}|\xi_{m}^{i}-\xi_{0}^{j}|^{2}\\
&+\sum_{j=1}^{N}\lambda_{j}\left(\sum_{i=1}^{N}\pi_{ij}-{1\over N}\right)+\sum_{i=1}^{N}\tilde{\lambda}_{i}\left(\sum_{j=1}^{N}\pi_{ij}-{1\over N}\right),
\end{align}
  The optimal values of $\pi^{\ast}_{ij}$ can be derived explicitly using the KKT condition
\begin{align}
\pi_{ij}^{\ast}=v_{i}\exp\left(-\dfrac{|\xi_{m}^{i}-\xi_{0}^{j}|^{2}}{\varepsilon}\right)u_{j},
\end{align}
where $v_{i}\df \exp(-{\tilde{\lambda}_{i}\over \varepsilon}-{1\over 2})$, and $u_{j}\df  \exp(-{\lambda_{j}\over \varepsilon}-{1\over 2})$. Alternatively, the transportation matrix is given by
\begin{align}
\label{Eq:Wasserstein_distance_1}
\bm{\pi}_{m}^{\ast}\df \mathrm{diag}(\bm{v}^{\ast})e^{-{\bm{D}_{m}\over \varepsilon}}\mathrm{diag}(\bm{u}^{\ast}),
\end{align}
where $\bm{v}^{\ast}\df (v^{\ast}_{i})_{1\leq i\leq N}$, $\bm{u}^{\ast}\df (u^{\ast}_{i})_{1\leq i\leq N}$, $\bm{D}_{m}=(d^{m}_{ij})_{ij}, d^{m}_{ij}\df |\xi_{m}^{i}-\xi_{0}^{j}|^{2}$,  and $\mathrm{diag}(\cdot)$ denotes the diagonal element of the matrix. The vectors $\bm{v}^{\ast}$ and $\bm{u}^{\ast}$ can be computed  efficiently via the Sinkhorn-Knopp matrix scaling algorithm \cite{sinkhorn1967concerning}
\begin{align}
(\bm{u}_{k+1},\bm{v}_{k+1})=\mathscr{F}(\bm{u}_{k},\bm{v}_{k}), \quad  k\in \integer,
\end{align}
where $\mathscr{F}:\real^{N}\times \real^{N}\rightarrow \real^{N}\times \real^{N}$ is the following contraction mapping
\begin{align}
\mathscr{F}(\bm{u},\bm{v})={1\over N}\left(\mathrm{diag}^{-1}\left(e^{-{\bm{D}\over \varepsilon}}\bm{v}_{m'-1}\right)\bm{1}_{N},\mathrm{diag}^{-1}\left(e^{-{\bm{D}\over \varepsilon}}\bm{u}_{m'-1}\right)\bm{1}_{N} \right),
\end{align} 
where $\bm{1}_{N}\df (1,1,\cdots,1)\in \real^{N}$. Moreover,
\begin{align}
\label{Eq:Wasserstein_distance_3}
W_{\varepsilon,2}^{2}(\widehat{\mu}^{N}_{m},\widehat{\nu}^{N})=\mathrm{Tr}\left(\bm{D}\bm{\pi}_{m}^{\ast}\right).
\end{align}

To optimize Eq. \eqref{Eq:Tends_to}, we consider the  projected noisy stochastic gradient descent optimization method (a.k.a. Langevin dynamics). In particular, at each iteration $m=0,1,\cdots,T$, two samples $\bm{z}_{m}=(y_{m},\bm{x}_{m})$ and $\tilde{\bm{z}}_{m}=(\tilde{y}_{m},\tilde{\bm{x}}_{m})$ uniformly and randomly are drawn from the training data-set.  Then, we analyze the following iterative optimization method
\begin{align}
\label{Eq:following_iterations}
	\bm{\xi}_{m+1}=\mathscr{P}_{\Xi^{N}}\left(\bm{\xi}_{m}-\eta_{m} \nabla \widehat{J}^{N}_{\varepsilon,h}(\bm{\xi}_{m};\bm{z}_{m},\tilde{\bm{z}}_{m})+\sqrt{\dfrac{2\eta}{\beta}} \bm{\zeta}_{m}\right),
\end{align}
where $(\xi_{0}^{k})_{1\leq k\leq N}\sim_{\mathrm{i.i.d}} \mu_{0}$ Above, $\beta>0$ is the temperature parameter, $\eta_{m}>0$ is the step-size, and $(\bm{\zeta}_{m})_{m\in \integer}\sim_{\text{i.i.d.}} \mathsf{N}(0,\bm{I}_{N\times N})$ is the isotropic Gaussian noise. Note that when $\gamma\rightarrow \infty$, the iterations in Eq. \eqref{Eq:following_iterations} correspond to the projected stochastic gradient descent method. Moreover, $\mathscr{P}_{\Xi^{N}}(\cdot)$ is the projection onto the sub-set $\Xi^{N}\subseteq \real_{+}^{N}$. Furthermore, $\nabla \widehat{J}^{N}_{\varepsilon,h}(\bm{\xi}_{m};\bm{z}_{m},\tilde{\bm{z}}_{m})\df \left(\nabla_{k} \widehat{J}^{N}_{\varepsilon,h}(\bm{\xi}_{m};\bm{z}_{m},\tilde{\bm{z}}_{m}) \right)_{1\leq k\leq m}$ is the stochastic gradient that has the following elements
\begin{align}
\label{Eq:Nabla_RN}
\nabla_{k} \widehat{J}^{N}_{\varepsilon,h}(\bm{\xi}_{m};\bm{z}_{m},\tilde{\bm{z}}_{m})&\df {\partial E_{\gamma}(\bm{\xi}_{m};\bm{z}_{m},\tilde{\bm{z}}_{m}) \over \partial \xi_{m}^{k}}+{h\over 2}\dfrac{\partial W^{2}_{\varepsilon,2}(\widehat{\mu}^{N}_{m},\widehat{\nu}^{N})}{\partial \xi_{m}^{k}}.
\end{align}
where
\small{\begin{subequations}
\begin{align}
{\partial E_{\gamma}(\bm{\xi}_{m};\bm{z}_{m},\tilde{\bm{z}}_{m}) \over \partial \xi_{m}^{k}}&={1\over N}\left(\gamma y_{m}\tilde{y}_{m}- {1\over N}\sum_{\ell=1}^{N}e^{-\xi_{m}^{\ell}\|\bm{x}_{m}-\tilde{\bm{x}}_{m} \|_{2}^{2}}\right)e^{-\xi_{m}^{k}\|\bm{x}_{m}-\tilde{\bm{x}}_{m} \|_{2}^{2}}\|\bm{x}_{m}-\tilde{\bm{x}}_{m} \|_{2}^{2}\\
\dfrac{\partial W^{2}_{\varepsilon,2}(\widehat{\mu}^{N}_{m},\widehat{\nu}^{N})}{\partial \xi_{m}^{k}}&=\sum_{\ell=1}^{N}u_{k}^{\ast}v_{\ell}^{\ast}\left({|\xi_{m}^{k}-\zeta^{\ell}|^{2}\over \varepsilon}-1\right)e^{-{|\xi_{m}^{k}-\zeta^{\ell}|^{2}\over \varepsilon}}(\xi_{m}^{k}-\zeta^{\ell}),
\end{align}
\end{subequations}}\normalsize
for all $k=1,2,\cdots,N$. To update the Lagrange multiplier $h$, we note that Eq. \eqref{Eq:Surrogate_Loss} is concave in $h>0$, and $h$ is a scaler. Therefore, a bisection method can be applied to optimize the Lagrange multiplier $h$. In Algorithm \ref{Algoirthm:1}, we summarize the main steps of the proposed kernel learning approach. To analyze the complexity of Algorithm \ref{Algoirthm:1}, we note that the Sinkhorn's divergence can be computed in $O(N^{2})$ time \cite{cuturi2013sinkhorn}. Since the Euclidean projection in Eq. \eqref{Eq:following_iterations}  is onto the hyper-cube $\Xi^{N}$, it can be computed efficiently in $O(N)$ time by computing $\min\{\xi^{k},\xi_{u}\}$ and $\min\{\xi^{k},\xi_{l}\}$ for each particle $k=1,2,\cdots,N$. Overall, the $\epsilon$-optimal solution to problem \eqref{Eq:Surrogate_Loss} can be reached in $O(N^{2}\log(1/\epsilon))$.

\begin{algorithm}[t!]\scriptsize{
		\caption{\small{Distributionally Robust Optimization Method for Learning the Radial Kernels}} 	    \label{Algoirthm:1}}

		\begin{algorithmic}
			\State {\bfseries Inputs:} {The learning rate $\eta>0$,the radius $R>0$, samples $(\bm{x}_{i},y_{i})_{1\leq i\leq n}$, parameters $\beta>0$ and $\gamma>0$, samples $\bm{\zeta}=(\zeta^{1},\cdots,\zeta^{N})\sim \nu$, divergence parameter $\varepsilon>0$}
			\State{\bfseries Output:}{ The samples $\bm{\xi}\in \real_{+}^{N}$ that is $\epsilon$-solution to Eq. \eqref{Eq:Surrogate_Loss}}
			\State{\bfseries Initialize:}{ The particles $\bm{\xi}\leftarrow \bm{\zeta}$}
			\State{Set $h_{\mathrm{u}}\leftarrow\infty, h_{\mathrm{l}}\leftarrow 0, h_{\mathrm{s}}\leftarrow 1$}
			\While{$h_{u}=\infty$}
			\State{Draw $\bm{z},\tilde{\bm{z}}\sim_{\mathrm{i.i.d.}}\mathrm{Uniform}{(y_{i},\bm{x}_{i})}_{i=1}^{n}$ and $\bm{\zeta}\sim \mathsf{N}(0,\bm{I}_{N\times N})$}
			\State{Update the particles } 
			\begin{align}
			\bm{\xi}\leftarrow\mathscr{P}_{\Xi^{N}}\left(\bm{\xi}-\eta\nabla \widehat{J}^{N}_{\varepsilon,h_{s}}(\bm{\xi};\bm{z},\tilde{\bm{z}})+\sqrt{2\eta\over \beta}\bm{\zeta} \right).
			\end{align}
		    \IfThenElse{$W_{2,\varepsilon}\bigg({1\over N}\sum_{k=1}^{N}\delta_{\xi^{k}},{1\over N}\sum_{k=1}^{N}\delta_{\zeta^{k}}\bigg)\leq R$}  {$h_{\mathrm{u}}\leftarrow h_{s}$}{$h_{\mathrm{s}}\leftarrow 2h_{\mathrm{s}}$} \Comment{Use Eqs. \eqref{Eq:Wasserstein_distance_1}-\eqref{Eq:Wasserstein_distance_3}}
			\EndWhile
						\While{$h_{\mathrm{u}}-h_{\mathrm{l}}\geq \epsilon h_{\mathrm{s}}$}
			\State{$h_{s}\leftarrow {(h_{\mathrm{u}}+h_{\mathrm{l}})/2}$}		
			\State{Draw $\bm{z},\tilde{\bm{z}}\sim_{\mathrm{i.i.d.}}\mathrm{Uniform}{(y_{i},\bm{x}_{i})}_{i=1}^{n}$ and $\bm{\zeta}\sim \mathsf{N}(0,\bm{I}_{N\times N})$}
			\State{Update the particles } 
			\begin{align}
			\bm{\xi}\leftarrow\mathscr{P}_{\Xi^{N}}\left(\bm{\xi}-\eta\nabla \widehat{J}^{N}_{\varepsilon,h_{s}}(\bm{\xi};\bm{z},\tilde{\bm{z}})+\sqrt{2\eta\over \beta}\bm{\zeta} \right).
			\end{align}
			\IfThenElse{$W_{2,\varepsilon}\bigg({1\over N}\sum_{k=1}^{N}\delta_{\xi^{k}},{1\over N}\sum_{k=1}^{N}\delta_{\zeta^{k}}\bigg)\leq R$}  {$h_{\mathrm{u}}\leftarrow h_{s}$}{$h_{\mathrm{l}}\leftarrow h_{\mathrm{s}}$} \Comment{Use Eqs. \eqref{Eq:Wasserstein_distance_1}-\eqref{Eq:Wasserstein_distance_3}}
			\EndWhile
		\end{algorithmic}
\end{algorithm}\normalsize

\section{Main Results}
\label{Section:Main_Results}

In this section, we state our main theoretical results regarding the performance of the proposed kernel learning procedure. The proof of theoretical results is presented in Appendix.

\subsection{Assumptions}
To state our theoretical results, we first state the technical assumptions underlying our theoretical results:
\begin{itemize}
	\item[\textbf{(A.1)}] The feature space $\mathcal{X}\subset \real^{d}$ has a finite diameter, \textit{i.e.},
	\begin{align}
	\nonumber
	K\df\sup_{\bm{x},\bm{y}\in \mathcal{X}} \|\bm{x}-\bm{y}\|_{2}<\infty.
	\end{align}
	\item[\textbf{(A.2)}] The slater condition holds for the empirical loss functions. In particular, we suppose there exists $\bm{\xi}_{s}\in \real_{+}^{N}$ and $\bar{\bm{\xi}}_{s}\in \real_{+}^{N}$ such that $\bm{\xi}_{s}\in \mathrm{relint}(\mathcal{P}^{N})$, and ${\bar{\bm{\xi}}}_{s}\in \mathrm{relint}(\mathcal{P}_{\varepsilon}^{N})$.\footnote{The relative interior of a convex set $C$, abbreviated $\relint(C)$, is defined as $\relint(C)\df \{\bm{x}\in C:\exists \epsilon, \ball_{\epsilon}(\bm{x})\cap \text{aff}(C)\subseteq C \}$, where $\text{aff}(C)$ denotes the affine hull of the set $C$, and $\ball_{\epsilon}(\bm{x})$ is a ball of the radius 
		$\epsilon$ centered on $\bm{x}$.}
	
	\item[\textbf{(A.3)}] The Langevin particles are confined to a compact sub-set of $\real_{+}$, \textit{i.e.}, $\Xi=[\xi_{l},\xi_{u}]$, for some $0\leq \xi_{l}<\xi_{u}<+\infty$. 	
	
	\item[\textbf{(A.4)}] The Langevin particles are initialized by sampling from a distribution $\mu_{0}\in \mathcal{M}_{+}(\Xi)$ whose Lebesgue density exists. Let $q_{0}(\xi)\df \mathrm{d}\mu_{0}/\mathrm{d}\xi$ denotes the associated Lebesgue density.
\end{itemize}

We remark that while \textbf{(A.3)} is essential to establish theoretical performance guarantees for Algorithm \ref{Algoirthm:1} in this paper, in practice, the same algorithm can be employed to optimize distributions with unbounded support.

\subsection{Non-asymptotic consistency}

In this part, we prove that the value of the population optimum evaluated at the solution of the finite sample optimization problem in Eq. \eqref{Eq:Tends_to} approaches its population value as the number of training data $(n)$, the number of Langevin particles $(N)$, the regularization parameter $(\gamma)$, and the inverse of the parameter of the Sinkhorn transport $(\varepsilon)$ tend to infinity.

\begin{theorem}\textsc{{(Non-asymptotic Consistency of Finite-Sample Estimator)}}
	\label{Thm:Consistency of Monte Carlo Estimation}
Further, consider the distribution balls $\mathcal{P}$ and $\mathcal{P}_{\varepsilon}^{N}$ of the radius $R$ that are defined with respect to the $2$-Wasserstein and Sinkhorn's divegences, respectively. Define the optimal  values of the population optimization and its finite sample estimate 
\begin{align}
\nonumber
\mu_{\ast}(\xi)\df \inf_{\mu\in \mathcal{P}} E_{0}(\mu(\xi))\df \expect_{P^{\otimes 2}_{\bm{x},y}}[\widehat{E}_{0}(\mu(\xi))],   \quad
\widehat{\mu}^{N}_{\ast}(\xi;\gamma,\varepsilon)\df  \inf_{\widehat{\mu}^{N}\in \mathcal{P}_{\varepsilon}^{N}} \widehat{E}_{\gamma}(\widehat{\mu}^{N}).
\end{align}
respectively. Then, with the probability of (at least) $1-3\varrho$ over the training data samples $\{(\bm{x}_{i},y_{i})\}_{i=1}^{n}$ and the Langevin particles $\{\xi_{0}^{k}\}_{k=1}^{N}$, the following non-asymptotic bound holds 
\begin{align}
\nonumber
\Big|E_{0}(\mu_{\ast}(\xi))-E_{0}(\widehat{\mu}_{\ast}^{N}(\xi;\varepsilon,\gamma)) \Big|&\leq  \dfrac{2\sqrt{2}(\xi_{u}-\xi_{l})}{R^{2}\sqrt{N}}\left(1+{\log^{1\over 2}\left(\dfrac{4\sqrt{2N}(\xi_{u}-\xi_{l})}{\varrho } \right)}\right)\\ \nonumber
&+2\max\left\{\dfrac{c_{1}L^{2}}{n} \ln^{1\over 2}\left(\dfrac{4}{\varrho}\right),\dfrac{c_{2}RL^{4} }{n^{2}}\ln\left(\dfrac{4e^{L^{4}\over 9}}{\varrho}\right) \right\}+\dfrac{3}{\gamma}+c_{3}e^{-{1\over \varepsilon}},
\end{align}
where $c_{1}=3^{1\over 4}\times 2^{4}$, and $c_{2}=9\times 2^{11}$, and $c_{3}$ is a constant independent of $\varepsilon$.
\end{theorem}
The proof of Theorem \ref{Thm:Consistency of Monte Carlo Estimation} is presented in Appendix \ref{Proof:Consistency_of_Monte_Carlo}.

\subsection{Scaling limits for the unconstrained optimization}

In this part, we provide a mean-field analysis of the projected particle optimization in Eq. \eqref{Eq:following_iterations} for the special case of the unconstrained optimization. In particular, throughout this part, we assume $R=\infty$ in the distributional ball (thus $h=0$). We then analyze the following Langevin optimization method for unconstrained distributional optimization:
\begin{align}
\label{Eq:Dump}
\bm{\xi}_{m+1}=\mathscr{P}_{\Xi^{N}}\left(\bm{\xi}_{m}-\eta_{m} \nabla \widehat{J}^{N}(\bm{\xi}_{m};\bm{z}_{m},\tilde{\bm{z}}_{m})+\sqrt{\dfrac{2\eta}{\beta}} \bm{\zeta}_{m}\right),
\end{align}
where $\widehat{J}^{N}(\bm{\xi}_{m};\bm{z}_{m},\tilde{\bm{z}}_{m})\df \widehat{J}_{\varepsilon,h}^{N}(\bm{\xi}_{m};\bm{z}_{m},\tilde{\bm{z}}_{m})\Big\vert_{h=0}$. The first result of this paper is concerned with the scaling limits of the Langevin optimization method in Eq. \eqref{Eq:Dump} in the mean-field regime:

\begin{theorem}\textsc{(Mean-Field Partial Differential Equation)}
\label{Theorem:Mean-Field Partial Differential Equation}	
Suppose Assumptions $\mathbf{(A.1)}$-$\mathbf{(A.4)}$ are satisfied. Let $\mu^{N}_{t}$ denotes the continuous-time embedding of the empirical measure $(\widehat{\mu}^{N}_{k})_{k\in \integer}$ associated with the projected Langevin particles in Eq. \eqref{Eq:Dump}, \textit{i.e.},
\begin{align}
\mu^{N}_{t}(\xi)\df \widehat{\mu}^{N}_{\lfloor {t\over \eta} \rfloor}= \dfrac{1}{N}\sum_{k=1}^{N}\delta(\xi-\xi^{k}_{\lfloor {t\over \eta}\rfloor}), \quad 0\leq t\leq T,
\end{align}
associated with the projected particle Langevin optimization in Eq. \eqref{Eq:Dump}.\footnote{For a real $x\in \real$, $\lfloor x \rfloor$ stands for the largest integer not exceeding $x$.} Suppose the step size $\eta=\eta_{N}$ satisfies $\eta_{N}\rightarrow 0$, $N/\log(\eta_{N}/N)\rightarrow \infty$, and $\eta_{N}/\log(\eta_{N}/N)\rightarrow 0$ as $N\rightarrow \infty$. Furthermore, suppose the Lebesgue density $q_{0}(\xi)=\mathrm{d}\mu_{0}/\mathrm{d}\xi$ exists. Then, for any fixed $t\in [0,T]$, $\mu^{N}_{t}\stackrel{\text{weakly}}{\rightarrow}\mu_{t}$ as $N\rightarrow\infty$, where the Lebesgue density of the limiting measure $p^{\ast}_{t}(\xi)=\mathrm{d}\mu^{\ast}_{t}/\mathrm{d}\xi$ is a solution to the following distributional dynamics with Robin boundary conditions as well as an initial datum\footnote{The notion of weak convergence of probability measures is formally defined in Appendix \ref{Appendix:Formally}.}
	\begin{subequations}
	\label{Eq:distributional_dynamics}
	\begin{align}
\dfrac{\partial p_{t}(\xi)}{\partial t}&=\dfrac{\partial}{\partial \xi}\left(p_{t}(\xi)\dfrac{\partial}{\partial \xi}J(\xi,p_{t}(\xi))\right)+\dfrac{1}{\beta} \dfrac{\partial^{2}}{\partial\xi^{2}}\left(p_{t}(\xi)\right),\quad \forall(t,\xi)\in [0,T]\times (\xi_{l},\xi_{u}), \\  \label{Eq:Robin_Boundary}
{\partial p_{t}(\xi)\over \partial \xi}&+\beta p_{t}(\xi)\dfrac{\partial}{\partial \xi}J(\xi,p_{t}(\xi))\Big\vert_{\xi=\xi_{l}}=0, \quad\forall t\in [0,T],\\
{\partial p_{t}(\xi)\over \partial \xi}&+\beta p_{t}(\xi)\dfrac{\partial}{\partial \xi}J(\xi,p_{t}(\xi))\Big\vert_{\xi=\xi_{u}}=0,\quad\forall t\in [0,T], \\ 
p_{0}(\xi)&=q_{0}(\xi), \quad \forall \xi\in [\xi_{l},\xi_{u}] \\
p_{t}(\xi)&\geq 0,\quad \forall \xi\in [\xi_{l},\xi_{u}],\quad \int_{\Xi}p_{t}(\xi)\mathrm{d}\xi=1, \quad \forall t\in [0,T],
	\end{align}
\end{subequations}
where the functional $J(\xi,p_{t}(\xi))$ is defined 
\begin{align}
\nonumber
J(\xi,p_{t}(\xi))=\expect\Big[y\hat{y}\exp(-\xi \|\bm{x}-\hat{\bm{x}}\|_{2}^{2}) \Big]+\dfrac{1}{\gamma}\int_{0}^{\infty}\expect\Big[\exp(-(\xi+\xi') \|\bm{x}-\hat{\bm{x}}\|_{2}^{2})\Big] p_{t}(\xi')\mathrm{d}\xi'.
\end{align}	
Above, the expectations are taken with respect to the tensor product of the joint distribution $P^{\otimes 2}_{\bm{x},y}$ of the features and class labels, and its marginal $P^{\otimes 2}_{\bm{x}}$, respectively.
\end{theorem}
The proof of Theorem \ref{Theorem:Mean-Field Partial Differential Equation} is presented in Appendix \ref{Appendix:Proof_of_Theorem_Mean_Field}.

Let us make several remarks regarding the distributional dynamics in Theorem \ref{Theorem:Mean-Field Partial Differential Equation}. 

First, we notice that despite some resemblance between Eq. \eqref{Eq:distributional_dynamics} and the Fokker-Plank (a.k.a. forward Kolmogorov) equations for the diffusion-drift processes, they differ in that the drift $J(\xi,p_{t}(\xi))$ of Eq. \eqref{Eq:distributional_dynamics} is a functional of the Lebesgue density of the law of the underlying process. In the context of kinetic gas theory of statistical physics, such distributional dynamics are known as the Mckean-Vlasov partial differential equations; see, \textit{e.g.}, \cite{chan1994dynamics}.  

Second, the proof of Theorem \ref{Theorem:Mean-Field Partial Differential Equation} is based on the notion of \textit{propogation of chaos} or \textit{asymptotic freedom} \cite{sznitman1991topics,khuzani2019mean}, meaning that the dynamics of particles are decoupled in the asymptotic of infinitely many particles $N\rightarrow \infty$. To formalize this notion, we require the following definitions:

\begin{definition}\textsc{(Exchangablity)}
	Let $\nu$ be a probability measure on a Polish space $\mathcal{S}$ and. For $N\in \integer$, we say that $\nu^{\otimes N}$ is an exchangeable probability measure on the product space $\mathcal{S}^{n}$ if it is invariant under the permutation $\bm{\pi}\df (\pi(1),\cdots,\pi(N))$ of indices. In particular,
	\begin{align}
	\label{Eq:Exchangable}
	\nu^{\otimes N}(\bm{\pi}\cdot B)=\nu^{\otimes N}(B),
	\end{align}
	for all Borel subsets $B\in \mathcal{B}(\mathcal{S}^{n})$.	
\end{definition}

An interpretation of the exchangablity condition \eqref{Eq:Exchangable} can be provided via De Finetti's representation theorem which states that the joint distribution of an infinitely exchangeable sequence of random variables is as if a random parameter were drawn from some distribution and then the random variables in question were independent and identically distributed, conditioned on that parameter.

Next, we review the mathematical definition of chaoticity, as well as the propagation of chaos in the product measure spaces:
\begin{definition}\textsc{(Chaoticity)} 
	\label{Def:Chaoticity}
	Suppose $\nu^{\otimes N}$ is exchangeable. Then, the sequence $\{\nu^{\otimes N}\}_{N\in \integer}$ is $\nu$-chaotic if, for any natural number $\ell \in \integer$ and any test function $f_{1},f_{2},\cdots,f_{k}\in C_{b}^{2}(\mathcal{S})$, we have
	\begin{align}
	\label{Eq:Chaoticity}
	\lim_{N\rightarrow \infty}\left\langle\prod_{k=1}^{\ell} f_{k}(s^{k}),\nu^{\otimes N}(\mathrm{d}s^{1},\cdots,\mathrm{d}s^{N})\right\rangle= \prod_{k=1}^{\ell}\langle f_{k},\nu\rangle
	\end{align}
\end{definition}

According to Eq. \eqref{Eq:Chaoticity} of Definition \ref{Def:Chaoticity}, a sequence of probability measures on the product spaces $\mathcal{S}$ is $\nu$-chaotic if, for fixed $k$ the joint probability measures for the first $k$ coordinates tend to the product measure $\nu(\mathrm{d}s_{1})\nu(\mathrm{d}s_{2})\cdots \nu(\mathrm{d}s_{k})=\nu^{\otimes k}$ on $\mathcal{S}^{k}$. If the measures $\nu^{\otimes N}$ are thought of as giving the joint distribution of $N$ particles residing in the space $\mathcal{S}$, then $\{\nu^{\otimes N}\}$ is $\nu$-chaotic if $k$ particles out of $N$ become more and more independent as $N$ tends to infinity, and each particle’s distribution tends to $\nu$. A sequence of symmetric probability measures on $\mathcal{S}^{N}$ is chaotic if it is $\nu$-chaotic for some probability measure $\nu$ on $\mathcal{S}$.

If a Markov process on $\mathcal{S}^{N}$ begins in a random state with the distribution $\nu^{\otimes N}$, the distribution of the state after $t$ seconds of Markovian random motion can be expressed in terms of the transition function $\mathcal{K}^{N}$ for the Markov process. The distribution at time $t>0$ is the probability measure $U_{t}^{N}\nu^{\otimes N}$ is defined by the kernel
\begin{align}
\label{Eq:Definition_of_Ut}
U_{t}^{N}\nu^{\otimes N}(B)\df \int_{\mathcal{S}^{N}}\mathcal{K}^{N}(s,B,t)\nu^{\otimes N}(\mathrm{d}s).
\end{align}

\begin{definition}\textsc{(Propogation of Chaos)} A sequence functions 
	\begin{align}
	\Big\{\mathcal{K}^{N}(s,B,t)\Big\}_{N\in \integer}
	\end{align}	
	whose $N$-th term is a Markov transition function on $\mathcal{S}^{N}$ that satisfies the  permutation condition
	\begin{align}
	\mathcal{K}^{N}(s,B,t)=\mathcal{K}^{N}(\bm{\pi}\cdot s,\bm{\pi}\cdot B,t),
	\end{align}
	propagates chaos if whenever $\{\nu^{\otimes N} \}_{N\in \integer}$ is chaotic, so is $\{U_{t}^{N} \}$ for any $t\geq 0$, where $U_{t}^{N}$ is defined in Eq. \eqref{Eq:Definition_of_Ut}.
\end{definition}

In the context of this paper, the propagation of chaos 
simply means that the time-scaled empirical measures $\widehat{\mu}^{N}_{t}$
(say, on the path space) converge weakly in probability to
the deterministic measure $\mu_{t}$, or equivalently that the law of $(\xi_{\lfloor {t\over \eta}\rfloor}^{1},\cdots, \xi_{\lfloor{t\over \eta}\rfloor}^{N})$ converges weakly to the product measure $\mu_{t}^{\otimes N}$ for any fixed $N$, where $\mu_{t}$ is the law of the following \textit{reflected} diffusion-drift It\^{o} process 
\begin{align}
\label{Eq:Lulu}
\theta_{t}=\theta_{0}-\int_{0}^{t}{\partial\over \partial \xi}J(\xi_{s},\mu_{s}(\xi))\mathrm{d}s +\dfrac{1}{\beta}W_{t}+Z_{t}^{-}-Z_{t}^{+}, \quad \mu_{t}\stackrel{\mathrm{law}}{=} \theta_{t}.
\end{align}
In Eqs. \eqref{Eq:Lulu}, $Z_{t}^{-}$ and $Z_{t}^{+}$ are the non-negative reflection processes from the boundaries $\xi_{l}$ and $\xi_{u}$, respectively. In particular, $Z_{t}^{-}$ and $Z_{t}^{+}$ are non-decreasing, c\'{a}dl\'{a}g, with the initial values $Z_{0}^{+}=Z_{0}^{-}=0$, and 
\begin{align}
\int_{0}^{\infty}(X_{t}-\xi_{\mathrm{l}})\mathrm{d}Z^{-}_{t}=0, \quad \int_{0}^{\infty}(\xi_{\mathrm{u}}-X_{t})\mathrm{d}Z^{+}_{t}=0,
\end{align}
where the integrals are in the Stieltjes sense. In fact, we show that the law of the Langevin particles $(\xi_{\lfloor {t\over \eta}\rfloor}^{1},\cdots, \xi_{\lfloor {t\over \eta}\rfloor}^{N})$ converges in the $2$-Wasserstein distance to $\mu_{t}^{\otimes N}$. Alternatively, due to the bounded equivalence of the Wasserstein distance and the total variation distance on compact metric spaces (cf. \textbf{(A.3)}), the sense in which $\widehat{\mu}^{N}_{t}$ converges in probability to $\mu_{t}$ can be strengthened; rather than working with the usual weak-$^\ast$ topology induced by duality with bounded continuous test functions (see, \textit{e.g.}, \cite{khuzani2019mean}), we work with the stronger topology induced by duality with bounded measurable test functions.

Third, the Robin boundary conditions in Eq. \eqref{Eq:Robin_Boundary} captures the effect of an \textit{elastic reflection} of the particles in the projected particle Langevin optimization method. Alternatively, the Robin boundary conditions asserts that the flux of particles in and out of the boundaries is zero. In particular, the particles are constrained to stay in $(\xi_{l},\xi_{u})$ by barriers at $\xi=\xi_{u}$ and $\xi=\xi_{l}$, in such a way that when the particle hits the barrier with incoming velocity $v<0$, it will instantly bounce back with velocity $−\rho v\geq 0$, where $\rho \geq 0$ is a parameter called the \textit{velocity restitution coefficient}. When $\rho=1$ we say that the reflection is perfectly elastic; when $\rho=0$ it is said to be a \textit{rigid reflection} resulting in the Neumann boundary conditions instead of the Robin boundary conditions. By definition of the Euclidean projection $\mathscr{P}_{\Xi^{N}}(\cdot)$ in the projected Langevin particle of Eq. \eqref{Eq:Dump}, the particles are bounced off of the boundary of the projection space with the elastic parameter $\rho=\beta$.  

Fourth, we note that the solution of the distributional dynamics in Eq. \eqref{Eq:distributional_dynamics} is absolutely continuous with respect to $\xi$, and lacks atoms at the boundaries $\xi=\xi_{l}$ and $\xi=\xi_{u}$. This is due to the fact that the boundaries are non-absorbing.

In the next proposition, we establish the existence and uniqueness of the steady-state solutions of the mean-field PDE \eqref{Eq:distributional_dynamics}:

\begin{proposition}\textsc{(Existence and Uniqueness of the Steady-State Solution)}
	\label{Proposition:Uniqueness and Fixed Points}
Suppose Assumptions $\mathbf{(A.1)}$-$\mathbf{(A.4)}$ are satisfied. Then, the following assertions hold:
\begin{itemize}
\item There exists a unique non-negative steady-state solution for the mean-field PDE in Eq. \eqref{Eq:distributional_dynamics} in the weak sense. 

\item  The fixed point of the restricted Gibbs-Boltzmann measure
\begin{align}
p_{\ast}(\xi)= \dfrac{\exp(-\beta J(\xi,p_{\ast}(\xi)) )}{Z_{\beta}}\bm{1}_{\xi\in \Xi}, \quad Z_{\beta}\df \int_{\Xi} \exp(-\beta J(\xi,p_{\ast}(\xi)))\mathrm{d}\xi.
\end{align}
is a steady-state solution to the mean-field PDE in Eq. \eqref{Eq:distributional_dynamics}. 
\end{itemize}
\end{proposition}

The proof of Proposition \ref{Proposition:Uniqueness and Fixed Points} is deferred to Appendix \ref{Appendix:Proof_of_Uniqueness}.

\section{Performance Evaluation on Synthetic and Benchmark Data-Sets}
\label{Section:Performance Evaluation on Synthetic and Benchmark Data-Sets}

For simulations of this section, we consider kernel SVMs in conjunction with the random feature model of Rahimi and Recht \cite{rahimi2008random,rahimi2009weighted}. Specifically, we consider the following model random feature model:
\begin{align}
\label{Eq:Random_Feature_Model}
K(\bm{x},\tilde{\bm{x}})=\int_{\Omega}\varphi(\bm{x};\bm{\omega})\varphi(\tilde{\bm{x}};\bm{\omega})\nu(\mathrm{d}\bm{\omega}),
\end{align}
where $\varphi: \mathcal{X}\times \Omega\rightarrow \real$ is the random feature. For translation invariant kernels, we have $\varphi(\bm{x};\bm{\omega})=\sqrt{2}\cos(\langle \bm{x},\bm{\omega} \rangle+b)$, where $b\sim \mathrm{Uniform}[-\pi,\pi]$ is a random bias term.  Given the i.i.d. samples $\bm{\zeta}^{1},\cdots,\bm{\zeta}^{D}\sim_{\text{i.i.d.}} \nu$, then we minimize the following risk function
 \begin{align}
\min_{\bm{\alpha}\in \real^{D}} \dfrac{1}{n}\sum_{i=1}^{n}\max\left\{0,1-\sum_{k=1}^{D}y_{i}\alpha_{i}\varphi(\bm{x}_{i};\bm{\omega}_{k}) \right\}+\dfrac{\lambda}{2}\|\bm{\alpha}\|_{2}^{2}.
\end{align}
To compute the random features associated with the radial kernels,  we note that the probability measure of the random feature has the following Lebesgue's density
\begin{align}
\nonumber
\dfrac{\mathrm{d}\nu}{\mathrm{d}\bm{\omega}}&=\int_{\Xi} \left(\dfrac{1 }{2\xi}\right)^{{d\over 2}} \exp\left(-{ \|\bm{\omega}\|_{2}^{2}\over 4\xi}\right)\mu(\mathrm{d}\xi).
\end{align}
More specifically, due to the fact that
\begin{align}
\int_{0}^{\infty}\left(\int_{\real^{d}}\left|e^{-i\langle \bm{\omega},\bm{x}-\tilde{\bm{x}} \rangle}\left(\dfrac{1 }{2\xi}\right)^{{d\over 2}} \exp\left(-{ \|\bm{\omega}\|_{2}^{2}\over 4\xi}\right) \right|\mathrm{d}\bm{\omega}\right)\mu(\mathrm{d}\xi) < \infty,
\end{align}
the Tonelli-Fubini theorem is admissible. Therefore, we may change the order of integral to derive
\small{\begin{align}
\nonumber
\int_{\real^{d}}\dfrac{e^{-i\langle \bm{\omega},\bm{x} \rangle}}{(2\pi)^{d\over 2}}\left(\int_{0}^{\infty}\left(\dfrac{1 }{2\xi}\right)^{{d\over 2}} \exp\left(-{ \|\bm{\omega}\|_{2}^{2}\over 4\xi}\right)\mu(\mathrm{d}\xi) \right)\mathrm{d}\bm{\omega}&=\int_{0}^{\infty}\left(\int_{\real^{d}}\dfrac{e^{-i\langle \bm{\omega},\bm{x} \rangle}}{(4\xi\pi)^{d\over 2}} \exp\left(-{ \|\bm{\omega}\|_{2}^{2}\over 4\xi} \right)\mathrm{d}\bm{\omega}\right)\mu(\mathrm{d}\xi)\\ \nonumber
&=\int_{0}^{\infty}e^{-\xi \|\bm{x}-\tilde{\bm{x}}\|_{2}^{2}},
\end{align}}\normalsize
where we recognize the last term as the integral representation of the radial kernel in Eq. \eqref{Eq:Sch}. Let $\bm{\xi}=(\xi^{k})_{1\leq k\leq N}\sim_{\text{i.i.d.}} \mu$ denotes the vector outputs of Algorithm \ref{Algoirthm:1}. We draw the random features,
\begin{align}
\bm{\varphi}(\bm{x})\df (\varphi(\bm{x};\bm{\omega}_{1}),\cdots,\varphi(\bm{x};\bm{\omega}_{D})),
\end{align}
where $\bm{\omega}_{1},\cdots,\bm{\omega}_{D}\sim_{\mathrm{i.i.d.}} \widehat{\nu}^{N}$.

We compare our method with three alternative kernel learning techniques, namely, the importance sampling of Sinha and Duchi \cite{sinha2016learning}, the Gaussian bandwidth optimization via $k$ nearest neighbor ($k$-NN) \cite{chen2017simple}, and our proposed particle SGD in \cite{khuzani2019mean}.  In the sequel, we provide a brief description of each method:
\begin{itemize}[leftmargin=*]
	\item \textit{Importance Sampling}: the importance sampling of Sinha and Duchi \cite{sinha2016learning} which proposes to assign a weight $w^{m}$ to each sample $\bm{\omega}^{m}\in \real^{d_{0}}$ for $m=1,2,\cdots,N$. The weights are then optimized via the following standard optimization procedure
	\begin{align}
	\label{Eq:Importance}
	\max_{w^{1},\cdots,w^{N_{\bm{\omega}}}\in \mathcal{Q}_{N}}\dfrac{1}{2n(n-1)} \sum_{1\leq i<j\leq n}y_{i}y_{j}\sum_{k=1}^{N_{\bm{\omega}}}w^{m}\varphi(\bm{x}_{i};\bm{\omega}^{k})\varphi(\bm{x}_{j};\bm{\omega}^{k}),
	\end{align}    
	where $\mathcal{Q}_{N_{\bm{\omega}}}\df \{\bm{w}\in \real_{+}^{N_{\bm{\omega}}}: \langle \bm{w},\bm{1}\rangle=1, \chi(\bm{w}|| \bm{1}/N_{\bm{\omega}})\leq R\}$ is the distributional ball, and $\chi(\cdot||\cdot)$ is the $\chi^{2}-$divergence. Notice that the complexity of solving the optimization problem in Eq. \eqref{Eq:Importance} is insensitive to the dimension of random features $\bm{x}_{i}\in \real^{d}$. 
	
	\item \textit{Gaussian Kernel with $k$-NN Bandwidth Selection Rule}: In this method, we fix the Gaussian kernel $K(\bm{x},\bm{y})=\exp\left(-{\|\bm{x}-\bm{y}\|_{2}^{2}\over \sigma^{2}}\right)$. To select a good bandwidth for the kernel, we use the deterministic $k$ nearest neighbor ($k$NN) approach of \cite{chen2017simple}. In particular, we choose the bandwidth according to the following rule
	\begin{align}
	\sigma^{2}= \dfrac{1}{n}\sum_{i=1}^{n}\|\bm{x}_{i}-\bm{x}_{i}^{k\mathrm{NN}}\|_{2}^{2},
	\end{align}
	where $\bm{x}_{i}^{k\mathrm{NN}}$ is defined as $k$ nearest neighbor of $\bm{x}_{i}$. In our experiments, we let $k=3$. We then generate random Fourier features $\varphi(\bm{x};\bm{\omega}^{k})=\sqrt{{2\over N_{\bm{\omega}}}}\cos(\langle \bm{w},\bm{\omega}^{k} \rangle+b^{k})$, with $(\bm{\omega}^{k})_{1\leq k\leq N}\sim_{\mathrm{i.i.d.}} {1\over \sqrt{2\pi\sigma^{2}}}\exp(-\sigma^{2}\|\bm{\omega}\|_{2}^{2})$ and $(b^k)_{1\leq k\leq N}\sim_{\mathrm{i.i.d.}} \mathrm{Uniform}[-\pi,\pi]$.
	
	\item \textit{Particle Stochastic Gradient Descend Method}: This is the method that we proposed in the earlier work \cite{khuzani2019mean}, where we optimized the samples in the random feature model of Eq. \eqref{Eq:Random_Feature_Model}. In particular,
		\small{\begin{align}
					\label{Eq:SGD}
			\bm{\omega}^{k}_{m+1}&=\bm{\omega}^{k}_{m}-{\eta\over N_{\bm{\omega}}}\left(y_{m}\widetilde{y}_{m}- {1\over \alpha N}\sum_{k=1}^{N_{\bm{\omega}}}\varphi(\bm{x}_{m};\bm{\omega}_{m}^{k})\varphi(\widetilde{\bm{x}}_{m};\bm{\omega}_{m}^{k})\right)\nabla_{\bm{\omega}}\Big(\varphi(\bm{x}_{m};\bm{\omega}_{m}^{k})\varphi(\widetilde{\bm{x}}_{m};\bm{\omega}_{m}^{k})\Big),
			\end{align}}\normalsize
	for $k=1,2,\cdots,N$. Notice that in the random feature model of \cite{rahimi2008random,rahimi2009weighted}, the explicit feature map is given by $\varphi(\bm{x};\bm{\omega})=\sqrt{2}\cos(\langle \bm{x},\bm{\omega}\rangle+b)$, where the dimension of the particles in Eq. \eqref{Eq:SGD} is the same as the dimension of the feature vectors $\bm{x}\in \real^{d}$.  As a result, for high dimensional data-sets $(d\gg 1)$, the computational complexity per iterations of Eq. \eqref{Eq:SGD} is potentially prohibitive.
\end{itemize}

\subsection{Empirical Results on the Synthetic Data-Set}

For experiments with the synthetic data, we use the setup of \cite{khuzani2019mean}. The synthetic data-set we consider is as follows:

\begin{itemize}
	\item The distribution of training data is $P_{\bm{V}}=\mathsf{N}(\bm{0},(1+\lambda)\bm{I}_{d\times d})$,
	
	\item The  distribution of generated data is $P_{\bm{W}}=\mathsf{N}(\bm{0},(1-\lambda)\bm{I}_{d\times d})$.
\end{itemize}
To reduce the dimensionality of data, we consider the embedding $\iota:\real^{d}\mapsto \real^{d_{0}}, \bm{x}\mapsto \iota(\bm{x})=\bm{\Sigma}\bm{x}$, where $\bm{\Sigma}\in \real^{d_{0}\times d}$ and $d_{0}<d$. In this case, the distribution of the embedded features are $P_{\bm{X}|Y=+1}=\mathsf{N}(\bm{0},(1+\lambda)\bm{\Sigma}\bm{\Sigma}^{T})$, and $P_{\bm{X}|Y=-1}=\mathsf{N}(\bm{0},(1-\lambda)\bm{\Sigma}\bm{\Sigma}^{T})$. 

Note that $\lambda\in [0,1]$ is a parameter that determines the separation of distributions. In particular, the Kullback-Leibler divergece of the two multi-variate Gaussian distributions is controlled by $\lambda\in [0,1]$,
\begin{align}
D_{\mathrm{KL}}(P_{\bm{X}|Y=-1},P_{\bm{X}|Y=+1})=\dfrac{1}{2}\left[ \log \left(\dfrac{1-\lambda}{1+\lambda}\right)-d_{0}+d_{0}(1-\lambda^{2})\right].
\end{align}

\begin{figure}[!t]
	\begin{center}
		\hspace*{-5mm}		\subfigure{
			\includegraphics[trim={.2cm .2cm .2cm  .6cm},width=.35\linewidth]{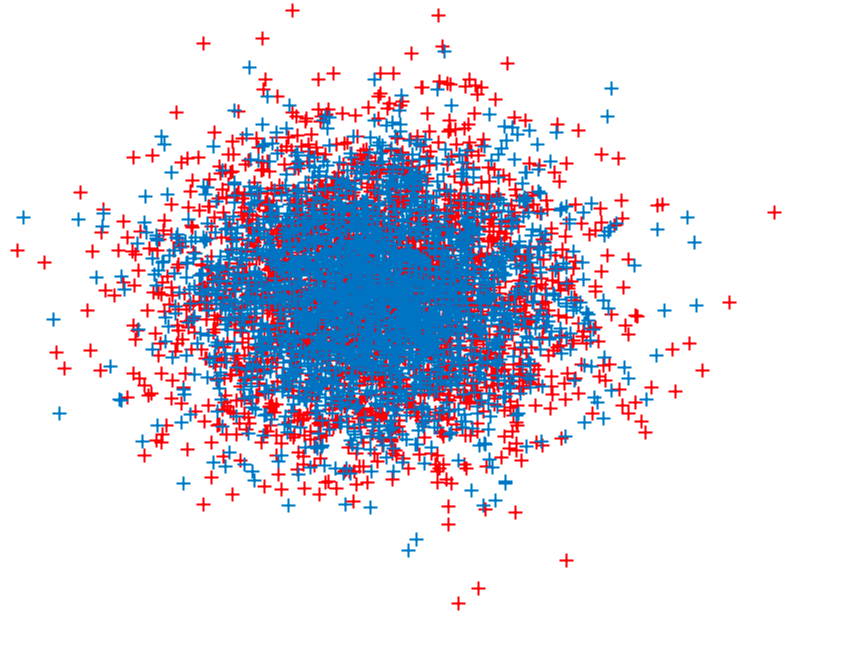} 
			\includegraphics[trim={.2cm .2cm .2cm  .6cm},width=.35\linewidth]{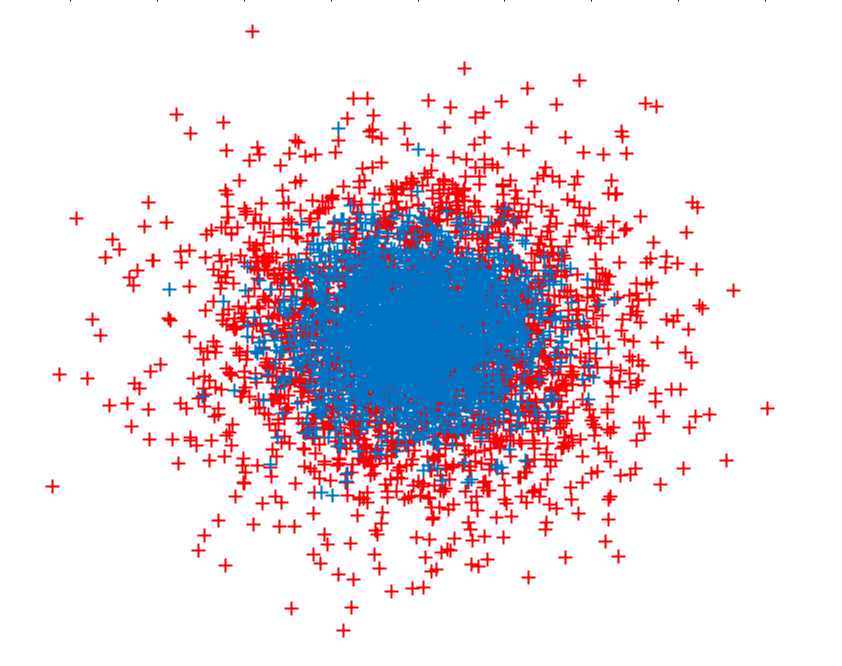} 
			\includegraphics[trim={.2cm .2cm .2cm  .2cm},width=.35\linewidth]{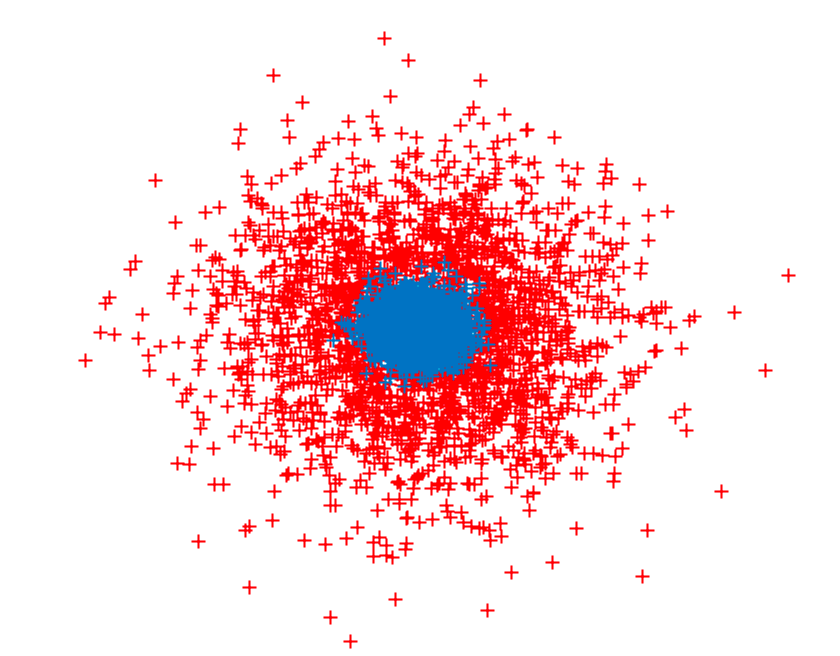}}
		\\ \vspace{-2mm}
		\subfigure{\footnotesize{\hspace{5mm} (a) \hspace{45mm} (b) \hspace{45mm} (c)} }
		\vspace{-2mm}
		\caption{\footnotesize{Visualization of data-points using the synthetic data generation models with multivariate Gaussian distributions $P_{\bm{V}}=\mathsf{N}(\bm{0},(1+\lambda)\bm{I}_{d\times d})$ and $P_{\bm{W}}=\mathsf{N}(\bm{0},(1-\lambda)\bm{I}_{d\times d})$ for $d=2$. Panel (a): $\lambda=0.1$, Panel (b): $\lambda=0.5$, and Panel (c): $\lambda=0.9$.}}
		\label{Fig:4} 
	\end{center}
\end{figure}

In Figure \ref{Fig:4}, we show the \textit{i.i.d.} samples from the distributions $P_{\bm{V}}$ and $P_{\bm{W}}$ for different choices of variance parameter of $\lambda=0.1$, $\lambda=0.5$, and $\lambda=0.9$.  Notice that for larger $\lambda$ the divergence is reduced and thus performing the two-sample test is more difficult. From Figure \ref{Fig:4}, we clearly observe that for large values of $\lambda$, the data-points from the two distributions $P_{\bm{V}}$ and $P_{\bm{W}}$ have a large overlap and conducting a statistical test to distinguish between these two distributions is more challenging.

\subsubsection{Kernel Learning Approach} 
\label{subsubsection:Kernel Learning Approach}
Figure \ref{Fig:5} depicts our two-phase kernel learning procedure. The kernel learning approach consists of training the auto-encoder for the dimensionality reduction and the kernel optimization sequentially, \textit{i.e.}, 
\begin{align}
\label{Eq:solution_of_the_optimization}
\sup_{\widehat{\mu}^{N}\in\mathcal{P}^{N}}\sup_{\iota \in \mathcal{Q}}\widehat{\mathrm{MMD}}_{K_{\widehat{\mu}^{N}}\circ \iota}^{\alpha}[P_{\bm{V}},P_{\bm{W}}].
\end{align}
where the function class is defined $\mathcal{Q}\df \{\iota(\bm{z})=\sigma(\bm{\Sigma}\bm{z}+\bm{b}),\bm{\Sigma}\in \real^{d_{0}\times d},\bm{b}\in \real^{d_{0}}\}$, and $(K_{\widehat{\mu}^{N}}\circ \iota)(\bm{x}_{1},\bm{x}_{2})=K_{\widehat{\mu}^{N}}(\iota(\bm{x}_{1}),\iota(\bm{x}_{2}))$. Here, $\bm{\sigma}(\cdot)$ is the sigmoid non-linearity. Now, we consider a two-phase optimization procedure: 
\begin{itemize}
	\item \textsc{Phase (I)}: we fix the kernel function, and optimize the auto-encoder to compute a co-variance matrix $\bm{\Sigma}$ and the bias term $\bm{b}$ for the dimensionality reduction.
	
	\item \textsc{Phase (II)}: we optimize the kernel based from the learned embedded features $\iota(\bm{x})$. 
\end{itemize}
This two-phase procedure significantly improves the computational complexity of SGD as it reduces the dimensionality of random feature samples $\bm{\xi}\in \real^{D}$, $D=d_{0}\ll d$.

\subsubsection{Statistical Hypothesis Testing with the Kernel MMD}

Let $\bm{V}_{1},\cdots,\bm{V}_{m}\sim_{\text{i.i.d.}} P_{\bm{V}}=\mathsf{N}(\bm{0},(1+\lambda)\bm{I}_{d\times d})$, and $\bm{W}_{1},\cdots,\bm{W}_{n}\sim_{\text{i.i.d.}} P_{\bm{W}}=\mathsf{N}(\bm{0},(1-\lambda)\bm{I}_{d\times d})$. Given these i.i.d. samples, the statistical test $\mathcal{T}(\{\bm{V}_{i}\}_{i=1}^{m},\{\bm{W}_{i}\}_{j=1}^{n}):\mathcal{V}^{m}\times \mathcal{W}^{n}\rightarrow \{0,1\}$ is used to distinguish between these hypotheses:
\begin{itemize}
	\item \textsc{Null hypothesis} $\mathsf{H}_{0}:P_{\bm{V}}=P_{\bm{W}}$ (thus $\lambda=0$),
	\item \textsc{Alternative hypothesis} $\mathsf{H}_{1}:P_{\bm{V}}\not= P_{\bm{W}}$ (thus $\lambda>0$).
\end{itemize}
To perform hypothesis testing via the kernel MMD, we require that $\mathcal{H}_{\mathcal{X}}$ is a universal RKHS, defined on a compact metric space $\mathcal{X}$. Universality requires that the kernel $K(\cdot,\cdot)$ be continuous and, $\mathcal{H}_{\mathcal{X}}$ be dense in $C(\mathcal{X})$.
Under these conditions, the following theorem establishes that the kernel MMD is indeed a metric:

\begin{theorem}\textsc{(Metrizablity of the RKHS)}
	Let $\mathcal{F}$ denotes a unit ball in a universal RKHS $\mathcal{H}_{\mathcal{X}}$ defined on a compact metric space $\mathcal{X}$ with the associated continuous kernel $K(\cdot,\cdot)$. Then,  the kernel MMD is a metric in the sense that $\mathrm{MMD}_{K}[P_{\bm{V}},P_{\bm{W}}]=0$ if and only if $P_{\bm{V}}=P_{\bm{W}}$.
\end{theorem}

A radial kernel is universal if the support of the measure in the integral representation in Eq. \eqref{Eq:Sch} excludes the origin, i.e., $\mathrm{supp}(\mu)\not=\{0\}$, see \cite{sriperumbudur2011universality}.

To design a test, let $\widehat{\mu}^{N}_{m}(\bm{\xi})={1\over N}\sum_{k=1}^{N}\delta(\bm{\xi}-\bm{\xi}_{m}^{k})$ denotes the solution of SGD in \eqref{Eq:SGD} for solving the optimization problem. Consider the following MMD estimator consisting of two $U$-statistics and an empirical function
\begin{align}
\nonumber
\widehat{\mathrm{MMD}}_{K_{\widehat{\mu}_{m}^{N}}\circ \iota}\big[\{\bm{V}_{i}\}_{i=1}^{m},\{\bm{W}_{i}\}_{i=1}^{n}\big]&=\dfrac{1}{m(m-1)}\sum_{k=1}^{N}\sum_{i\not =j}\varphi(\iota(\bm{V}_{i}),\bm{\omega}_{m}^{k})\varphi(\iota(\bm{V}_{j}),\bm{\omega}_{m}^{k})\\ \nonumber
&\hspace{4mm}+\dfrac{1}{n(n-1)}\sum_{k=1}^{N}\sum_{i\not =j}\varphi(\iota(\bm{W}_{i}),\bm{\omega}_{m}^{k})\varphi(\iota(\bm{W}_{j}),\bm{\omega}_{m}^{k})\\ \label{Eq:Negative_Estimator}
&\hspace{4mm}-\dfrac{2}{nm}\sum_{k=1}^{N}\sum_{i=1}^{m}\sum_{j=1}^{n}\varphi(\iota(\bm{W}_{i}),\bm{\omega}_{m}^{k})\varphi(\iota(\bm{V}_{j}),\bm{\omega}^{k}_{m}).
\end{align}
Given the samples $\{\bm{V}_{i}\}_{i=1}^{m}$ and $\{\bm{W}_{i}\}_{i=1}^{n}$, we design a test statistic as below
\begin{align}
\label{Eq:Test_statistics}
\mathcal{T}(\{\bm{V}_{i}\}_{i=1}^{m},\{\bm{W}_{i}\}_{i=1}^{n})\df \begin{cases}
\mathsf{H}_{0}  & \text{if } \widehat{\mathrm{MMD}}_{K_{\widehat{\mu}_{m}^{N}}\circ \iota}\big[\{\bm{V}_{i}\}_{i=1}^{m},\{\bm{W}_{i}\}_{i=1}^{n}\big]\leq \tau \\
\mathsf{H}_{1}  & \text{if } \widehat{\mathrm{MMD}}_{K_{\widehat{\mu}_{m}^{N}}\circ \iota}\big[\{\bm{V}_{i}\}_{i=1}^{m},\{\bm{W}_{i}\}_{i=1}^{n}\big]> \tau,
\end{cases}.
\end{align}
where $\tau\in \real$ is a threshold. Notice that the unbiased MMD estimator of  \eqref{Eq:Negative_Estimator} can be negative despite the fact that the population MMD is non-negative. Consequently, negative values for the statistical threshold $\tau$  \eqref{Eq:Test_statistics} are admissible. Nevertheless, int our simulations, we only consider non-negative values for the threshold $\tau$. 

A Type I error is made when $\mathsf{H}_{0}$ is rejected based on the observed samples, despite the null hypothesis having generated the data. Conversely, a Type II error occurs when $\mathsf{H}_{0}$ is accepted despite the alternative hypothesis $\mathsf{H}_{1}$ being true. The \textit{significance level} $\alpha$ of a test is an upper bound on the probability of a Type I error: this is a design parameter of the test which must be set in advance, and is used to determine the threshold to which we compare the test statistic. The \textit{power of a test} is the probability of rejecting the null hypothesis $\mathsf{H}_{0}$ when it is indeed incorrect. In particular,
\begin{align}
\mathrm{Power}\df \prob(\text{reject}\ \mathsf{H}_{0}| \mathsf{H}_{1} \ \text{is true}).
\end{align}

In this sense, the statistical power controls the probability of making Type II errors.

\subsection{Empirical results on benchmark data-sets}

In Figure \ref{Fig:5}, we evaluate the power of the test for $100$ trials of hypothesis test using the test statistics of \eqref{Eq:Test_statistics}. To obtain the result, we used an autoencoder to reduce the dimension from $d=100$ to $p=50$. Clearly, for the trained kernel in Panel (a) of Figure \ref{Fig:5}, the threshold $\tau$ for which $\mathrm{Power}=1$ increases after learning the kernel via the two phase procedure described earlier. In comparison, in Panel (b), we observe that training an auto-encoder only with a fixed standard Gaussian kernel $K(\bm{x},\bm{y})=\exp(-\|\bm{x}-\bm{y}\|_{2}^{2})$ attains lower thresholds compared to our two-phase procedure. In Panel (c), we demonstrate the case of a fixed Gaussian kernel without an auto-encoder. In this case, the threshold is significantly lower due to the large dimensionality of the data. We also observer that the particle SGD for optimization we proposed in \cite{khuzani2019mean} to optimize the translation invariant kernels provides the highest statistical power for a given threshold value. Nevertheless, as we observe in Figure \ref{Fig:5}, the run-time of the particle SGD for a given number of particles is significantly higher than Algorithm \ref{Algoirthm:1}. This is due to the fact that the dimension of the particles in the particle SGD is dependent on the number of hidden layers of auto-encoder. In contrast, the particles in Algorithm \ref{Algoirthm:1} are one dimensional.

\begin{figure}[!t]	
	
	\begin{center}
		\subfigure{	\hspace{-15mm}
		\includegraphics[trim={.2cm .2cm .2cm  .6cm},width=.30\linewidth]{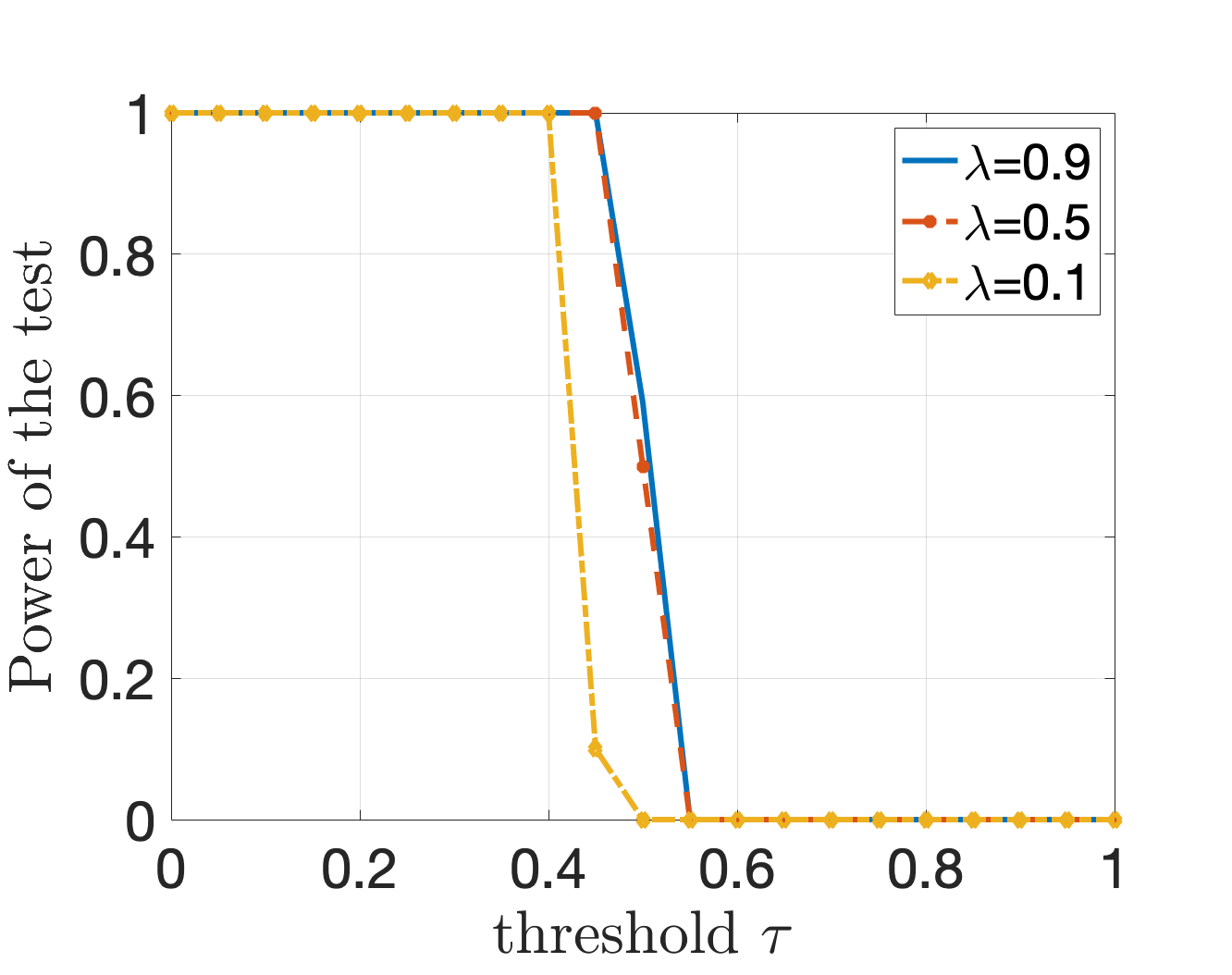} \hspace{-3mm}
					\includegraphics[trim={.2cm .2cm .2cm  .6cm},width=.30\linewidth]{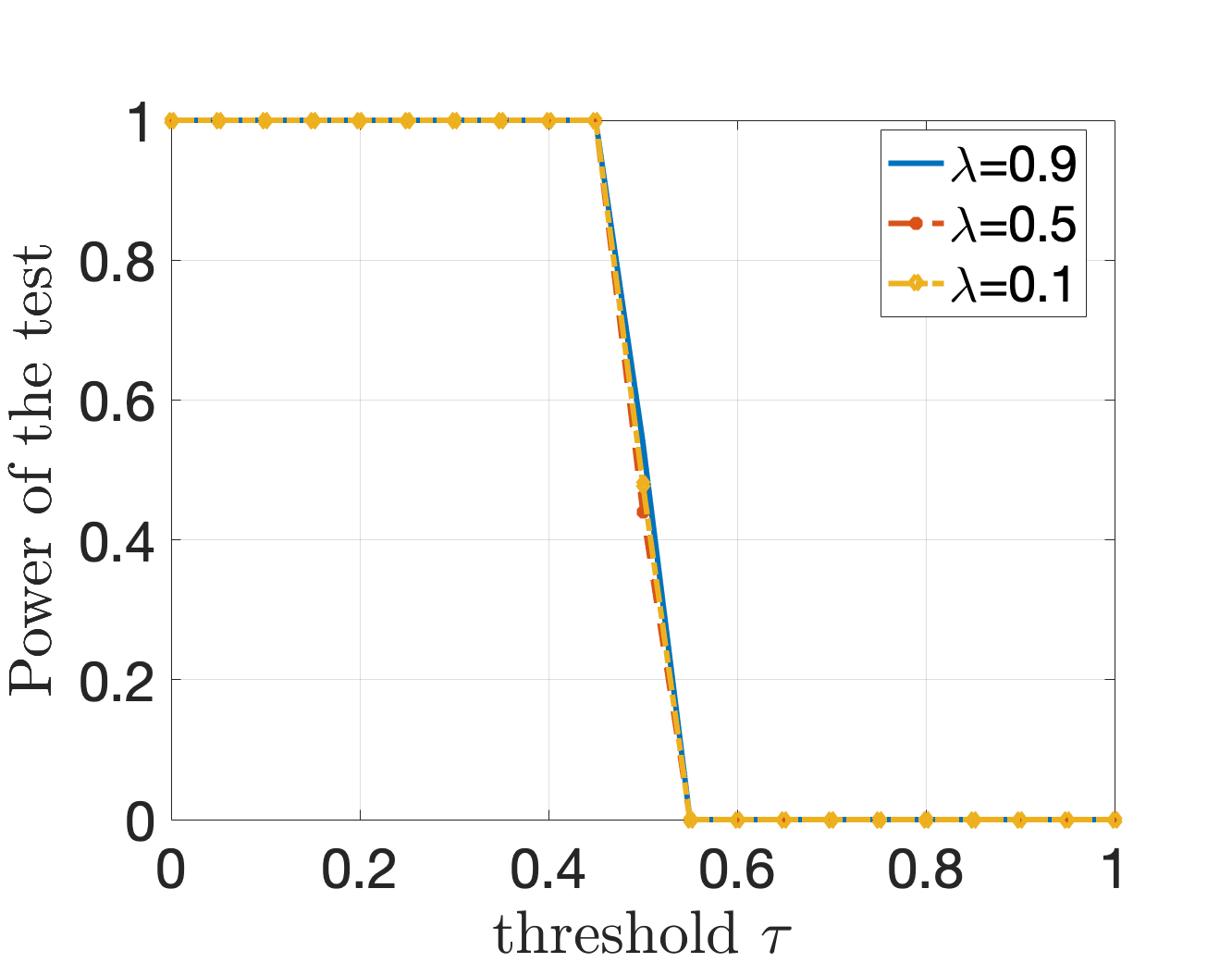} \hspace{-3mm}
			\includegraphics[trim={.2cm .2cm .2cm  .6cm},width=.30\linewidth]{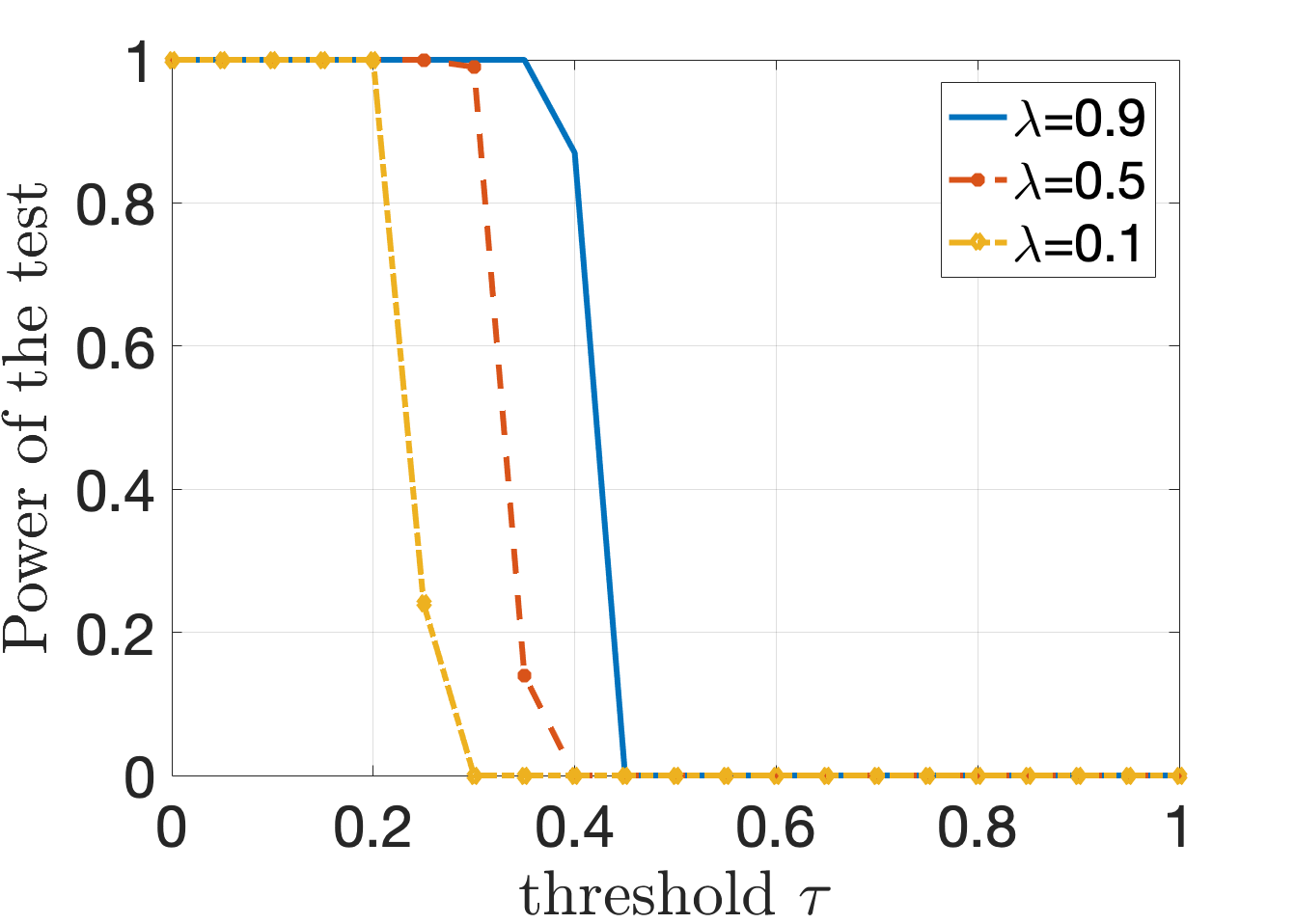}\hspace{-3mm}
			\includegraphics[trim={.2cm .2cm .2cm  .6cm},width=.30\linewidth]{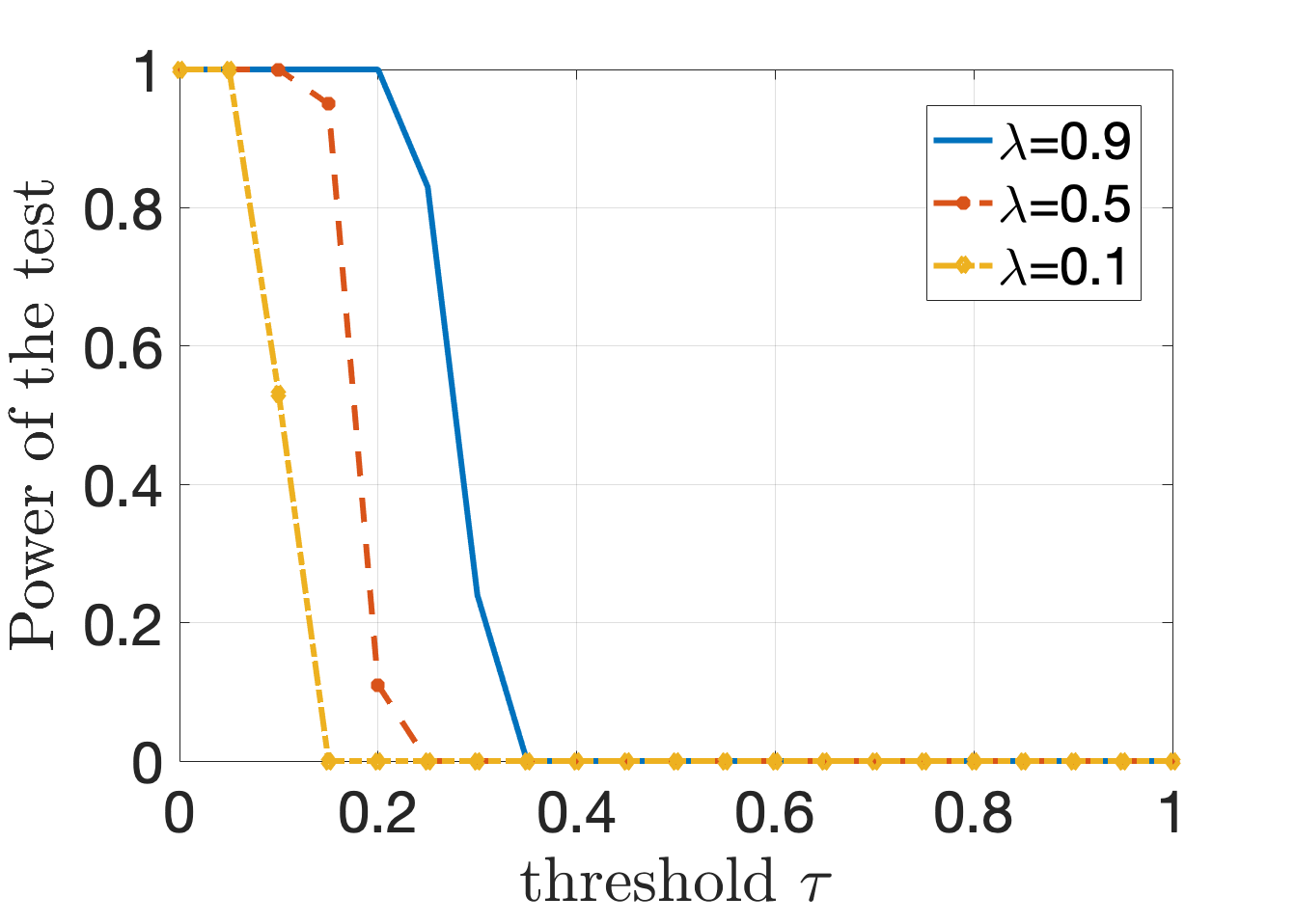}	
		}
		\\ \vspace{-2mm}
		\subfigure{\footnotesize{\hspace{5mm} (a) \hspace{36mm} (b)\hspace{36mm} (c)\hspace{36mm} (d) }}
		\vspace{-2mm}
		\caption{\footnotesize{The statistical power versus the threshold $\tau$ for the binary hypothesis testing via the unbiased estimator of the kernel MMD. The parameters for this simulations are $\lambda\in \{0.1,0.5,0.9\}$, $d=100$, $d_{0}=50$, $n+m=100$, and $N=5000$.  Panel (a): Trained radial kernel using the two-phase procedure with the particle Langevin optimization in Algorithm \ref{Algoirthm:1}, and an auto-encoder, Panel (b): Trained shift-invariant kernel using particle SGD in Eq. \eqref{Eq:SGD} and an auto-encoder, Panel (c): Trained kernel with an auto-encoder and a fixed Gaussian kernel whose the bandwidth is $\sigma=1$, Panel (d): Untrained kernel without an auto-encoder.}}
		\label{Fig:5} 
	\end{center}
\end{figure}

\begin{figure}[!t]	
	
	\begin{center}
	\subfigure{	\hspace{-15mm}
	\includegraphics[trim={.2cm .2cm .2cm  .6cm},width=.30\linewidth]{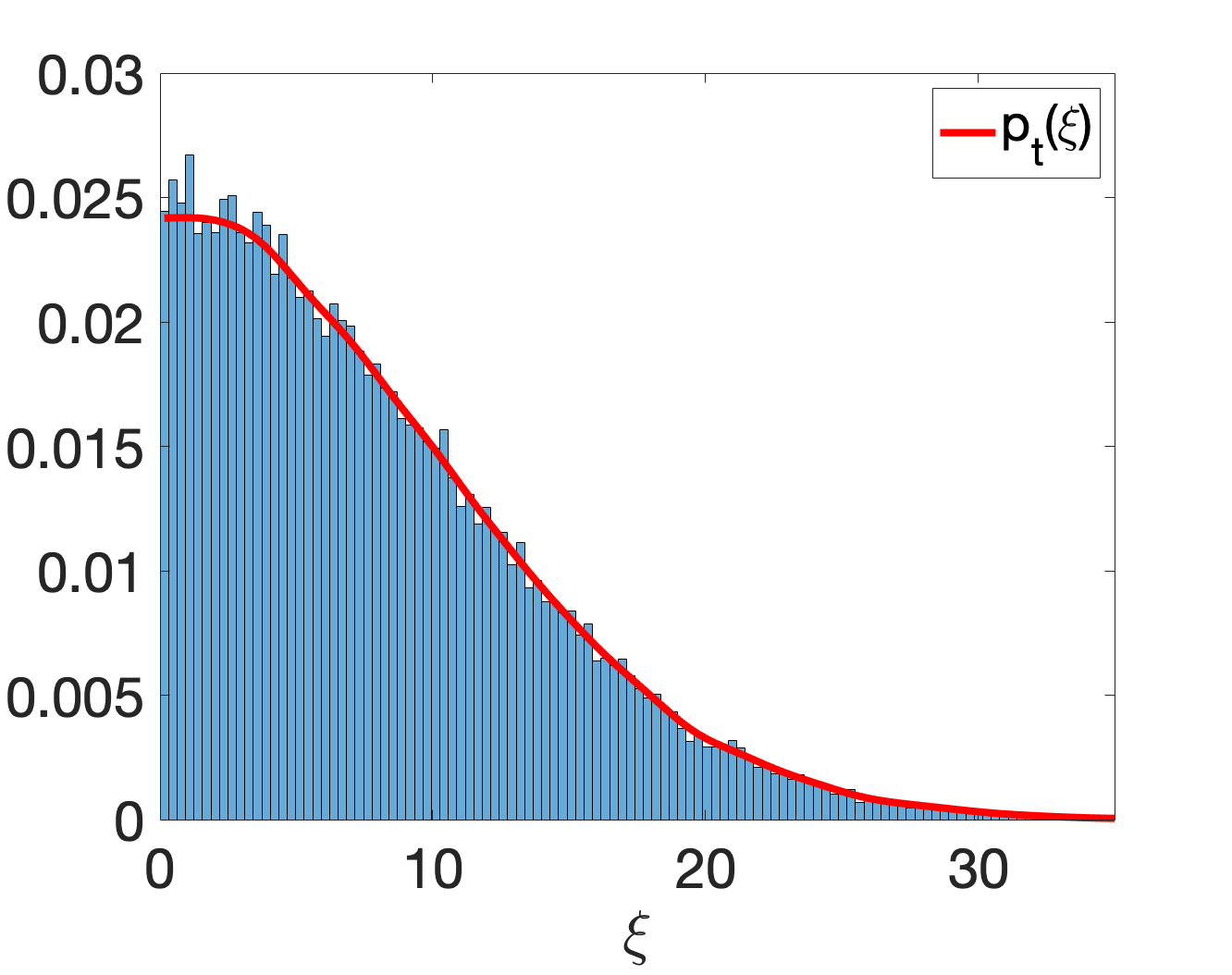} \hspace{-3mm}
	\includegraphics[trim={.2cm .2cm .2cm  .6cm},width=.30\linewidth]{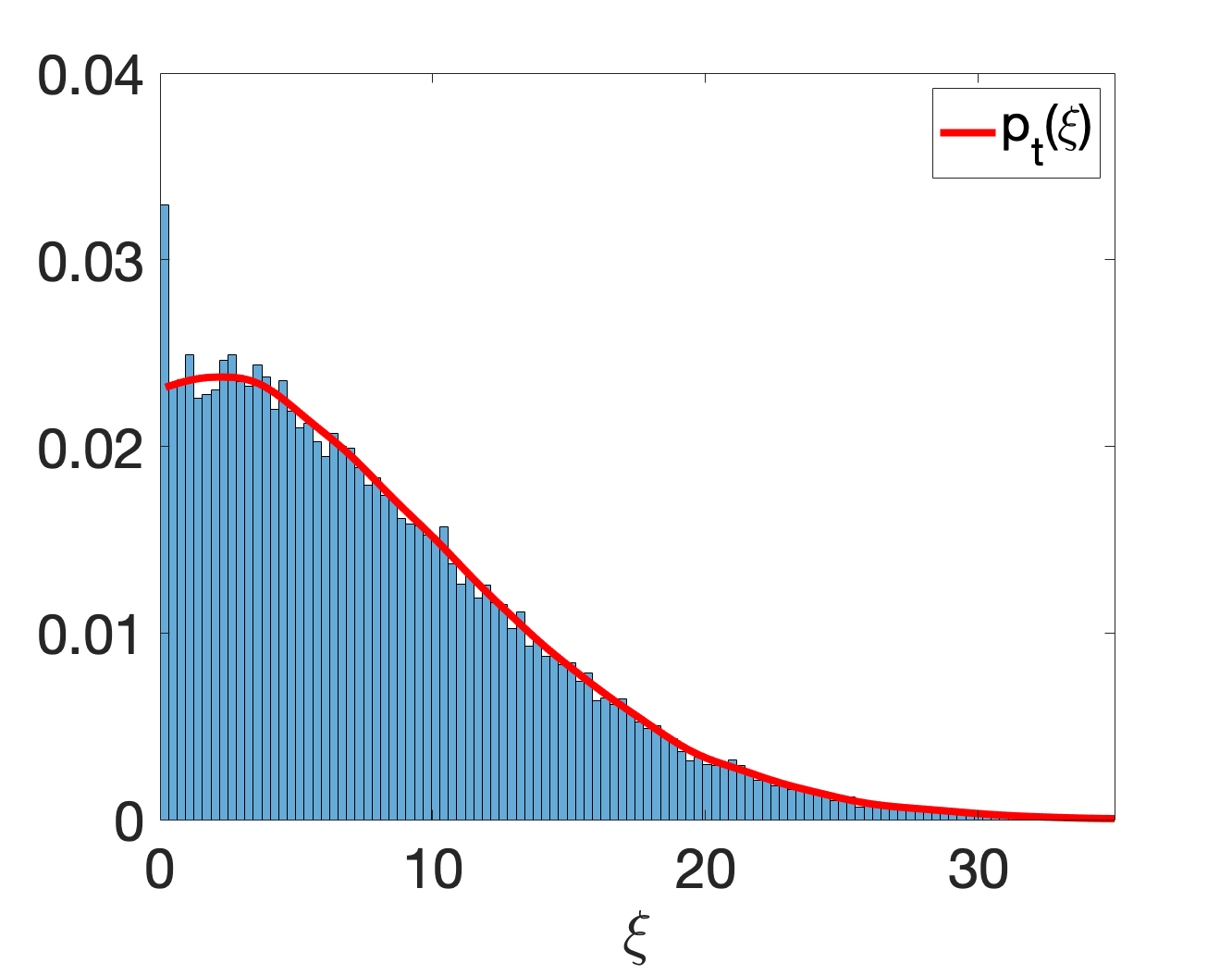} \hspace{-3mm}
	\includegraphics[trim={.2cm .2cm .2cm  .6cm},width=.30\linewidth]{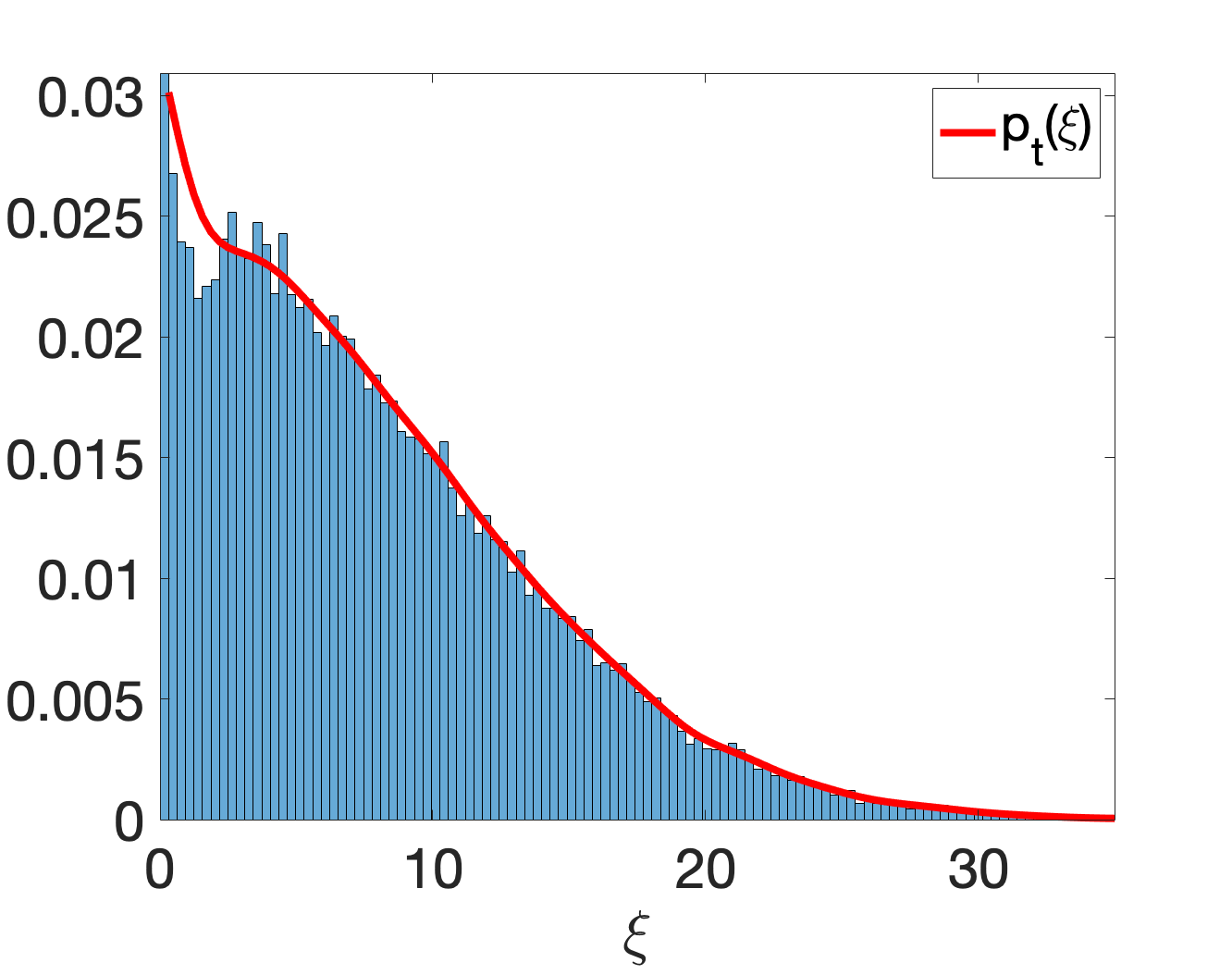}\hspace{-3mm}
	\includegraphics[trim={.2cm .2cm .2cm  .6cm},width=.30\linewidth]{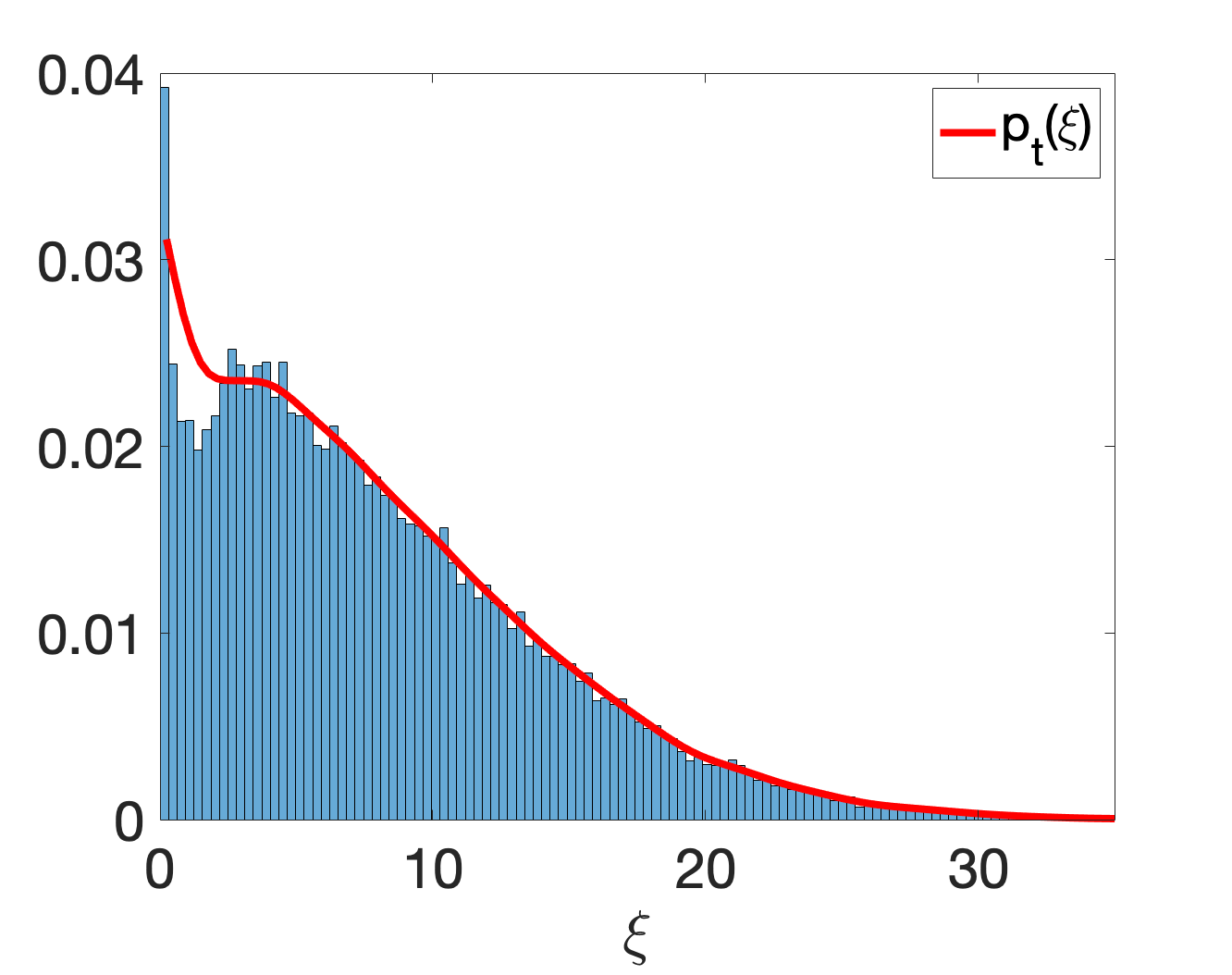}	
}
\\ \vspace{-2mm}
\subfigure{\footnotesize{\hspace{5mm} (a) \hspace{36mm} (b)\hspace{36mm} (c)\hspace{36mm} (d) }}
\vspace{-2mm}\\
	\subfigure{	\hspace{-15mm}
	\includegraphics[trim={.2cm .2cm .2cm  .6cm},width=.30\linewidth]{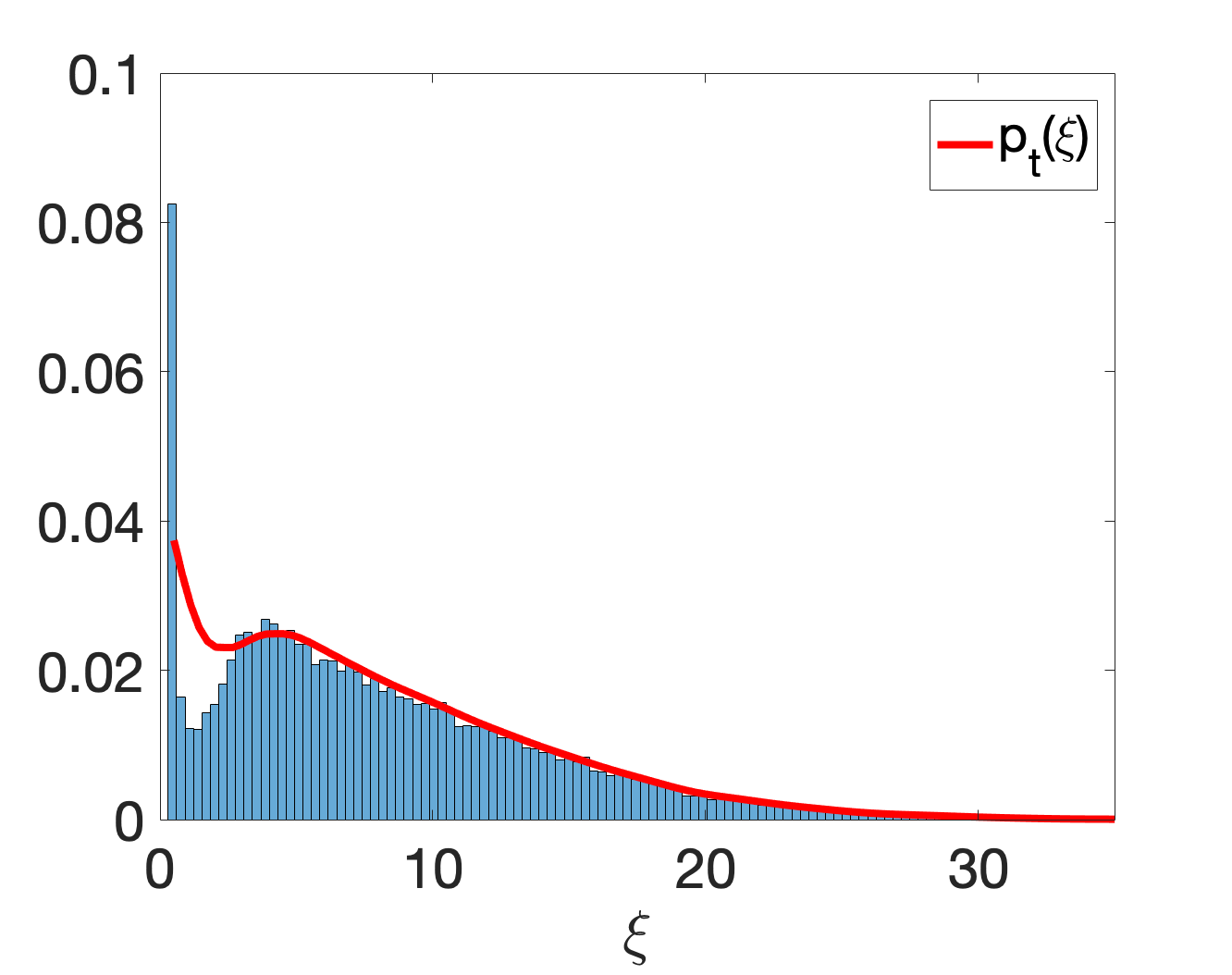} \hspace{-3mm}
	\includegraphics[trim={.2cm .2cm .2cm  .6cm},width=.30\linewidth]{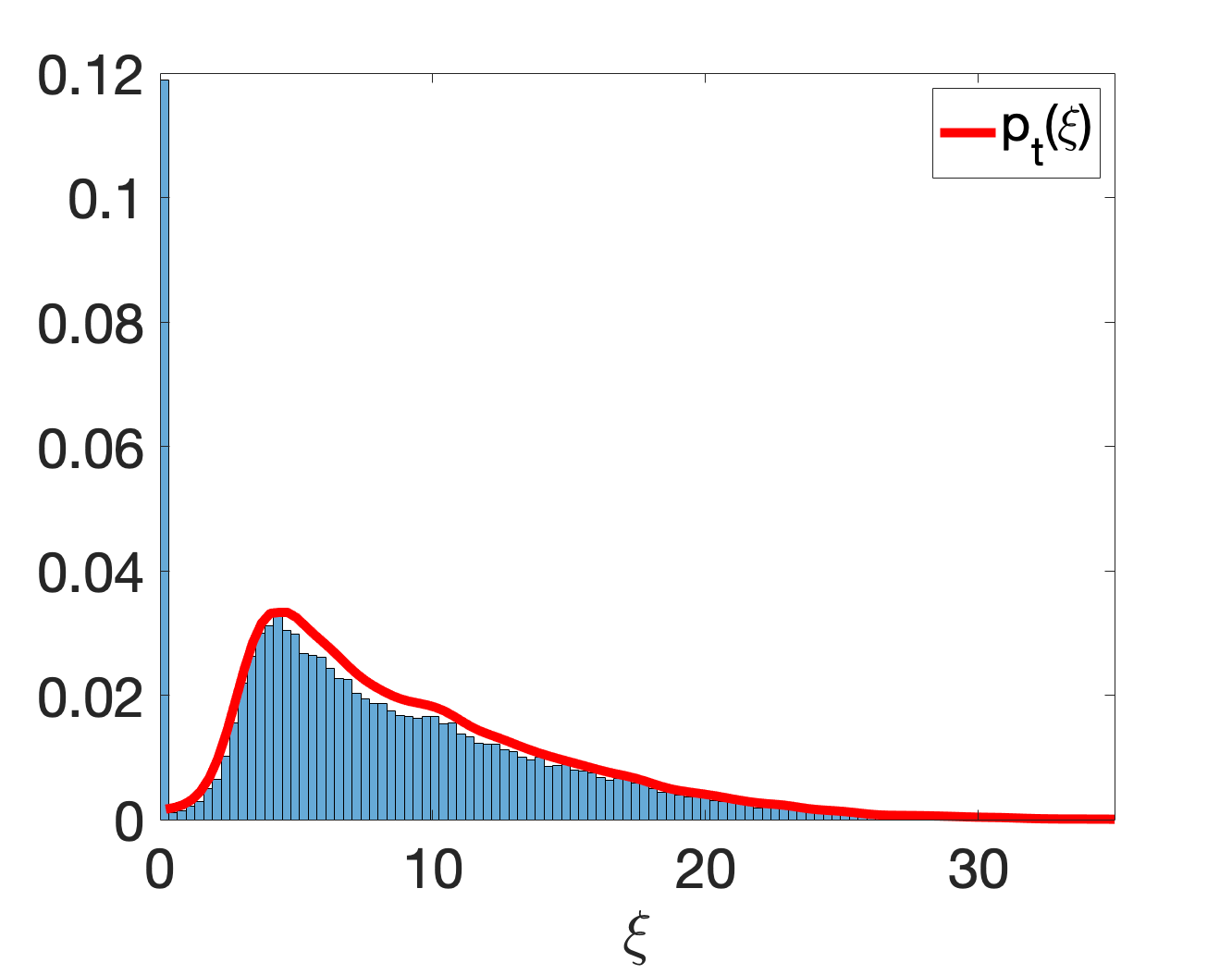} \hspace{-3mm}
	\includegraphics[trim={.2cm .2cm .2cm  .6cm},width=.30\linewidth]{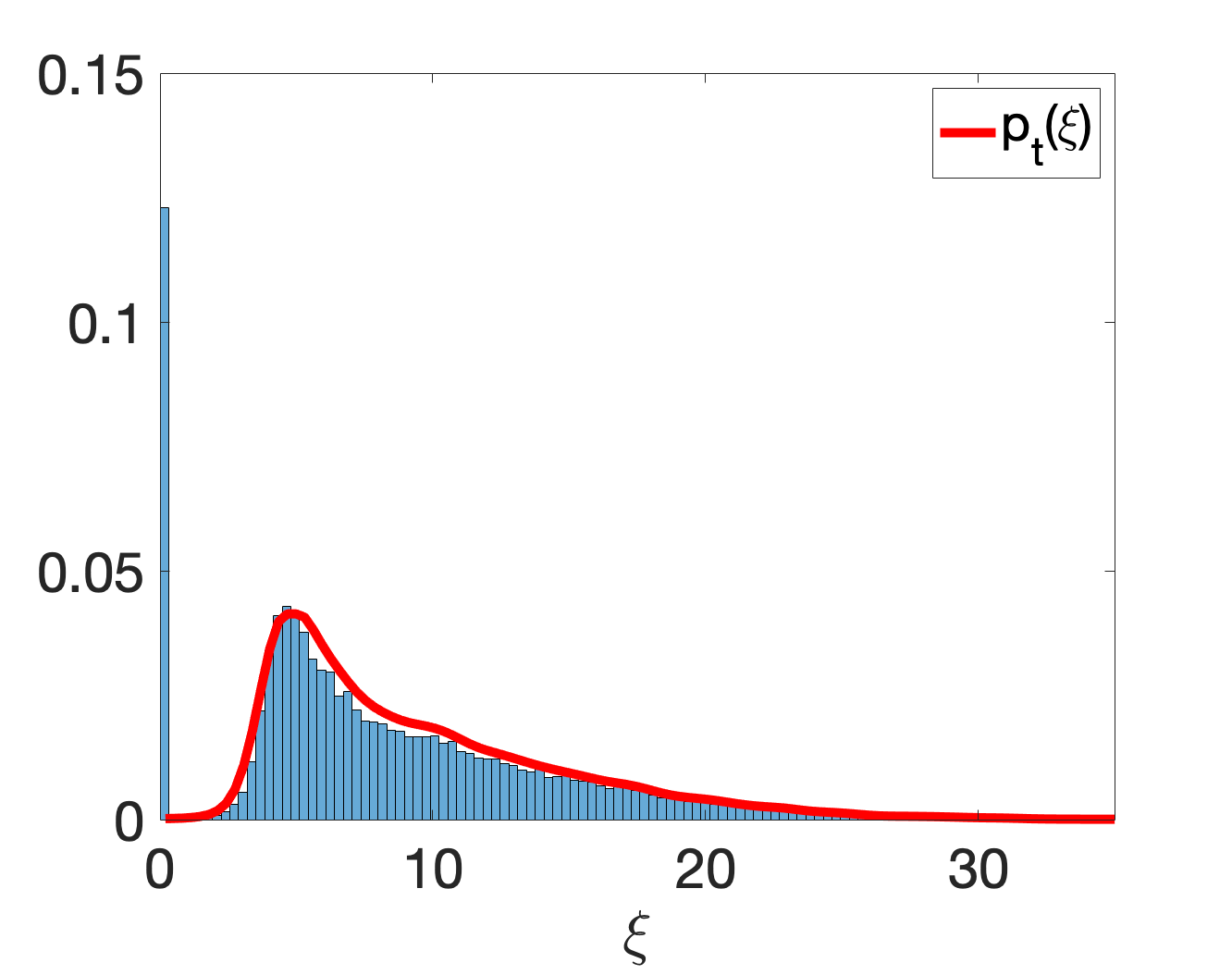}\hspace{-3mm}
	\includegraphics[trim={.2cm .2cm .2cm  .6cm},width=.30\linewidth]{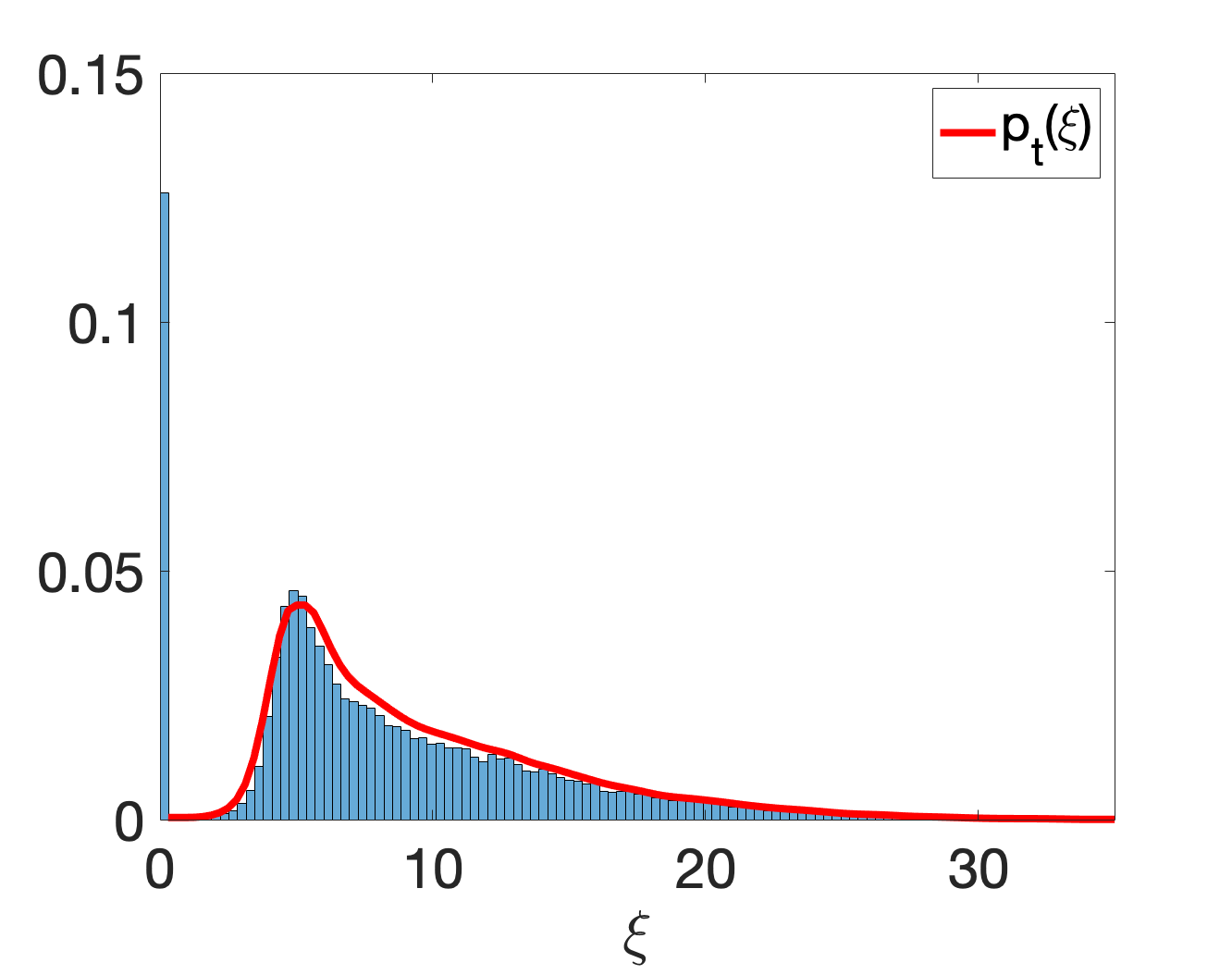}	
}
\\ \vspace{-2mm}
\subfigure{\footnotesize{\hspace{5mm} (e) \hspace{36mm} (f)\hspace{36mm} (g)\hspace{36mm} (h) }}
\vspace{-2mm}
		\caption{\footnotesize{The evolution of the histogram of Langevin particles $\xi^{1},\cdots,\xi^{N}\in \Xi=\real_{+}$ and predictions from the theory (Theorem \ref{Theorem:Mean-Field Partial Differential Equation}) with $\gamma=10000$,  $R=10000$, $\eta=10^{-5}$. Panel (a): $m=10\ (t=10^{-4}s)$, Panel (b): $m=1000 (t=0.1s)$, Panel (c):$m=5000\ (t=0.5s)$, Panel (d): $m=10000 (t=1s)$, Panel (e): $m=50000\ (t=5s)$, Panel (f):$m=100000\ (t=10 s)$, Panel (g): $m=200000\ (t=20s)$, Panel (h): a$m=300000\ (t=30s)$.}}
		\label{Fig:6} 
	\end{center}
\end{figure}

\section{Application to Locally Sensitive Hashing}

\begin{figure}[!t]	
	
	\begin{center}
		\subfigure{	
			\includegraphics[trim={.2cm .2cm .2cm  .6cm},width=.31\linewidth]{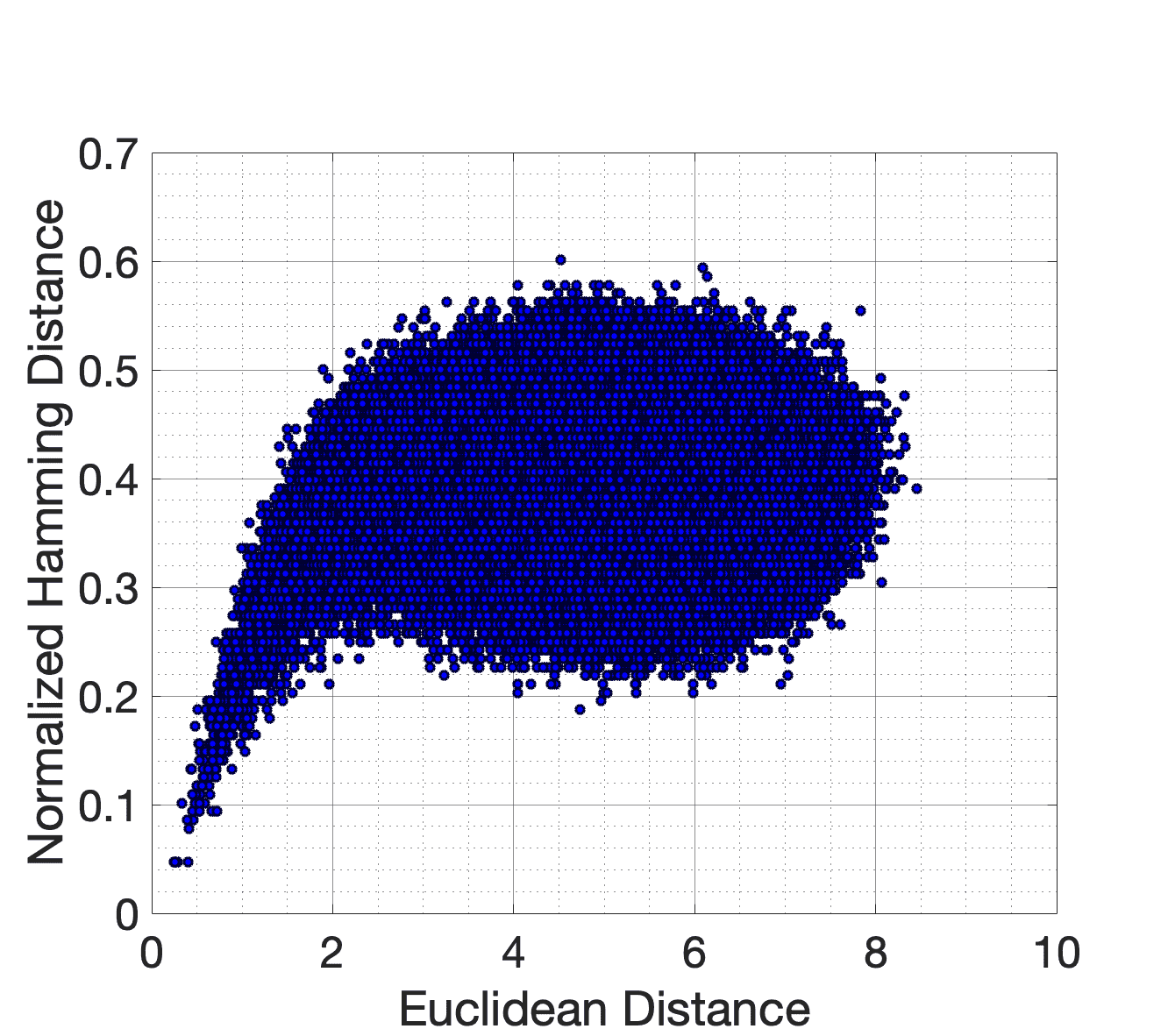} \hspace{-3mm}
			\includegraphics[trim={.2cm .2cm .2cm  .6cm},width=.31\linewidth]{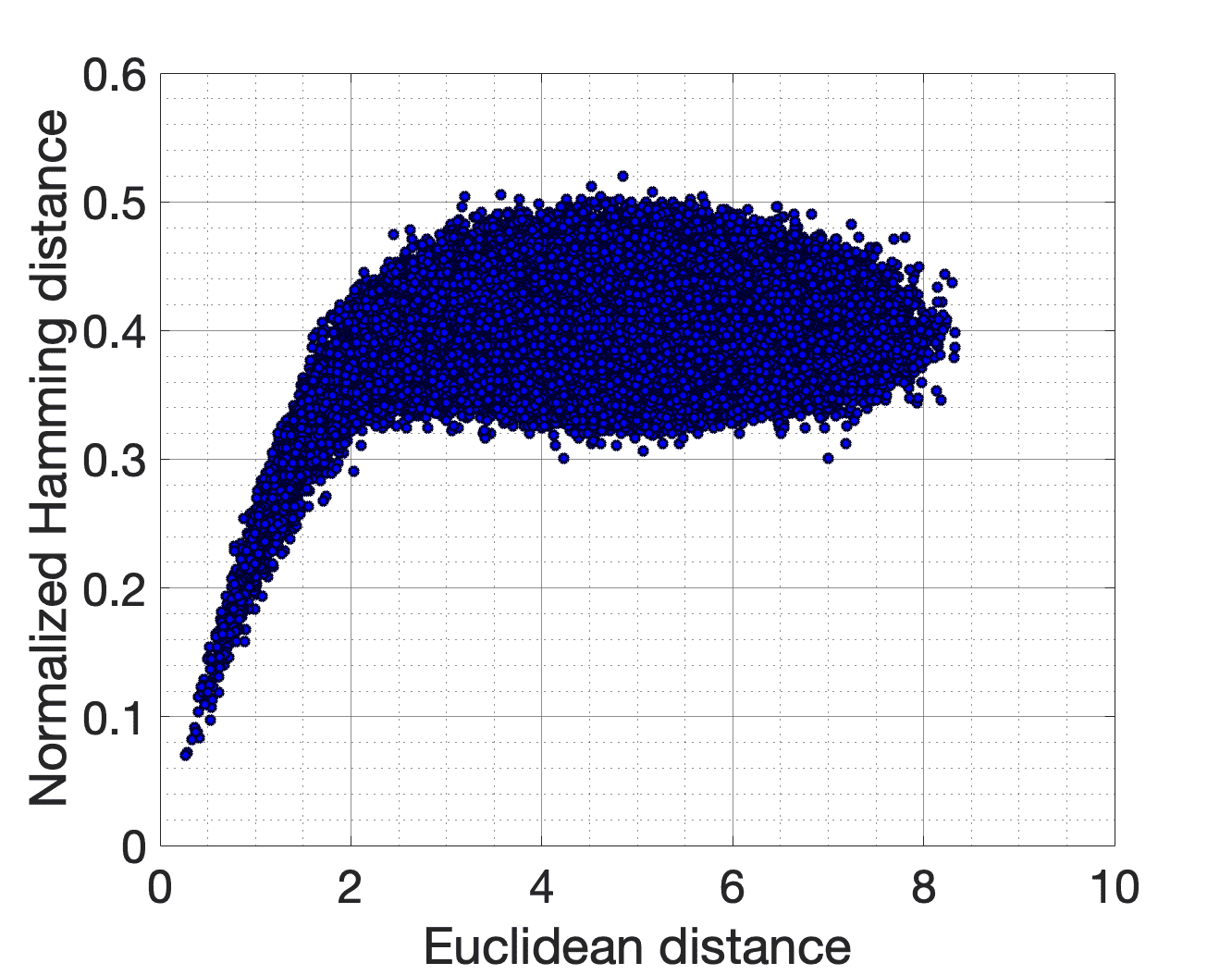} \hspace{-3mm}
			\includegraphics[trim={.2cm .2cm .2cm  .6cm},width=.31\linewidth]{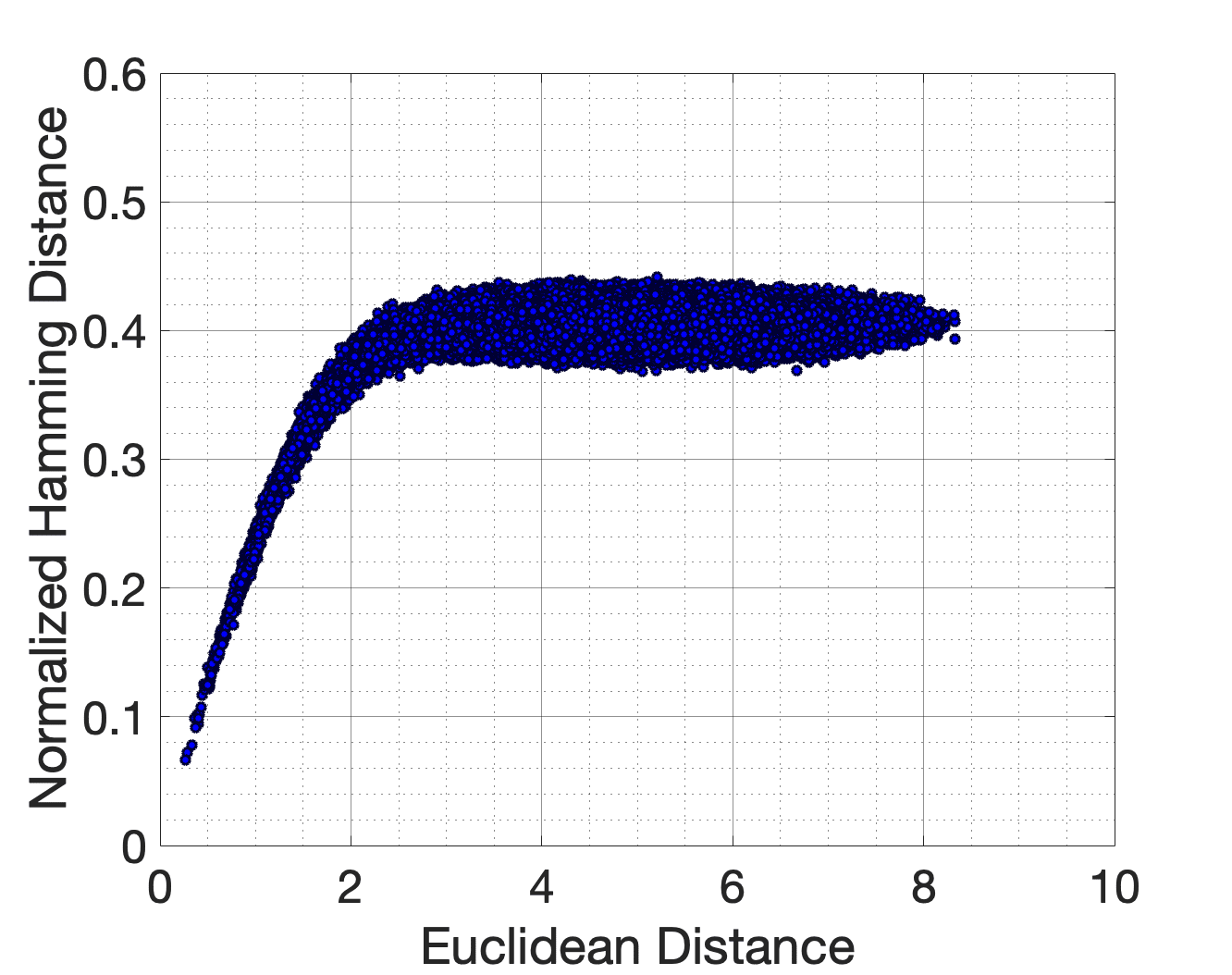}\hspace{-3mm}
		}
		\\ \vspace{-2mm}
		\subfigure{\footnotesize{\hspace{5mm} (a) \hspace{35mm} (b)\hspace{45mm} (c) }}\\
			\subfigure{
			\includegraphics[trim={.2cm .2cm .2cm  .6cm},width=.31\linewidth]{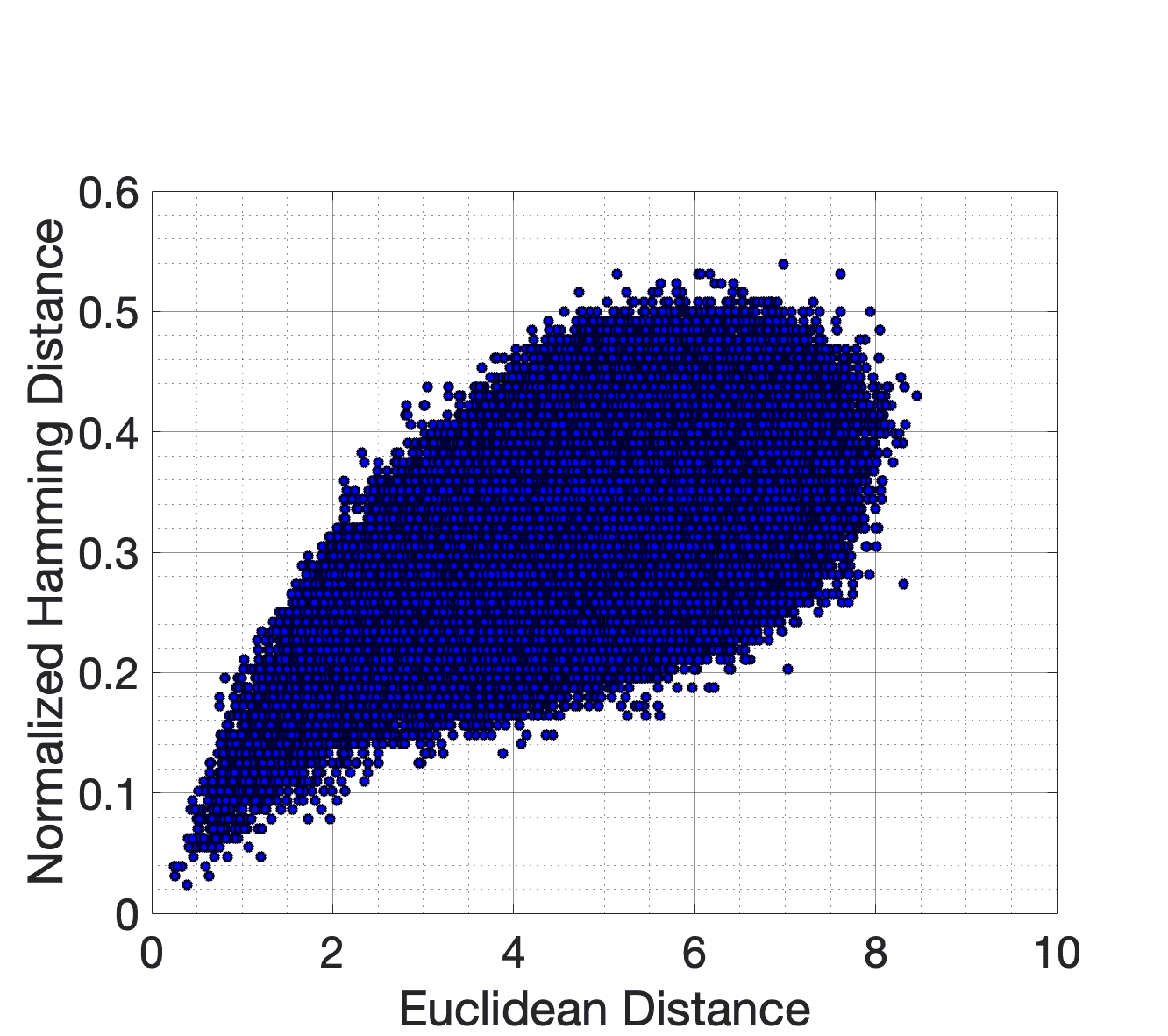} \hspace{-3mm}
			\includegraphics[trim={.2cm .2cm .2cm  .6cm},width=.31\linewidth]{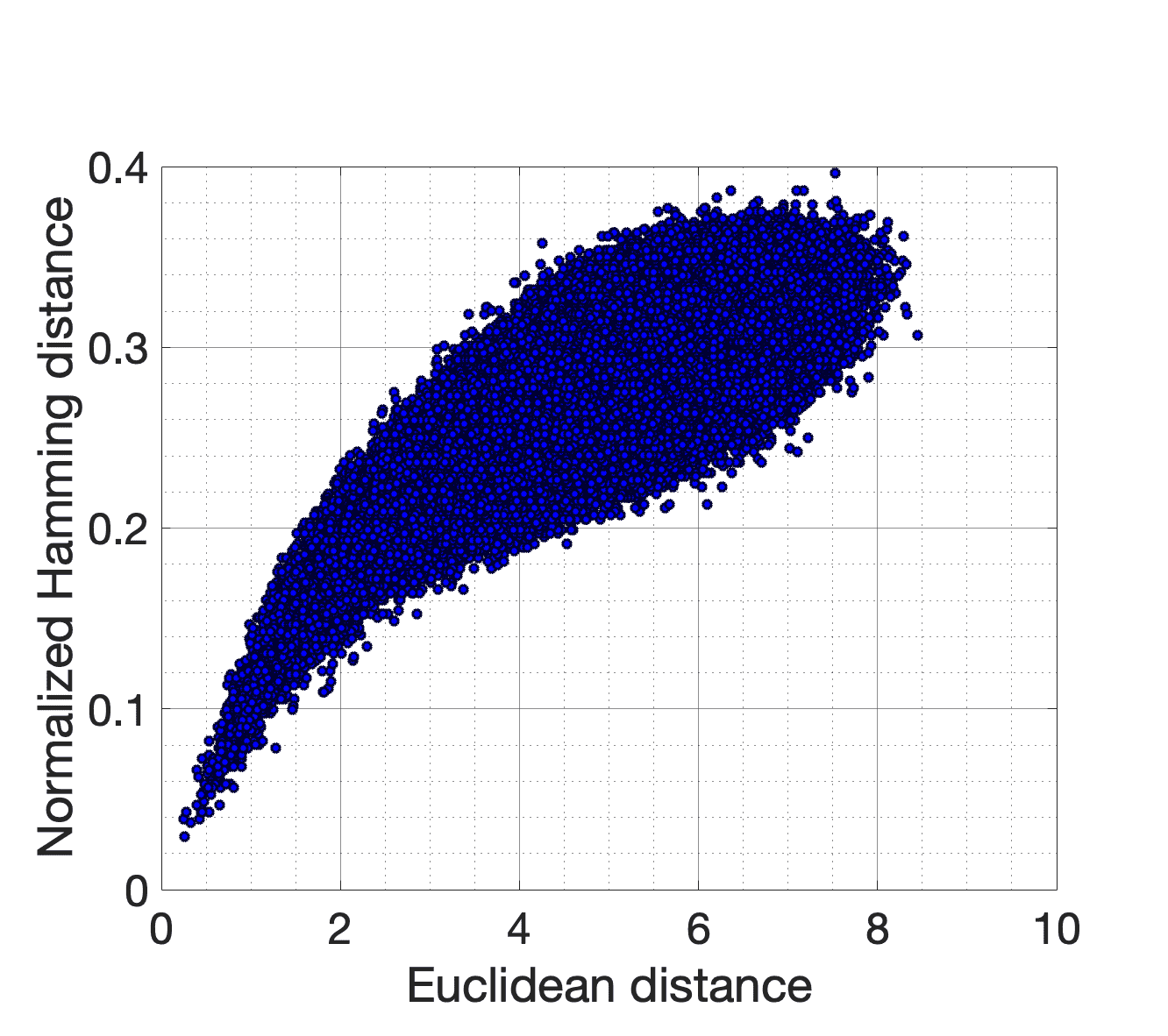} \hspace{-3mm}
			\includegraphics[trim={.2cm .2cm .2cm  .6cm},width=.31\linewidth]{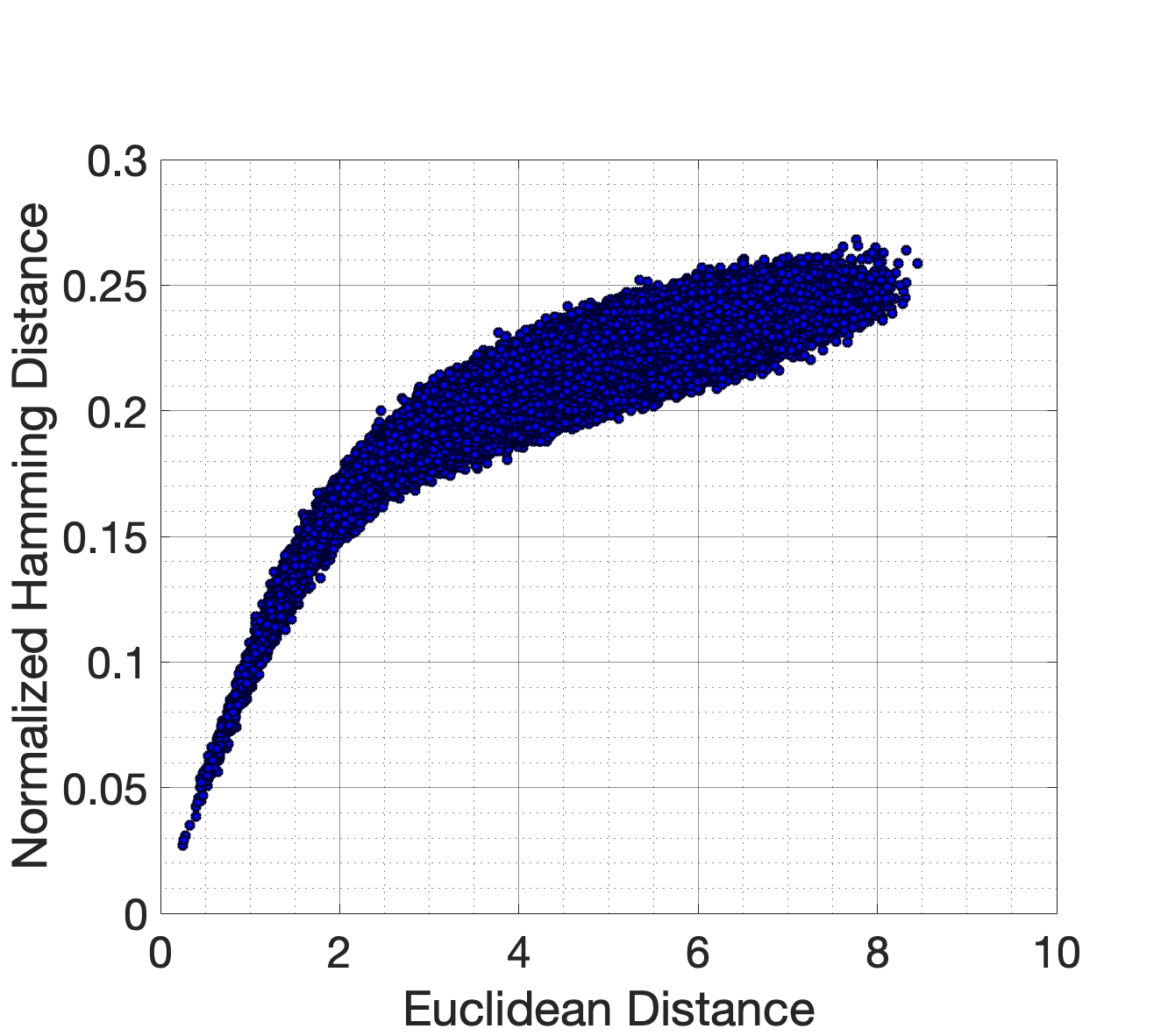}\hspace{-3mm}
		}
		\\ \vspace{-2mm}
		\subfigure{\footnotesize{\hspace{5mm} (d) \hspace{35mm} (e)\hspace{35mm} (f) }}
			\caption{\footnotesize{The scatter plots of the normalized Hamming distance versus Euclidean distance for the locally sensitive hash function of Lemma \ref{Lemma:KHF}. The hash functions are generated with the random features associated with an untrained Gaussian kernel with the bandwidth parameter of $\sigma=1$ (top row), and a trained kernel with Algorithm \ref{Algoirthm:1} in conjunction with a $k$-mean clustering method (bottom row). Panels (a),(d): 128 bits, Panels (b),(e): 512 bits, Panels (c),(f): 4096 bits.}}
		\label{Fig:5} 
	\end{center}
\end{figure}

\begin{figure}[!t]
	\begin{center}
	\subfigure{	\hspace{-4mm}
		\includegraphics[trim={.2cm .2cm .2cm  .6cm},width=.35\linewidth]{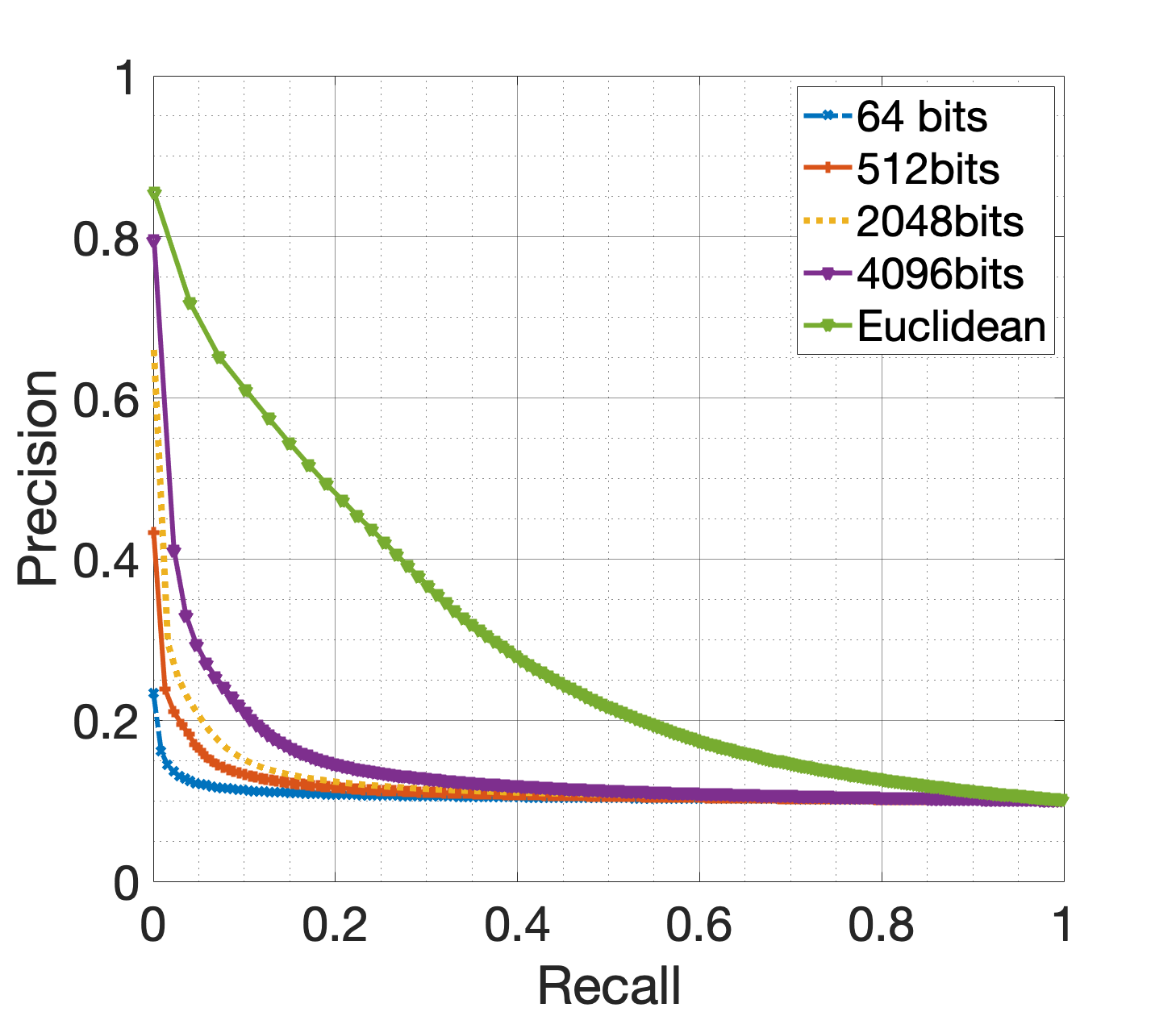} \hspace{-3mm}
		\includegraphics[trim={.2cm .2cm .2cm  .6cm},width=.35\linewidth]{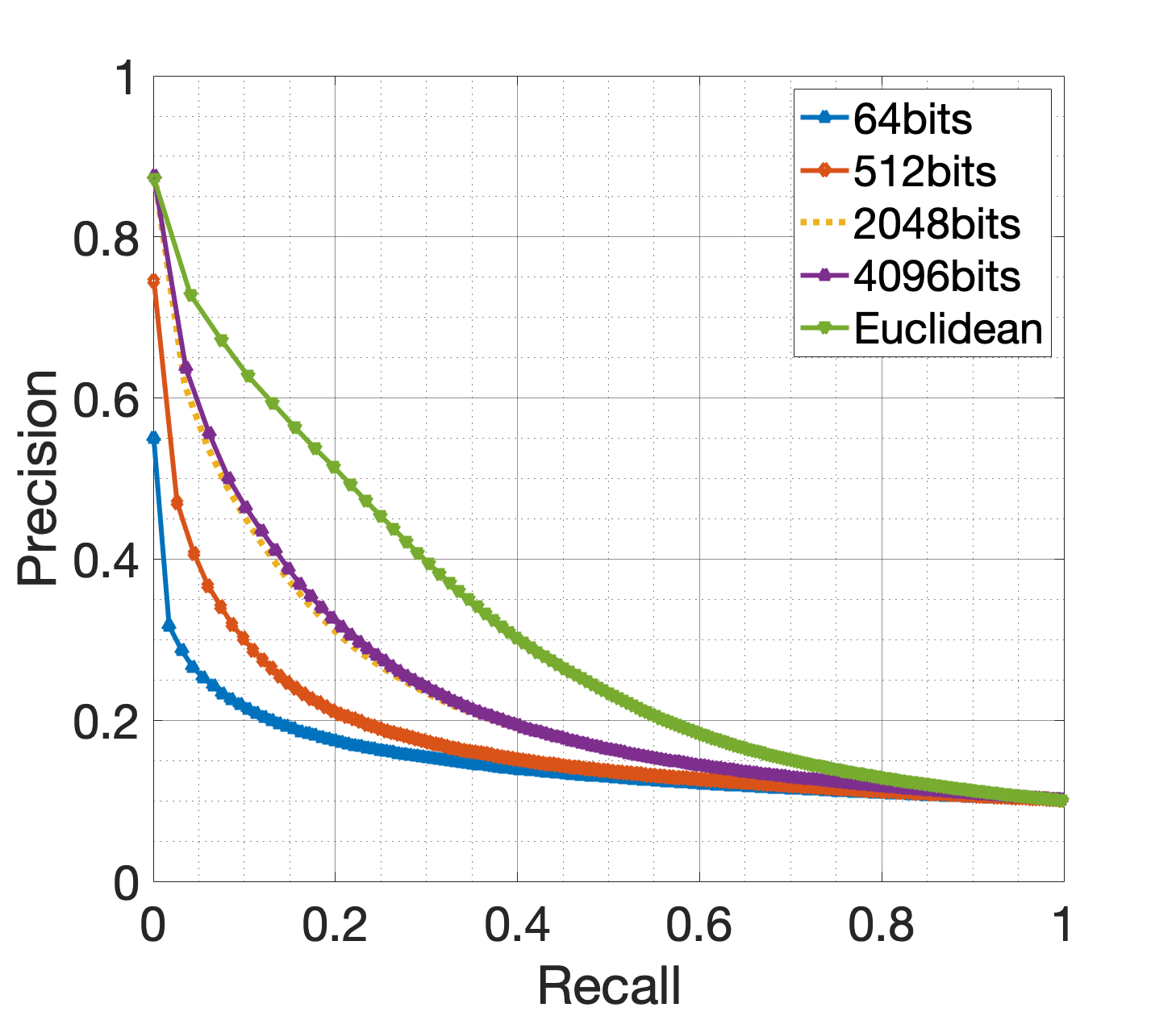} \hspace{-3mm}
\includegraphics[trim={.2cm .2cm .2cm  .6cm},width=.35\linewidth]{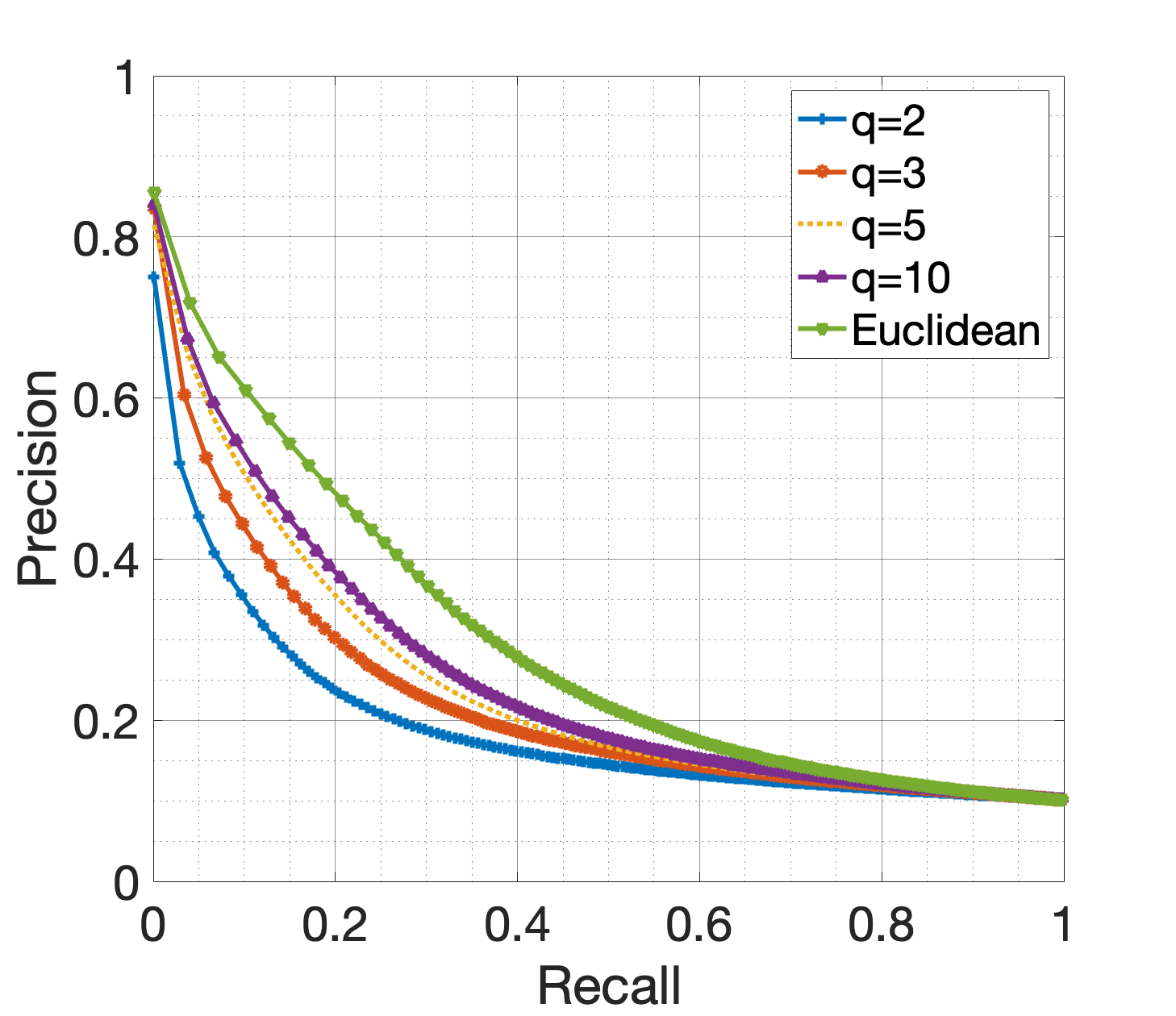} \hspace{-3mm}}
	\\ \vspace{-2mm}
	\subfigure{\footnotesize{\hspace{5mm} (a) \hspace{40mm} (b)  \hspace{45mm} (c)  }}
	\caption{\footnotesize{The precision-recall curves for the kernel locally sensitive hashing. Panel (a): hash function $F^{n}:\mathcal{X}\rightarrow \mathbb{Z}_{2}^{n}$ of Lemma \ref{Lemma:KHF} with an untrained Gaussian kernel, Panel (b): hash function $F^{n}:\mathcal{X}\rightarrow \mathbb{Z}_{2}^{n}$ of Lemma \ref{Lemma:KHF} with the trained kernel with Algorithm \ref{Algoirthm:1}, where the labeled data for kernel training is generated via the $k$-mean clustering with $k=10$, and one versus all rule. Panel (c): hash function $G^{n}:\mathcal{X}\rightarrow \mathbb{Z}_{q}^{n}$ of Lemma \ref{Lemma:Kernel_LSH_via_Embedding} with the quantization levels $q\in \{2,3,5,10\}$, and the fixed code-word length of $n=512$.}}
	
	\label{Fig:6} 
\end{center}
\end{figure}

\begin{figure}[!t]
	\begin{center}
		\subfigure{		
			\includegraphics[trim={.2cm .2cm .2cm  .6cm},width=.3\linewidth]{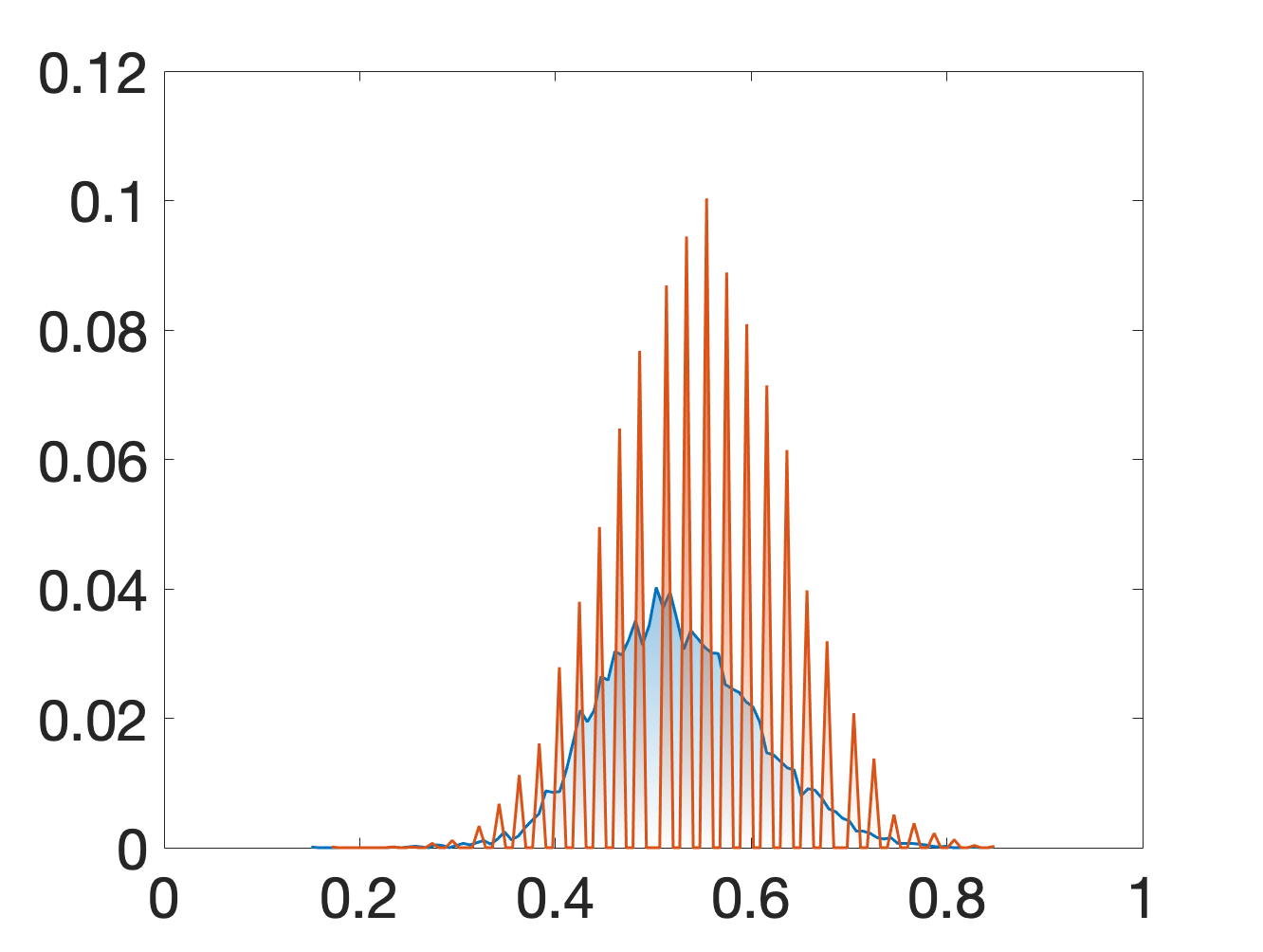} \hspace{-3mm}
			\includegraphics[trim={.2cm .2cm .2cm  .6cm},width=.3\linewidth]{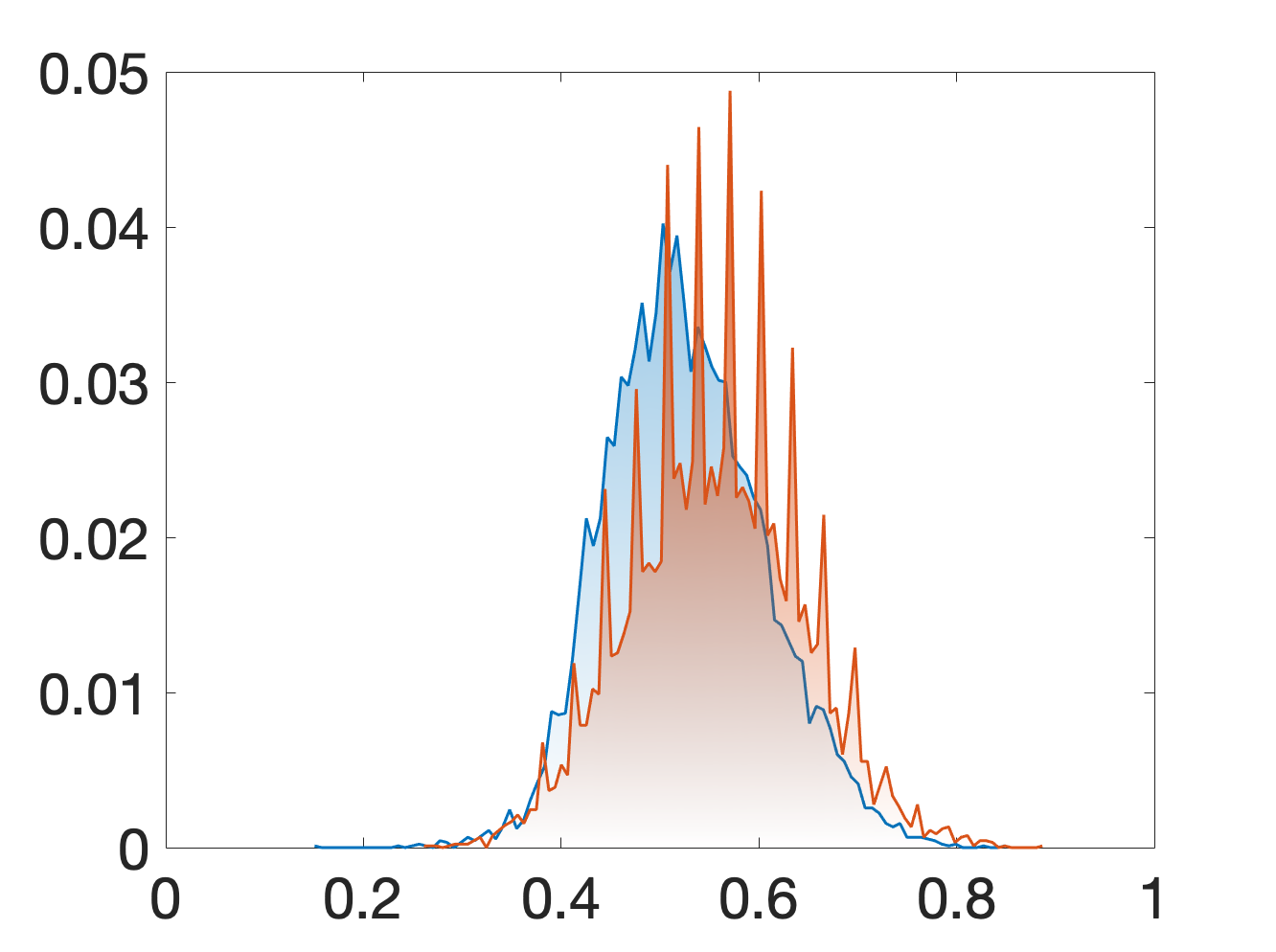} \hspace{-3mm}
				\includegraphics[trim={.2cm .2cm .2cm  .6cm},width=.3\linewidth]{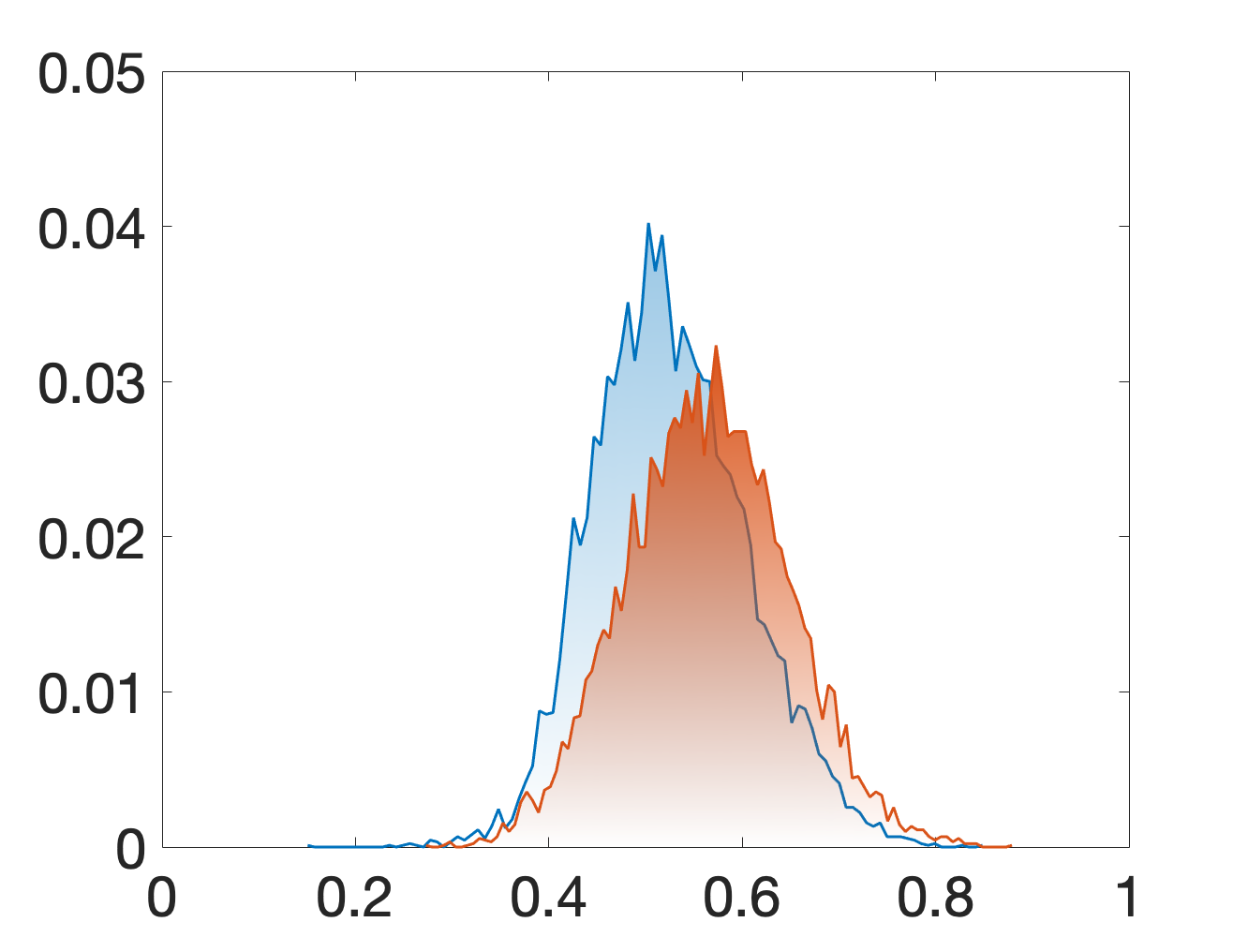} \hspace{-3mm}
		}
		\\ \vspace{-2mm}	\subfigure{\hspace{4mm}(a)\footnotesize{\hspace{36mm} (b) \hspace{36mm} (c) }}
		\caption{\footnotesize{The (normalized) histogram of the (rescaled) Lee distance (red color) and the (rescaled) Euclidean distance (blue color) of a sample query point from each point of the MNIST data-base for different quantization level $q$ in Lemma \ref{Lemma:Kernel_LSH_via_Embedding}, Panel (a): $q=2$, Panel (b): $q=20$, Panel (c): $q=50$. Increasing the quantization level $q$ in the Lee distance yields a more fine grained approximation for the Euclidean distance. (best viewed in color)}}
		
		\label{Fig:7} 
	\end{center}
\end{figure}

\begin{figure}[!t]	
	\begin{center}
		\subfigure{	
			\includegraphics[trim={.2cm .2cm .2cm  .6cm},width=.25\linewidth]{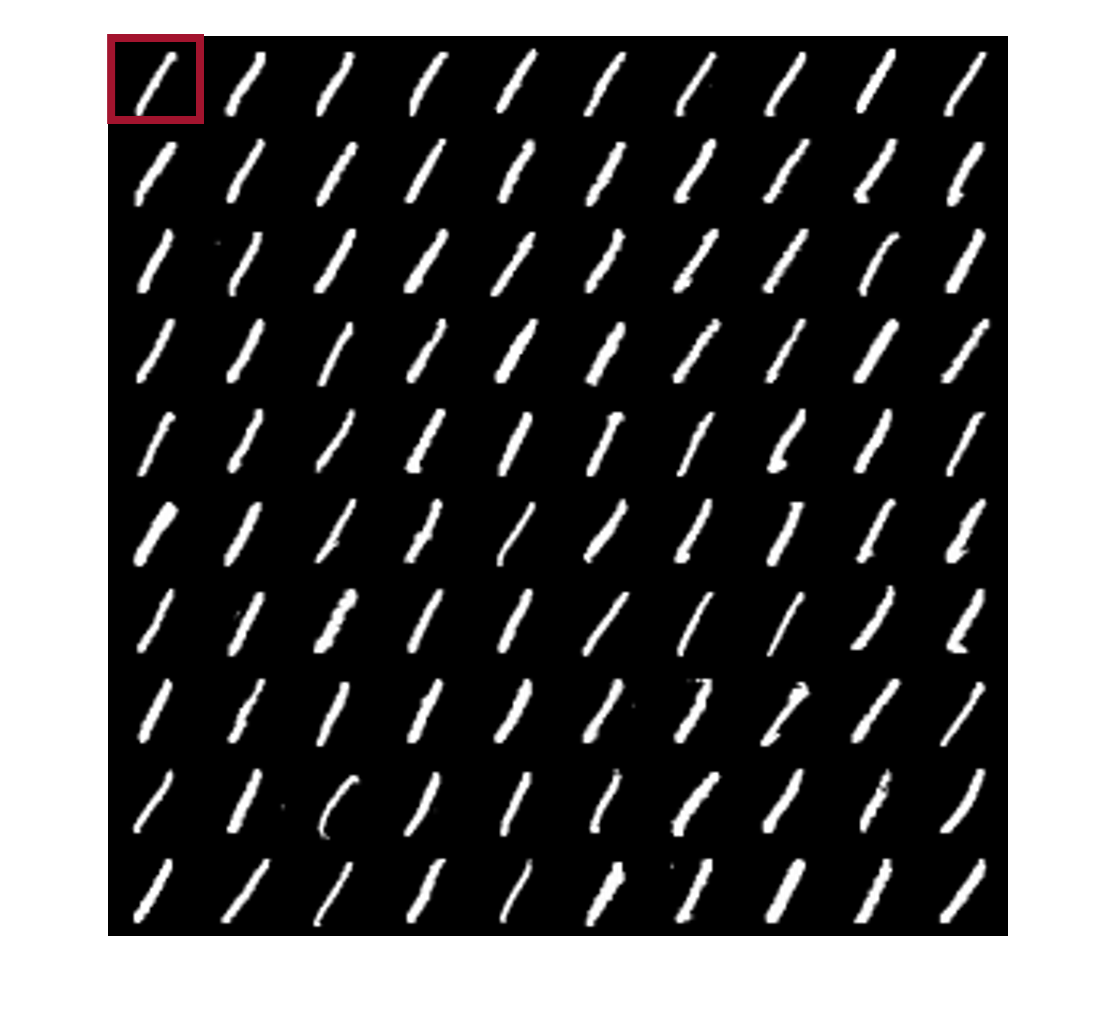} \hspace{-3mm}
			\includegraphics[trim={.2cm .2cm .2cm  .6cm},width=.25\linewidth]{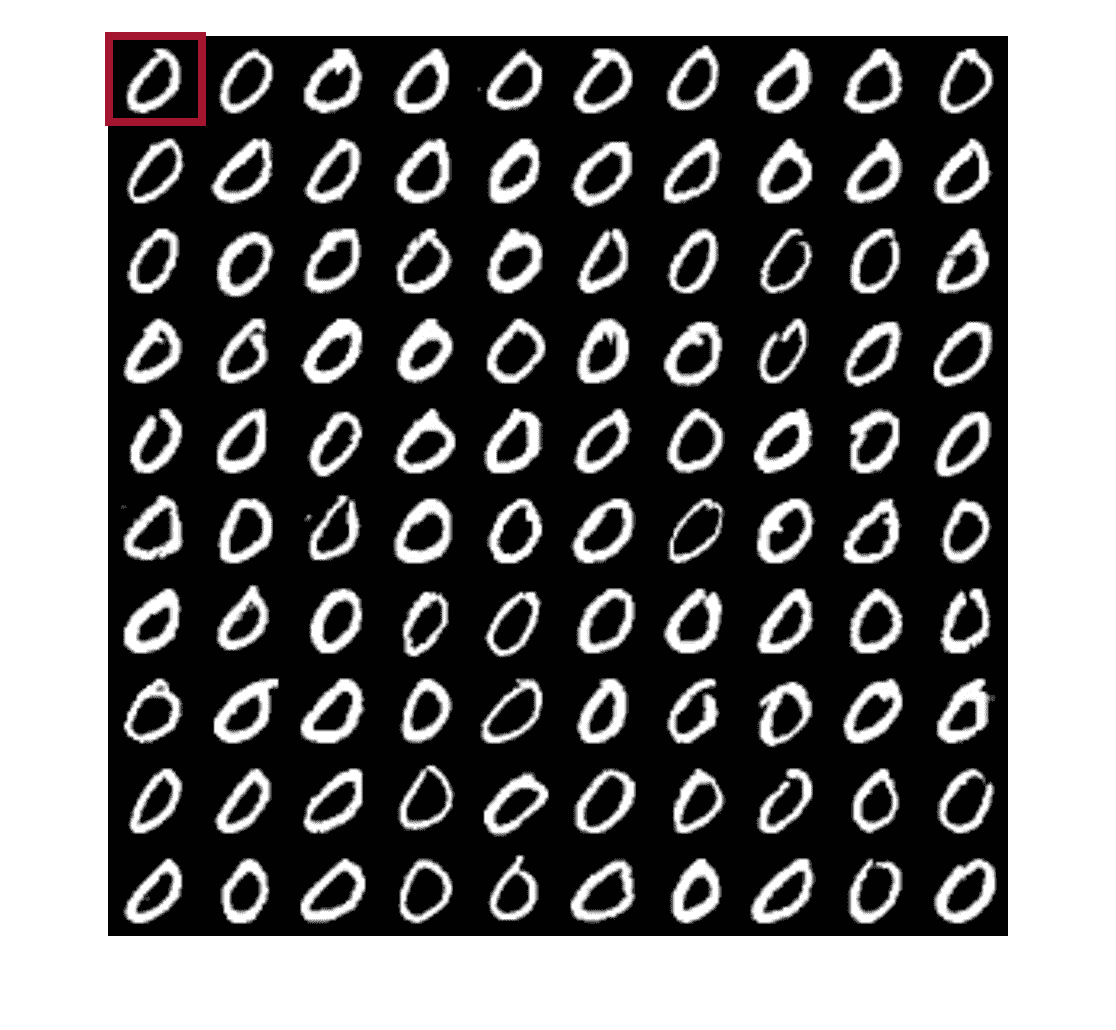} \hspace{-3mm}
			\includegraphics[trim={.2cm .2cm .2cm  .6cm},width=.25\linewidth]{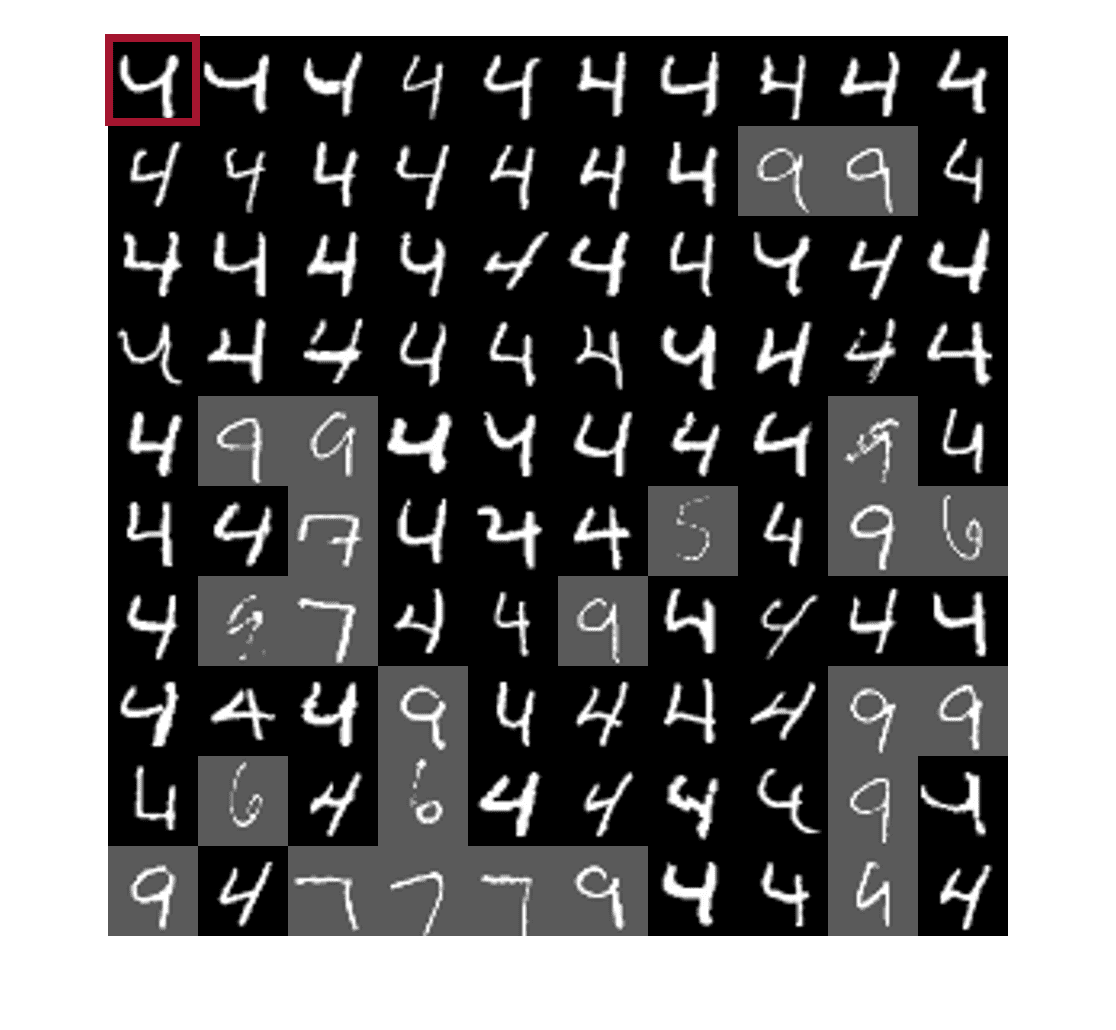}\hspace{-3mm}
		}
		\\ \vspace{-3mm}
		\subfigure{\footnotesize{\hspace{5mm} (a) \hspace{27mm} (b)\hspace{27mm} (c) }}
		\vspace{-3mm}
		\\
		\subfigure{	
		\includegraphics[trim={.2cm .2cm .2cm  .6cm},width=.25\linewidth]{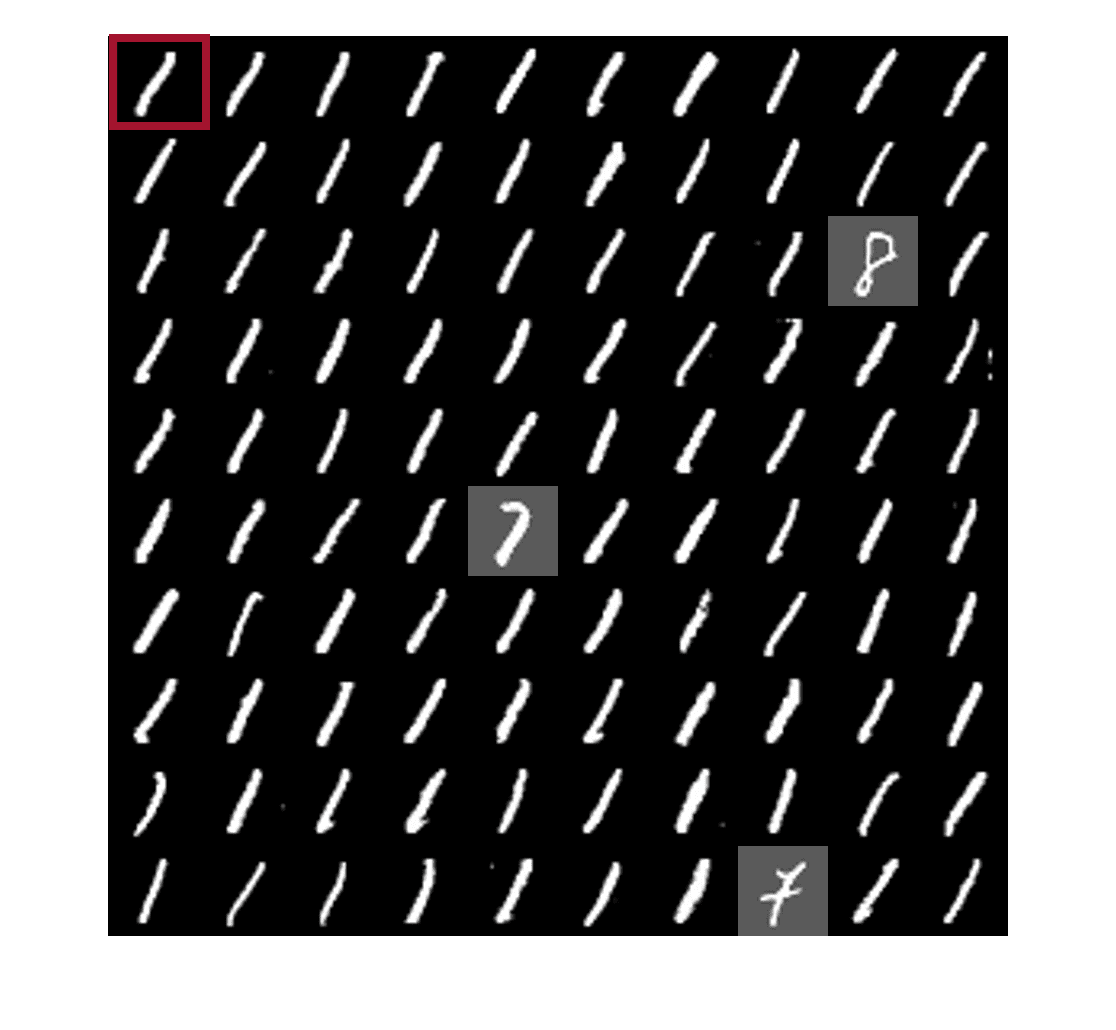} \hspace{-3mm}
		\includegraphics[trim={.2cm .2cm .2cm  .6cm},width=.25\linewidth]{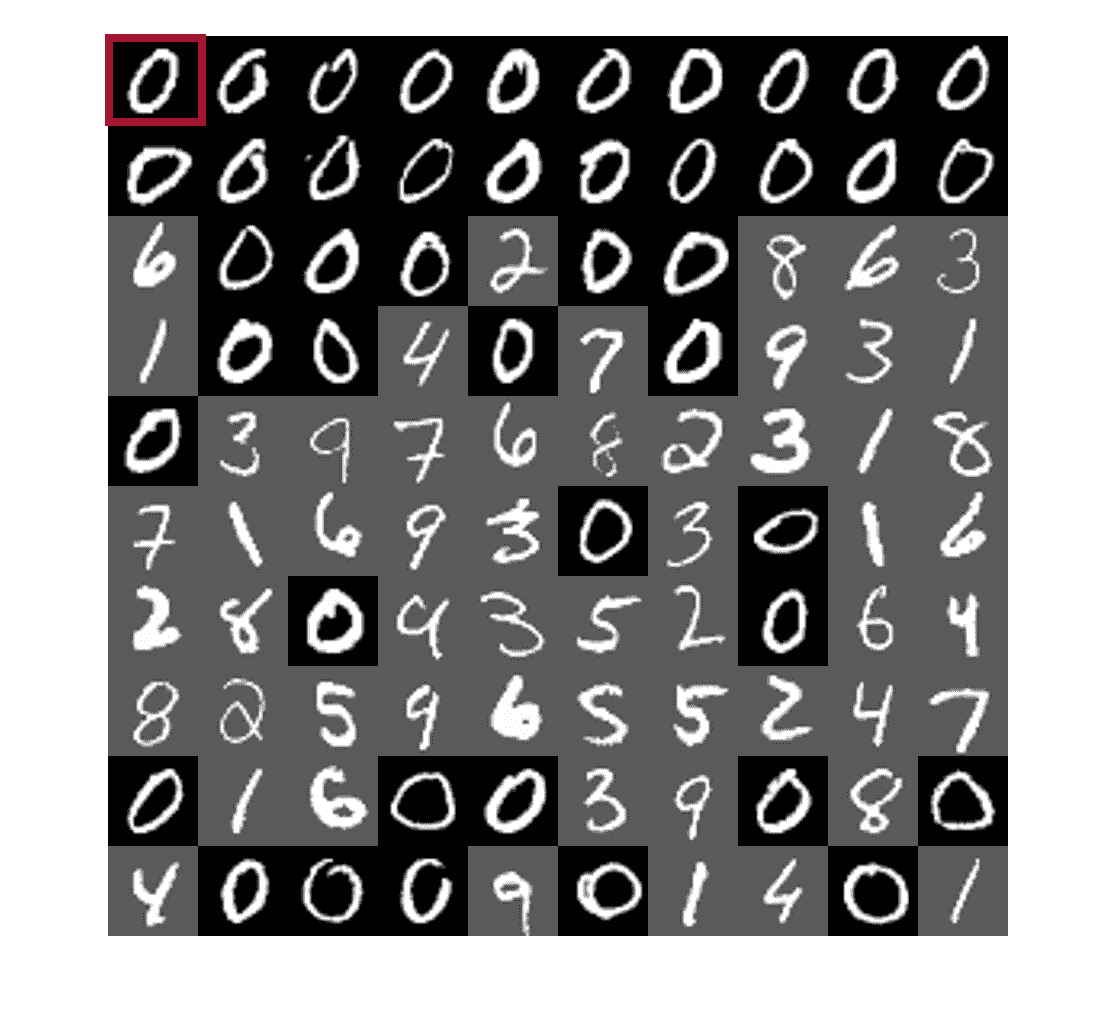} \hspace{-3mm}
		\includegraphics[trim={.2cm .2cm .2cm  .6cm},width=.25\linewidth]{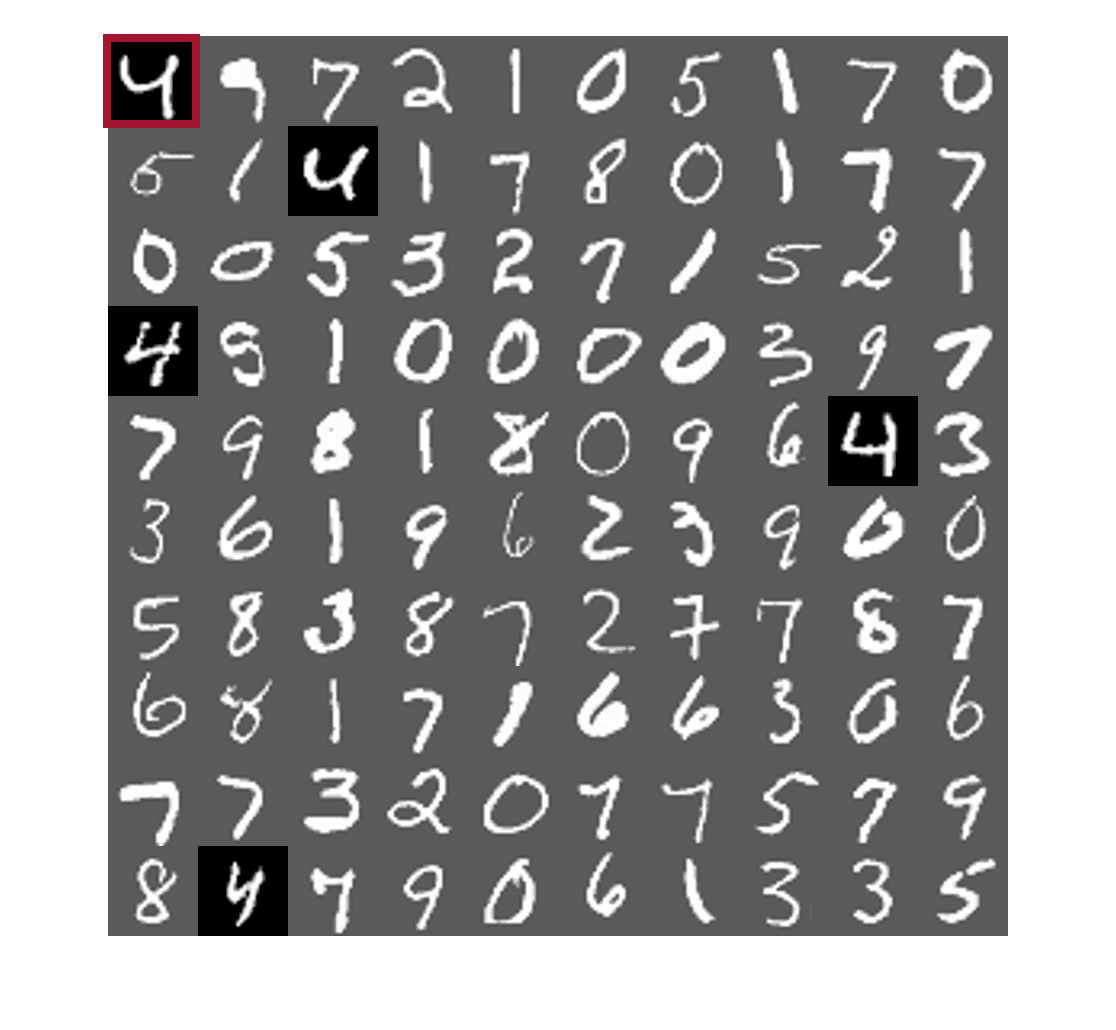}\hspace{-3mm}}
			\\ \vspace{-4mm}
				\subfigure{\footnotesize{\hspace{5mm} (d) \hspace{27mm} (e)\hspace{27mm} (f) }}
\\			
\vspace{-2mm}
		\subfigure{	
			\includegraphics[trim={.2cm .2cm .2cm  .6cm},width=.25\linewidth]{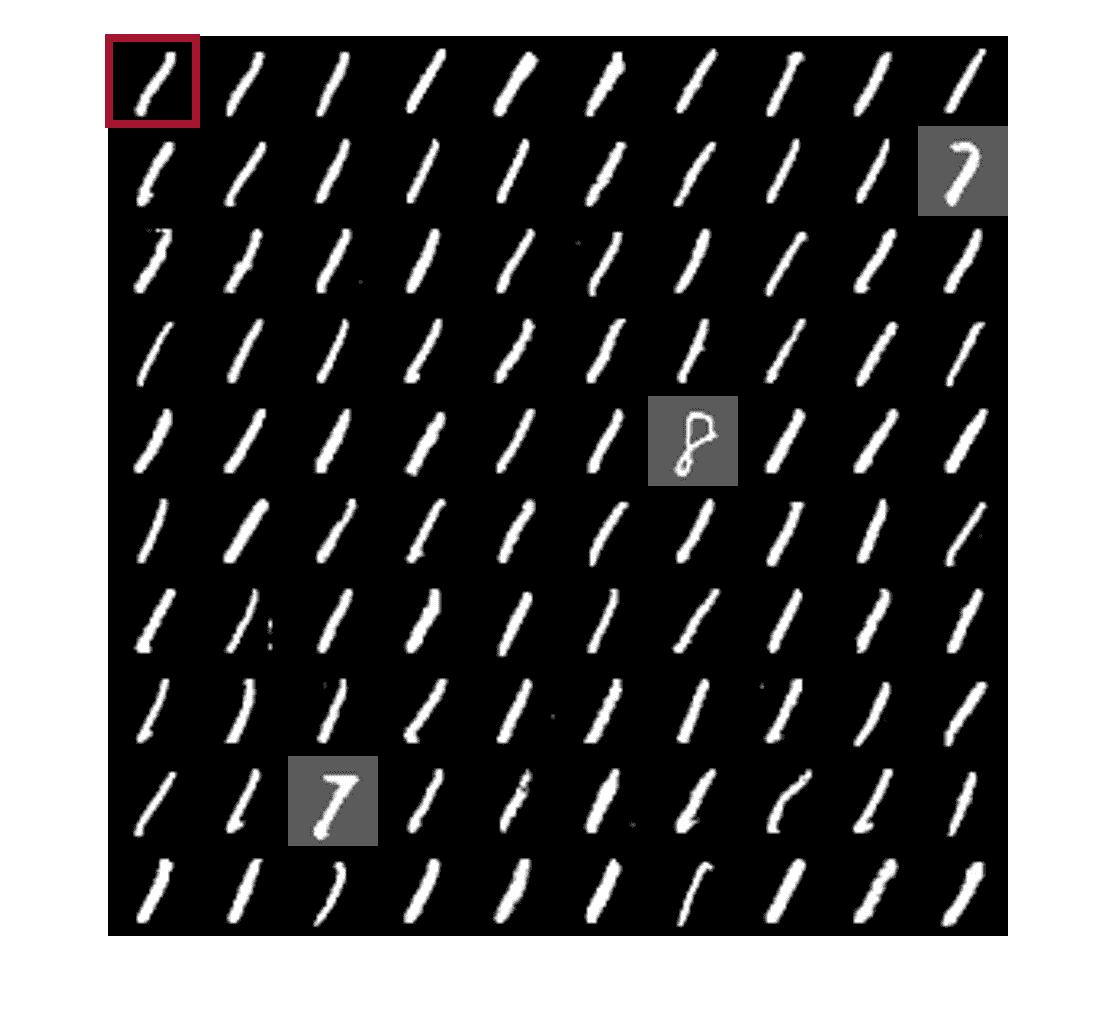} \hspace{-3mm}
			\includegraphics[trim={.2cm .2cm .2cm  .6cm},width=.25\linewidth]{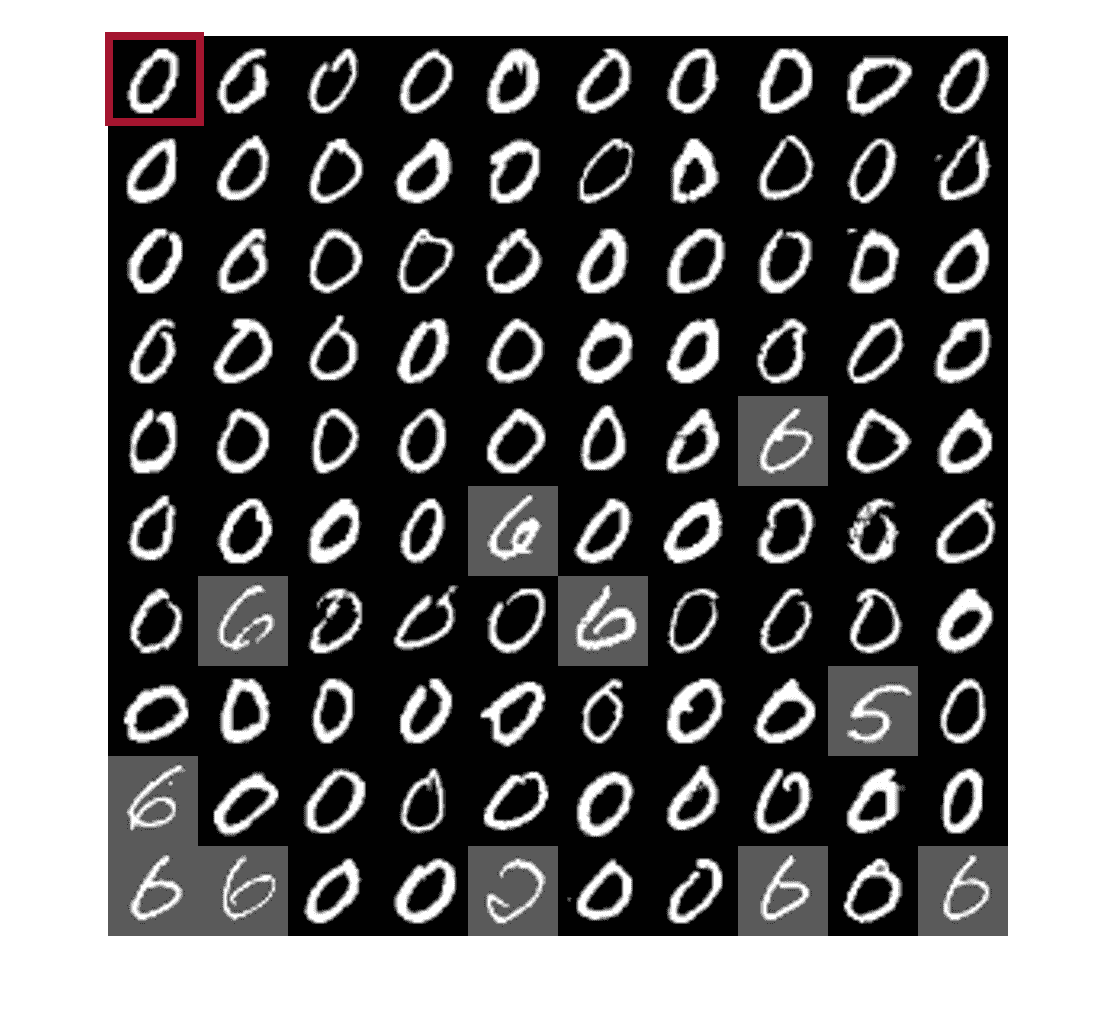} \hspace{-3mm}
			\includegraphics[trim={.2cm .2cm .2cm  .6cm},width=.25\linewidth]{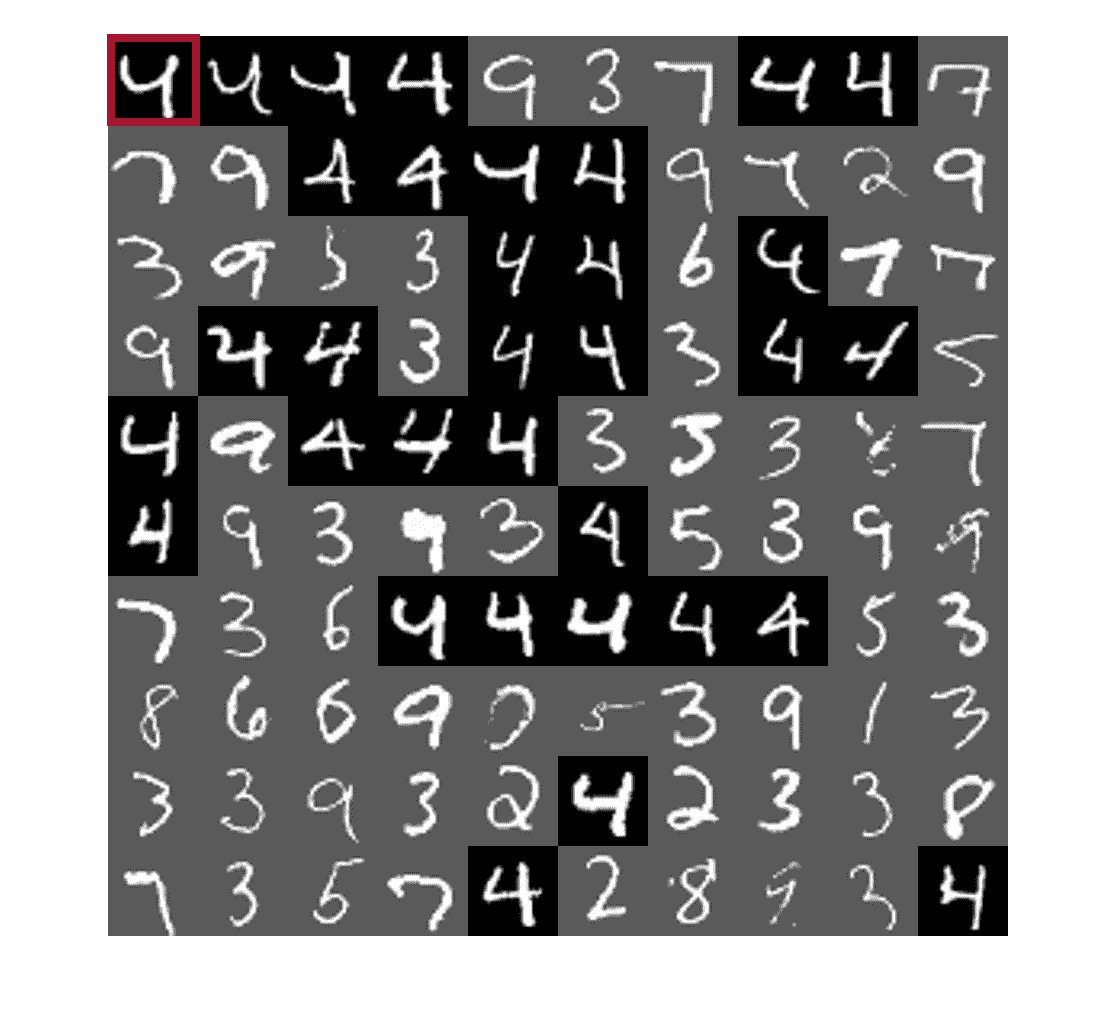}\hspace{-3mm}
	}\\
\vspace{-4mm}
	\subfigure{\footnotesize{\hspace{5mm} (g) \hspace{27mm} (h)\hspace{27mm} (i) }}
\subfigure{	
	\includegraphics[trim={.2cm .2cm .2cm  .6cm},width=.25\linewidth]{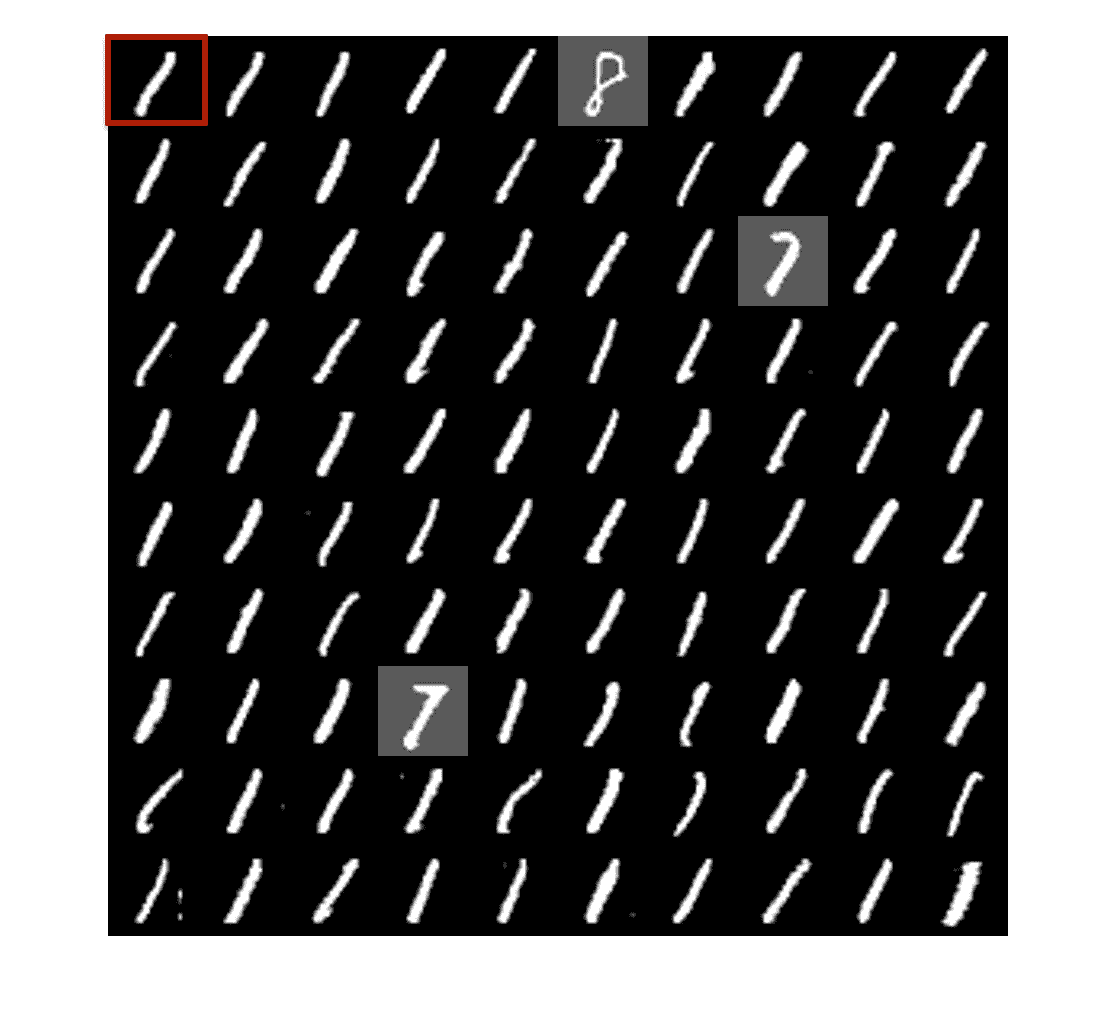} \hspace{-3mm}
	\includegraphics[trim={.2cm .2cm .2cm  .6cm},width=.25\linewidth]{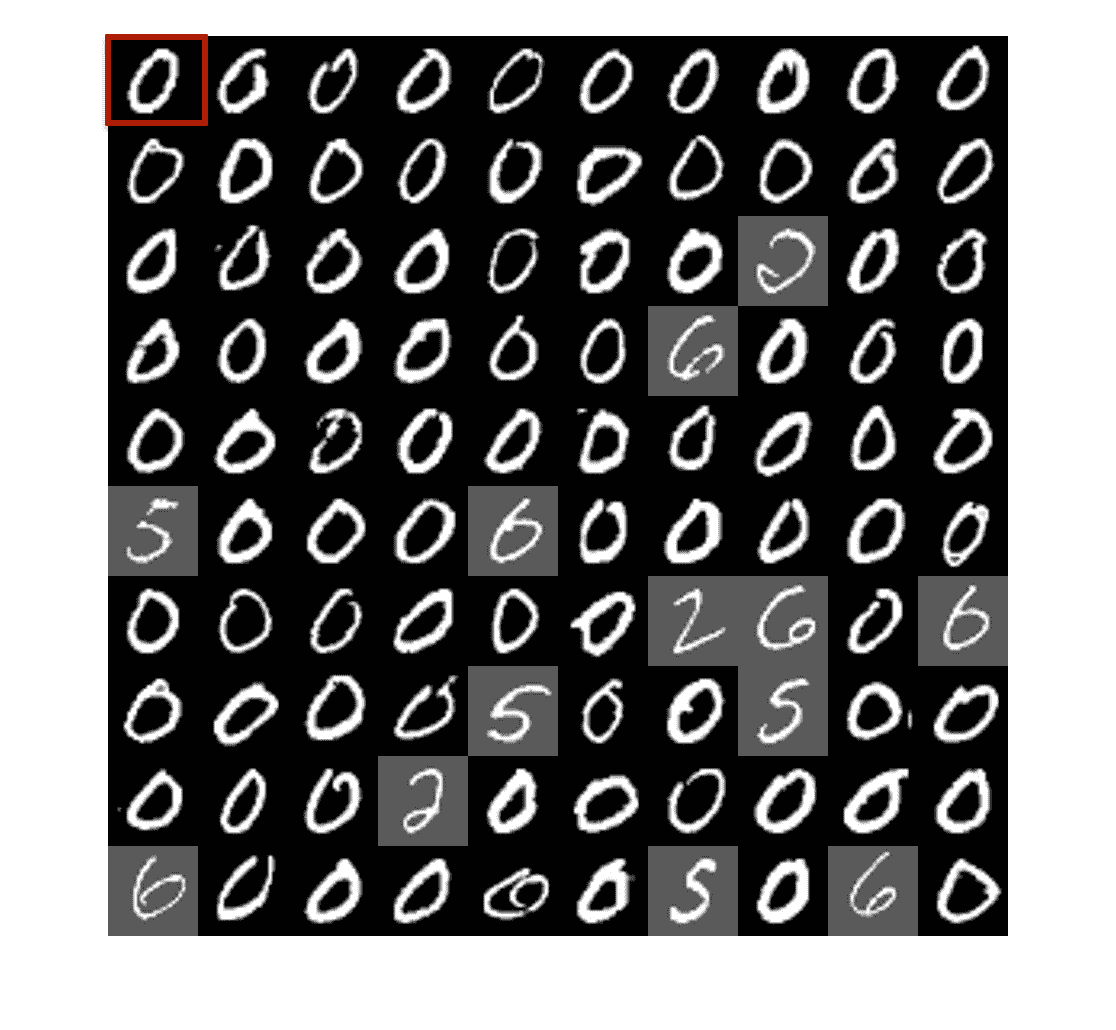} \hspace{-3mm}
	\includegraphics[trim={.2cm .2cm .2cm  .6cm},width=.25\linewidth]{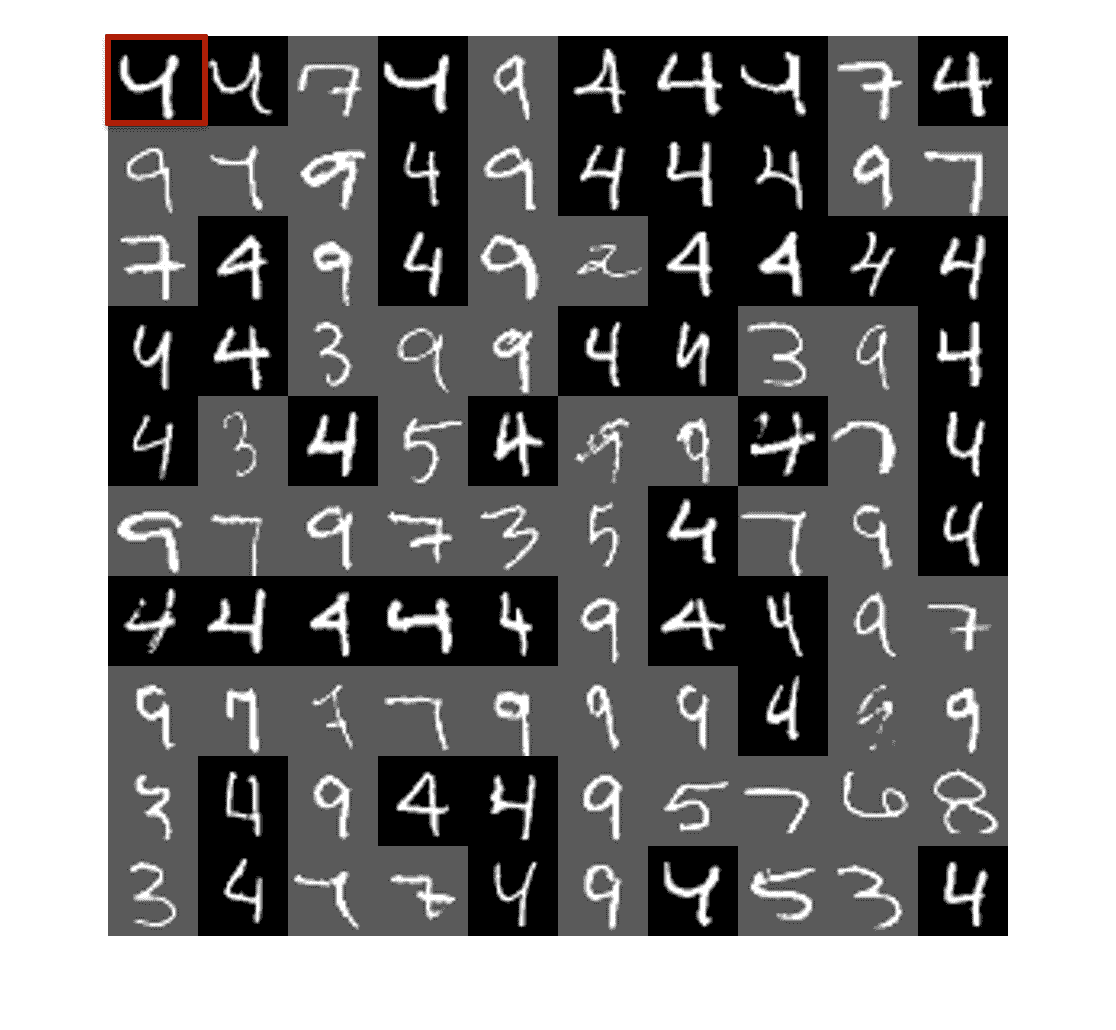}\hspace{-3mm}
} \\ \vspace{-4mm}
	\subfigure{\footnotesize{\hspace{5mm} (j) \hspace{27mm} (k)\hspace{27mm} (l) }}
	\vspace{-4mm}
	\caption{\footnotesize{Image retrieval for three query digits $\{1,0,4\}$ on the MNIST database. The query image is in the top left of the collage (red box), and the incorrect retrieved images are shaded with the gray color. The Panels (a)-(c): Euclidean distance. Panels (d)-(f): $\mathbb{Z}_{2}^{n}$ code with an untrained kernel and $n=4096$ bits. Panels (g)-(i): $\mathbb{Z}_{2}^{n}$ code with a trained kernel with Algoirthm \ref{Algoirthm:1} and $n=4096$. Panels (j)-(l): $\mathbb{Z}_{q}^{n}$ codeword length with $q=10$ and $n=1024$. (best viewed in color)}}
	\label{Fig:8} 
	\end{center}
\end{figure}

The research on computer systems and web search in the mid-nineties was focused on designing ``hash" functions that were sensitive to the topology on the input domain, and preserved distances approximately during hashing. Constructions of such hash functions led to efficient methods to detect similarity of files in distributed file systems and proximity of documents on the web. We describe the problem and results below. Recall that a metric space $(\mathcal{C},d)$ is given by a set $\mathcal{C}$ (here the set of cells) and a distance measure $d:\mathcal{X}\times \mathcal{X}\rightarrow 
\real_{+}$ which satisfies the axioms of being a metric.

\begin{definition}\textsc{(Locally Sensitive Hash Function)}
	Given a metric space $(\mathcal{C},\rho)$ and the set $\mathcal{L}$, a
	family of function $H\subseteq \{h :\mathcal{C}\rightarrow \mathcal{L}\}$ is said to be a basic locality sensitive hash (LSH) family if there exists an increasing invertible function $\alpha:\real_{+}\rightarrow [0,1] $ such that for all $\bm{x},\tilde{\bm{x}}\in \mathcal{X}$, we have
	\begin{align}
\prob_{h\in \mathcal{H}}[h(\bm{x})=h(\tilde{\bm{x}})]\leq \alpha(d(\bm{x},\tilde{\bm{x}})).
	\end{align}
\end{definition}
To contrast the locally sensitive hash functions with the standard hash function families note
that in the latter, the goal is to map a domain $\mathcal{C}$ to a range $\mathcal{L}$ such that the probability of a collision among any pair of elements $\bm{x}\not=\tilde{\bm{x}}\in \mathcal{C}$ is small. In contrast, with LSH families, we wish for the probability of a collision to be small only when the pairwise distance $d(\bm{x},\tilde{\bm{x}})$ is large, and we do want a high probability of collision when $d(\bm{x},\tilde{\bm{x}})$ is small. Indyk and Motwani \cite{indyk1998approximate} showed that such a hashing scheme facilitates the construction of efficient data structures for answering approximate nearest-neighbor queries on the collection of objects.

In the sequel, we describe a family of locally sensitive hash functions due to Raginsky, et al. \cite{raginsky2009locality}, where the measure $d(\bm{x},\tilde{\bm{x}})$ is determined by a kernel function. We draw a random threshold $t\sim \mathrm{Uniform}[-1,1]$ and define the quantizer $Q_{t}(u)=\mathrm{sgn}(t+u)$. 
The following theorem is due to Raginsky, \textit{et al.} \cite{raginsky2009locality}:
\begin{lemma}\textsc{(Kernel Hash Function, Raginsky, et al. \cite{raginsky2009locality})}
	\label{Lemma:KHF}
Consider a translation invariant kernel $K:\mathcal{X}\times \mathcal{X}\rightarrow \real, K(\bm{x},\tilde{\bm{x}})=\phi(\bm{x}-\tilde{\bm{x}})$. Define the following hamming $\mathbb{Z}_{2}^{n}$-code
\begin{align}
F^{n}(\bm{x})\df (F_{t_{1},b_{1},\bm{\omega}_{1}}(\bm{x}),\cdots,F_{\tau_{n},b_{n},\bm{\omega}_{n}}(\bm{x})),
\end{align}
where $(b_{i})_{1\leq i\leq n}\sim_{\text{i.i.d}} \mathrm{Uniform}[-\pi,\pi]$, $(t_{i})_{1\leq i\leq n}\sim_{\text{i.i.d.}} \mathrm{Uniform}[-1,1]$, and $(\bm{\omega}_{i})_{1\leq i\leq n}\sim_{\text{i.i.d.}} \nu$, where $\nu$ is the distribution in the Rahimi and Recht random feature model
\begin{align}
\phi(\bm{x}-\tilde{\bm{x}})=\int_{\real^{d}}\varphi(\bm{x};\bm{\omega})\varphi(\tilde{\bm{x}};\bm{\omega})\mathrm{d}\nu(\bm{\omega}).
\end{align}
Furthermore, $F_{t,b,\bm{\omega}}:\real^{d}\rightarrow \mathbb{Z}_{2}$ is a mapping defined as follows
\begin{align}
\label{Eq:Ftb}
F_{t,b,\bm{\omega}}(\bm{x})\df \dfrac{1}{2}\Big[1+Q_{t}(\cos(\langle \bm{\omega},\bm{x}\rangle+b)) \Big].
\end{align}
Fix $\epsilon,\delta\in (0,1)$. For any finite data set $D = \{\bm{x}_{1},\cdots,\bm{x}_{N}\}\subset \real^{d}$, $F^{n}$ is such that
\begin{align}
h_{K}(\bm{x}_{i}-\bm{x}_{j})-\delta\leq {1\over n}d_{\mathrm{H}}(F^{n}(\bm{x}_{i}),F^{n}(\bm{x}_{j}))\leq h_{K}(\bm{x}_{i}-\bm{x}_{j})+\delta,
\end{align}
where $d_{H}(\cdot,\cdot):\mathbb{Z}^{n}_{2}\times \mathbb{Z}^{n}_{2}\rightarrow \{0,1\}$ is the Hamming distance. Furthermore,
\begin{align}
h_{K}(\bm{x}_{i}-\bm{x}_{j})\df \dfrac{8}{\pi^{2}} \sum_{m=0}^{\infty}\dfrac{1-K(m\bm{x}_{i}-m\tilde{\bm{x}}_{i})}{4m^{2}-1} .
\end{align}
\end{lemma}
The construction of LSH in Eq. \eqref{Eq:Ftb} of Lemma \eqref{Lemma:KHF}  is based on the random feature model of Rahimi and Recht \cite{rahimi2008random,rahimi2009weighted}, and is restricted to the Hamming codes defined on a hyper-cube.  In the sequel, we present an alternative kernel LSH using the embedding of functions in $L^{p}(\Omega,\nu)$ space into $\ell^{p}$-spaces, generalizing the result of Lemma \ref{Lemma:KHF} to any finite integer alphabet $\mathbb{Z}_{q}\df \{0,1,\cdots,q-1\}$.  To describe the result, we define the Lee distance between two code-words $\bm{x}=(x_{1},\cdots,x_{n}),\bm{y}=(y_{1},\cdots,y_{n})\in \mathbb{Z}_{q}^{n}$ of length $n$ as follows
\begin{subequations}
\begin{align}
d_{\mathrm{Lee}}(\bm{x},\bm{y})&\df \sum_{i=1}^{n}\min\{(x_{i}-y_{i})\ \mathrm{mod}\ q, (y_{i}-x_{i})\ \mathrm{mod}\ q\}\\
&=\sum_{i=1}^{n}\min\{|y_{i}-x_{i}|, q-|y_{i}-x_{i}|\}.
\end{align}
\end{subequations}
In the special cases of $q=2$ and $q=3$, the Lee distance corresponds to the Hamming distance as distances are $0$ for two single equal symbols and 1 for two single non-equal symbols. For $q>3$ this is not the case anymore, and the Lee distance can become larger than one.

\begin{lemma}\textsc{(Generalized Kernel LSH)}
	\label{Lemma:Kernel_LSH_via_Embedding}
	Consider a translation invariant kernel $K:\real^{d}\times \real^{d}\rightarrow \real, K(\bm{x},\tilde{\bm{x}})=\psi(\|\bm{x}-\tilde{\bm{x}}\|_{2}^{2})$. Define a random map $G_{t,\bm{\omega}}:\real^{d}\rightarrow \mathbb{Z}_{q}$ through the following function
	\begin{align}
	G_{t,\bm{w}}(\bm{x})\df \left\lceil q^{-1}({\langle \bm{w},\bm{\varphi}_{N}(\bm{x})  \rangle+t}) \right\rceil \mathrm{mod}\ q,
	\end{align}
	where  $\bm{w}\sim \mathsf{N}(0,\bm{I}_{N\times N})$, $t\sim \mathrm{Uniform}[0,q]$. Moreover, 
	$\bm{\varphi}_{N}(\bm{x})\df (\cos(\langle\bm{\omega}_{k},\bm{x} \rangle+b_{k}))_{1\leq k\leq N}$, where $b_{1},\cdots,b_{N}\sim_{\text{i.i.d.}}\nu_{0}\df \mathrm{Uniform}[-\pi,\pi]$, and $\bm{\omega}_{1},\cdots,\bm{\omega}_{N}\sim_{\text{i.i.d.}}\nu$. Then, the probability of collision is bounded as follows
	\begin{align}
	\nonumber
	&\left|\prob\Big[G_{t,\bm{w}}(\bm{x})=G_{t,\bm{\omega}}(\tilde{\bm{x}})\Big]-\Psi_{q}(K(\bm{x},\tilde{\bm{x}}))\right|=\mathcal{O}\left(\dfrac{\ln(N)}{N}\right).
	\end{align}
	where $\Psi_{q}:\real_{+}\mapsto \real_{+}$ is defined
	\begin{align}
	\Psi_{q}(u)\df \int_{0}^{q}{1\over \sqrt{\pi (1-u)}}e^{-{s^{2}\over 4(1-u)}} \left( 1-\dfrac{s}{q}\right)\mathrm{d}s.
	\end{align}
\end{lemma}

The proof of Lemma \ref{Lemma:Kernel_LSH_via_Embedding} is presented in Appendix \ref{Appendix:Proof_of_General_Kernel_LSH}, and uses the property of the stable distributions.

\subsection{Description of the experiment}

We present image retrieval results for 10000 images from MNIST databases [10]. The images are 28 by 28 pixels and are compressed via an auto-encoder with 500 hidden layers, which have proven to be effective at learning useful features. For this experiment, we randomly select 1,000 images to serve as queries, and the remaining 9,000 images make up the ``database''. To distinguish true positives from false positives for the performance evaluation retrieval performance, we select a ``nominal" neighborhood. 

\subsection{Results} 

In Figure \ref{Fig:5}, we illustrate the scatter plots of the hash codes $F^{n}(\bm{x})$ of length $n\in \{128,512,4096\}$ bits. In Figure \ref{Fig:5} (a)-(c), we show the scatter plots for Hamming codes constructed by Lemma \ref{Lemma:KHF} using an untrained Gaussian kernel $K(\bm{x},\widehat{\bm{x}})=\exp(-\|\bm{x}-\widehat{\bm{x}}\|_{2}^{2})$. From Figure \ref{Fig:5} (a)-(c) we observe that as the number of bits of the Hamming code increase, the scatter plots concentrate around a curve that has a flat region. Evidently, this curve deviates from the ideal curve which is a straight line passing through the origin. In particular, in \ref{Fig:5} (a)-(c), as the Euclidean distance changes on the interval $[2,8]$ of $x$-axis, the normalized Hamming distance of the constructed codes remains constant. 

To attain proportionality between the Hamming distance and the Euclidean distance, we train the kernel using Algorithm \ref{Algoirthm:1}. To generate the labels for kernel training, we first perform a $k$-NN on 1000 images from 9,000 images of the data-base, where here $k=10$. Then, we generate the binary class labels by assigning $y=+1$ to the points in a randomly chosen cluster, and $y=-1$ to the remaining data-points in other clusters (\textit{i.e.}, the one-versus-all rule). After training the kernel, the kernel is used to generate the Hash codes using the construction of Lemma \ref{Lemma:KHF}. The resulting scattering plots are depicted in Figure \ref{Fig:5} (d)-(f), and the Hamming distances are more proportional to the Euclidean distance.

Central to our study is the precision-recall curve, which captures the trade-off between precision and recall for different threshold. A high area under the curve represents both high recall and high precision, where high precision relates to a low false positive rate, and high recall relates to a low false negative rate. High scores for both show that the classifier is returning accurate results (high precision), as well as returning a majority of all positive results (high recall). Specifically,
\begin{subequations}
\begin{align}
{\displaystyle {\text{Precision}}\df {\frac {|\{{\text{relevant images}}\}\cap \{{\text{retrieved images}}\}|}{|\{{\text{retrieved images}}\}|}}},\\
{\displaystyle {\text{Recall}}\df {\frac {|\{{\text{relevant images}}\}\cap \{{\text{retrieved images}}\}|}{|\{{\text{relevant images}}\}|}}}.
\end{align}
\end{subequations}

In Figure \ref{Fig:6} we illustrate the precision-recall curves for Hamming codes of different lengths $n\in \{64,512,2048,4096\}$. We also depict the precision-recall curve for the Euclidean distance. We recall the notions of the recall and precision as follows
\
From Figure \ref{Fig:6}(a), we depict the precision-recall curves for Hamming codes constructed from the Gaussian kernel. We observe that even Hamming codes of length $4096$ bits attain the precision-recall curve that is significantly lower than the Euclidean curve. In Figure \ref{Fig:6}(b), we illustrate the precision-recall curves after training kernels using Algorithm \ref{Algoirthm:1}. The curves from Hamming codes are markedly closer to the curve associated with the Euclidean distance. In Figure \ref{Fig:6}(c), we depict the precision-recall curves using the $q$-ary codes in Lemma \ref{Lemma:Kernel_LSH_via_Embedding} for different quantization levels $q\in \{2,3,5,10\}$. Clearly, increasing $q$ leads to a better performance initially while larger values of $q$ has a diminishing return. This can be justified using Figure \ref{Fig:7}, where we plot the normalized histogram of the distance of a sample query point from each point of the MNIST data-base consisting of 9,000 of data-points. The blue color histogram depict the histogram using a Euclidean distance, and the orange color histogram with spikes depicts the histogram of the Lee distance. Figure \ref{Fig:7}(a)-(c) are plotted with the quantization levels $q=2$, $q=20$, and $q=50$, respectively. As the quantization  level increases, the histogram of the Lee distance provides a better approximation for the Euclidean distance histogram. However, even for large $q$, there is a disparity between the two histograms, as the length of the $q$-ary codes are fixed at $n=512$ and is finite. 

In Figure \ref{Fig:8}, we depict the retrieval results for three query digits $\{1,0,4\}$ on MNIST data-set. In Fig. \ref{Fig:8}(a)-(c), we show the performance of retrieval system using a Euclidean distance for searching nearest neighbor. In particular, for each query digit, the retrieval system returns 99 similar images among 9,000 images in the data-base using nearest neighborhood search. The retrieved images for less challenging digits such as 1 and 0 are all correct, while for more challenging query points such as 4, some erroneous images are retrieved. In Fig. \ref{Fig:8}(d)-(f), we depict the retrieved images using a Hamming distance for the nearest neighborhood search, where the Hamming codes for each feature vector are generated via the kernel LSH in Lemma \ref{Lemma:KHF}, where a Gaussian kernel $K(\bm{x},\widehat{\bm{x}})=\exp(-\|\bm{x}-\widehat{\bm{x}}\|_{2}^{2})$ is employed. We observe that compared to the Euclidean distance in Fig. \ref{Fig:8}(a)-(c), the performance of Hamming based nearest neighbor search is significantly less accurate. However, we point out that from a computation complexity perspective, computing the Hamming distance is significantly more efficient than the Euclidean distance which is the main motivation in using LSH for nearest neighborhood search.

In Fig. \ref{Fig:8} (g)-(i), we show the retrieval results after training the kernel using Algorithm \ref{Algoirthm:1}.  Clearly, the retrievals are more accurate on the digits $0$ and $4$. However, the length of the Hamming code required to achieve such performance is quite large $n=4096$. In Fig. \ref{Fig:8} (j)-(l), we show the retrieval performance using the $q$-ary code of length $n=1024$ and the quanitization level of $q=10$. Although the length of the $q$-ary code is significantly smaller than that of the Hamming code, the retrieval error is slightly less. However, the computational cost is slightly higher than that of the Hamming code.

In the future, we will test our method on data-sets consisting of millions of data points. At present, our promising initial results on MNIST data-set, combined with our comprehensive theoretical analysis, convincingly demonstrate the potential usefulness of our kernel learning scheme for large-scale indexing and search applications.

\section{Classification on Benchmark Data-Sets}

We now apply our kernel learning method for classification and regression tasks on real-world data-sets. 

\subsection{Data-Sets Description}

We apply our kernel learning approach to classification and regression tasks of real-world data-sets. In Table \ref{Table:Benchmark_Data_Descrpition}, we provide the characteristics of each data-set. All of these datasets are publicly available at UCI repository.\footnote{\url{https://archive.ics.uci.edu/ml/index.php}}
\subsubsection{Online news popularity}

This data-set summarizes a heterogeneous set of features about articles published by Mashable in a period of two years. The goal is to predict the number of shares in social networks (popularity).

\subsubsection{Buzz in social media dataset}

This data-set contains examples of buzz events from two different social networks: Twitter, and Tom's Hardware, a forum network focusing on new technology with more conservative dynamics

\subsubsection{Adult}

Adult data-set contains the census information of individuals including education, gender, and capital gain. The assigned classification task is to predict whether a person earns over 50K annually. The train and test sets are two separated files consisting of roughly 32000 and 16000 samples respectively.

\subsubsection{Epileptic Seizure Detection}

The epileptic seizure detection data-set consists of a recording of brain activity for 23.6 seconds. The corresponding time-series is sampled into 4097 data points. Each data point is the value of the EEG recording at a different point in time. So we have total 500 individuals with each has 4097 data points for 23.5 seconds. The 4097 data points are then divided and shuffled every into 23 segments, each segment contains 178 data points for 1 second, and each data point is the value of the EEG recording at a different point in time.

\subsection{Quantitative Comparison}

\begin{figure}[!t]
	\begin{center}
		\subfigure{	\hspace{-20mm}
			\includegraphics[trim={.2cm .2cm .2cm  .6cm},width=.3\linewidth]{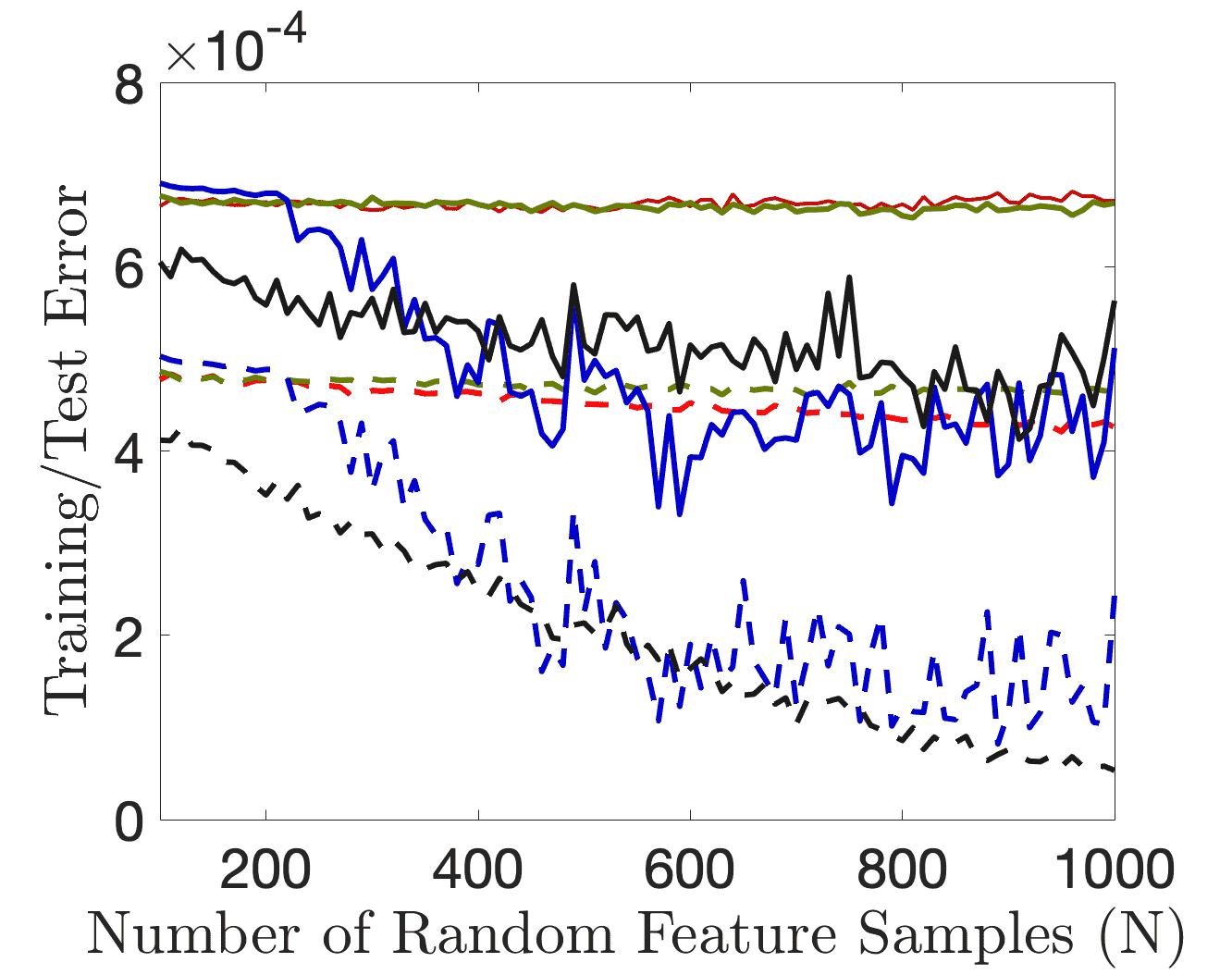} \hspace{-2mm}
			\includegraphics[trim={.2cm .2cm .2cm  .6cm},width=.3\linewidth]{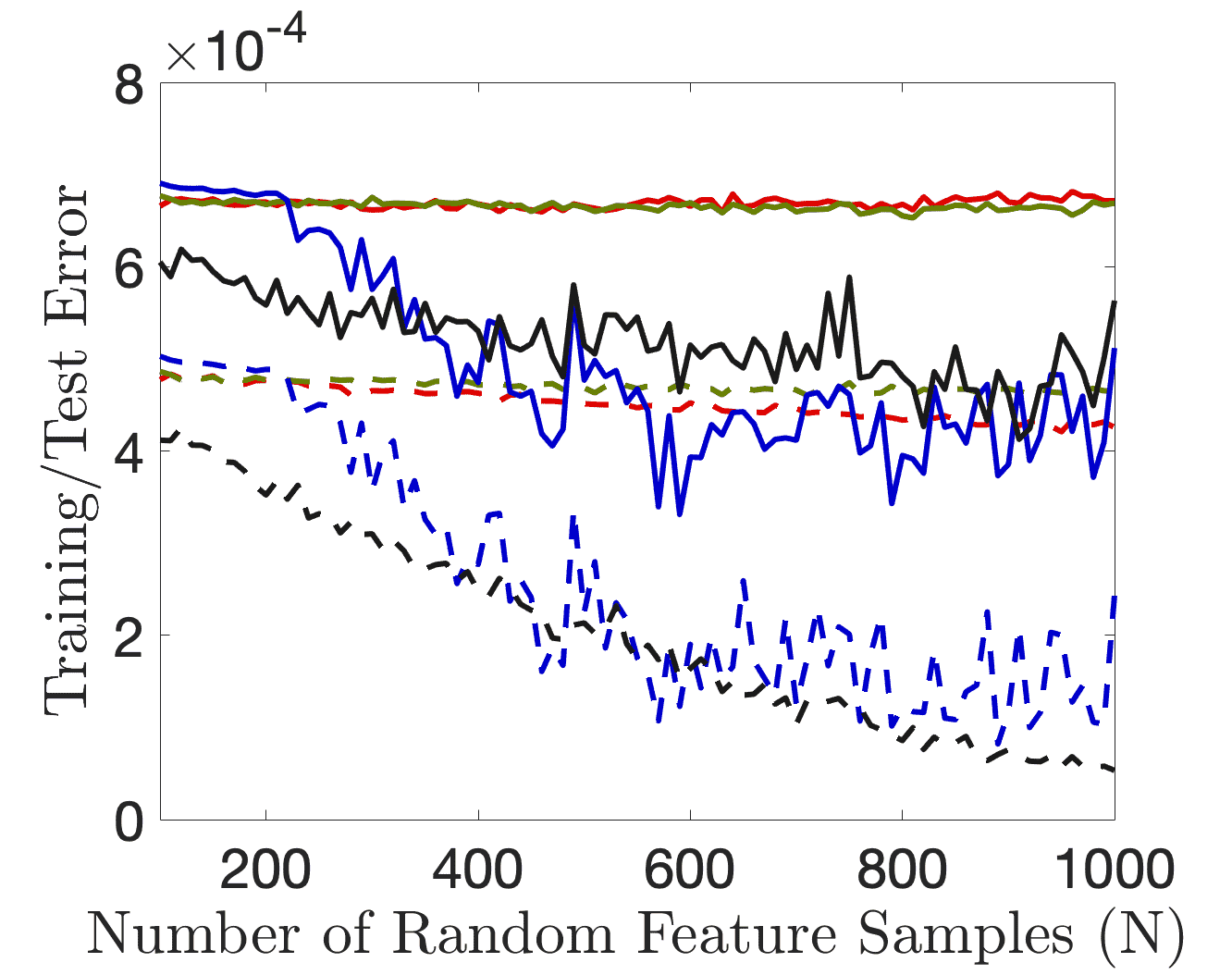} \hspace{-2mm}
			\includegraphics[trim={.2cm .2cm .2cm  .6cm},width=.3\linewidth]{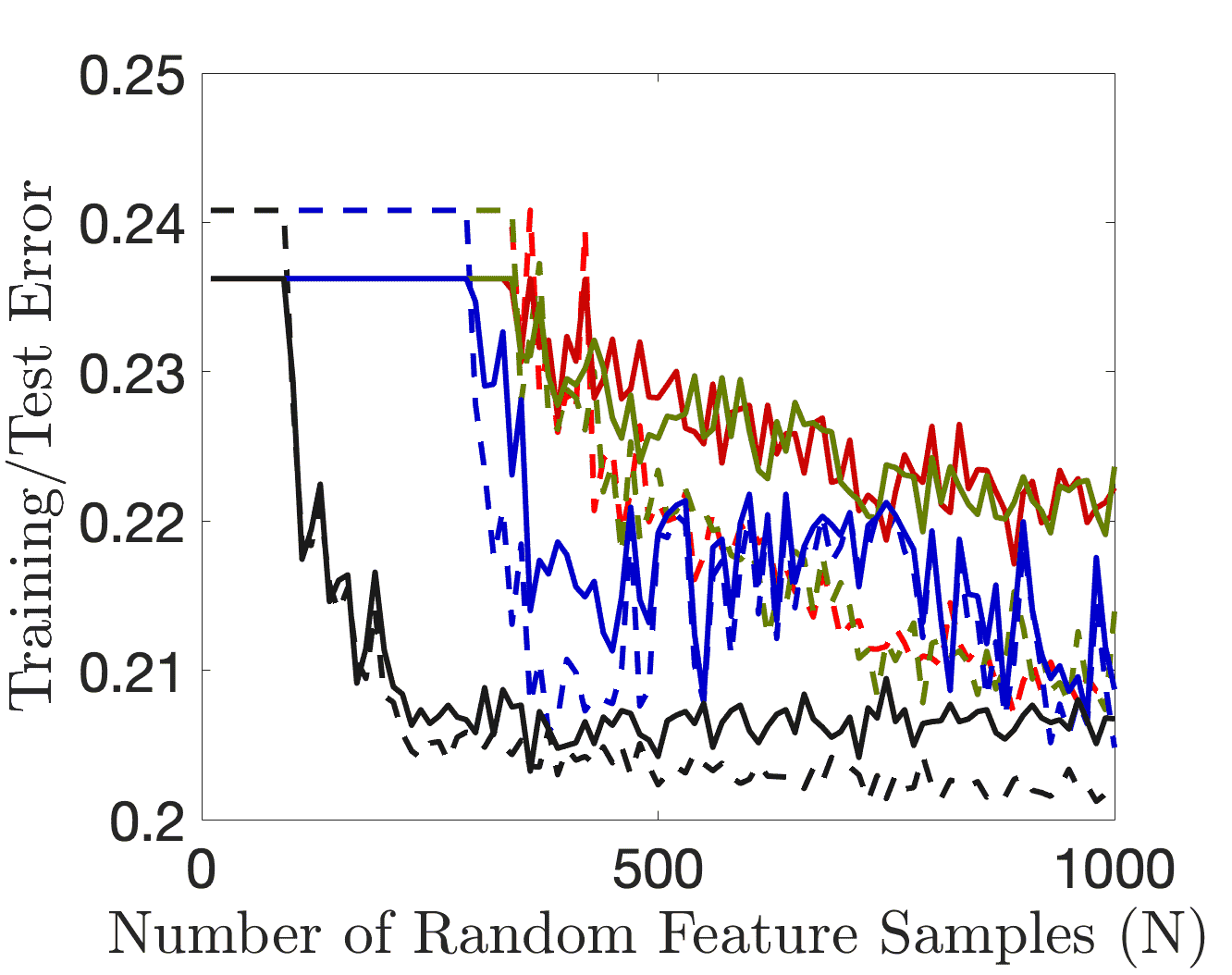}  \hspace{-2mm}
			\includegraphics[width=.3\linewidth]{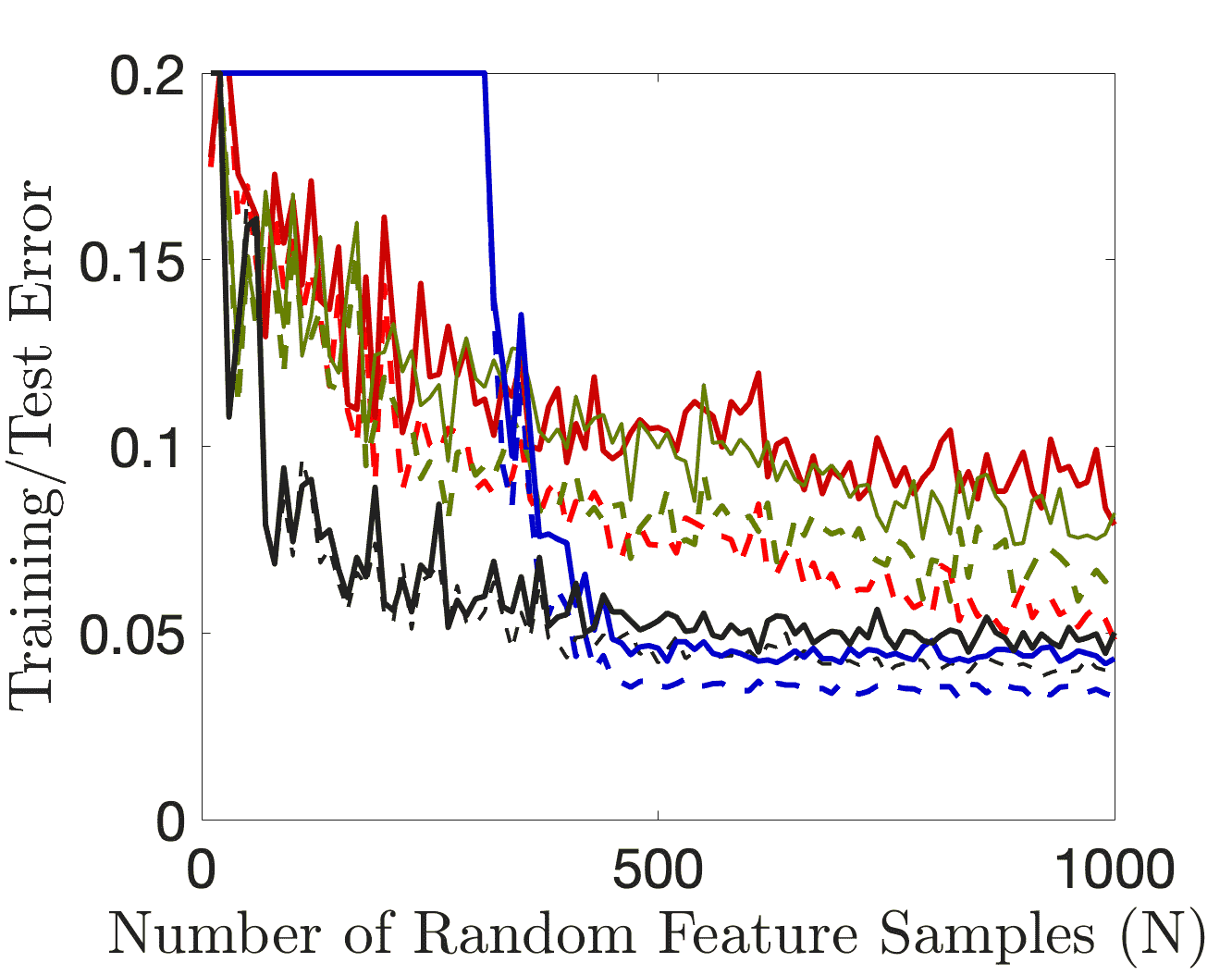}}
		\\ 
		\subfigure{\footnotesize{(a)\hspace{40mm} (b)  \hspace{40mm} (c) \hspace{40mm} (d)}} 
		\vspace{-3mm}
		\caption{\footnotesize{Training and test errors of Algorithm \ref{Algoirthm:1} (black line), the SGD optimization (blue lines), the importance sampling (green lines), and a Gaussian kernel with the \textit{optimized} bandwidth (red lines) for classification and linear regression tasks. The training and test errors are depicted with dashed and solid lines, respectively.  Panel (a): \textsc{Buzz}, Panel (b): \textsc{Online news popularity}, Panel (c): \textsc{Adult}, Panel (d): \textsc{Seizure}.} (Best viewed in color)}
		\label{Fig:Seizure} 
	\end{center}
	\begin{center}
		\subfigure{
			\hspace{-20mm}	\includegraphics[trim={.2cm .2cm .2cm  .6cm},width=.3\linewidth]{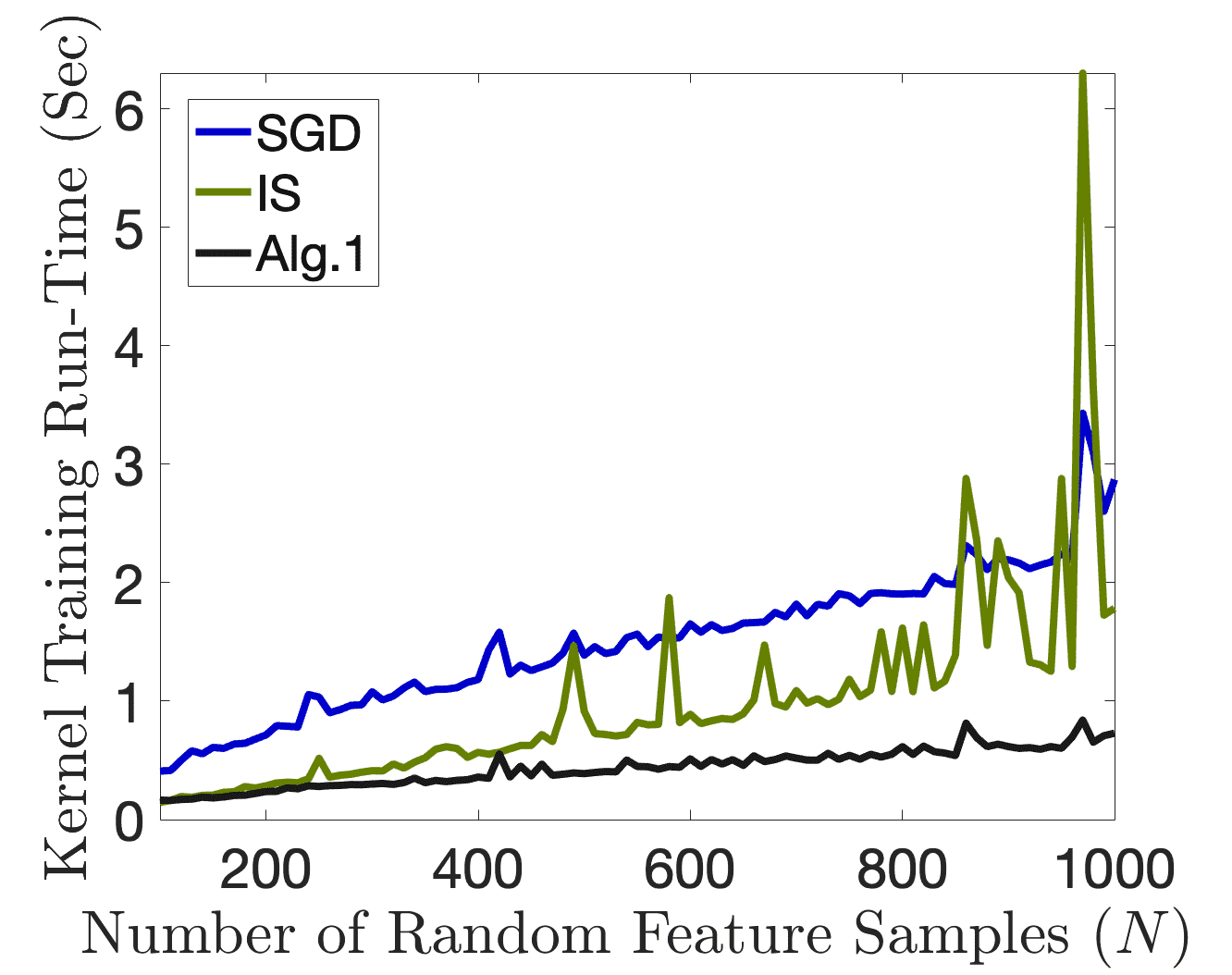} \hspace{-2mm}
			\includegraphics[trim={.2cm .2cm .2cm  .6cm},width=.3\linewidth]{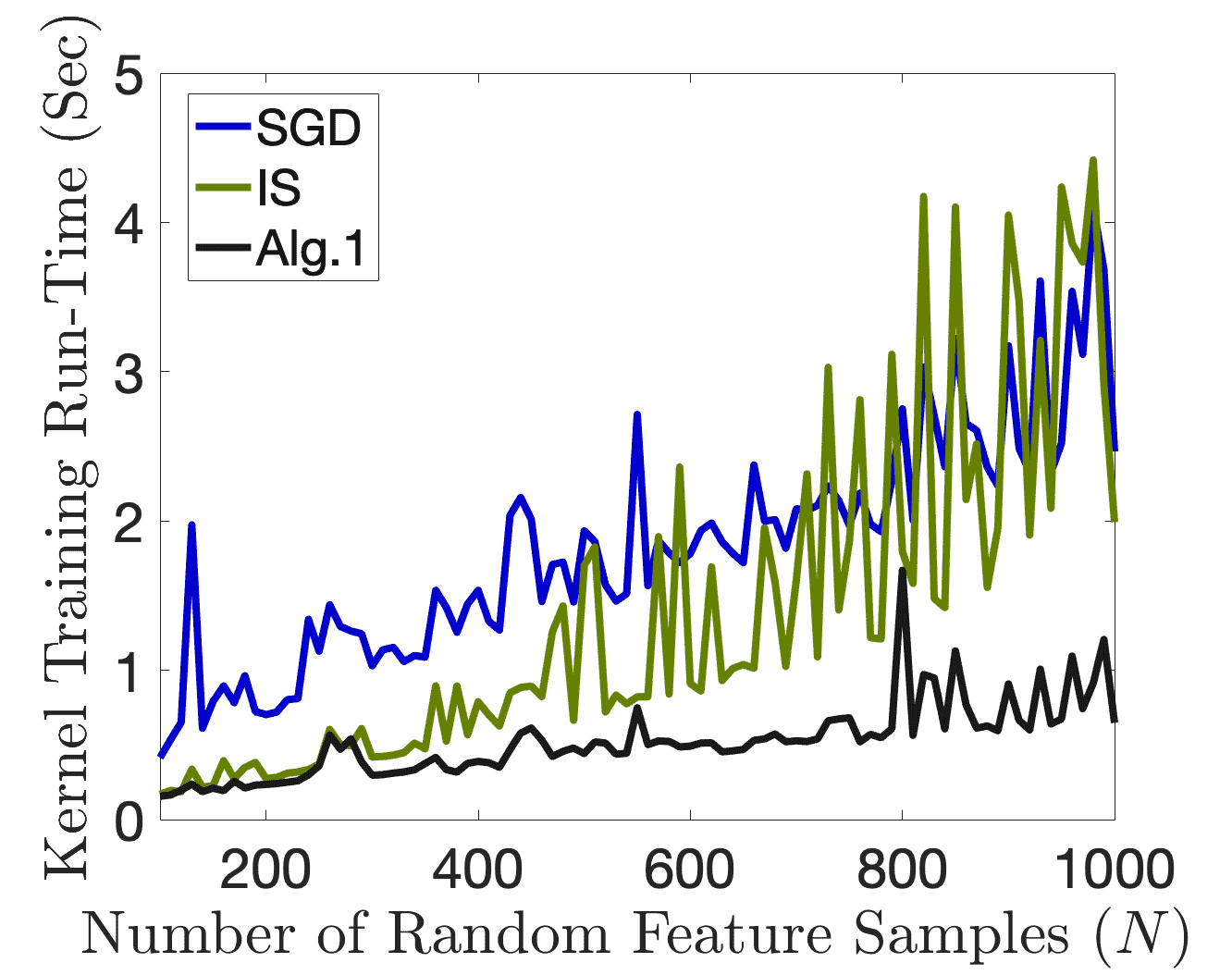} 
			\includegraphics[trim={.2cm .2cm .2cm  .6cm},width=.3\linewidth]{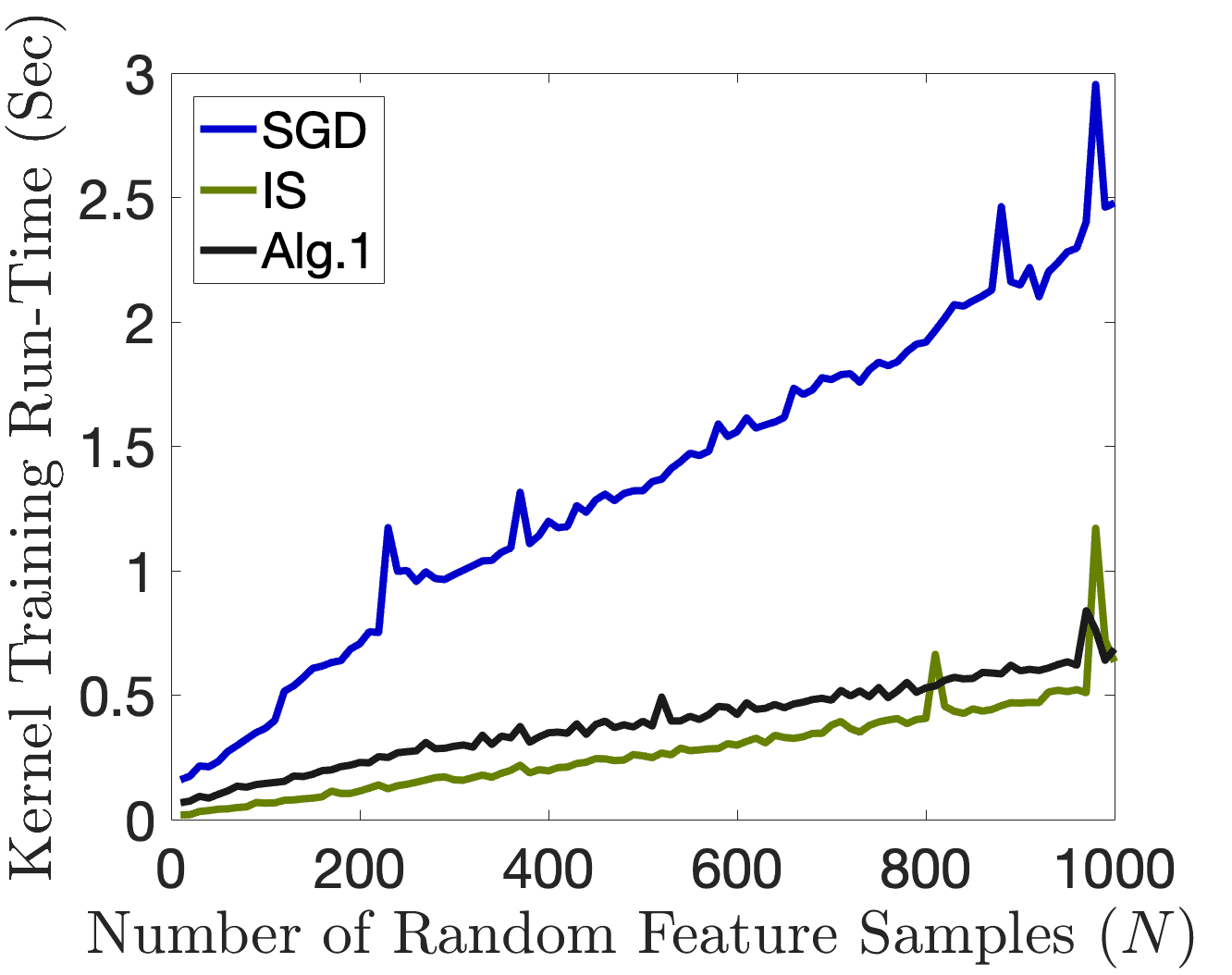}  \hspace{-2mm}
			\includegraphics[trim={.2cm .2cm .2cm  .6cm},width=.3\linewidth]{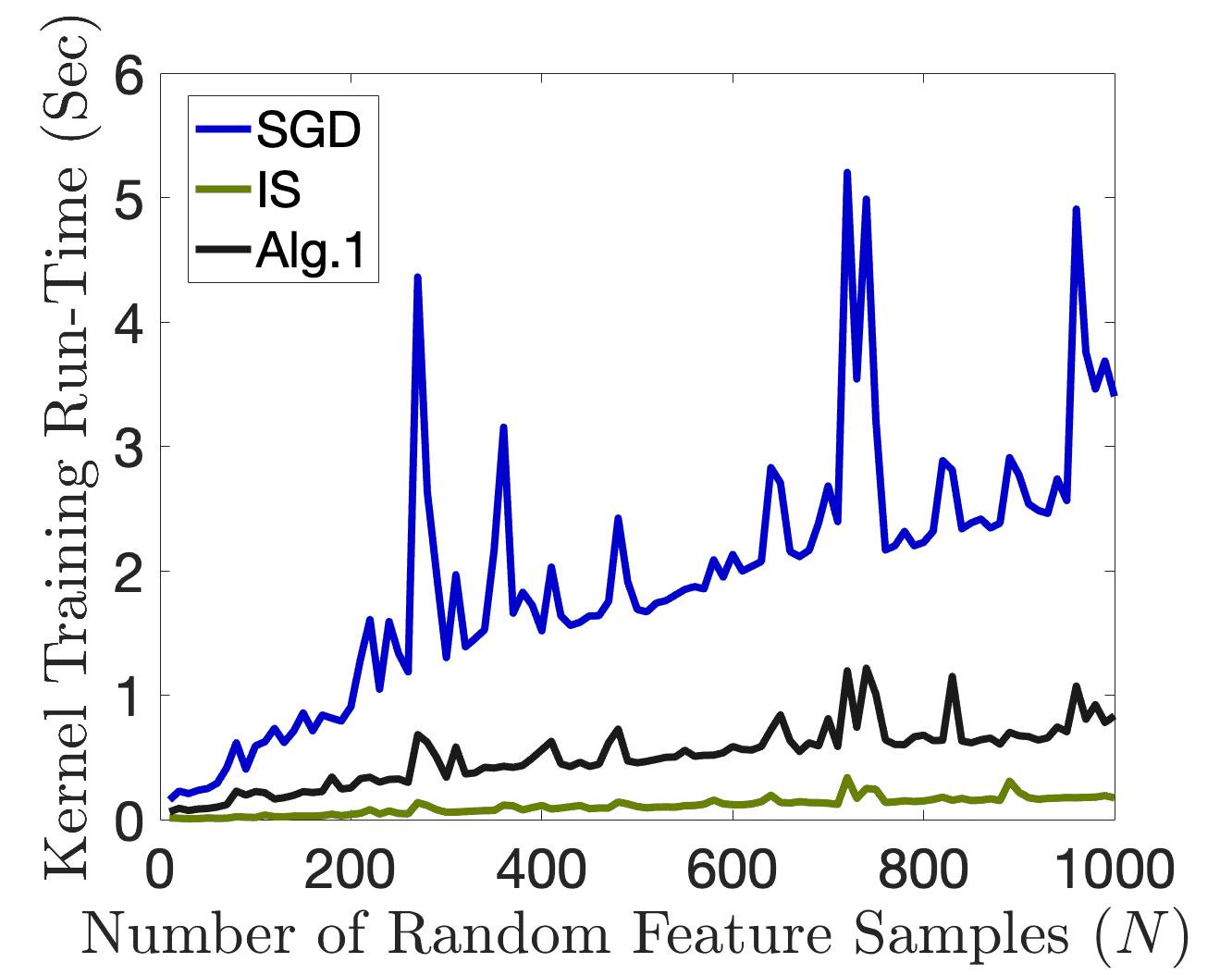}  }
		\\ 
		\subfigure{\footnotesize{(a)\hspace{40mm} (b)  \hspace{40mm} (c) \hspace{40mm} (d)}} 
		\vspace{-3mm}
		\caption{\footnotesize{The run-times of kernel optimization algorithms using Algorithm \ref{Algoirthm:1}, the SGD optimization, and Importance Sampling (IS).  Panel (a): \textsc{Buzz}, Panel (b): \textsc{Online news popularity}, Panel (c): \textsc{Adult}, Panel (d): \textsc{Seizure} (Best viewed in color).}}
		\label{Fig:Seizure_1} 
	\end{center}
\end{figure}

In Figure \ref{Fig:Seizure}, we present the training and test results for regression and classification tasks on benchmark data-sets, using top $d=35$ features from each data-set and for different number of random feature samples $N$. In all the experiments, the SGD method provides a better accuracy in both the training and test phases. Nevertheless, in the case of seizure detection, we observe that for a small number of random feature samples, the importance sampling and Gaussian kernel with $k$NN for bandwidth outperforms SGD. 

In Figure \ref{Fig:Seizure_1}, we illustrate the time consumed for training the kernel using SGD and importance sampling methods versus the number of random features $N$. For the linear regression on \textsc{Buzz}, and \textsc{Online news popularity} the difference between the run-times are negligible. However, for classification tasks using \textsc{Adult} and \textsc{Seizure}, the difference in run-times are more pronounced.

Those visual observations are also repeated in a tabular form in Table \ref{Table:long_table}, where the training and test errors as well algorithmic efficiencies (run-times) are presented for $N=100$, $N=1000$, $N=2000$, and $N=3000$ number of features.

\section{Conclusions}

We have proposed and analyzed a distributionally robust optimization method to learn the distribution in the Sch\"{o}nberg integral representation of radial basis functions. In particular, we analyzed a projected particle Langevin dynamics to optimize the samples of the distribution in the integral representation of the radial kernels. We established theoretical performance guarantees for the proposed distributional optimization procedure. Specifically, we derived a non-asymptotic bound for the consistency of the finite sample approximations. Furthermore, using a mean-field analysis, we derived the scaling limits of a particle stochastic optimization method. We showed that in the scaling limits, the projected particle Langevin optimization converges to a reflected diffusion-drift process. We then derived a partial differential equation, describing the evolution of the law of the underlying reflected process. We evaluated the performance of the proposed kernel learning approach for classification on benchmark data-sets. We also used the kernel learning approach in conjunction with a $k$-mean clustering method to train the kernel in the kernel locally sensitive hash functions.

\begin{center}
	\centering
	\begin{table}[t]
		\begin{tabular}{||c|| c| c| c|c||}
			\hline
			\rowcolor{Gray}
		\color{white}{	Data-set}      & \color{white}{\text{Task}}   & \color{white}{$d$}  &  \color{white}{$n_{\text{training}}$}   &\color{white}{$n_{\text{test}}$}  \\
			\hline
			\hline
			Buzz& Regression &77  &93800 &46200 \\
			\hline
			Online news popularity& Regression & 58  &26561  &13083 \\ 
			\hline
			Adult & Classification  & 122 &32561  &16281  \\
			\hline
			Seizure &Classification  & 178  & 8625  &2875 \\
			\hline 
		\end{tabular}
		\caption{\footnotesize{Description of the benchmark data-sets used in this paper}.}
		\label{Table:Benchmark_Data_Descrpition}
	\end{table}
\end{center}

\begin{center}
\footnotesize{\centering
	\begin{table}[!htbp]
		\vspace{-20mm}
		\scriptsize{ \begin{tabular}{||c|| c| c| c|c||}
				\hline	
				\multicolumn{5}{c}{{\textbf{Buzz}}}\\ \hline
				\rowcolor{Gray}
				\hline
				Algorithm 1 & $N$=100 & $N$=1000 & $N$=2000 &$N$=3000 \\ [0.5ex] 
				\hline\hline
				Training Error   &0.448e$-3$   &0.228e$-3$  &   0.781e$-4$    &0.288e$-3$  \\ 
				\hline
				Test Error   & 0.635e$-3$    &0.554e$-3$    & 0.512e$-3$    &0.397e$-3$ \\
				\hline
				Run Time (sec)   &  0.160&  0.748 &1.454 &1.867\\
				\hline
				\hline
			
				\rowcolor{Gray}
				\hline
				Particle SGD & $N$=100 & $N$=1000 & $N$=2000 &$N$=3000 \\ [0.5ex] 
				\hline\hline
				Training Error   &0.501e$-3$   &0.190e$-3$   & 0.180e$-3$    &0.151e$-3$  \\ 
				\hline
				Test Error   & 0.501e$-3$    &0.190e$-3$    &0.180e$-3$    &0.151e$-3$ \\
				\hline
				Run Time (sec)   &  0.316&  1.971 &4.617 &14.621\\
				\hline
				\hline
				\rowcolor{Gray}
				Importance Sampling & $N$=100 & $N$=1000 & $N$=2000 &$N$=3000\\ [0.5ex] 
				\hline\hline
				Training Error  & 0.486e$-3$    &0.466e$-3$    &0.460e$-3$    &0.455e$-3$ \\ 
				\hline
				Test Error &  0.677e$-3$   & 0.661e$-3$    &0.662e$-3$    &0.661e$-3$ \\
				\hline
				Run Time (sec)     &  0.188  &  1.527 & 3.592 &4.955\\
				\hline
				\hline
				\rowcolor{Gray}
				Gaussian Kernel & $N$=100 & $N$=1000 & $N$=2000 &$N$=3000\\ [0.5ex] 
				\hline\hline
				Training Error   &0.484e$-3$    &0.429e$-3$    &0.379e$-3$   & 0.327e$-3$ \\ 
				\hline
				Test Error  &  0.673e$-3$    &0.673$-3$    &0.709e$-3$    &0.744e$-3$
				\\				\hline

				\multicolumn{5}{c}{{\textbf{Online news popularity}}}\\ 
					\rowcolor{Gray}
				\hline
				Algorithm 1 & $N$=100 & $N$=1000 & $N$=2000 &$N$=3000 \\ [0.5ex] 
				\hline\hline
				Training Error   &0.463e$-3$   &0.239e$-3$   & 0.676e$-3$    &0.288e$-4$  \\ 
				\hline
				Test Error   & 0.654e$-3$    &0.537e$-3$    &0.474e$-3$    &3.974e$-4$ \\
				\hline
				Run Time (sec)   & 0.165&  0.686 &1.264 &1.736\\
				\hline
				\hline
				\hline
				\rowcolor{Gray}	
				Particle SGD & $N$=100 & $N$=1000 & $N$=2000 &$N$=3000 \\ [0.5ex] 
				\hline\hline
				Training Error   &0.501e$-3$    &0.190e$-3$   & 0.180e$-3$  &  0.151e$-3$  \\ 
				\hline
				Test Error   & 0.687e$-3$    &0.524e$-3$    &0.567e$-3$   & 0.588e$-3$ \\
				\hline
				Run Time (sec)   &  0.309    &2.071    &3.152    &5.054\\
				\hline
				\hline
				\rowcolor{Gray}
				Importance Sampling & $N$=100 & $N$=1000 & $N$=2000 &$N$=3000\\ [0.5ex] 
				\hline\hline
				Training Error  & 0.486e$-3$    &0.466e$-3$    &0.460e$-3$    &0.455e$-3$ \\ 
				\hline
				Test Error &  0.677e$-3$     &0.662e$-3$     &0.662e$-3$     &0.661e$-3$  \\
				\hline
				Run Time (sec)     & 0.222    &1.640    &4.254   &10.782\\
				\hline
				\hline
				\rowcolor{Gray}
				Gaussian Kernel & $N$=100 & $N$=1000 & $N$=2000 &$N$=3000\\ [0.5ex] 
				\hline\hline
				Training Error   &0.484e$-3$    &0.429e$-3$    &0.379e$-3$    &0.327e$-3$ \\ 
				\hline
				Test Error  & 0.673e$-3$    &0.673e$-3$   &0.709e$-3$    &0.744e$-3$	
				\\				\hline
				\multicolumn{5}{c}{{\textbf{Adult}}}\\ \hline
					\hline
				\rowcolor{Gray}	
				Algorithm 1 & $N$=100 & $N$=1000 & $N$=2000 &$N$=3000 \\ [0.5ex] 
				\hline\hline
				Training Error   &0.221   &0.201   &  0.199    &0.199  \\ 
				\hline
				Test Error   & 0.219    &0.207   & 0.205    &0.205 \\
				\hline
				Run Time (sec)   &  0.155&  0.645 &1.426 &1.802\\
				\hline
				\hline
				\rowcolor{Gray}	
				Particle SGD & $N$=100 & $N$=1000 & $N$=2000 &$N$=3000 \\ [0.5ex] 
				\hline\hline
				Training Error   &
				0.240    &0.204    &0.199    &0.197 \\ 
				\hline
				Test Error   & 0.236    &0.209   &0.208   &0.210
				\\
				\hline
				Run Time (sec)   &  0.274    &1.784    &9.272   &11.295\\
				\hline
				\hline
				\rowcolor{Gray}
				Importance Sampling & $N$=100 & $N$=1000 & $N$=2000 &$N$=3000\\ [0.5ex] 
				\hline\hline
				Training Error  &0.240   & 0.207    &0.214    &0.220
				\\ 
				\hline
				Test Error &  0.236    &0.219    &0.222   &0.223 \\
				\hline
				Run Time (sec)     &  0.058    &0.511    &2.661    &8.068\\
				\hline
				\hline
				\rowcolor{Gray}
				Gaussian Kernel & $N$=100 & $N$=1000 & $N$=2000 &$N$=3000\\ [0.5ex] 
				\hline\hline
				Training Error   & 0.240    &0.206    &0.201    &0.197 \\ 
				\hline
				Test Error  &  0.236    &0.222    &0.216    &0.214
				\\				\hline
				\multicolumn{5}{c}{{\textbf{Seizure}}}\\ \hline
				\rowcolor{Gray}	
				Algorithm 1 & $N$=100 & $N$=1000 & $N$=2000 &$N$=3000 \\ [0.5ex] 
				\hline\hline
				Training Error   &0.0814  &0.0408   & 0.0364    & 0.0343  \\ 
				\hline
				Test Error   & 0.0873    & 0.0480    &0.0449    &0.0470 \\
				\hline
				Run Time (sec)   & 0.134  &0.701 &1.229 &1.882\\
				\hline
				\hline
				\rowcolor{Gray}	
				Particle SGD & $N$=100 & $N$=1000 & $N$=2000 &$N$=3000 \\ [0.5ex] 
				\hline\hline
				Training Error   &0.200    &0.034   & 0.031    &0.032  \\ 
				\hline
				Test Error   & 0.200    &0.043    &0.043    &0.043 \\
				\hline
				Run Time (sec)   &  0.350    &1.747    &3.833    &5.222\\
				\hline
				\hline
				\rowcolor{Gray}
				Importance Sampling & $N$=100 & $N$=1000 & $N$=2000 &$N$=3000\\ [0.5ex] 
				\hline\hline
				Training Error  & 0.135    &0.078    &0.055    &0.051
				\\ 
				\hline
				Test Error &  0.144    &0.099    &0.063    &0.056 \\
				\hline
				Run Time (sec)     &  0.038    &0.138    &0.352    &0.404\\
				\hline
				\hline
				\rowcolor{Gray}
				Gaussian Kernel & $N$=100 & $N$=1000 & $N$=2000 &$N$=3000\\ [0.5ex] 
				\hline\hline
				Training Error   & 0.146    &0.056    &0.033    &0.030 \\ 
				\hline
				Test Error  & 0.154    &0.093   & 0.076    &0.069 \\
				\hline
		\end{tabular}}\normalsize
		\caption{\footnotesize{Performance comparison between kernel learning method (Algorithm \ref{Algoirthm:1}), the particle SGD \cite{khuzani2019mean}, the importance sampling \cite{sinha2016learning}, and a regular Gaussian kernel in conjunction with the kernel SVMs for different random feature samples $N$ on benchmark data-sets.}}
		\label{Table:long_table}
	\end{table}
}
\end{center}

\appendix

\section{Proof of Main Results}
\label{Appendix:Formally}
\textbf{Notations}.  We define the following notion of distances between two measures $\mu,\nu\in \mathcal{M}(\mathcal{X})$ on the metric space $\mathcal{X}$:
\begin{itemize}
	\item \textit{$f$-divergence}:
	\begin{align}
	D_{f}(\mu||\nu)=\int_{\mathcal{X}} f\left(\dfrac{\mathrm{d}\mu}{\mathrm{d}\nu} \right)\mathrm{d}\nu.
	\end{align}
	In the case of $f(x)=x\log x$, the corresponding distance is the Kullback-Leibler divergence which we denote by $D_{\mathrm{KL}}(\cdot||\cdot)$.
	\item \textit{Wasserstein distance}: 
	\begin{align}
W_{p}(\mu,\nu)=\inf_{\pi\in \Pi(\mu,\nu)}\left(\int_{\mathcal{X}\times \mathcal{X}}d^{p}(\bm{x},\tilde{\bm{x}})\pi(\bm{x},\tilde{\bm{x}}) \right)^{1\over p},
\end{align}
where $\Pi(\mu,\nu)$ is the set of all measure $\mu$ and $\nu$. For $p=1$, the Kantorovich-Rubinstein duality yields 
\begin{align}
W_{1}(\mu,\nu)=\sup_{f\in \mathcal{F}_{\mathrm{L}}} \int_{\mathcal{X}}f(\bm{x})(\mathrm{d}\mu(\bm{x})-\mathrm{d}\nu(\bm{x})),
\end{align}
where $\mathcal{F}_{\mathrm{L}}\df \{f:\mathcal{X}\rightarrow \real:\|f\|_{\mathrm{Lip}}\leq 1, f\in C^{0}(\mathcal{X})\}$, where $C^{0}(\mathcal{X})$ is the class of continuous functions.

	\item \textit{Bounded Lipschitz metric:} 
	\begin{align}
	D_{\mathrm{BL}}(\mu,\nu)=\sup_{f\in \mathcal{F}_{\mathrm{BL}}}\left|\int_{\mathcal{X}}f(\bm{x})\mathrm{d}\mu(\bm{x})-\int_{\mathcal{X}}f(\bm{x})\mathrm{d}\nu(\bm{x}) \right|,
	\end{align}
	where $\mathcal{F}_{\mathrm{BL}}\df \{f:\real^{d}\rightarrow \real: \|f\|_{\mathrm{Lip}}\leq 1, \|f\|_{\infty}\leq 1\}$.
   	
\end{itemize}

We denote vectors by lower case bold letters, \textit{e.g.} $\bm{x}=(x_{1},\cdots,x_{n})\in \real^{n}$, and matrices by the upper case bold letters, \textit{e.g.}, $\bm{M}=[M_{ij}]\in \real^{n\times m}$. For a real $x\in \real$, $\lfloor x \rfloor$ stands for the largest integer not exceeding
$x$. Let $\ball_{r}(\bm{x})\df \{\bm{y}\in \real^{d}: \|\bm{y}-\bm{x}\|_{2}\leq r \}$ denote the Euclidean ball of radius $r$ centered at $\bm{x}$. For a sub-set $\mathcal{S}\subset \real^{d}$ of the Euclidean distance, we define its closure $\bar{\mathcal{S}}\df \{s\in \mathcal{S}: s=\lim_{n\rightarrow\infty}x_{n}, x_{n}\in \mathcal{S},\forall n\in \integer \}$, and its boundary $\partial \mathcal{S}\df \bar{\mathcal{S}}\backslash \mathcal{S}$. Given a random variable $\bm{x}$, we denote its law with $\mathbb{P}_{\bm{x}}\stackrel{\mathrm{law}}{=}\bm{x}$. We use $W^{k,p}(\Omega)$ for $1\leq p\leq \infty$ to denote the Sobolev space of functions $f\in L^{p}(\Omega)$ whose weak derivatives up to order $k$ are in $L^{p}(\Omega)$. The Sobolev space $W^{1,2}(\Omega)$ is denoted by $H^{1}(\Omega)$, and $H_{0}^{1}(\Omega)$ is the space of functions in $H^{1}(\Omega)$ that vanishes at the boundary. The dual space of $H_{0}^{1}(\Omega)$ is denoted by $H^{-1}(\Omega)$.

To establish the concentration results in this paper, we require the following two definitions:
\begin{definition}\textsc{(Sub-Gaussian Norm)}
	\label{Definition:Sub-Gaussian Norm}
	The sub-Gaussian norm of a random variable $Z$, denoted by $\|Z\|_{\psi_{2}}$, is defined as
	\begin{align}
	\|Z\|_{\psi_{2}}= \sup_{q\geq 1} q^{-1/2}(\expect|Z|^{q})^{1/q}.
	\end{align}
	For a random vector $\bm{Z}\in \real^{n}$, its sub-Gaussian norm is defined as  follows
	\begin{align}
	\label{Eq:Sub_Gaussian_random_vector}
	\|\bm{Z}\|_{\psi_{2}}=\sup_{\bm{x}\in \mathrm{S}^{n-1}}\|\langle  \bm{x},\bm{Z}\rangle \|_{\psi_{2}}.
	\end{align}
\end{definition}

\begin{definition}\ \textsc{(Sub-exponential Norm)}
	The sub-exponential norm of a random variable $Z$, denoted by $\|Z\|_{\psi_{1}}$, is defined as follows
	\begin{align}
	\|Z\|_{\psi_{1}}=\sup_{q\geq 1} q^{-1}(\expect[|Z|^{q}])^{1/q}.
	\end{align}
	For a random vector $\bm{Z}\in \real^{n}$, its sub-exponential norm is defined below
	\begin{align}
	\|\bm{Z}\|_{\psi_{1}}= \sup_{\bm{x}\in \mathrm{S}^{n-1}} \|\langle \bm{Z},\bm{x}\rangle\|_{\psi_{1}}.
	\end{align}
\end{definition}

\subsection{Proof of Theorem \ref{Thm:Consistency of Monte Carlo Estimation}}
\label{Proof:Consistency_of_Monte_Carlo}

We begin the proof by recalling the following definitions from the main text:
\begin{subequations}
	\label{Eq:Kernel_1}
	\begin{align}
	\label{Eq:222}
	\widehat{E}_{0}(\mu)&=  \dfrac{2}{n(n-1)} \sum_{1\leq i< j\leq n}y_{i}y_{j}K_{\mu}(\bm{x}_{i},\bm{x}_{j}),   \\
	E_{0}(\mu)&=  \expect\left[y\widehat{y}K_{\mu}(\bm{x},\widehat{\bm{x}})\right].
	\end{align}
\end{subequations}
Furthermore, we recall the definition of the regularized risk function:
\begin{subequations}
	\label{Eq:Kernel_2}
	\begin{align}
	\widehat{E}_{\gamma}(\mu)&=\widehat{E}_{0}(\mu)-\dfrac{2}{n(n-1)\gamma}\sum_{1\leq i< j\leq n}K_{\mu}^{2}(\bm{x}_{i},\bm{x}_{j}),\\
	E_{\gamma}(\mu)&=E_{0}(\mu)-{1\over \gamma}\expect\left[K^{2}_{\mu}(\bm{x},\widehat{\bm{x}})\right].
	\end{align}
\end{subequations}
In Eqs. \eqref{Eq:Kernel_1} and \eqref{Eq:Kernel_2}, the kernel $K_{\mu}(\bm{x},\hat{\bm{x}})$ has the following integral form
\begin{align}
K_{\mu}(\bm{x},\widehat{\bm{x}})=\int_{0}^{\infty}e^{-\xi \|\bm{x}-\widehat{\bm{x}}\|_{2}^{2}}\mu(\mathrm{d}\xi).
\end{align}

Lemmas \ref{Lemma:Consistency with Respect to Training Data} and \ref{Lemma:Consistency with Respect to Langevin_Particles} provide
consistency guarantees with respect to the training data $(\bm{x}_{i},y_{i})_{1\leq i\leq n}\sim_{\text{i.i.d.}} P_{\bm{x},y}$  and the particles $\xi_{0},\cdots,\xi_{N}\sim_{\text{i.i.d.}} \mu$.

\begin{lemma}\textsc{(Consistency with Respect to the Training Data)}
	\label{Lemma:Consistency with Respect to Training Data}
	Suppose  $\bf{(A.1)}-\bf{(A.4)}$ holds. Consider the distribution ball $\mathcal{P}=\{\mu\in \mathcal{M}_{+}(\Xi): W_{2}(\mu,\mu_{0})\leq R \}$, where $\mu_{0}\in \mathcal{M}_{+}(\Xi)$ is an arbitrary distribution. Then,  	
	\begin{align}
	\sup_{\mu\in \mathcal{P}}\left|E_{0}(\mu)-\widehat{E}_{0}(\mu)\right|\leq \max\left\{\dfrac{c_{0} K^{2}}{n} \sqrt{\ln\left(\dfrac{4}{\rho}\right)},\dfrac{c_{1}RK^{2} }{n^{2}}\ln\left(\dfrac{4e^{\sqrt{2}K^{2}}}{\rho}\right) \right\},
	\end{align}	
	with the probability of at least $1-\rho$ over the draw of the training data $(\bm{x}_{i},y_{i})_{1\leq i\leq n}\sim_{\text{i.i.d.}} P_{\bm{x},y}$, where $c_{0}=3^{1\over 4}\times 2^{7\over 2}$ and $c_{1}=\sqrt{3}\times 2^{7}$.
\end{lemma}

The proof of Lemma \ref{Lemma:Consistency with Respect to Training Data} is provided in Section \ref{Appendix:Proof_Lemma:Consistency with Respect to Training Data}. Lemma \ref{Lemma:Consistency with Respect to Training Data} asserts that when the number of training data $n$ tends to infinity, the error due to the finite sample approximation becomes negligible.  In the next lemma, we provide a similar consistency type result for the sample average approximation with respect to the Langevin particles:

\begin{lemma}\textsc{(Consistency with Respect to the Langevin Particles)}
	\label{Lemma:Consistency with Respect to Langevin_Particles}
	Suppose  $\bf{(A.1)}-\bf{(A.4)}$ holds. Consider the distributional ball $\mathcal{P}$ in Lemma \ref{Lemma:Consistency with Respect to Training Data} and consider its empirical approximation $\mathcal{P}^{N}=\{\widehat{\mu}^{N}:W_{2}(\widehat{\mu}^{N},\widehat{\mu}_{0}^{N})\leq R\}$. Then,
	\begin{align}
	\left|\sup_{\mu\in \mathcal{P}}\widehat{E}_{0}(\mu)-\sup_{\widehat{\mu}^{N}\in \mathcal{P}^{N}}\widehat{E}_{0}(\mu)\right|\leq \dfrac{2\sqrt{2}(\xi_{u}-\xi_{l})}{R^{2}\sqrt{N}}\left(1+\ln^{1\over 2}\left(\dfrac{4\sqrt{2N}(\xi_{u}-\xi_{l})}{\rho } \right)\right),
	\end{align}
	with the probability of at least $1-\rho$ over the draw of the initial particles $\xi_{0}^{1},\cdots,\xi_{0}^{N}\sim_{\text{i.i.d.}} \mu_{0}$.
\end{lemma}
The proof of Lemma \ref{Lemma:Consistency with Respect to Langevin_Particles} is presented in Section \ref{Appendix:Proof_Lemma:Consistency with Respect to Langevin_Particles}.

In the next lemma, we quantify the amount of error due to including the regularization $\gamma$ to the risk function:
\begin{lemma}\textsc{(Regularization Error)} 
	\label{Lemma:Regularization_Error}	
	Suppose  $\bf{(A.1)}-\bf{(A.4)}$ holds. Consider the empirical distribution ball $\mathcal{P}^{N}$ of Lemma \ref{Lemma:Consistency with Respect to Langevin_Particles}, and define
	\begin{align}
	\widehat{\mu}^{N}_{\ast}(\gamma)\df \arg\sup_{\widehat{\mu}^{N}\in \mathcal{P}^{N}} E_{\gamma}(\widehat{\mu}^{N}).
	\end{align} 
	Then, for any $\gamma>0$, the following upper bound holds
	\begin{align}
	\left|\sup_{\widehat{\mu}^{N}\in \mathcal{P}^{N}}E_{0}(\widehat{\mu}^{N})-E_{0}(\widehat{\mu}_{\ast}^{N}(\gamma))\right|\leq \dfrac{1}{\gamma}.
	\end{align}
\end{lemma}
The proof of Lemma \ref{Lemma:Regularization_Error} is presented in Section \ref{Appendix:Proof_of_Lemma_Regularization_Error}.

The last lemma of this part is concerned with the approximation error due to using the Sinkhorn divergence in lieu of the Wasserstein distance:  

\begin{lemma}\textsc{(Sinkhorn Divergence Approximation Error)} 
	\label{Lemma:Sinkhorn Divergence Approximation Error}	
	Suppose  $\bf{(A.1)}-\bf{(A.4)}$ holds. Let  $\widehat{\mu}^{N}_{\ast}(\gamma)$ denotes the empirical measure of Lemma \ref{Lemma:Regularization_Error}.  Furthermore, consider the distribution ball $\mathcal{P}_{\varepsilon}^{N}=\{\widehat{\mu}^{N}:W_{2,\varepsilon}(\widehat{\mu}^{N},\widehat{\mu}^{N}_{0})\leq R\}$, and
	\begin{align}
	\widehat{\mu}^{N}_{\ast}(\gamma,\varepsilon)\df \arg\sup_{\widehat{\mu}^{N}\in \mathcal{P}_{\varepsilon}^{N}} E_{\gamma}(\widehat{\mu}^{N})
	\end{align} 
	Then, for any $\varepsilon>0$, the following upper bound holds
	\begin{align}
	\left|E_{0}(\widehat{\mu}_{\ast}^{N}(\gamma))-E_{0}(\widehat{\mu}_{\ast}^{N}(\gamma,\varepsilon))\right|\leq c\exp\left(-{2\over \varepsilon}\right)+\dfrac{2}{\gamma},
	\end{align}
	where $c>0$ is a universal constant independent of $\varepsilon$.
\end{lemma}
The proof of Lemma \ref{Lemma:Sinkhorn Divergence Approximation Error} is presented in Section \ref{Appendix:Proof_Lemma_Sinkhorn Divergence Approximation Error}.

Equipped with Lemmas \ref{Lemma:Consistency with Respect to Training Data}-\ref{Lemma:Sinkhorn Divergence Approximation Error}, we are now in position to proof main inequality of Theorem \ref{Thm:Consistency of Monte Carlo Estimation}.

We use the triangle inequality to decompose the error term into four different components:
\begin{align}
\label{Eq:A1-A4}
&\Big|E_{0}(\mu_{\ast})-E_{0}(\widehat{\mu}^{N}_{\ast}(\gamma,\varepsilon))\Big|\leq \mathsf{E}_{1}+\mathsf{E}_{2}+\mathsf{E}_{3}+\mathsf{E}_{4}+\mathsf{E}_{5},
\end{align}
where the error terms $\mathsf{E}_{i},i=1,2,3,4,5$ are defined as follows
\begin{align*}
\mathsf{E}_{1}&\df \Big| E_{0}(\mu_{\ast})-\sup_{\mu\in \mathcal{P}}\widehat{E}_{0}(\mu) \Big| \\
\mathsf{E}_{2}&\df \Big| \sup_{\mu\in \mathcal{P}}\widehat{E}_{0}(\mu)- \sup_{\widehat{\mu}^{N}\in \mathcal{P}^{N}}\widehat{E}_{0}(\widehat{\mu}^{N}) \Big|\\
\mathsf{E}_{3}&\df \Big|\sup_{\widehat{\mu}^{N}\in \mathcal{P}^{N}}\widehat{E}_{0}(\widehat{\mu}^{N})-E_{0}(\widehat{\mu}^{N}_{\ast}) \Big|\\
\mathsf{E}_{4}&\df \Big|E_{0}(\widehat{\mu}_{\ast}^{N})-E_{0}(\widehat{\mu}_{\ast}^{N}(\gamma))\Big|\\
\mathsf{E}_{5}&\df \Big| E_{0}(\widehat{\mu}_{\ast}^{N}(\gamma))-E_{0}(\widehat{\mu}_{\ast}^{N}(\gamma,\varepsilon))\Big|.
\end{align*}
The error terms $\mathsf{E}_{1}$ and $\mathsf{E}_{3}$ can be bounded using Lemma \ref{Lemma:Consistency with Respect to Training Data}, $\mathsf{E}_{2}$ using Lemma \ref{Lemma:Consistency with Respect to Langevin_Particles}, $\mathsf{E}_{4}$ using Lemma \ref{Lemma:Regularization_Error}, and $\mathsf{E}_{5}$ using Lemma \ref{Lemma:Sinkhorn Divergence Approximation Error}.

\subsection{Proof of Theorem  \ref{Theorem:Mean-Field Partial Differential Equation}}
\label{Appendix:Proof_of_Theorem_Mean_Field}

Consider the projected particle Langevin dynamics in Eq. \eqref{Eq:Tends_to} which we repeat here for the convenience of the reader
\begin{align}
\label{Eq:Projected_L}
\bm{\xi}_{m}=\mathscr{P}_{\Xi^{N}}\left(\bm{\xi}_{m-1}-\eta \nabla \widehat{J}^{N}(\bm{\xi}_{m-1};\bm{z}_{m},\tilde{\bm{z}}_{m})+\sqrt{2\eta\over  \beta}\bm{\zeta}_{m}\right),
\end{align}
for all $m\in \left[0,T\eta^{-1}\right]\cap \integer$, where $\bm{\zeta}_{m}\sim\mathsf{N}(\bm{0},\bm{I}_{N\times N})$. Associated with the discrete-time process $(\bm{\xi}_{m})_{m\in \integer}$, we define the continuous-time c\'{a}dl\'{a}g process $(\bar{\bm{\xi}}_{t})_{0\leq t\leq T}$ that is constant on the interval $\bar{\bm{\xi}}_{t}=\bar{\bm{\xi}}_{m\eta},\forall t\in [m\eta,(m+1)\eta)$ and satisfies the following recursion
\begin{align}
\label{Eq:Recursion_1}
\bar{\bm{\xi}}_{\eta m}=\mathscr{P}_{\Xi^{N}}\left(\bar{\bm{\xi}}_{\eta (m-1)}-\eta \nabla \widehat{J}^{N}(\bar{\bm{\xi}}_{\eta (m-1)};\bm{z}_{m},\tilde{\bm{z}}_{m})+\sqrt{2\over \beta}\bm{\zeta}_{\eta m} \right),
\end{align}
for all $m\in \left[0,T\eta^{-1}\right]\cap \integer$, where $\bm{\zeta}_{\eta m}=\bm{W}_{\eta m}-\bm{W}_{\eta (m-1)}$, and $\bm{W}_{t}$ is a $\mathcal{F}_{t}$-adapted Wiener process with the initial value $\bm{W}_{0}={0}$,  and independent increments 
\begin{align}
\expect\left[e^{i \langle \bm{x},\bm{W}_{t}-\bm{W}_{s}\rangle }|\mathcal{F}_{s}\right]=e^{-{1\over 2}(t-s)\|\bm{x}\|^{2}}, \quad \forall \bm{x}\in \real^{N}.
\end{align}
Therefore, by comparing the dynamics in Eqs.  \eqref{Eq:Projected_L} and \eqref{Eq:Recursion_1}, we observe that $\bar{\bm{\xi}}_{\eta m}=\bm{\xi}_{m}$ for all $m\in [0,T\eta^{-1}]\cap\integer$.

We compare the iterations in Eq. \eqref{Eq:Projected_L} with the following \textit{decoupled} $N$-dimensional reflected It\^{o} stochastic differential
equation
\begin{subequations}
	\label{Eq:Discrete_Time_Version}
	\begin{align}
	\label{Eq:local_time}
	&\bm{\theta}_{t}=\bm{\theta}_{0}-\int_{0}^{t}\nabla J(\bm{\theta}_{s};\mu_{s})\mathrm{d}s+\sqrt{2\over \beta} \bm{W}_{t}+\int_{0}^{t}\bm{n}_{s}L(\mathrm{d}s),\quad 0\leq t\leq T\\
	&\bm{\theta}_{t}\in \Xi^{N},\quad \mathrm{d}L(s)\geq 0, \quad \int_{0}^{t}\bm{1}_{\bar{\Xi}^{N}\backslash\partial\Xi^{N}}(\bm{\theta}_{t})L(\mathrm{d}s)=0,
	\end{align}
\end{subequations}
where 
\begin{itemize}[leftmargin=.25in]
	\item $\bm{\theta}_{0}=(\theta^{1}_{0},\cdots,\theta^{N}_{0})\sim\mu^{\otimes N}_{0}$,
	\item $\mu^{\otimes N}_{s}=\mathbb{P}_{\bm{\theta}_{s}}$, with $\mathbb{P}_{\bm{\theta}_{s}}\stackrel{\mathrm{law}}{{=}}\bm{\theta}_{s}$,
	\item $J(\bm{\theta}_{s};\mu_{s})$ is the drift process defined as follows
	\begin{align}
	J(\bm{\theta}_{s};\mu_{s})\df \dfrac{1}{N}\sum_{k=1}^{N}\expect\left[y\widehat{y}e^{-\theta_{s}^{k} \|\bm{x}-\widehat{\bm{x}} \|^{2}_{2}}\right]+\dfrac{1}{\gamma N}\sum_{k=1}^{N}\int_{0}^{\infty}\expect\left[e^{-(\theta_{s}^{k}+\theta) \|\bm{x}-\widehat{\bm{x}}\|^{2}_{2} }\right]\mu_{s}(\mathrm{d}\theta).
	\end{align}
	Moreover, $\nabla J(\bm{\theta}_{s};\mu_{s})=\left(\nabla_{k} J(\bm{\theta}_{s};\mu_{s})\right)_{1\leq k\leq N}$ is the gradient vector with the following elements
	\begin{align}
	\nonumber
	\nabla_{k} J(\bm{\theta}_{s};\mu_{s})&\df \dfrac{\partial}{\partial \theta_{s}^{k}}J(\bm{\theta}_{s};\mu_{s})\\ 
	&={1\over N}{\partial \over \partial \theta_{s}^{k}} \expect\left[y\widehat{y}e^{-\theta_{s}^{k} \|\bm{x}-\widehat{\bm{x}} \|^{2}_{2}}\right]+{1\over \gamma N}{\partial \over \partial \theta_{s}^{k}}  \int_{0}^{\infty}\expect\left[e^{-(\theta_{s}^{k}+\theta) \|\bm{x}-\widehat{\bm{x}}\|^{2}_{2} }\right]\mu_{s}(\mathrm{d}\theta),
	\end{align}
	
	\item $L$ is the \textit{local time} of $\bm{\theta}_{t}$ at the boundary of the projection space $\Xi^{N}$. In particular, $L$ is a measure on $[0,T]$ that is non-negative, non-decreasing, and whose support is defined
	\begin{align}
	\mathrm{supp}(L)\subseteq \{t\geq 0:\bm{\theta}_{t}\in \partial \Xi^{N} \},
	\end{align}
	\item and, $\bm{n}_{t}$ is the normal vector to the boundary $\partial\Xi^{N}$ at the time $t\in [0,T]$.
\end{itemize}

Let us make some remarks about the stochastic differential equation in Eq. \eqref{Eq:Discrete_Time_Version}.  First, let $\bm{Z}_{t}=\int_{0}^{t}\bm{n}_{s} L(\mathrm{d}s)$. Then, due to Skorokhod \cite{skorokhod1961stochastic} and Tanka's \cite{tanaka2002stochastic} theorems, it is known that the processes $\bm{\theta}_{t}$ and $\bm{Z}_{t}$ are uniquely defined. The uniqueness property is also referred to as the \textit{Skorohod's reflection mapping principle} in stochastic calculus (cf. Definition \ref{Definition:Skorokhod problem}). Second, notice that at each given time $t\in [0,T]$, the drift of the $k$th coordinate $\nabla_{k}J(\bm{\theta}_{s},\mu_{s})$ only depends on the location of the $k$th particle $(\theta_{s}^{k})_{0\leq s\leq t}$, and is independent of other particles $(\theta^{j}_{s})_{0\leq s\leq t},j\not=k$.

Now, consider the following c\`{a}dl\`{a}g process $(\bar{\bm{\theta}}_{t})_{0\leq t\leq T}$, where $\bar{\bm{\theta}}_{t}=\bar{\bm{\theta}}_{m\eta},\forall t\in [m\eta,(m+1)\eta)$, and 
\begin{align}
\label{Eq:Standard_Wiener}
\bar{\bm{\theta}}_{\eta m }=\mathscr{P}_{\Xi^{N}}\left(\bar{\bm{\theta}}_{\eta (m-1)}-\eta \nabla J(\bar{\bm{\theta}}_{\eta(m-1)};\rho_{\eta(m-1)})+\sqrt{2\over \beta}\bm{\zeta}_{\eta m}\right),
\end{align}
for all $m\in \left[0,T\eta^{-1}\right]\cap \integer$, where $\rho^{\otimes N}_{s}=\mathbb{P}_{\bar{\bm{\theta}}_{s}}$, and $\bar{\bm{\theta}}_{0}=\bm{\theta}_{0}$.  We also consider the following  c\`{a}dl\`{a}g process $(\tilde{\bm{\theta}}_{t})_{0\leq t\leq T}$ that is constant on the interval $\tilde{\bm{\theta}}_{t}=\tilde{\bm{\theta}}_{m\eta },\forall t\in [m\eta,(m+1)\eta)$, and has the following recursions
\begin{align}
\label{Eq:Compared}
\tilde{\bm{\theta}}_{\eta m}=\mathscr{P}_{\Xi^{N}}\left(\tilde{\bm{\theta}}_{\eta (m-1)}-\eta \nabla J(\tilde{\bm{\theta}}_{\eta (m-1)};\nu_{\eta (m-1)})+\sqrt{2\over \beta}\tilde{\bm{\zeta}}_{\eta m}\right),
\end{align}
for all $m\in \left[0,T\eta^{-1}\right]\cap \integer$, where $\nu^{\otimes N}_{s}=\mathbb{P}_{\tilde{\bm{\theta}}_{s}}$, $\tilde{\bm{\theta}}_{\eta m}=\bm{\theta}_{0}$, and $\tilde{\bm{\zeta}}_{\eta m}\df \tilde{\bm{W}}_{\eta m}-\tilde{\bm{W}}_{\eta (m-1)}$ is the Wiener process with a drift. In particular,
\begin{align}
\tilde{\bm{W}}_{t}=\bm{W}_{t}+\sqrt{{\beta\over 2}}\int_{0}^{t}(\nabla J(\tilde{\bm{\theta}}_{s},\nu_{s})-\nabla J(\bm{\theta}_{s},\mu_{s}))\mathrm{d}s.
\end{align}
We notice that the processes $(\bar{\bm{\theta}}_{t})_{0\leq t\leq T}$ and $(\tilde{\bm{\theta}}_{t})_{0\leq t\leq T}$ in Eqs. \eqref{Eq:Standard_Wiener} and \eqref{Eq:Compared}, respectively, only differ in the definition of the Wiener process.

In the sequel, we establish three propositions to derive bounds on the Wasserstein distances between the laws of c\'{a}dl\'{a}g processes in Eqs. \eqref{Eq:Recursion_1},\eqref{Eq:Discrete_Time_Version},\eqref{Eq:Standard_Wiener}, and \eqref{Eq:Compared}:  

\begin{proposition}\textsc{(Wasserstein Distance between the Laws of C\'{a}dl\'{a}g Processes)}
	\label{Lemma:X1}	
	Consider the empirical measure $\mu^{N}_{t}\df\widehat{\mu}^{N}_{\lfloor{t\over \eta}\rfloor}\df {1\over N}\sum_{k=1}^{N}\delta_{\xi_{\lfloor{t\over \eta}\rfloor}^{k}}(\xi)$ associated with the Langevin particles in Eq. \eqref{Eq:Recursion_1}. Furthermore, let $\rho_{\eta \lfloor{t\over \eta}\rfloor}^{\otimes N}=\mathbb{P}_{\bar{\bm{\theta}}_{\eta \lfloor{t\over \eta}\rfloor}}$ denotes the law of the random variable $\bar{\bm{\theta}}_{\eta \lfloor{t\over \eta}\rfloor}$ defined by the recursion in Equation \eqref{Eq:Recursion_1}. Then, the following upper bound holds
	\begin{align}
	\sup_{0\leq t\leq T} W_{2}^{2}\left(\mu_{t}^{N},\rho_{\eta \lfloor{t\over \eta}\rfloor}\right)\leq {{{4(\xi_{u}-\xi_{l})\over N}(1+\gamma^{-1})^{2}K^{2}}}\sqrt{{\eta T}\log \left(\dfrac{T}{\eta\rho}\right)}\exp\left( \dfrac{2K^{4}(1+\gamma^{-1})T}{\sqrt{N}}\right),
	\end{align}
	with the probability of (at least) $1-2\rho$.
\end{proposition}

The proof of Proposition \ref{Lemma:X1} is presented in Appendix \ref{Appendix:Proof_of_Lemma_X1}.

\begin{proposition}\textsc{(Change of Measure)}
	\label{Lemma:X2}	
	Consider the laws $\rho_{\eta \lfloor {t\over \eta} \rfloor }^{\otimes N}=\mathbb{P}_{\bar{\bm{\theta}}_{\eta \lfloor {t\over \eta} \rfloor }}$ and $\nu_{\eta \lfloor {t\over \eta} \rfloor }^{\otimes N}=\mathbb{P}_{\tilde{\bm{\theta}}_{\eta \lfloor {t\over \eta} \rfloor }}$, where $\bar{\bm{\theta}}_{\eta \lfloor {t\over \eta} \rfloor }$ and $\tilde{\bm{\theta}}_{\eta \lfloor {t\over \eta} \rfloor }$ evolve according to the recursions in Eqs. \eqref{Eq:Standard_Wiener} and \eqref{Eq:Compared}, respectively. Then,
	\begin{align}
	\nonumber
	\sup_{0\leq t\leq T}W_{2}^{2}\left(\rho_{\eta \lfloor {t\over \eta} \rfloor },\nu_{\eta \lfloor {t\over \eta} \rfloor }\right)\leq \dfrac{8\beta (\xi_{u}-\xi_{l})^{2} K^{4}(1+\gamma^{-1})^{2}}{N}\left(\int_{0}^{\eta \lfloor {T\over \eta} \rfloor}\expect\Big[\|\tilde{\bm{\theta}}_{s}-\bm{\theta}_{s}\|_{2}^{4}\Big]\mathrm{d}s\right)^{1\over 2}.
	\end{align}
\end{proposition}
The proof of Proposition \ref{Lemma:X2} is presented in Appendix \ref{Appendix:Proof_of_Lemma_X2}.

\begin{proposition}\textsc{(Time Discreteization Error)}
	\label{Proposition:X3}	
	Consider the laws $\mu_s^{\otimes N}=\mathbb{P}_{\bm{\theta}_{s}}$ and $\nu_{\eta \lfloor {t\over \eta} \rfloor }^{\otimes N}=\mathbb{P}_{\tilde{\bm{\theta}}_{\eta \lfloor {t\over \eta} \rfloor }}$, where $\bm{\theta}_{t}$ and $\tilde{\bm{\theta}}_{\eta \lfloor {t\over \eta} \rfloor }$ are governed by the dynamics of Eqs. \eqref{Eq:Standard_Wiener} and \eqref{Eq:Compared}, respectively. Then,	 the following inequality holds.
	\begin{align}
	\sup_{0\leq t\leq T}W_{2}^{2}\left(\nu_{\eta \lfloor {t\over \eta} \rfloor },\mu_{t}\right) \leq \expect\left[\sup_{0\leq t\leq T}\|\bm{\theta}_{s}-\tilde{\bm{\theta}}_{s}\|_{2}^{2}\right],
	\end{align}	
	Moreover, for any $p\in \integer$, the following maximal inequality holds
	\begin{align}
	\expect\left[\sup_{0\leq t\leq T}\|\bm{\theta}_{s}-\tilde{\bm{\theta}}_{s}\|_{2}^{2p}\right]	 \leq \sqrt{2B_{p}}\left(\dfrac{2^{4p-1}}{\beta^{2p}}A_{p}+2^{2p-1}C_{p}\right) e^{2^{2p-1}TD_{p}\sqrt{2B_{p}}}.
	\end{align}
	where the constants $A_{p},B_{p},C_{p}$ and $D_{p}$ are defined as follows
	\begin{subequations}
		\begin{align}
		\label{Eq:Recall_A}
		A_{p}&\df  T \left(\dfrac{2p}{2p-1} \right)^{2p}2^{p-1}\eta^{2p-1}N\dfrac{\Gamma\left({N+2p\over 2}\right)}{\Gamma\left({N+2\over 2}\right)},\\
		\label{Eq:Recall_B}
		B_{p}&\df  c_{p}\left(\mathrm{dist}(\bm{\theta}_{0},\partial \Xi^{N}) \right)^{-2p}\left(\left({2\over \beta}\right)^{p}N^pT^{p} +\dfrac{K^{4p}}{N^{p}}T^{2p} \right), \\
		\label{Eq:Recall_C}
		C_{p}&\df 2^{2p-1}\eta \dfrac{K^{4p}}{N^{p}}(1+\gamma^{-1})^{2p},\\
		\label{Eq:Recall_D}
		D_{p}&\df \left(\dfrac{K^{4p}(1+\gamma^{-1})^{2p}}{N^{p}}\log\left(\dfrac{2^{2p}K^{4p}T(1+\gamma^{-1})^{2p}}{N^{p}}\right)+\dfrac{2^{4p-2}K^{4p}}{N^{4p}} \right).
		\end{align}
		Above, $c_{p}>0$ is a constant independent of $N$ and $T$, and $\mathrm{dist}(\bm{\theta}_{0},\partial \Xi^{N})=\min_{\bm{\theta}\in \partial \Xi^{N}} \|\bm{\theta}_{0}-\bm{\theta}\|_{2}$ is the distance of the initial value from the boundary.
	\end{subequations}
\end{proposition}

From the triangle inequality of the 2-Wasserstein distance, we obtain that
\begin{align}
\nonumber
\sup_{0\leq t\leq T}W^{2}_{2}\left(\widehat{\mu}^{N}_{\eta \lfloor {t\over \eta} \rfloor},\mu_{t}\right)
&\leq \sup_{0\leq t\leq T}W^{2}_{2}\left(\widehat{\mu}^{N}_{\eta \lfloor {t\over \eta} \rfloor},\rho_{\eta \lfloor {t\over \eta} \rfloor} \right)+\sup_{0\leq t\leq T}W_{2}^{2}\left(\rho_{\eta \lfloor {t\over \eta} \rfloor},\nu_{\eta \lfloor {t\over \eta} \rfloor} \right)\\  \label{Eq:We_Leverage}
&\hspace{4mm}+\sup_{0\leq t\leq T}W_{2}^{2}\left(\nu_{\eta \lfloor {t\over \eta} \rfloor},\mu_{t}\right).
\end{align}
We now leverage the upper bounds in Propositions \ref{Lemma:X1}, \ref{Lemma:X2}, and \ref{Proposition:X3} to bound each term on the right hand side of Eq. \eqref{Eq:We_Leverage}. Taking the limit $N\rightarrow \infty$ and choosing the prescribed step-size $\eta=\eta_{N}$ in Theorem \ref{Theorem:Mean-Field Partial Differential Equation} yields
\begin{align}
\label{Eq:Due_to_useless}
\lim_{N\rightarrow\infty}\sup_{0\leq t\leq T}W^{2}_{2}\left(\widehat{\mu}^{N}_{\eta \lfloor {t\over \eta} \rfloor},\mu_{t}\right)=0.
\end{align}
For any pseudo-Lipschitz function $\psi:\real_{+}
\rightarrow \real$ and empirical measures $\{\mu_{n}\}_{n\in\integer},\mu_{n}\in \mathcal{M}_{+}(\real_{+})$, it holds that $\int \psi\mu_{n}\rightarrow \int \psi\mu^{\ast}$ as $n\rightarrow \infty$ if and only if $W_{2}(\mu_{n},\mu)\rightarrow 0$, see
\cite[Thm. 6.9]{villani2008optimal}. Therefore, due to Eq. \eqref{Eq:Due_to_useless}, we conclude that $\widehat{\mu}^{N}_{\lfloor {t\over \eta} \rfloor}\stackrel{\text{weakly}}{\rightarrow}  \mu_{t}$ uniformly on the interval $0\leq t\leq T$.

In the sequel, we characterize a distributional dynamics describing the evolution of the Lebesgue density of the limiting measure $\mu_{t}$:
\begin{proposition}
	\label{Proposition:McKean-Vlasov}
	\textsc{(McKean-Vlaso Mean-Field Equation)} Let $(\Omega, \mathcal{F},\mathcal{F}_{t},\prob)$ denotes a probability space, and consider the  $\mathcal{F}_{t}$-adapted reflected diffusion-drift process  $X:\Omega\times [0,T] \mapsto [\xi_{l},\xi_{u}]$ described as follows
	\begin{subequations}
		\begin{align}
		\label{Eq:Reff}
		&X_{t}=X_{0}+\int_{0}^{t} g(X_{s};\mu_{s})\mathrm{d}s+\dfrac{1}{\beta}W_{t}-Z^{+}_{t}+Z^{-}_{t},\\
		&\theta_{0}\sim \mu_{0}, \quad \mathrm{supp}(\mu_{0})\subset \Xi=[\xi_{l},\xi_{u}],\quad ,0\leq t\leq T,
		\end{align}
	\end{subequations}
	where $\mu_{s}=\mathbb{P}_{X_{s}}$ is the law of the underlying process, and $Z^{+}_{t}$ and $Z^{-}_{t}$ are the reflection processes from the boundaries $\xi_{\mathrm{u}}$ and $\xi_{\mathrm{l}}$, respectively.  In particular, $Z_{0}^{+}=Z_{0}^{-}=0$, non-decreasing, c\'{a}dl\'{a}g, and 
	\begin{align}
	\int_{0}^{\infty}(X_{t}-\xi_{\mathrm{l}})\mathrm{d}Z^{-}_{t}=0, \quad \int_{0}^{\infty}(\xi_{\mathrm{u}}-X_{t})\mathrm{d}Z^{+}_{t}=0.
	\end{align}
	Suppose the Lebesgue density $q_{0}(\xi)={\mathrm{d}\mu_{0}\over \mathrm{d}\xi}$ exists. Then, the Lebesgue density of the law of the stochastic process $(p_{t}(\xi)=\mathrm{d}\mu_{t}/\mathrm{d\xi})_{0\leq t}$ at subsequent times is governed by the following one dimensional partial differential equation with Robin boundary conditions
	\begin{subequations}
		\label{Eq:DDD}
		\begin{align}
		&\dfrac{\partial p(t,\xi)}{\partial t}=-\dfrac{\partial}{\partial \xi}(p_{t}(\xi)g(\xi,p_{t}(\xi)))+\dfrac{1}{\beta} \dfrac{\partial^{2}}{\partial\xi^{2}}p_{t}(\xi),\quad \forall t\in [0,T],\forall \xi\in (\xi_{l},\xi_{u}) \\
		&{\partial p_{t}(\xi)\over\partial \xi}+\beta p_{t}(\xi)g(\xi,\mu_{t})\Big\vert_{\xi=\xi_{l}}=0, \quad 	{\partial p_{t}(\xi)\over\partial \xi}+\beta p_{t}(\xi)g(\xi,\mu_{t})\Big\vert_{\xi=\xi_{u}}=0,\quad \forall t\in [0,T] \\
		&p_{0}(\xi)=q_{0}(\xi), \quad \forall \xi\in \Xi.
		\end{align}
	\end{subequations}
\end{proposition}
The proof of Theorem \ref{Proposition:McKean-Vlasov} is a special case of the proof provided by Harrison and Reiman \cite{harrison1981distribution} for general multi-dimensional diffusion-drift processes that solve the Skorokhod problem for a general reflection matrix (cf. Definition \ref{Definition:Skorokhod problem}). More specifically, the proof of \cite{harrison1981distribution} is based on an extension of It\^{o}'s formula proved by Kunita and Watanabe \cite{kunita1967square}, and Grisanov's change of measure technique in Theorem \ref{Thm:Grisanov_Theorem}.

To finish the proof, we apply the result of Proposition \ref{Proposition:McKean-Vlasov} to each coordinate of the stochastic process in Eqs. \eqref{Eq:Discrete_Time_Version}.

$\hfill$ $\blacksquare$

\subsection{Proof of Proposition \ref{Appendix:Proof_of_Uniqueness}}
\label{Appendix:Proof_of_Uniqueness}

We establish the existence and uniqueness via the standard technique of Lax-Milgram theorem in the theory of elliptic PDEs, see \cite[Chapters 5,6]{evans1998}, and \cite{burger2018fokker,burger2016flow}.  In particular, we closely follow the work of \cite{burger2018fokker}, and provide the proof for a slightly more general case, namely we consider the Poisson equation
\begin{subequations}
	\label{Eq:distributional_dynamics_general}
\begin{align}
\mathrm{div} F(\bm{\xi})&=g(\bm{\xi}), \quad  F(\bm{\xi})\df \dfrac{1}{\beta} \nabla p_{\ast}(\bm{\xi})+p_{\ast}(\bm{\xi})\nabla V(\bm{\xi}), \quad \forall \bm{\xi}\in \Omega, \\
\langle F(\bm{\xi}),\bm{n}\rangle&=0, \quad \bm{\xi}\in  \partial \Omega.
\end{align}
\end{subequations}
where $\bm{n}$ is the outward normal vectors at the boundary $\partial \Omega$. Furthermore, we suppose the potential satisfies $V\in W^{1,\infty}(\Omega)$ for every $t\in [0,T]$. In the steady-state regime (\textit{i.e.} when $\partial_{t}p_{t}=0$), the mean-field PDE in \eqref{Eq:distributional_dynamics} is a special case of the general form described in Eq. \eqref{Eq:distributional_dynamics_general} with $\Omega=\Xi\subset \real_{+}$, $\partial \Omega=\{\xi_{l},\xi_{u}\}$, $g(\xi)=0$ for all $\xi\in\Omega$, and $V(\xi)=J(\xi,p_{\ast}(\xi))$, where $p_{\ast}(\xi)$ is a steady-state solution of the mean-field PDE in \eqref{Eq:distributional_dynamics}. Equation \eqref{Eq:distributional_dynamics_general} describes a probability in the space of probability measures in $\mathcal{M}_{+}(\Omega)$, and is to be interpreted in the weak sense. The main difficulty arises from assigning boundary values along $\partial \Omega$ to a function $p_{\ast}\in W^{1,p}(\Omega)$ as it is not in general continuous. The following trace theorem attempts to address this issue:

\begin{definition}(\textsc{Trace Theorem, see, \textit{e.g.}, \cite[Thm. 1]{evans1998}})
	\label{Definition:Trace}
Suppose $\Omega$ is bounded and $\partial \Omega$ is in $C^{1}$. Then, there exists a bounded linear operator $\mathcal{T}:W^{1,p}(\Omega)\rightarrow L^{p}(\partial \Omega)$ such that
\begin{itemize}
	\item $\mathcal{T}p_{\ast}=p_{\ast}\vert_{\partial \Omega}$ if $p_{\ast}\in W^{1,p}(\Omega)\cap C(\Omega)$
	\item $\|\mathcal{T}p_{\ast}\|_{L^{p}(\partial \Omega)}\leq c_{p,\Omega}\|p_{\ast}\|_{W^{1,p}(\Omega)}$,
\end{itemize}
for each $p_{\ast}\in W^{1,p}(\Omega)$, with the constant $c_{p,\Omega}$ depending only on $p$ and $\Omega$. Furthermore, $\mathcal{T}p_{\ast}$ is called \textit{the trace} of $p_{\ast}$ on $\partial \Omega$.
\end{definition}
Equipped with Definition \ref{Definition:Trace}, we can now provide the form of weak solutions of the PDE in Eq. \eqref{Eq:distributional_dynamics_general}:

\begin{definition}(\textsc{Weak Solutions of the Steady-State Poisson Equation})
We say that a function $p_{\ast}\in H^{1}(\Omega)$ is a weak solution to equation \eqref{Eq:distributional_dynamics_general} supplemented with the boundary condition if for all the test functions $\psi\in H^{1}(\Omega)$ the following identity holds
\begin{align}
\nonumber
\dfrac{1}{\beta}\int_{\Omega} \langle \nabla p_{\ast}(\bm{\xi}),\nabla \psi(\bm{\xi})\rangle \mathrm{d}\bm{\xi}&+\int_{\Omega}\langle \nabla V(\bm{\xi}),\nabla \psi(\bm{\xi}) \rangle p_{\ast}(\bm{\xi})\mathrm{d}\bm{\xi}  \\ \label{Eq:weak_formulation}
&-\int_{\partial \Omega}\langle \nabla V(\bm{\sigma}),\bm{n}\rangle  \mathcal{T}p_{\ast}(\bm{\sigma})\mathcal{T}\psi(\bm{\sigma})\mathrm{d}\bm{\sigma} =\langle  g,\psi\rangle_{L^{2}(\Omega)}.
\end{align}
almost everywhere in time $t\in [0,T]$.
\end{definition}

To establish the uniqueness of the steady-state solution, we reformulate the weak solution in terms of the \textit{Slotboom variable} $\rho_{\ast}(\bm{\xi})=p_{\ast}(\bm{\xi})e^{\beta V(\bm{\xi})}$. The transformed flux is given by $F(\bm{\xi})=\dfrac{1}{\beta}e^{-\beta V(\bm{\xi})}\nabla \rho_{\ast}(\bm{\xi})$. Then, the weak formulation in Eq. \eqref{Eq:weak_formulation} can be rewritten as follows
\begin{align}
\label{Eq:Transformed}
\dfrac{1}{\beta}\int_{\Omega}e^{-\beta V(\bm{\xi})} \langle \nabla \rho_{\ast}(\bm{\xi}),\nabla \psi(\bm{\xi})\rangle \mathrm{d}\bm{\xi}-\int_{\partial \Omega}\langle \nabla V(\bm{\sigma}),\bm{n}\rangle  \mathcal{T}\rho_{\ast}(\bm{\sigma})^{-\beta V(\bm{\sigma})}\mathcal{T}\psi\mathrm{d}\bm{\sigma} =\langle  g,\psi\rangle_{L^{2}(\Omega)}.
\end{align}
Define the following bilinear form
\begin{align}
\label{Eq:bilinear_form}
&B:H^{1}(\Omega)\times H^{1}(\Omega)\rightarrow \real,\\
&(\phi,\psi)\mapsto B[\phi,\psi]\df \dfrac{1}{\beta}\int_{\Omega}e^{-\beta V(\bm{\xi})} \langle \nabla \phi,\nabla \psi\rangle \mathrm{d}\bm{\xi}-\int_{\partial \Omega}\langle \nabla V,\bm{n}\rangle  \mathcal{T}\phi e^{-\beta V}\mathcal{T}\psi\mathrm{d}\bm{\sigma}
\end{align}
The bilinear form $B$ is continuous on $H^{1}(\Omega)$ due to the Cauchy-Schwarz inequality and the fact that the Sobolev norm is controlled by the $L^{2}(\Omega)$-norm of the gradient. Moreover, the bilinear form is non-coercive.\footnote{A bilinear form $B:H^{1}(\Omega)\times H^{1}(\Omega)\rightarrow \real$ is coercive if $B[\psi,\psi]\geq c\|\psi\|_{H^{1}(\Omega)}$ for some constant $c>0$.} 

A weak solution $\rho_{\ast}\in H^{1}(\Omega)$ for Eq. \eqref{Eq:Transformed} satisfies the following condition
\begin{align}
B[\rho_{\ast},\psi]=\langle g,\psi \rangle_{L^{2}(\Omega)}, \quad \forall \psi\in H^{1}(\Omega).
\end{align}
We claim $B[\phi,\psi]$ satisfies the hypotheses of Lax-Milgram on $H^{1}(\Omega)$. In particular, we must show that there exist constants $c_{0},c_{1},c_{2}>0$ such that the following energy estimates are satisfied for all $\psi,\phi\in H^{1}(\Omega)$ (see \cite[6.2.2 Theorem 2]{evans1998})
\begin{subequations}
	\label{Eq:Gar}
\begin{align}
\label{Eq:Ff}
\text{(Boundedness):} \hspace{4mm}|B[\phi,\psi]|&\leq c_{0}\|\phi\|_{H^{1}(\Omega)} \|\psi\|_{H^{1}(\Omega)},\\ \label{Eq:Ss}
\text{(G\r{a}rding Inequality):}\hspace{4mm} B[\phi,\phi] &\geq c_{1}\|\phi\|^{2}_{H^{1}(\Omega)}-c_{2}\|\phi\|^{2}_{L^{2}(\Omega)}.
\end{align} 
\end{subequations}
The first condition follows by boundedness of the trace operator $\mathcal{T}$ in Definition \ref{Definition:Trace} as well as boundedness of $V(\bm{\xi})$ and its gradient. In particular,
\begin{align}
\nonumber
|B[\phi,\psi]|&\leq \dfrac{1}{\beta}c_{V}\|\phi\|_{H^{1}(\Omega)}\|\psi\|_{H^{1}(\Omega)}+c'_{V}\left\|\mathcal{T}\phi e^{-\beta V}\right\|_{L^{2}\left(\partial \Omega\right)}\left\|\mathcal{T}\psi \right\|_{L^{2}\left(\partial \Omega\right)}\\
&\leq \dfrac{1}{\beta}c_{V}\|\phi\|_{H^{1}(\Omega)}\|\psi\|_{H^{1}(\Omega)}+c^{2}_{2,\Omega}c'_{V}c_{V}\|\phi\|_{H^{1}(\Omega)}\|\psi\|_{H^{1}(\Omega)},
\end{align}
where $c_{V}\df \sup_{\bm{\xi}\in \Omega} e^{-\beta V(\bm{\xi})}$ and $c'_{V}=\sup_{\bm{\xi}\in \Omega}|\langle \nabla V(\bm{\xi}),\bm{n} \rangle|$.

To establish the G\r{a}rding Inequality, we invoke the following variation of Sobolev type inequality for boundary value problems:
\begin{theorem}\textsc{(Sobolev Type Inequality for Bounded Domains, \cite[p. 4]{roland2015poincar})}
	\label{Eq:Sobolev Type Inequality for Bounded Domains}
For a bounded domain $\Omega \subset \real^{d}$ and for all the functions $\phi\in H^{1}(\Omega)$, there exists a constant $K_{\Omega}>0$ depending on $\Omega$ only, such that the following inequality holds
\begin{align}
\label{Eq:Follow_my_lead}
\|\phi\|^{2}_{H^{1}(\Omega)}\leq K_{\Omega}\left(\|\nabla \phi\|^{2}_{H^{1}(\Omega)}+\|\phi\|^{2}_{L^{2}(\Omega)} \right).
\end{align}
\end{theorem}

Now, for the bilinear form in Eq. \eqref{Eq:bilinear_form} we obtain the following lower bound
\begin{align}
\nonumber
B[\psi,\psi]&=\dfrac{1}{\beta}\int_{\Omega}e^{-\beta V(\bm{\xi})}\|\nabla \psi\|_{2}^{2}\mathrm{d}\bm{\xi}-\int_{\partial \Omega}\langle \nabla V,\bm{n}\rangle \mathcal{T}\psi e^{-\beta V}\mathcal{T}\psi \mathrm{d}\bm{\sigma}\\ \nonumber
&\geq \dfrac{1}{\beta}\tilde{c}_{V}\int_{\Omega}\|\nabla \psi\|_{2}^{2}\mathrm{d}\bm{\xi}-c_{V}c'_{V}\int_{\partial \Omega}|\mathcal{T}\psi|^{2}\mathrm{d}\bm{\sigma}\\ \nonumber
&= \dfrac{1}{\beta}\tilde{c}_{V}\|\nabla \psi\|^{2}_{H^{1}(\Omega)}-c_{V}c'_{V}\| \psi\|^{2}_{L^{2}(\Omega)}   \\ \nonumber
&=\dfrac{1}{\beta K_{\Omega}}\tilde{c}_{V}\|\psi\|_{H^{1}(\Omega)}^{2}-\left(\dfrac{1}{\beta K_{\Omega}}\tilde{c}_{V}+c_{V}c_{V}' \right)\|\psi\|^{2}_{L^{2}(\Omega)},
\end{align}
where $\tilde{c}_{V}\df \inf_{\bm{\xi}\in \Omega}e^{-\beta V(\bm{\xi})}$, and the last inequality follows from Inequality \eqref{Eq:Follow_my_lead} in Theorem \ref{Eq:Sobolev Type Inequality for Bounded Domains}. Having verified the criteria described in Eqs. \eqref{Eq:Gar}, the existence and uniqueness follows by invoking the Lax-Milgram theorem.

\subsection{Proof of Lemma \ref{Lemma:Consistency with Respect to Training Data}}
\label{Appendix:Proof_Lemma:Consistency with Respect to Training Data}

We begin the proof by writing the following basic inequalities 
\begin{align}
\nonumber
\Big|\sup_{\mu\in \mathcal{P}}E_{0}(\mu)-\sup_{\mu\in \mathcal{P}}\widehat{E}_{0}(\mu)\Big|&\leq \sup_{\mu\in \mathcal{P}}\big|E_{0}(\mu)-\widehat{E}_{0}(\mu)\big|\\ \nonumber
&=\sup_{\mu\in \mathcal{P}}\Bigg|{1\over n(n-1)}\sum_{i\not= j}y_{i}y_{j}K_{\mu}(\bm{x}_{i},\bm{x}_{j})-\expect_{P_{y,\bm{x}}^{\otimes 2}}\left[y\widehat{y}K_{\mu}(\bm{x},\widehat{\bm{x}})\right]\Bigg|\\ \nonumber
&=4\sup_{\mu\in \mathcal{P}} \big|\expect_{\mu}[e_{n}(\xi)]\big|\\ \nonumber
&\leq 4\left|\sup_{\mu\in \mathcal{P}} \expect_{\mu}[e_{n}(\xi)]\right|,
\end{align}
where the error term is defined as follows
\begin{align}
\label{Eq:Error_Term}
e_{n}(\xi)
&\df \dfrac{2}{n(n-1)}\sum_{1\leq i< j\leq n}y_{i}y_{j}e^{-\xi \| \bm{x}_{i}-\bm{x}_{j}\|_{2}^{2}}-\expect_{P^{\otimes 2}_{\bm{x},y}}\left[y\widehat{y}e^{-\xi \|\bm{x}-\widehat{\bm{x}}\|_{2}^{2}}\right],
\end{align}
where the last equality follows by using the integral representation of the kernel function in Equation \eqref{Eq:Sch}.

Now, we invoke the following strong duality theorem \cite{gao2016distributionally}:

\begin{theorem}\textsc{(Strong Duality for Robust Optimization, \cite[Theorem 1]{gao2016distributionally})}
	\label{Thm:Strong_Duality}
	Consider the general metric space $(\Xi,d)$, and any normal distribution $\nu\in \mathcal{M}_{+}(\Xi)$, where $\mathcal{M}(\Xi)$ is the set of Borel probability measures on $\Xi$. Then, 
	\begin{align}
	\sup_{\mu\in \mathcal{M}_{+}(\Xi)} \Big\{\expect_{\mu}[\Psi(\xi)]: W_{p}(\mu,\nu)\leq R\Big\}=\min_{\lambda\geq 0}\left\{\lambda R^{p}-\int_{\Xi}\inf_{\xi\in \Xi}[\lambda d^{p}(\xi,\zeta)-\Psi(\xi)] \nu(\mathrm{d}\zeta) \right\},
	\end{align}
	provided that $\Psi$ is upper semi-continuous in $\xi$.
\end{theorem}

Under the strong duality of Theorem \ref{Thm:Strong_Duality}, we obtain that
\begin{align}
\label{Eq:Returning_To_Strong_Duality}
\Big|\inf_{\mu\in \mathcal{P}}E_{0}(\mu)&-\inf_{\mu\in \mathcal{P}}\widehat{E}_{0}(\mu) \Big|\leq 4\left|\min_{\lambda\geq 0}\left\{\lambda R^{p}-\int_{\real^{D}}\inf_{\bm{\zeta}\in \real^{D}}\big[\lambda \|\xi-\zeta\|_{2}^{p}-e_{n}(\zeta)\big] \mu_{0}(\mathrm{d}\xi) \right\}\right|.
\end{align}
In the sequel, let $p=2$. The \textit{Moreau's envelope}  \cite{parikh2014proximal} of a function $f:\mathcal{X}\rightarrow \real$ is defined as follows
\begin{align}
\label{Eq:Moreau's envelope}
M_{f}^{\alpha}(\bm{y})\df\inf_{\bm{x}\in \mathcal{X}} \left\{\dfrac{1}{2\alpha}\|\bm{x}-\bm{y}\|_{2}^{2}+f(\bm{x})\right\},\quad \forall \bm{y}\in \mathcal{X},
\end{align}
where $\alpha>0$ is the regularization parameter. We also define the proximal operator as follows:
\begin{align}
\mathrm{prox}_{f}^{\alpha}(\bm{y})\df\arg\inf_{\bm{x}\in \mathcal{X}} \left\{\dfrac{1}{2\alpha}\|\bm{x}-\bm{y}\|_{2}^{2}+f(\bm{x})\right\},\quad \forall \bm{y}\in \mathcal{X}.
\end{align}

When the function $f$ is differentiable, the following lemma is established in \cite{khuzani2019mean}:

\begin{lemma}\textsc{(Moreau's envelope of Differentiable Functions)}
	\label{Lemma:M_envelope}
	Suppose the function $f:\mathcal{X}\rightarrow \real$ is differentiable. Then, 	the Moreau's envelope defined in Eq. \eqref{Eq:Moreau's envelope} has the following upper bound and lower bounds
	\begin{align}
	\label{Eq:THe_Moreau_lower}
	f(\bm{y})-\dfrac{\alpha}{2} \int_{0}^{1}\sup_{\bm{x}\in \mathcal{X}} \|\nabla f(\bm{y}+s(\bm{x}-\bm{y}))\|_{2}^{2}\mathrm{d}s \leq M_{f}^{\alpha}(\bm{y})\leq f(\bm{y}).
	\end{align}
	In particular, when $f$ is $L_{f}$-Lipschitz, we have
	\begin{align}
	\label{Eq:M_envelope}
	f(\bm{y})-\dfrac{\alpha L^{2}_{f}}{2}\leq M_{f}^{\alpha}(\bm{y})\leq f(\bm{y}).
	\end{align}
\end{lemma}

We now require the following result due to \cite{rockafellar2009variational}:
\begin{proposition}\textsc{(Basic Properties of Moreau's envelope,\cite{rockafellar2009variational})} 
	\label{Proposition:Rockafellar}	
	Let $f:\real\rightarrow \real$ be a lower semi-continuous, proper and convex. The following statements hold for any $\alpha>0$:
	
	\begin{itemize}[leftmargin=*]
		\item The proximal operator $\mathrm{prox}_{f}^{\alpha}(x)$ is unique and continuous, in the sense that $\mathrm{prox}_{f}^{\alpha'}(x')\rightarrow \mathrm{prox}_{f}^{\alpha}(x)$ whenever $(x',\alpha')\rightarrow (x,\alpha)$ with $\alpha>0$.
		\item The value of $M^{\alpha}(x)$ is finite and depends continuously on $(\alpha, x)$, with $M_{\alpha}(x)\rightarrow f(x)$ for all $x\in \real$ as
		$\alpha\rightarrow 0_{+}$.
		\item 	The Moreau envelope function is differentiable with respect to $x$ and the regularization parameter $\alpha$. Specifically, for all $x\in \real$,
		the following properties are true:
		\begin{subequations}
			\begin{align}
			\dfrac{\mathrm{d}}{\mathrm{d}x}M_{f}^{\alpha}(x)&=\dfrac{1}{\alpha} (x-\mathrm{prox}_{f}^{\alpha}(x)),\\ \label{Eq:Extreme}
			\dfrac{\mathrm{d}}{\mathrm{d}\alpha}M_{f}^{\alpha}(x)&=-\dfrac{1}{2\alpha^{2}} (x-\mathrm{prox}_{f}^{\alpha}(x))^{2}.
			\end{align}
		\end{subequations}
	\end{itemize}
\end{proposition}

Now, we return to Equation \eqref{Eq:Returning_To_Strong_Duality}. We leverage the lower bound on Moreau's envelope in Eq. \eqref{Eq:THe_Moreau_lower} of Lemma \ref{Lemma:M_envelope} as follows
\begin{align}
\nonumber
&\Big|\sup_{\mu\in \mathcal{P}}E_{0}(\mu)-\sup_{\mu\in \mathcal{P}}\widehat{E}_{0}(\mu)\Big|\\ \nonumber
&\leq 4\left|\min_{\lambda\geq 0}\left\{\lambda R^{2}-\int_{\real^{D}}M_{-e_{n}}^{{1\over 2\lambda}}(\xi) \mu_{0}(\mathrm{d}\xi) \right\}\right|\\ \nonumber
&\leq  4\left|\min_{\lambda\geq 0}\left\{\lambda R^{2}+\expect_{\mu_{0}}[e_{n}(\xi)]+{1\over 4\lambda}  \expect_{\mu_{0}}\left[\int_{0}^{1}\sup_{\bm{\zeta}\in \real^{D}}\|\nabla e_{n}((1-s)\xi+s\zeta)\|_{2}^{2}\mathrm{d}s\right]  \right\}\right|\\  \label{Eq:To_proceed_From_1}
&\leq 4|\expect_{\mu_{0}}[e_{n}(\xi)]|+4R\expect_{\mu_{0}}\left[\int_{0}^{1}\sup_{\bm{\zeta}\in \real^{D}}\|\nabla e_{n}((1-s)\xi+s\zeta)\|_{2}^{2}\mathrm{d}s \right].
\end{align}
Let $\zeta_{\ast}=\zeta_{\ast}(\xi,s)=\arg\sup_{\zeta\in \real^{D}}\|\nabla e_{n}((1-s)\xi+s\zeta)\|_{2}$. Then, applying the union bound in conjunction with Inequality \eqref{Eq:To_proceed_From_1} yields
\begin{align}
\nonumber
\prob\Bigg(\Big|\sup_{\mu\in \mathcal{P}}E_{0}(\mu)&-\sup_{\mu\in \mathcal{P}}\widehat{E}_{0}(\mu)\Big|\geq \delta \Bigg)\leq \prob\Bigg(\Bigg|\int_{\real^{D}}e_{n}(\xi)\mu_{0}(\mathrm{d}\xi) \Bigg|\geq {\delta\over 8}  \Bigg)\\
&+\prob\Bigg(\int_{\real^{D}}\hspace{-1mm}\int_{0}^{1} \|\nabla e_{n}((1-s)\xi+s\zeta_{\ast})\|^{2}_{2}\mathrm{d}s\mu_{0}(\mathrm{d}\xi)\geq {\delta\over 8R} \Bigg). \label{Eq:union_bound}
\end{align}

Now, we state the following lemma:
\begin{lemma}\textsc{(Tail Bounds for the Finite Sample Estimation Error)}
	\label{Lemma:Tail Bounds for the Finite Sample Estimation Error}
	Consider the estimation error $e_{n}$ defined in Eq. \eqref{Eq:Error_Term}. Then, the following statements hold:
	\begin{itemize}[leftmargin=*]
		\item $Z=\|\nabla e_{n}(\xi)\|_{2}^{2}$ is a sub-exponential random variable with the Orlicz norm of $\|Z\|_{\psi_{1}}\leq 64\sqrt{3}K^{2}/n^{2}$ for every $\xi\in \real_+$. Moreover, 
		\begin{align}
		\label{Eq:probability_bound_1}
		\hspace{10mm}\prob\Bigg( \int_{\real^{+}}\int_{0}^{1}\| \nabla e_{n}((1-s)\xi+s\zeta_{\ast}) \|_{2}^{2}\mathrm{d}s\mu_{0}(\mathrm{d}\xi)\geq \delta \Bigg) \leq 2 e^{-{n^{2}\delta\over 16\sqrt{3}K^{2}}+{\sqrt{2}K^{2}}},
		\end{align}
		\item $e_{n}(\xi)$ is a zero-mean sub-Gaussian random variable with the Orlicz norm of $\|e_{n}(\xi)\|_{\psi_{2}}\leq {16\sqrt{3}K^{2}\over n}$ for every $\xi\in \real_{+}$. Moreover,
		\begin{align}
		\label{Eq:probability_bound_2}
		\prob\left(\left|\int_{\Xi} e_{n}(\xi)\mu_{0}(\mathrm{d}\xi)\right|\geq \delta \right)\geq 2e^{-{n^{2}\delta^{2} \over 16\sqrt{3}K^{2}} }.
		\end{align}
	\end{itemize}
\end{lemma}
The proof of Lemma \ref{Lemma:Tail Bounds for the Finite Sample Estimation Error} is deferred to Appendix \ref{Appendix:Proof_of_Thm_Tail Bounds}. 

Now, we leverage the concentration bounds \eqref{Eq:probability_bound_1} and \eqref{Eq:probability_bound_2} of Lemma \ref{Lemma:Tail Bounds for the Finite Sample Estimation Error} to upper bound the terms on the right hand side of Eq. \eqref{Eq:union_bound} as below
\begin{align}
\nonumber
\prob\Bigg(\Big|\sup_{\mu\in \mathcal{P}}E_{0}(\mu)-\sup_{\mu\in \mathcal{P}}E_{0}(\mu)\Big|\geq \delta \Bigg)&\leq 2e^{-{n^{2}\delta^{2} \over \sqrt{3}\times 2^{7}\times K^{2}} }+2 e^{-{n^{2}\delta\over \sqrt{3}\times 2^{7}\times RL^{4}}+{\sqrt{2}K^{2}}}\\ 
&\leq 4\max\left\{ e^{-{n^{2}\delta^{2} \over \sqrt{3}\times 2^{7}\times K^{2}} },e^{-{n^{2}\delta\over \sqrt{3}\times 2^{7}\times RK^{2}}+{\sqrt{2}K^{2}}}\right\},
\end{align}
where the last inequality comes from $a+b\leq 2\max\{a,b\}$.
Therefore, with the probability of (at least) $1-\varrho$, we have that
\begin{align}
\nonumber
\Big|\sup_{\mu\in \mathcal{P}}E_{0}(\mu)-\sup_{\mu\in \mathcal{P}}E_{0}(\mu)\Big|\leq  \max\left\{\dfrac{c_{0} K^{2}}{n} \sqrt{\ln\left(\dfrac{4}{\varrho}\right)},\dfrac{c_{1}RK^{2} }{n^{2}}\ln\left(\dfrac{4e^{\sqrt{2}K^{2}}}{\varrho}\right) \right\},
\end{align}
where $c_{0}=3^{1\over 4}\times 2^{7\over 2}$, and $c_{1}=\sqrt{3}\times 2^{7}$. \hfill $\blacksquare$

\subsection{Proof of Lemma \ref{Lemma:Consistency with Respect to Langevin_Particles}}
\label{Appendix:Proof_Lemma:Consistency with Respect to Langevin_Particles}

We recall the definitions of the population and empirical distributional balls
\begin{subequations}
	\begin{align}
	\mathcal{P}&= \{\mu\in \mathcal{M}(\real_{+}): W_{2}(\mu,\mu_{0})\leq R \}.  \\
	\mathcal{P}^{N}&= \{\widehat{\mu}^{N}\in \mathcal{M}(\real_{+}):W_{2}(\widehat{\mu}^{N},\widehat{\mu}^{N}_{0})\leq R \}.
	\end{align}
\end{subequations}
Due to the strong duality result of Theorem \ref{Thm:Strong_Duality}, the following identity holds
\begin{align}
\nonumber
\sup_{\mu\in \mathcal{P}}\widehat{E}_{0}(\mu)&=\sup_{\mu\in \mathcal{P}}\dfrac{2}{n(n-1)}\sum_{1\leq i<j\leq n}y_{i}y_{j}\expect_{\mu}\left[e^{-\xi \|\bm{x}_{i}-\bm{x}_{j}\|^{2}_{2}}\right]\\  \nonumber
&= \inf_{\lambda\geq 0}\left\{ \lambda R^{2}-\int_{\Xi}\hspace{-.2mm}\inf_{\zeta \in \Xi}\left\{ \lambda(\xi-\zeta)^{2}-\phi_{n}(\zeta)\right\}\mu_{0}(\mathrm{d}\xi) \right\}\\ \label{Eq:Lambda_1}
&=\inf_{\lambda\geq 0}\left\{ \lambda R^{2}-\int_{\Xi}M^{1\over 2\lambda}_{\phi_{n}}(\xi) \mu_{0}(\mathrm{d}\xi) \right\}.
\end{align}
Similarly, we have 
\begin{align}
\nonumber
\sup_{\widehat{\mu}^{N}\in \mathcal{P}^{N}} \widehat{E}_{0}(\mu)&=\sup_{\widehat{\mu}^{N}\in \mathcal{P}^{N}} \dfrac{2}{n(n-1)}\sum_{1\leq i<j\leq n}y_{i}y_{j}\expect_{\widehat{\mu}^{N}}\left[e^{-\xi \|\bm{x}_{i}-\bm{x}_{j} \|_{2}^{2}}\right]   \\ \nonumber
&= \inf_{\lambda\geq 0}\left\{ \lambda R^{2}-\dfrac{1}{N}\sum_{k=1}^{N}\inf_{\zeta\in \Xi}\left\{ \lambda (\xi_{0}^{k}-\zeta)^{2}-\phi_{n}(\zeta)\right\} \right\}\\ \label{Eq:Lambda_2}
&=\inf_{\lambda\geq 0}\left\{\lambda R^{2}-\dfrac{1}{N}\sum_{k=1}^{N}M_{\phi_{n}}^{1\over 2\lambda}(\xi_{0}^{k}) \right\}.
\end{align}
In Eqs. \eqref{Eq:Lambda_1} and \eqref{Eq:Lambda_2}, $\phi_{n}(\zeta)$ has the following definition
\begin{align}
\phi_{n}(\zeta)\df \dfrac{2}{n(n-1)} \sum_{1\leq i<j\leq n}y_{i}y_{j}e^{-\zeta \|\bm{x}_{i}-\bm{x}_{j}\|^{2}_{2}}.
\end{align}

Let $\lambda_{\ast}$ and $\widehat{\lambda}^{N}_{\ast}$ denote the solution of the optimization problems in Eqs. \eqref{Eq:Lambda_1} and \eqref{Eq:Lambda_2}, respectively. We claim that 
\begin{align}
\lambda_{\ast},\widehat{\lambda}^{N}_{\ast}\in \Lambda\df \left[0,\dfrac{2}{R^{2}} \right].
\end{align}
To establish the upper bound on $\widehat{\lambda}^{N}_{\ast}$, we first derive a crude upper bound on the objective value
\begin{align}
\nonumber
\inf_{\lambda\geq 0}\left\{ \lambda R^{2}-\dfrac{1}{N}\sum_{k=1}^{N}\inf_{\zeta\in \Xi}\left\{ \lambda (\xi_{0}^{k}-\zeta)^{2}-\phi_{n}(\zeta)\right\} \right\} &\nonumber\stackrel{\rm{(a)}}{\leq} \sup_{\zeta\in \Xi} \phi_{n}(\zeta)\\ \nonumber
&\leq  \sup_{\zeta\in \Xi}|\phi_{n}(\zeta) |\\ \label{Eq:Baba_1}
&\stackrel{\rm{(b)}}{\leq} 1,
\end{align}
where $\rm{(a)}$ follows by evaluating the objective function at $\lambda=0$, and $\rm{(b)}$ is due to the following inequality
\begin{align}
\nonumber
|\phi_{n}(\zeta)|&=\left| \dfrac{2}{n(n-1)} \sum_{1\leq i<j\leq n}y_{i}y_{j}e^{-\zeta \|\bm{x}_{i}-\bm{x}_{j}\|^{2}_{2}} \right|\\ \nonumber
&\leq \dfrac{2}{n(n-1)}\sum_{1\leq i<j\leq n}  \left|y_{i}y_{j}\right|e^{-\zeta \|\bm{x}_{i}-\bm{x}_{j}\|_{2}^{2}}\\ \label{Eq:Due_to}
&\leq 1,
\end{align}
which holds for all $\zeta\in \Xi$. To derive the last inequality, we used the fact that $y_{i},y_{j}\in \{-1,1\}$ and $\exp(-\zeta \|\bm{x}_{i}-\bm{x}_{j}\|_{2}^{2})\leq 1$ for all $\zeta\in \Xi$. Now, for $\lambda={\delta\over R^{2}},\delta>2$, we obtain that
\begin{align}
\nonumber
\dfrac{\delta}{R^{2}}\cdot R^{2}-\dfrac{1}{N}\sum_{k=1}^{N}\inf_{\zeta\in \Xi}\left\{\dfrac{\delta}{R^{2}}(\xi_{0}^{k}-\zeta)-\phi_{n}(\zeta)\right\}
&\stackrel{\rm{(c)}}{\geq} \delta-\dfrac{1}{N}\sum_{k=1}^{N}\phi_{n}(\xi_{0}^{k})\\ \nonumber
&\geq \delta-\dfrac{1}{N}\sum_{k=1}^{N}|\phi_{n}(\xi_{0}^{k})|\\ \label{Eq:Baba_2}
&\stackrel{\rm{(d)}}{\geq} \delta-1>1,
\end{align}
where $\mathrm{(c)}$ is due to the upper bound on Moreau's envelop in Lemma \ref{Lemma:M_envelope}, and $\mathrm{(d)}$ is due to the upper bound in Eq. \eqref{Eq:Due_to}. From Inequalities \eqref{Eq:Baba_1} and \eqref{Eq:Baba_2} we conclude that the objective value of the minimization for $\lambda=\delta/R^{2},\delta>2$ is strictly larger than the objective value evaluated at $\lambda=0$. Thus, necessarily $\widehat{\lambda}_{\ast}^{N}\leq {2\over R^{2}}$. Using a similar argument (\textit{mutatis mutandis} for the upper bound \eqref{Eq:Due_to}), we can prove that $\lambda_{\ast}\leq {2\over R^{2}}$.

Then, the following inequality holds
\begin{align}
\nonumber
&\Bigg|\sup_{\mu\in \mathcal{P}}\widehat{E}_{0}(\mu)-\sup_{\widehat{\mu}^{N}\in \mathcal{P}^{N}} \widehat{E}_{0}(\widehat{\mu}^{N})\Bigg|\\ \nonumber
&= \Bigg|\inf_{\lambda\in \Lambda}\left\{\lambda R^{2}-\int_{\Xi}\inf_{\zeta\in \Xi} M^{1\over 2\lambda}_{\phi_{n}}(\xi)\mu_{0}(\mathrm{d}\xi) \right\} -\inf_{\lambda\in \Lambda}\left\{\lambda R^{2}-\dfrac{1}{N}\sum_{k=1}^{N}\inf_{\zeta\in \Xi}M_{\phi_{n}}^{1\over 2\lambda}(\xi_{0}^{k}) \right\}   \Bigg|\\ 
&\leq \sup_{\lambda\in \Lambda}\left| \dfrac{1}{N}\sum_{k=1}^{N}\inf_{\zeta\in \Xi}M_{\phi_{n}}^{1\over 2\lambda}(\xi_{0}^{k})-\int_{\Xi}\inf_{\zeta\in \Xi}M_{\phi_{n}}^{1\over 2\lambda}(\xi)\mu_{0}(\mathrm{d}\xi) \right|, \label{Eq:From_T}
\end{align}
where the last inequality follows by the fact that for two bounded functions $f,g:\mathcal{X}\rightarrow \real$, we have $|\inf_{\mathcal{X}} f-\inf_{\mathcal{X}} g|\leq \sup_{\mathcal{X}} |f-g|$.

For any given $\lambda\in \Lambda$, define the following function
\small{\begin{align}
	\nonumber
	T_{\lambda}:\mathcal{\real}^{N}_{+}&\rightarrow \real\\ \nonumber
	\hspace{-4mm}\bm{\xi}_{0}\df (\xi_{0}^{1},\cdots,\xi_{0}^{N})&\mapsto T_{\lambda}(\bm{\xi}_{0})=\dfrac{1}{N}\sum_{k=1}^{N}M_{\phi_{n}}^{{1\over 2\lambda}}(\xi^{k}_{0})-\int_{\Xi}M_{\phi_{n}}^{1\over 2\lambda}(\xi)\mu_{0}(\mathrm{d}\xi),
	\end{align}}\normalsize
Then, we rewrite Inequality \eqref{Eq:From_T} as follows
\begin{align}
\label{Eq:damaghoo}
\Bigg|\sup_{\mu\in \mathcal{P}}\widehat{E}_{0}(\mu)-\sup_{\widehat{\mu}^{N}\in \mathcal{P}^{N}} \widehat{E}_{0}(\widehat{\mu}^{N})\Bigg|\leq \sup_{\lambda\in \Lambda} |T_{\lambda}(\bm{\xi}_{0})| .
\end{align}

Let $\bm{\xi}_{0}=(\xi_{0}^{1},\cdots,\xi_{0}^{k},\cdots,\xi_{0}^{N})\in \Xi^{N}_{+}$ and $\widetilde{\bm{\xi}}_{0}=(\xi_{0}^{1},\cdots,\widetilde{\xi}_{0}^{k},\cdots,\xi_{0}^{N})\in \Xi^{N}_{+}$ denote two sequences that differs in the $k$-th coordinate for $1\leq k\leq N$. Then,
\begin{align}
\label{Eq:customer}
\big|T_{\lambda}(\bm{\xi})-T_{\lambda}(\widetilde{\bm{\xi}}_{0})\big|
&\leq \dfrac{1}{N}\Big|M_{\phi_{n}}^{{1\over 2\lambda}}(\xi^{k}_{0})-M_{\phi_{n}}^{{1\over 2\lambda}}(\widetilde{\xi}^{k}_{0})\Big|\\ \nonumber
&\leq \dfrac{1}{N} \Bigg|\inf_{\zeta\in \Xi}\left\{\lambda|\zeta-\xi_{0}^{k}|^{2}-\dfrac{2}{n(n-1)}\sum_{1\leq i<j\leq n}y_{i}y_{j}e^{-\zeta \|\bm{x}_{i}-\bm{x}_{j}\|_{2}^{2}} \right\} \\ \nonumber
&\hspace{4mm}-\inf_{\zeta\in \Xi}\left\{\lambda|\zeta-\widetilde{\xi}_{0}^{k}|^{2}-\dfrac{2}{n(n-1)}\sum_{1\leq i<j\leq n}y_{i}y_{j}e^{-\zeta \|\bm{x}_{i}-\bm{x}_{j}\|_{2}^{2}} \right\}\Bigg| \\ \nonumber
&\stackrel{\rm{(a)}}{\leq}\dfrac{1}{N} \sup_{\zeta\in \Xi} \left| \lambda |\zeta-\xi_{0}^{k}|^{2}-\lambda|\zeta-\tilde{\xi}^{k}_{0}|^{2} \right| \\ \label{Eq:Mystery_1}
&\leq \dfrac{4}{NR^{2}}(\xi_{u}-\xi_{l})^{2},
\end{align}
where last inequality follows from the fact that $\xi_{0}^{k},\tilde{\xi}_{0}^{k}\in \Xi=[\xi_{l},\xi_{u}]$, and that $\lambda\in \Lambda=[0,(2/R^{2})]$. From McDiarmid's martingale inequality \cite{mcdiarmid1989method} we obtain that
\begin{align}
\prob\left(|T_{\lambda}(\bm{\xi}_{0})|\geq \delta \right)\leq \exp\left(-\dfrac{NR^{4}\delta^{2}}{8(\xi_{u}-\xi_{l})^{2}}\right). 
\end{align}
for any fixed $\lambda\in \Lambda$, and $\delta>0$. Consider an $\varepsilon$-net covering of the interval $\Lambda$ denoted by $\mathcal{N}(\varepsilon,\Lambda)=\{\lambda_{1},\cdots,\lambda_{N(\varepsilon,\Lambda)} \}$, where $N(\varepsilon, \Lambda)\leq {2\over R^{2}\varepsilon}$.

The mapping $\lambda\mapsto T_{\lambda}(\bm{\xi}_{0})$ is Lipschitz on the domain $\lambda\in \Lambda$. Indeed, for any $\lambda_{1},\lambda_{2}\in \Lambda$, we obtain 
\begin{align}
\nonumber
|T_{\lambda_{1}}(\bm{\xi}_{0})-T_{\lambda_{2}}(\bm{\xi}_{0})|&\leq \dfrac{1}{N}\sum_{k=1}^{N}\left|M_{\phi_{n}}^{1\over 2\lambda_{1}}(\xi_{0}^{k})-M_{\phi_{n}}^{1\over 2\lambda_{2}}(\xi_{0}^{k}) \right|\\ \nonumber
&\hspace{4mm}+\int_{\Xi} \left|M_{\phi_{n}}^{1\over 2\lambda_{1}}(\xi)-M_{\phi_{n}}^{1\over 2\lambda_{2}}(\xi) \right|\mu(\mathrm{d}\xi).
\end{align}

We leverage Eq. \eqref{Eq:Extreme} of Proposition \ref{Proposition:Rockafellar} to obtain
\begin{align}
{\mathrm{d}\over \mathrm{d}\lambda}M^{{1\over 2\lambda}}_{f}(\xi)=(\xi-\mathrm{prox}_{f}^{1\over 2\lambda}(\xi))^{2}\leq 4(\xi_{u}0\xi_{l})^{2}, 
\end{align}
where the last inequality is due to the fact that $\xi,\mathrm{prox}_{f}^{1\over 2\lambda}(\xi)\in \Xi=[\xi_{l},\xi_{u}]$. Consequently, $M_{f}^{1\over 2\lambda}$ is $4\xi_{u}^{2}$-Lipschitz, and hence $T_{\lambda}(\bm{\xi}_{0})$ is $8\xi_{u}^{2}$-Lipschitz.
\begin{align}
\label{Eq:Travel_Bound}
\sup_{\lambda \in \Lambda}|T_{\lambda}(\bm{\xi}_{0})|\leq \max_{i\in \mathcal{N}(\varepsilon,\Lambda)}|T_{\lambda_{i}}(\bm{\xi}_{0})|+ 8(\xi_{u}-\xi_{l})^{2}\varepsilon.
\end{align}
Using the union bound yields
\begin{align}
\nonumber
\prob\left(\max_{i\in \mathcal{N}(\varepsilon,\Lambda)} |T_{\lambda_{i}}(\bm{\xi}_{0})|\geq \delta \right)&\leq \cup_{i=1}^{N(\varepsilon,\Lambda)} \prob(|T_{\lambda_{i}}(\bm{\xi}_{0})|\geq \delta)\\
&\leq \dfrac{2}{R^{2}\varepsilon}\exp\left(-\dfrac{NR^{4}\delta^{2}}{8(\xi_{u}-\xi_{l})^{2}}\right)
\end{align}
Therefore, with the probability of (at least) $1-\rho$, we have
\begin{align}
\label{Eq:Save}
\max_{i\in \mathcal{N}(\varepsilon,\Lambda)} |T_{\lambda_{i}}(\bm{\xi}_{0})|\leq \dfrac{2\sqrt{2}(\xi_{u}-\xi_{l})}{R^{2}\sqrt{N}}\sqrt{\log\left(\dfrac{2}{\varepsilon R^{2} \rho } \right)}.
\end{align}
We plug Eq. \eqref{Eq:Save} into \eqref{Eq:Travel_Bound} 
\begin{align}
\sup_{\lambda \in \Lambda}|T_{\lambda}(\bm{\xi}_{0})|\leq  \dfrac{2\sqrt{2}(\xi_{u}-\xi_{l})}{R^{2}\sqrt{N}}\sqrt{\log\left(\dfrac{2}{\varepsilon R^{2} \rho } \right)}+8(\xi_{u}-\xi_{l})^{2}\varepsilon.
\end{align}
Since the size of the net $0<\varepsilon<{2\over R^{2}}$ is arbitrary, we let $\varepsilon={1\over 2\sqrt{2}(\xi_{u}-\xi_{l})R^{2}\sqrt{N}}$. Then, from Eq. \eqref{Eq:damaghoo} we obtain
\begin{align}
\Bigg|\sup_{\mu\in \mathcal{P}}\widehat{E}_{0}(\mu)-\sup_{\widehat{\mu}^{N}\in \mathcal{P}^{N}} \widehat{E}_{0}(\widehat{\mu}^{N})\Bigg|\leq \dfrac{2\sqrt{2}(\xi_{u}-\xi_{l})}{R^{2}\sqrt{N}}\left(1+\sqrt{\log\left(\dfrac{4\sqrt{2N}(\xi_{u}-\xi_{l})}{\rho } \right)}\right),
\end{align}
with the probability of at least $1-\rho$. \hfill $\blacksquare$

\subsection{Proof of Lemma  \ref{Lemma:Regularization_Error}}
\label{Appendix:Proof_of_Lemma_Regularization_Error}

We define the solution of the population objective function as follows
\begin{align}
\nonumber
\widehat{\mu}^{N}_{\ast}(\gamma)&\df \arg\sup_{\widehat{\mu}^{N}\in \mathcal{P}^{N}} E_{\gamma}(\widehat{\mu}^{N})\\  \label{Eq:optimization111}
&=\arg\sup_{\widehat{\mu}^{N}\in \mathcal{P}^{N}}E_{0}(\widehat{\mu}^{N}) 
-\dfrac{1}{\gamma}\expect_{P_{\bm{x}}^{\otimes 2}}\left[K^{2}_{\mu}(\bm{x},\widehat{\bm{x}})\right].
\end{align}
We also define the solution of the empirical kernel alignment as follows
\begin{align}
\label{Eq:Optimization222}
\widehat{\mu}^{N}_{\circ}&\df   \arg\sup_{\widehat{\mu}^{N}\in \mathcal{P}^{N}} E_{0}(\widehat{\mu}^{N}).
\end{align}
Due to the optimality of the empirical measure $\widehat{\mu}_{\ast}^{N}(\gamma)$ for the inner optimization in Eq. \eqref{Eq:optimization111}, the following inequality holds
\begin{align}
\label{Eq:rearrange}
E_{\gamma}(\widehat{\mu}_{\circ}^{N}) &\leq E_{\gamma}(\widehat{\mu}^{N}_{\ast})\leq E_{0}(\widehat{\mu}^{N}_{\ast}).
\end{align}
We now expand $E_{\gamma}(\widehat{\mu}_{\circ}^{N})$ to obtain the following inequality
\begin{align}
\nonumber
E_{0}(\widehat{\mu}^{N}_{\circ})-\dfrac{1}{\gamma}\expect_{P_{\bm{x}}^{\otimes 2}}[K_{\widehat{\mu}^{N}_{\circ}}^{2}(\bm{x},\widehat{\bm{x}})]\leq E_{0}(\widehat{\mu}^{N}_{\ast}).
\end{align}
After rearranging the terms in Eq. \eqref{Eq:rearrange}, we arrive at
\begin{align}
\label{Eq:Majani_1}
E_{0}(\widehat{\mu}_{\circ}^{N})-E_{0}(\widehat{\mu}_{\ast}^{N})
&\leq \dfrac{1}{\gamma}\expect_{P_{\bm{x}}^{\otimes 2}}[K_{\widehat{\mu}^{N}_{\circ}}^{2}(\bm{x},\widehat{\bm{x}})]\leq \dfrac{1}{\gamma},
\end{align}
where the last step follows by the fact that for radial kernels $K(\bm{x},\widehat{\bm{x}})\leq 1$ for all $\bm{x},\widehat{\bm{x}}\in \mathcal{X}$. Similarly, due to the optimality of the empirical measure $\widehat{\mu}^{N}_{\circ}$ for the optimization in Eq. \eqref{Eq:Optimization222} we have that
\begin{align}
\label{Eq:Majani_2}
E_{0}(\widehat{\mu}^{N}_{\ast})\leq E_{0}(\widehat{\mu}^{N}_{\circ}).
\end{align}
Combining Eqs. \eqref{Eq:Majani_1} and \eqref{Eq:Majani_2} now yields
\begin{align}
\big| E_{0}(\widehat{\mu}^{N}_{\ast})- E_{0}(\widehat{\mu}^{N}_{\circ})\big|\leq \dfrac{1}{\gamma}.
\end{align}
\hfill$\blacksquare$

\subsection{Proof of Lemma \ref{Lemma:Sinkhorn Divergence Approximation Error}}
\label{Appendix:Proof_Lemma_Sinkhorn Divergence Approximation Error}

To establish the proof, we state the following proposition due to \cite[Proposition 1]{luise2018differential}:
\begin{proposition}\textsc{(Approximation Error in Sikhorn's divergence, \cite{{luise2018differential}})}
	\label{Proposition:Approximation_Error}
	For any pair of discrete measures $\widehat{\mu}^{N},\widehat{\nu}^{N}\in \mathcal{M}(\mathcal{X})$, the following inequality holds
	\begin{align}
	\label{Eq:Inequality_Sikhorn}
	\left|W_{2,\varepsilon}(\widehat{\mu}^{N},\widehat{\nu}^{N})-W_{2}(\widehat{\mu}^{N},\widehat{\nu}^{N})\right|\leq ce^{-{1\over \varepsilon}},
	\end{align}
	where $c$ is a constant independent of $\varepsilon$, that depends on the support of $\widehat{\mu}^{N}$, and $\widehat{\nu}^{N}$.
\end{proposition}

Using Lagrange's multipliers for the distributional constraints yield the following saddle point problems
\begin{subequations}
	\begin{align}
	\label{Eq:BB_1}
	\sup_{\widehat{\mu}^{N}\in \mathcal{M}(\Xi)}\inf_{h\in \real_{+}}	J_{h}(\widehat{\mu}^{N}(\xi))&=E_{\gamma}(\widehat{\mu}^{N}(\xi))-\dfrac{h}{2}(W_{2}(\widehat{\mu}^{N},\widehat{\mu}_{0}^{N})-R)\\
	\label{Eq:BB_2}
	\sup_{\widehat{\mu}^{N}\in \mathcal{M}(\Xi)}\inf_{h\in \real_{+}}	J_{h}^{\varepsilon}(\widehat{\mu}^{N}(\xi))&= E_{\gamma}(\widehat{\mu}^{N}(\xi))-\dfrac{h}{2}(W_{2,\varepsilon}(\widehat{\mu}^{N},\widehat{\mu}_{0}^{N})-R).
	\end{align}
\end{subequations}
Let $(\widehat{\mu}^{N}_{\diamond},h_{\diamond})$ and $(\widehat{\mu}^{N}_{\ast},h_{\ast})$  denote the saddle points of \eqref{Eq:BB_1} and \eqref{Eq:BB_2}, respectively. Due to the approximation error in Eq. \eqref{Eq:Inequality_Sikhorn} of Proposition \ref{Proposition:Approximation_Error}, the following inequality holds for all $\varepsilon>0$,
\begin{align}
\label{Eq:Scilly}
\left|J_{h}(\widehat{\mu}^{N}(\xi))-J_{h}^{\varepsilon}(\widehat{\mu}^{N}(\xi))\right|\leq {1\over 2}hce^{-{1\over \varepsilon}}.
\end{align}
respectively. Due to the optimality of $(\widehat{\mu}^{N}_{\diamond},h_{\diamond})$ for the saddle point optimization \eqref{Eq:BB_1}, we have that
\begin{align}
\label{Eq:damaghoo_1} E_{\gamma}(\widehat{\mu}_{\diamond}^{N})&=J_{h_{\diamond}}(\widehat{\mu}_{\diamond}^{N})\geq J_{h_{\diamond}}(\widehat{\mu}_{\ast}^{N}).
\end{align}
Using the upper bound in Eq. \eqref{Eq:Scilly} yields
\begin{align}
\nonumber
J_{h_{\diamond}}(\widehat{\mu}^{N}_{\ast})&\geq J_{h_{\diamond}}^{\varepsilon}(\widehat{\mu}^{N}_{\ast})-\dfrac{1}{2}h_{\diamond}ce^{-{1\over \varepsilon}}\\ \nonumber
&\geq J_{h_{\ast}}^{\varepsilon}(\widehat{\mu}^{N}_{\ast})-\dfrac{1}{2}h_{\diamond}ce^{-{1\over \varepsilon}}\\ \label{Eq:damaghoo_2}
&=E_{\gamma}(\widehat{\mu}^{N}_{\ast})-\dfrac{1}{2}h_{\diamond}ce^{-{1\over \varepsilon}}.
\end{align}
Combining Eqs. \eqref{Eq:damaghoo_1} and \eqref{Eq:damaghoo_2} yields
\begin{align}
\label{Eq:Trickle_1}
E_{\gamma}(\widehat{\mu}^{N}_{\ast})-E_{\gamma}(\widehat{\mu}^{N}_{\diamond})\leq \dfrac{1}{2}h_{\diamond}ce^{-{1\over \varepsilon}}.
\end{align}
Similarly, it can be shown that
\begin{align}
\label{Eq:Trickle_2}
E_{\gamma}(\widehat{\mu}^{N}_{\diamond})-E_{\gamma}(\widehat{\mu}^{N}_{\ast})\leq \dfrac{1}{2}h_{\ast}ce^{-{1\over \varepsilon}}.
\end{align}
Putting together Eqs. \eqref{Eq:Trickle_1} and \eqref{Eq:Trickle_2} yields the following inequality
\begin{align}
\label{Equation:A51}
\left|E_{\gamma}(\widehat{\mu}_{\ast}^{N})-E_{\gamma}(\widehat{\mu}_{\diamond}^{N}) \right|\leq{1\over 2} \max\{h_{\ast},h_{\diamond} \}c\exp\left(-{1\over \varepsilon}\right).
\end{align}
From Eqs. \eqref{Eq:Kernel_1} and \eqref{Eq:Kernel_2}, we also have
\begin{align}
\label{Eq:shash}
\left|E_{0}(\widehat{\mu}_{\ast}^{N})-E_{0}(\widehat{\mu}_{\diamond}^{N})  \right|\leq \left|E_{\gamma}(\widehat{\mu}_{\ast}^{N})-E_{\gamma}(\widehat{\mu}_{\diamond}^{N})\right|+\dfrac{1}{\gamma}\expect\left[K^{2}_{\widehat{\mu}^{N}_{\ast}}(\bm{x},\widehat{\bm{x}})+K^{2}_{\widehat{\mu}^{N}_{\diamond}}(\bm{x},\widehat{\bm{x}})\right].
\end{align}
Plugging Eq. \eqref{Equation:A51} into Eq. \eqref{Eq:shash} and using the fact that $K_{\mu}(\bm{x},\widehat{\bm{x}})\leq 1$ for all $\mu\in \mathcal{M}_{+}(\Xi)$ yields
\begin{align}
\label{Eq:It_now_remains}
\left|E_{0}(\widehat{\mu}_{\ast}^{N})-E_{0}(\widehat{\mu}_{\diamond}^{N})  \right|\leq {1\over 2}\max\{h_{\ast},h_{\diamond} \}c\exp\left(-{2\over \varepsilon}\right)+\dfrac{2}{\gamma}.
\end{align}
It now remains to show that the optimal Lagrange multipliers $h_{\ast}$ and $h_{\star}$ in Equation \eqref{Eq:It_now_remains} are bounded. To this end, define the following Lagrangian dual function
\begin{subequations}
	\begin{align}
	Q(h)&\df \inf_{\widehat{\mu}^{N}\in \mathcal{M}(\real_{+})} J_{h}(\widehat{\mu}^{N}),\\
	Q^{\varepsilon}(h)&\df  \inf_{\widehat{\mu}^{N}\in \mathcal{M}(\real_{+})} J_{h}^{\varepsilon}(\widehat{\mu}^{N}).
	\end{align}
\end{subequations}
Furthermore, let $\widehat{\mu}_{sl,1}^{N}=\dfrac{1}{N}\sum_{k=1}^{N}\delta_{\xi^{k}_{sl}}$ and $\widehat{\mu}_{sl,2}^{N}=\dfrac{1}{N}\sum_{k=1}^{N}\delta_{\bar{\xi}^{k}_{sl}}$ are the empirical measures in conjunction with the slater vectors $\bm{\xi}_{sl},\bar{\bm{\xi}}_{sl}\in \real_{+}^{N}$ in Assumption \textbf{(A.2)}. We leverage \cite[Lemma 1]{nedic2009approximate}, to obtain following upper bounds 
\begin{subequations}
	\begin{align}
	h_{\diamond}&\leq  \dfrac{1}{R-W_{2}(\widehat{\mu}^{N}_{s,1},\widehat{\mu}^{N}_{0})}\left(E_{\gamma}(\widehat{\mu}_{s,1}^{N})-Q(\bar{h}) \right), \\
	h_{\ast}&\leq \dfrac{1}{R-W_{2}(\widehat{\mu}^{N}_{s,2},\widehat{\mu}^{N}_{0})}\left(E_{\gamma}(\widehat{\mu}_{s,2}^{N})-Q^{\varepsilon}(\bar{h}) \right),
	\end{align}
\end{subequations}
where $\bar{h}\in \real_{+}$ is arbitrary, \textit{e.g.}  $\bar{h}=0$.\hfill $\blacksquare$

\subsection{Proof of Proposition \ref{Lemma:X1}}
\label{Appendix:Proof_of_Lemma_X1}
To compute the 2-norm difference between the processes $(\bar{\bm{\xi}}_{t})_{0\leq t\leq T}$ and $(\bar{\bm{\theta}}_{t})_{0\leq t\leq T}$, we define the following auxiliary vectors
\begin{subequations}
	\begin{align}
	\bar{\bm{\xi}}^{\ast}_{\eta (m-1)}&\df \bar{\bm{\xi}}_{\eta (m-1)}-\eta \nabla \widehat{J}^{N}\big(\bar{\bm{\xi}}_{\eta (m-1)};\bm{z}_{m-1},\tilde{\bm{z}}_{m-1}\big)+\sqrt{2\over \beta}\bm{\zeta}_{\eta m},\\
	\bar{\bm{\theta}}^{\ast}_{\eta (m-1)}&\df \bar{\bm{\theta}}_{\eta (m-1)}-\eta \nabla J\big(\bar{\bm{\theta}}_{\eta(m-1)};\rho_{\eta(m-1)}\big)+\sqrt{2\over \beta}\bm{\zeta}_{\eta m}.
	\end{align}
\end{subequations}
The recursions in Eqs. \eqref{Eq:Recursion_1}-\eqref{Eq:Standard_Wiener} then take the following forms
\begin{subequations}
	\begin{align}
	\bar{\bm{\xi}}_{\eta m}&=\mathscr{P}_{\Xi^{N}}\left(\bm{\xi}^{\ast}_{\eta (m-1)}\right),\\
	\bar{\bm{\theta}}_{\eta m}&=\mathscr{P}_{\Xi^{N}}\left(\bm{\theta}^{\ast}_{\eta (m-1)} \right). 
	\end{align}
\end{subequations}
By the non-expansive property of the Euclidean projection onto a non-empty, closed, convex set $\Xi^{N}$ we obtain (see \cite{bertsekas2003convex})
\begin{align}
\Big\|\mathscr{P}_{\Xi^{N}}\big(\bm{\xi}^{\ast}_{\eta (m-1)}\big)-\mathscr{P}_{\Xi^{N}}\left(\bm{\theta}^{\ast}_{\eta (m-1)}\right)\Big\|_{2}\leq\left\|\bar{\bm{\xi}}_{\eta (m-1)}- \bar{\bm{\theta}}_{\eta (m-1)}\right\|_{2}.
\end{align}
Using a triangle inequality yields
\begin{align}
\nonumber
\left\|\bar{\bm{\xi}}_{\eta m}-\bar{\bm{\theta}}_{\eta m} \right\|_{2}&\leq \Big\|\bar{\bm{\xi}}_{\eta (m-1)}-\bar{\bm{\theta}}_{\eta (m-1)} \Big\|_{2}\\ \label{Eq:Recursive}
&\hspace{4mm}+\eta \left\|\nabla \widehat{J}^{N}\Big(\bar{\bm{\xi}}_{\eta(m-1)};\bm{z}_{m-1},\tilde{\bm{z}}_{m-1}\Big)-\nabla J\Big(\bar{\bm{\theta}}_{\eta(m-1)};\rho_{\eta(m-1)}\Big) \right\|_{2}.
\end{align}
Computing Eq. \eqref{Eq:Recursive} recursively yields
\begin{align}
\label{Eq:Rewrite}
\Big\|\bar{\bm{\xi}}_{\eta m}-\bar{\bm{\theta}}_{\eta m} \Big\|_{2}\leq \Big\|\bar{\bm{\xi}}_{0}-\bar{\bm{\theta}}_{0} \Big\|_{2}+\eta \sum_{\ell=0}^{m-1}\left\|\nabla \widehat{J}^{N}\Big(\bar{\bm{\xi}}_{\eta \ell};\bm{z}_{\ell},\tilde{\bm{z}}_{\ell}\Big)-\nabla J\Big(\bar{\bm{\theta}}_{\eta \ell};\rho_{\eta \ell}\Big) \right\|_{2}.
\end{align}
Therefore, using the triangle inequality and based on the initialization $\bar{\bm{\xi}}_{0}=\bm{\theta}_{0}$, we can rewrite Eq. \eqref{Eq:Rewrite} as follows
\begin{align}
\nonumber
\Big\|\bar{\bm{\xi}}_{\eta m}-\bar{\bm{\theta}}_{\eta m} \Big\|_{2} 
&\leq \eta \sum_{\ell=0}^{m-1}\left\|\nabla \widehat{J}^{N}\Big(\bar{\bm{\xi}}_{\eta \ell};\bm{z}_{\ell},\tilde{\bm{z}}_{\ell}\Big)-\nabla J\Big(\bar{\bm{\xi}}_{\eta \ell};\widehat{\mu}^{N}_{\eta \ell}\Big) \right\|_{2}\\ \nonumber
&\hspace{4mm}+\eta \sum_{\ell=0}^{m-1}\left\| \nabla J\Big(\bar{\bm{\xi}}_{\eta \ell};\widehat{\mu}^{N}_{\eta \ell}\Big)-\nabla J\Big(\bar{\bm{\theta}}_{\eta \ell};\widehat{\mu}^{N}_{\eta \ell}\Big) \right\|_{2}\\ \nonumber
&\hspace{4mm}+\eta \sum_{\ell=0}^{m-1}\left\|\nabla J\Big(\bar{\bm{\theta}}_{\eta \ell};\widehat{\mu}^{N}_{\eta \ell}\Big)-\nabla J\left(\bar{\bm{\theta}}_{\eta \ell};\rho_{\eta\ell}\right) \right\|_{2}\\ \label{Eq:Smolin}
&\df \mathsf{E}_{1}(\eta m)+\mathsf{E}_{2}(\eta m)+\mathsf{E}_{3}(\eta m).
\end{align}

In the sequel, we analyze each term separately:
\subsubsection{Upper Bound on $\mathsf{E}_{1}(\eta m)$}

Let $\mathcal{F}_{m}$ denotes the $\sigma$-algebra generated by the samples $(\bm{z}_{k},\tilde{\bm{z}}_{k})_{k\leq m}$, and the initial condition $\bm{\xi}_{0}$.  Let $\mathcal{F}_{0}=\emptyset$. Then, taking the expectation with respect to the joint distribution $P_{\bm{x},y}$ yields
\begin{align}
\nonumber
\expect_{P_{\bm{x},y}^{\otimes 2}}\left[\nabla \widehat{J}^{N}\big(\bar{\bm{\xi}}_{\eta \ell};\bm{z}_{\ell},\tilde{\bm{z}}_{\ell}\big)\Big|\mathcal{F}_{\ell-1}  \right]&=\nabla J\big(\bar{\bm{\xi}}_{\eta \ell};\widehat{\mu}^{N}_{\eta \ell}\big),
\end{align}
where $\nabla J\big(\bar{\bm{\xi}}_{\eta \ell};\widehat{\mu}^{N}_{\eta \ell}\big)=\big(\nabla_{k} J\big(\bar{\bm{\xi}}_{\eta \ell};\widehat{\mu}^{N}_{\eta \ell}\big)\big)_{1\leq k\leq N}$ has the following elements
\begin{align}
\label{Eq:Totally_0}
\nabla_{k} J(\bar{\bm{\xi}}_{\eta \ell};\widehat{\mu}^{N}_{\eta \ell})&=Q(\bar{\xi}_{\eta\ell}^{k})+\dfrac{1}{N}\sum_{m=1}^{N}R(\bar{\xi}_{\eta\ell}^{k},\bar{\xi}_{\eta\ell}^{m}),
\end{align}
where 
\begin{subequations}
	\label{Eq:Totally}
	\begin{align}
	Q(\bar{\xi}_{\eta\ell}^{k})&\df {1\over N} \expect_{P_{\bm{x},y}^{\otimes 2}}\left[ \|\bm{x}-\tilde{\bm{x}} \|_{2}^{2} y\tilde{y}e^{-\bar{\xi}^{k}_{\eta \ell} \|\bm{x}-\widetilde{\bm{x}} \|_{2}^{2}}\right],\\ 
	R(\bar{\xi}_{\eta\ell}^{k},\bar{\xi}_{\eta\ell}^{m})&\df \dfrac{1}{N\gamma}\expect_{P_{\bm{x}}^{\otimes 2}}\left[\|\bm{x}-\tilde{\bm{x}}\|^{2}_{2}e^{-(\bar{\xi}^{m}_{\eta \ell}+\bar{\xi}_{\eta \ell}^{k}) \|\bm{x}-\tilde{\bm{x}}\|^{2}_{2} }\right].
	\end{align}
\end{subequations}

Define the following random vector
\begin{subequations}
	\begin{align}
	\bm{Z}_{m}&\df \eta \sum_{\ell=0}^{m}\Bigg( \nabla \widehat{J}^{N}\Big(\bar{\bm{\xi}}_{\eta \ell};\bm{z}_{\ell},\tilde{\bm{z}}_{\ell}\Big)-\nabla J\Big(\bar{\bm{\xi}}_{\eta \ell};\widehat{\mu}^{N}_{\eta \ell}\Big)\Bigg) \\
	&=\eta  \sum_{\ell=0}^{m}\Bigg(\nabla \widehat{J}^{N}\Big(\bar{\bm{\xi}}_{\eta \ell};\bm{z}_{\ell},\tilde{\bm{z}}_{\ell}\Big)-\expect\left[\nabla \widehat{J}^{N}(\bar{\bm{\xi}}_{\eta \ell};\bm{z}_{\ell},\tilde{\bm{z}}_{\ell})\Big|\mathcal{F}_{\ell-1}  \right]\Bigg),
	\end{align}
\end{subequations}
with $\bm{Z}_{0}=\bm{0}$. Clearly, $\bm{Z}_{m}$ is a martingale $\expect[\bm{Z}_{m}|\mathcal{F}_{m-1}]=\bm{Z}_{m-1}$. Moreover, it has a bounded difference
\begin{align}
\label{Eq:How_Much}
\|\bm{Z}_{m}-\bm{Z}_{m-1}\|_{2}&\leq \eta \left\| \nabla \widehat{J}^{N}\Big(\bar{\bm{\xi}}_{\eta \ell};\bm{z}_{m},\tilde{\bm{z}}_{m}\Big)\right\|_{2}+\eta\left\|\nabla J\Big(\bar{\bm{\xi}}_{\eta m};\widehat{\mu}^{N}_{\eta m}\Big)\right\|_{2}.
\end{align}
Now, for all $k=1,2,\cdots,N$, the following inequalities can be established using Assumption \textbf{(A.1)},
\begin{subequations}
	\label{Eq:Proof_System}
	\begin{align}
	\left|\nabla_{k} J(\bar{\bm{\xi}}_{\eta \ell};\widehat{\mu}^{N}_{\eta \ell})\right|&\leq {1\over N}(1+\gamma^{-1})K^{2}\\
	\left|\nabla_{k} \widehat{J}^{N}(\bar{\bm{\xi}}_{\eta \ell};\widehat{\mu}^{N}_{\eta \ell})\right|&\leq {1\over N}(1+\gamma^{-1})K^{2},
	\end{align}
\end{subequations}
respectively. From Eq. \eqref{Eq:How_Much} we have
\begin{align}
\|\bm{Z}_{m}-\bm{Z}_{m-1}\|_{2}&\leq {2\over \sqrt{N}}(1+\gamma^{-1})K^{2}.
\end{align}

Therefore, $\bm{Z}_{m}-\bm{Z}_{m-1}$ is (conditionally) zero mean and bounded, and is thus (conditionally) sub-Gaussian with the Orlicz norm of $\|\bm{Z}_{m}-\bm{Z}_{m-1}\|_{\psi_{2}}\leq 2(1+\gamma^{-1})K^{2}\sqrt{N}.$, \textit{i.e.},
\begin{align}
\expect\left[e^{\langle \bm{u},\bm{Z}_{m}-\bm{Z}_{m-1}\rangle} \Big|\mathcal{F}_{m-1}\right]\leq \exp\left({2\over N}(1+\gamma^{-1})^{2}K^{4}\|\bm{u}\|_{2}^{2}\right) , \quad \forall \bm{u}\in \real^{N}.
\end{align}
We thus conclude that $\|\bm{Z}_{m}-\bm{Z}_{m-1}\|_{2}$ is sub-Gaussian with the Orlicz norm of ${2\over \sqrt{N}}(1+\gamma^{-1})K^{2}$. Now, we note that
\begin{align}
\mathsf{E}_{1}(\eta m)=\eta \sum_{\ell=0}^{m-1}\left\|\bm{Z}_{\ell}-\bm{Z}_{\ell-1} \right\|_{2}, \quad \bm{Z}_{-1}\df \bm{0}.
\end{align}
Applying the Bernestein inequality yields
\begin{align}
\label{Eq:Smolin_1}
\prob(\mathsf{E}_{1}(\eta m)\geq \varepsilon)&\leq \exp\left(-\dfrac{\varepsilon^{2}N}{2m\eta^{2}(1+\gamma^{-1})^{2}K^{4}}\right)\\
&\leq \exp\left(-\dfrac{\varepsilon^{2}N}{2\eta^{2}(1+\gamma^{-1})^{2}K^{2}\lfloor {T\over \eta} \rfloor} \right), \quad \forall m\in [0,T\eta^{-1}]\cap \integer.
\end{align}
Applying a union bound yields
\begin{align}
\prob\left(\sup_{m\in [0,{T\over \eta}]\cap \integer} \mathsf{E}_{1}(\eta m)\geq \varepsilon\right)
&\leq \left\lfloor {T\over \eta} \right\rfloor \exp\left(-\dfrac{\varepsilon^{2}N}{2\eta^{2}(1+\gamma^{-1})^{2}K^{2}\lfloor {T\over \eta} \rfloor} \right).
\end{align}
Therefore, with the probability of at least $1-\rho$, we obtain
\begin{align}
\sup_{m\in [0,{T\over \eta}]\cap \integer} \mathsf{E}_{1}(\eta m)\leq \sqrt{\dfrac{2\lfloor {T\over \eta} \rfloor \eta^{2}(1+\gamma^{-1})^{2}K^{2}}{N}\log\dfrac{\lfloor{T\over \eta}\rfloor}{\rho}}.
\end{align}

\subsubsection{Upper Bound on $\mathsf{E}_{2}(\eta m)$}

To characterize the upper bound on $\mathsf{E}_{2}(\eta m)$, we write
\begin{align}
\left\|\nabla J(\bar{\bm{\xi}}_{\eta\ell};\widehat{\mu}^{N}_{\eta \ell})-\nabla J(\bar{\bm{\theta}}_{\eta\ell};\widehat{\mu}^{N}_{\eta\ell})  \right\|_{2}&\leq \left\|\nabla J(\bar{\bm{\xi}}_{\eta\ell};\widehat{\mu}^{N}_{\eta \ell})-\nabla J(\bar{\bm{\theta}}_{\eta\ell};\widehat{\mu}^{N}_{\eta\ell})  \right\|_{1}
\\ \label{Eq:Gross} &= \sum_{k=1}^{N}\left|\nabla_{k}J\left(\bar{\bm{\xi}}_{\eta\ell};\widehat{\mu}_{\eta\ell}^{N}\right)-\nabla_{k} J\left(\bar{\bm{\theta}}_{\eta\ell};\widehat{\mu}_{\eta\ell}^{N} \right)  \right|.
\end{align}
Moreover, using Eqs. \eqref{Eq:Totally_0}-\eqref{Eq:Totally} and the triangle inequality yields
\begin{align}
\nonumber
\left|\nabla_{k}J\left(\bar{\bm{\xi}}_{\eta\ell};\widehat{\mu}_{\eta\ell}^{N}\right)-\nabla_{k} J\left(\bar{\bm{\theta}}_{\eta\ell};\widehat{\mu}_{\eta\ell}^{N} \right)  \right|&\leq  \left|Q(\bar{\xi}_{\eta\ell}^{k})-Q(\bar{\theta}_{\eta\ell}^{k})\right|\\
&\hspace{4mm}\label{Eq:Insane}+\dfrac{1}{N}\sum_{m=1}^{N}\left|R(\bar{\xi}_{\eta\ell}^{k},\bar{\xi}_{\eta\ell}^{m})-R(\bar{\theta}_{\eta\ell}^{k},\bar{\xi}_{\eta\ell}^{m})\right|.
\end{align}
The first term on the right hand side of Eq.\eqref{Eq:Insane} has the following upper bound
\begin{align}
\nonumber
\left|Q(\bar{\xi}_{\eta\ell}^{k})-Q(\bar{\theta}_{\eta\ell}^{k})\right|&\stackrel{\rm{(a)}}{\leq} \dfrac{1}{N}\expect_{P_{\bm{x},y}^{\otimes 2}}\left[\|\bm{x}-\tilde{\bm{x}} \|_{2}^{2}|y\tilde{y}| \left|e^{-\bar{\xi}_{\eta\ell}^{k}\|\bm{x}-\tilde{\bm{x}}\|_{2}^{2}}-e^{-\bar{\theta}_{\eta\ell}^{k}\|\bm{x}-\tilde{\bm{x}}\|_{2}^{2}} \right| \right]\\  \nonumber
&\stackrel{\rm{(b)}}{\leq} \dfrac{K^{2}|\bar{\xi}_{\eta\ell}^{k}-\bar{\theta}_{\eta\ell}^{k}|}{N} \expect_{P_{\bm{x},y}^{\otimes 2}}\left[\|\bm{x}-\tilde{\bm{x}} \|_{2}^{2}|y\tilde{y}|\right]\\ \label{Eq:Insane_1}
&\leq \dfrac{K^{4}|\bar{\xi}_{\eta\ell}^{k}-\bar{\theta}_{\eta\ell}^{k}|}{N},
\end{align}
where $\mathrm{(a)}$ we used the fact that the mapping $\xi\mapsto \exp(-\xi \|\bm{x}-\tilde{\bm{x}}\|_{2}^{2})$ is $K^{2}$-Lipschitz, and $\mathrm{(b)}$ follows by using Assumption $\textbf{(A.2)}$.

The second term on the right hand side of Eq. \eqref{Eq:Insane} has the following upper bound
\begin{align}
\nonumber
\left|R(\bar{\xi}_{\eta\ell}^{k},\bar{\xi}_{\eta\ell}^{m})-R(\bar{\theta}_{\eta\ell}^{k},\bar{\xi}_{\eta\ell}^{m})\right|&\leq \dfrac{1}{N\gamma}\expect_{P_{\bm{x}}^{\otimes 2}}\left[\|\bm{x}-\tilde{\bm{x}}\|^{2}_{2}e^{-\bar{\xi}^{m}_{\eta \ell} \|\bm{x}-\tilde{\bm{x}}\|^{2}_{2} }\left|e^{-\bar{\xi}^{k}_{\eta \ell} \|\bm{x}-\tilde{\bm{x}}\|^{2}_{2} }-e^{-\bar{\theta}^{k}_{\eta \ell} \|\bm{x}-\tilde{\bm{x}}\|^{2}_{2} } \right|\right]\\ \nonumber
&\leq  \dfrac{K^{2}|\bar{\xi}_{\eta\ell}^{k}-\bar{\theta}_{\eta\ell}^{k}|}{N\gamma}\expect_{P_{\bm{x}}^{\otimes 2}}\left[\|\bm{x}-\tilde{\bm{x}}\|^{2}_{2}e^{-\bar{\xi}^{m}_{\eta \ell} \|\bm{x}-\tilde{\bm{x}}\|^{2}_{2} }\right]\\ \label{Eq:Insane_2}
&\leq \dfrac{K^{4}|\bar{\xi}_{\eta\ell}^{k}-\bar{\theta}_{\eta\ell}^{k}|}{N\gamma}.
\end{align}
Plugging Eqs. \eqref{Eq:Insane_1}, \eqref{Eq:Insane_2} into Eq. \eqref{Eq:Insane} yields
\begin{align}
\nonumber
\left|\nabla_{k}J\left(\bar{\bm{\xi}}_{\eta\ell};\widehat{\mu}_{\eta\ell}^{N}\right)-\nabla_{k} J\left(\bar{\bm{\theta}}_{\eta\ell};\widehat{\mu}_{\eta\ell}^{N} \right)  \right| \label{Eq:Gross_2}
&\leq  \dfrac{K^{4}(1+\gamma^{-1})}{N}|\bar{\xi}_{\eta\ell}^{k}-\bar{\theta}_{\eta\ell}^{k}|,
\end{align}
where the last inequality follows from the fact that $|y\tilde{y}|\leq 1$, $\|\bm{x}-\tilde{\bm{x}}\|_{2}\leq K$, and $\exp(-\xi^{m}_{\eta\ell}\|\bm{x}-\tilde{\bm{x}} \|_{2}^{2})\leq 1$. Now, plugging Eq. \eqref{Eq:Gross_2} into Eq. \eqref{Eq:Gross} yields
\begin{align}
\left\|\nabla J(\bar{\bm{\xi}}_{\eta\ell};\widehat{\mu}^{N}_{\eta \ell})-\nabla J(\bar{\bm{\theta}}_{\eta\ell};\widehat{\mu}^{N}_{\eta\ell})  \right\|_{2}\leq \dfrac{K^{4}(1+\gamma^{-1})}{\sqrt{N}}\|\bar{\bm{\xi}}_{\eta\ell}-\bar{\bm{\theta}}_{\eta\ell}\|_{2}.
\end{align}  
Therefore, the error $\mathsf{E}_{2}(\eta m)$ has the following upper bound
\begin{align}
\mathsf{E}_{2}(\eta m)\leq \dfrac{K^{4}(1+\gamma^{-1})\eta}{\sqrt{N}} \sum_{\ell=0}^{m-1}\|\bar{\bm{\xi}}_{\eta\ell}-\bar{\bm{\theta}}_{\eta\ell}\|_{2},
\end{align}
for all $m\in [0,T\eta^{-1}]\cap \integer$. Therefore,
\begin{align}
\label{Eq:Smolin_2}
\sup_{m\in [0,T\eta^{-1}]\cap \integer}\mathsf{E}_{2}(\eta m)\leq \dfrac{K^{4}(1+\gamma^{-1})\eta}{\sqrt{N}} \sum_{\ell=0}^{\lfloor{T\over \eta}\rfloor -1}\|\bar{\bm{\xi}}_{\eta\ell}-\bar{\bm{\theta}}_{\eta\ell}\|_{2},
\end{align}

\subsubsection{Upper Bound on $\mathsf{E}_{3}(m)$}

To upper bound on $\mathsf{E}_{3}(\eta m)$, we bound the following norm
\begin{align}
\label{Eq:One_Player}
\left\|\nabla J\Big(\bar{\bm{\theta}}_{\eta \ell};\widehat{\mu}^{N}_{\eta \ell}\Big)-\nabla J\left(\bar{\bm{\theta}}_{\eta \ell};\rho_{\eta\ell}\right) \right\|_{2}\leq \sum_{k=1}^{N} \left|\nabla_{k}J\Big(\bar{\bm{\theta}}_{\eta \ell};\widehat{\mu}^{N}_{\eta \ell}\Big) -\nabla_{k}J\Big(\bar{\bm{\theta}}_{\eta \ell};\rho_{\eta \ell}\Big) \right|.
\end{align}
Each term inside the parenthesis on the right hand side of Eq. \eqref{Eq:One_Player}  has the following upper bound
\begin{align}
\nonumber
\left|\nabla_{k}J\Big(\bar{\bm{\theta}}_{\eta \ell};\widehat{\mu}^{N}_{\eta \ell}\Big) -\nabla_{k}J\Big(\bar{\bm{\theta}}_{\eta \ell};\rho_{\eta \ell}\Big) \right|&\leq \dfrac{1}{N} \Bigg| \dfrac{1}{N}\sum_{m=1}^{N} R(\bar{\theta}^{k}_{\eta\ell},\bar{\xi}_{\eta\ell}^{m})-\int_{\Xi}R(\bar{\theta}^{k}_{\eta\ell},\bar{\theta}) \rho_{\eta\ell}(\mathrm{d}\bar{\theta}) \Bigg|\\ \nonumber
&\leq  \dfrac{1}{N}\left| \dfrac{1}{N}\sum_{m=1}^{N} R(\bar{\theta}^{k}_{\eta\ell},\bar{\xi}_{\eta\ell}^{m})- \dfrac{1}{N}\sum_{m=1}^{N} R(\bar{\theta}^{k}_{\eta\ell},\bar{\theta}_{\eta\ell}^{m}) \right|\\ \label{Eq:Remnants}
&+\dfrac{1}{N} \Bigg| \dfrac{1}{N}\sum_{m=1}^{N} R(\bar{\theta}^{k}_{\eta\ell},\bar{\theta}_{\eta\ell}^{m})-\int_{\Xi}R(\bar{\theta}^{k}_{\eta\ell},\bar{\theta}) \rho_{\eta\ell}(\mathrm{d}\bar{\theta}) \Bigg|,
\end{align}
where the last step is due to the triangle inequality. The first term on the right hand side of the last inequality in \eqref{Eq:Remnants} has the following upper bound
\begin{align}
\label{Eq:Applaud}
\left| \dfrac{1}{N}\sum_{m=1}^{N} R(\bar{\theta}^{k}_{\eta\ell},\bar{\xi}_{\eta\ell}^{m})- \dfrac{1}{N}\sum_{m=1}^{N} R(\bar{\theta}^{k}_{\eta\ell},\bar{\theta}_{\eta\ell}^{m}) \right|&\leq \dfrac{K^{4}}{N^{2}}\sum_{m=1}^{N}\left|\bar{\xi}_{\eta\ell}^{m}-\bar{\theta}^{m}_{\eta\ell}\right|.
\end{align}
The second term on the right hand side of the last inequality in \eqref{Eq:Remnants} has the following upper bound
\begin{align}
\nonumber
\Bigg| \dfrac{1}{N}\sum_{m=1}^{N} R(\bar{\theta}^{k}_{\eta\ell},\bar{\theta}_{\eta\ell}^{m})-\int_{\Xi}R(\bar{\theta}^{k}_{\eta\ell},\bar{\theta}) \rho_{\eta\ell}(\mathrm{d}\bar{\theta}) \Bigg|
&=\left|\int_{\Xi}R(\bar{\theta}^{k}_{\eta\ell},\bar{\theta})\left(\widehat{\rho}^{N}_{\eta\ell}(\mathrm{d}\bar{\theta}) -\rho_{\eta\ell}(\mathrm{d}\bar{\theta})\right)\right| 
\\ \label{Eq:From_Inequality}
&\leq {K^{2}\over N}D_{\mathrm{BL}}(\rho_{\eta\ell},\widehat{\rho}^{N}_{\eta \ell}),
\end{align}
where the inequality follows by the fact that $|R(\bar{\theta}^{k}_{\eta\ell},\bar{\theta})|\leq K^{2}/N$, and $\widehat{\rho}_{\eta\ell}^{N}$ is the empirical measure associated with the samples $\bar{\theta}^{k}_{\eta\ell}$, \textit{i.e.},
\begin{align}
\widehat{\rho}_{\eta \ell}^{N}\df \dfrac{1}{N}\sum_{m=1}^{N}\delta_{\theta_{0}}(\bar{\theta}^{m}_{\eta\ell}).
\end{align}
We use the upper bounds Eqs. \eqref{Eq:From_Inequality} and \eqref{Eq:Applaud} in conjunction with Inequality \eqref{Eq:Remnants}. We derive
\begin{align}
\nonumber
\left\|\nabla J\Big(\bar{\bm{\theta}}_{\eta \ell};\widehat{\mu}^{N}_{\eta \ell}\Big)-\nabla J\left(\bar{\bm{\theta}}_{\eta \ell};\rho_{\eta\ell}\right) \right\|_{2}&\leq \dfrac{K^{2}}{N}\|\bar{\bm{\xi}}_{\eta\ell}-\bar{\bm{\theta}}_{\eta\ell}\|_{1}+\dfrac{K^{2}}{N} D_{\mathrm{BL}}(\rho_{\eta\ell},\widehat{\rho}_{\eta\ell}^{N})\\
&\leq \dfrac{K^{2}}{\sqrt{N}}\|\bar{\bm{\xi}}_{\eta\ell}-\bar{\bm{\theta}}_{\eta\ell}\|_{2}+\dfrac{K^{2}}{N} D_{\mathrm{BL}}(\rho_{\eta\ell},\widehat{\rho}_{\eta\ell}^{N}).
\end{align}
We apply McDiarmid's martingale inequality to obtain the following concentration inequality
\begin{align}
\label{Eq:Invoke}
\prob\Big(D_{\mathrm{BL}}(\rho_{\eta\ell},\widehat{\rho}^{N}_{\eta\ell})\geq \delta \Big)\leq \exp\left(-{2N\delta^{2}}\right).
\end{align}
Let $Q(m)=\sum_{\ell=0}^{m-1}D_{\mathrm{BL}}(\rho_{\eta\ell},\widehat{\rho}^{N}_{\eta\ell})$. Then, we obtain for all $m=0,1,\cdots, \lfloor {T\over \eta} \rfloor$ that
\begin{align}
\prob\left(Q(m)\geq \delta \right)&\leq \exp\left({-2N\delta^{2}\over m}\right)\\
&\leq \exp\left({-2N\delta^{2}\over \lfloor {T\over \eta} \rfloor}\right).
\end{align}
Applying a union bound yields
\begin{align}
\prob\left(\sup_{m\in [0,{T\over \eta}]\cap \integer}Q(m)\geq \delta \right)
&\leq  \left\lfloor{T\over \eta}\right\rfloor \exp\left(-{2N\delta^{2}\over \lfloor {T\over \eta} \rfloor}\right).
\end{align}
Therefore, with the probability of at least $1-\rho$, we have
\begin{align}
\sup_{m\in [0,{T\over \eta}]\cap \integer}\mathsf{E}_{3}(\eta m) \label{Eq:Smolin_3}
&\leq \dfrac{K^{2}\eta}{\sqrt{N}}\sum_{\ell=0}^{\lfloor {T\over \eta} \rfloor-1}\|\bar{\bm{\xi}}_{\eta\ell}-\bar{\bm{\theta}}_{\eta\ell}\|_{2}+\dfrac{K^{2}\eta }{N}\sqrt{\dfrac{\left\lfloor{T\over \eta}\right\rfloor}{2N}\log\left(\dfrac{\left\lfloor{T\over \eta}\right\rfloor}{\rho}\right)}.
\end{align}

\subsection{Combining the upper bounds}

We now leverage the upper bounds on $\mathsf{E}_{1}(\eta m)$, $\mathsf{E}_{2}(\eta m)$, and $\mathsf{E}_{3}(\eta m)$ in Eqs. \eqref{Eq:Smolin_1}, \eqref{Eq:Smolin_2}, and \eqref{Eq:Smolin_3}. Define
\begin{align}
S_{\eta}(n)\df \sup_{0\leq m\leq n}\Big\|\bar{\bm{\xi}}_{\eta m}-\bar{\bm{\theta}}_{\eta m} \Big\|_{2}.
\end{align}
Applying a union bound yields the following inequality from Eq. \eqref{Eq:Smolin}
\begin{align}
\nonumber
S_{\eta}\Big(\Big\lfloor {T\over \eta} \Big\rfloor\Big)
&\leq {4\eta(1+\gamma^{-1})^{2}K^{2}}\left({1\over \sqrt{N}}+{1\over N\sqrt{2N}}\right)\sqrt{\Big\lfloor {T\over \eta}\Big\rfloor\log \left(\dfrac{\lfloor {T\over \eta}\rfloor}{\rho}\right)}\\  \label{Eq:In_conjunction}
&\hspace{4mm}+\dfrac{2K^{4}(1+\gamma^{-1})\eta}{\sqrt{N}}\sum_{\ell=0}^{\lfloor {T\over \eta}\rfloor-1}S_{\eta}(\ell),
\end{align}
with the probability of at least $1-2\rho$. In deriving the last inequality, we assumed that $K\geq 1$.

Now, we invoke the discrete Gr\"{o}nwall's inequality \cite{holte2009discrete}:

\begin{lemma}\textsc{(Discrete Gr\"{o}nwall's inequality, \cite{holte2009discrete})}
	\label{Lemma:B1}	
	If $\{y_{m}\}_{m\in \integer}$, $\{x_{m} \}_{m\in \integer}$, and $\{z_{m}\}_{m\in \integer}$ are non-negative sequences, and
	\begin{align}
	\label{Eq:Discete_Gronwall_0}
	y_{m}\leq x+\sum_{\ell=0}^{m-1}z_{\ell}y_{\ell}, \quad m\in \integer,
	\end{align}	
	then,
	\begin{align}
	\label{Eq:Discete_Gronwall}
	y_{m}\leq x\prod_{0\leq \ell<n}(1+z_{\ell})\leq x\exp\left(\sum_{\ell=0}^{m-1}z_{\ell}\right).
	\end{align}
\end{lemma}
Employing the (discrete) Gr\"{o}nwall's inequalities \eqref{Eq:Discete_Gronwall_0}- \eqref{Eq:Discete_Gronwall} of Lemma \ref{Lemma:B1} in conjunction with Inequality \eqref{Eq:In_conjunction} yields
\begin{align}
\label{Eq:Trading}
\sup_{m\in [0,{T\over \eta}]\cap \integer}\Big\|\bar{\bm{\xi}}_{\eta m}-\bar{\bm{\theta}}_{\eta m} \Big\|_{2} \leq & {{4\eta(1+\gamma^{-1})^{2}K^{2}}\over \sqrt{N}}\sqrt{\Big\lfloor {T\over \eta}\Big\rfloor\log \left(\dfrac{\lfloor {T\over \eta}\rfloor}{\rho}\right)}\exp\left( \dfrac{2K^{4}(1+\gamma^{-1})\eta \Big\lfloor {T\over \eta}\Big\rfloor}{\sqrt{N}}\right),
\end{align}
with the probability of $1-2\rho$. Alternatively, since $\bar{\bm{\xi}}_{\eta\ell}=\bm{\xi}_{\ell},\forall\ell\in [0,{T/\eta}]\cap \integer$, we have
\begin{align}
\label{Eq:Trading}
\sup_{m\in [0,{T\over \eta}]\cap \integer}\Big\|\bm{\xi}_{m}-\bar{\bm{\theta}}_{\eta m} \Big\|_{2} \leq & {{4\sqrt{\eta}(1+\gamma^{-1})^{2}K^{2}}\over \sqrt{N}}\sqrt{T\log \left(\dfrac{T}{\eta \rho}\right)}\exp\left( \dfrac{2K^{4}(1+\gamma^{-1})T}{\sqrt{N}}\right),
\end{align}
Recall the definition of the Wasserstein distance between two measures $\nu,\mu\in \mathcal{M}(\mathcal{X})$,
\begin{align}
W_{p}(\mu,\nu)= \inf_{\bm{X},\bm{Y}} \left(\expect\left[\|\bm{X}-\bm{Y}\|^{p}_{p}\right]\right)^{1\over p},
\end{align}
where the infimum is over all pair of random variables $(\bm{X},\bm{Y})$ with the marginals $\bm{X}\sim \mu$ and $\bm{Y}\sim \nu$. Accordingly,
\begin{align}
\nonumber
\sup_{m\in [0,{T\over \eta}]\cap \integer} W_{1}\left(\left(\widehat{\mu}_{ m}^{N}\right)^{\otimes N},\rho^{\otimes N}_{\eta m}\right)&= \sup_{m\in [0,{T\over \eta}]\cap \integer}\inf_{\bm{\xi}_{m},\bar{\bm{\theta}}_{\eta m}} \expect\left[\Big\|\bm{\xi}_{m}-\bar{\bm{\theta}}_{\eta m} \Big\|_{2}\right]\\ \nonumber
&\leq \sup_{m\in [0,{T\over \eta}]\cap \integer}\expect\left[\Big\|\bm{\xi}_{ m}-\bar{\bm{\theta}}_{\eta m} \Big\|_{2}\right]     \\ \label{Eq:Na_Li}
&\leq \expect\left[\sup_{m\in [0,{T\over \eta}]\cap \integer}\Big\|\bm{\xi}_{ m}-\bar{\bm{\theta}}_{\eta m} \Big\|_{2}\right].
\end{align}
For two measures $\mu,\nu\in \mathcal{M}(\mathcal{X})$ on a compact metric space $(\mathcal{X},d)$ with the diameter $\mathrm{diam}(\mathcal{X})=\sup_{\bm{x}_{1},\bm{x}_{2}\in \mathcal{X}}\|\bm{x}_{1}-\bm{x}_{2}\|_{2}$, the following inequalities can be shown for all $p,q\in [1,\infty),p\leq q$, (see, \textit{e.g.}, \cite{panaretos2019statistical})
\begin{subequations}
	\begin{align}
	\label{Eq:SDP_0}
	W_{p}(\mu,\nu)&\leq W_{q}(\mu,\nu), \\ \label{Eq:SDP}
	W^{q}_{q}(\mu,\nu)&\leq \left(\mathrm{diam}(\mathcal{X})\right)^{q-p}W_{p}^{p}(\mu,\nu).
	\end{align}
\end{subequations}
Specifically, Inequalities \eqref{Eq:SDP_0} and \eqref{Eq:SDP} are due to Jensen's and H\"{o}lder's inequalities, respectively. From Eqs. \eqref{Eq:SDP} on the metric space $(\Xi^{N},\|\cdot\|_{2})$ with $p=1$ and  $q=2$ and \eqref{Eq:Na_Li} we obtain that
\begin{align}
\label{Eq:Finance}
\sup_{m\in [0,{T\over \eta}]\cap \integer} W_{2}^{2}\left(\left(\widehat{\mu}_{ m}^{N}\right)^{\otimes N},\rho^{\otimes N}_{\eta m}\right)\leq \sqrt{N}(\xi_{u}-\xi_{l}) \expect\left[\sup_{m\in [0,{T\over \eta}]\cap \integer}\Big\|\bm{\xi}_{m}-\bar{\bm{\theta}}_{\eta m} \Big\|_{2}\right],
\end{align}
where we used the fact that $\mathrm{diam}(\Xi^{N})=\sqrt{N}(\xi_{u}-\xi_{l})$.

Now, we use the following tensorization property of the Wasserstein distances:

\begin{theorem}\textsc{(Tensorization, \cite[Lemma 3]{mariucci2018wasserstein})}
	\label{Thm:Tensorization}	
	Consider the metric measure space $(\real^{n},\|\cdot\|_{r})$, and let $\mu=\bigotimes_{i=1}^{n}\mu_{i}$ and $\nu=\bigotimes_{i=1}^{n}\nu_{i}$ denotes two probability measures on $\real^{n}$. Then,
	\begin{align}
	W^{p}_{p}(\mu,\nu)\leq \max\{1,n^{{p\over r}-1} \}\sum_{i=1}^{n}W^{p}_{p}(\mu_{i},\nu_{i}).
	\end{align}
	for any $p\geq 1$. In the particular case of $p=2,r=2$, the following exact identity holds
	\begin{align}
	\label{Eq:Suck_0}
	W_{2}^{2}(\mu,\nu)=\sum_{i=1}^{n}W_{2}^{2}(\mu_{i},\nu_{i}).
	\end{align}
	In particular, when $\mu_{1}=\mu_{2}=\cdots=\mu_{n}=\mu_{0}$ and $\nu_{1}=\nu_{2}=\cdots=\nu_{n}=\nu_{0}$, then 
	\begin{align}
	W^{p}_p(\mu,\nu)\leq n^{p\over r}W_{p}^{p}(\mu_{0},\nu_{0}), \quad p\geq 1,
	\end{align}
	and for $p=2,r=2$,  
	\begin{align}
	\label{Eq:Suck}
	W^{2}_2(\mu,\nu)=nW_{2}^{2}(\mu_{0},\nu_{0}).
	\end{align}
\end{theorem}

Let us remark that the special case of $p=2,r=2$ is not stated in \cite[Lemma 3]{mariucci2018wasserstein}, and indeed is due to \cite{panaretos2019statistical}.

Applying the identity \eqref{Eq:Suck} of Theorem \ref{Thm:Tensorization} to Eq. \eqref{Eq:Finance} yields
\begin{align}
\label{Eq:Save_Me}
\sup_{m\in [0,{T\over \eta}]\cap \integer} W_{2}^{2}\left(\widehat{\mu}_{ m}^{N},\rho_{\eta m}\right)\leq {(\xi_{u}-\xi_{l})\over \sqrt{N}}\expect\left[\sup_{m\in [0,{T\over \eta}]\cap \integer}\Big\|\bm{\xi}_{m}-\bar{\bm{\theta}}_{\eta m} \Big\|_{2}\right].
\end{align}
Combining \eqref{Eq:Save_Me} with the upper bound \eqref{Eq:Trading} yields the desired result. \hfill $\blacksquare$

\subsection{Proof of Proposition \ref{Lemma:X2}}
\label{Appendix:Proof_of_Lemma_X2}
Let $(\Omega,\mathcal{F}, (\mathcal{F})_{0\leq t\leq T}, \mathbb{P})$ denotes a filtered probability space, and let $(\bm{W}_{t})_{0\leq t\leq T}$ denotes the Wiener processes that is $\mathcal{F}_{t}$-adapted. To establish the proof, we consider the mapping $\mathcal{T}:C([0,\eta m],\real^{N}) \rightarrow \real^{N\times m}$ from the space of sample paths such that $\mathcal{T}((\bm{W}_{\eta k})_{k\leq m})=\bm{\theta}_{\eta m}$, where $\bm{\theta}_{\eta m}$ is defined recursively as follows
\begin{align}
\bm{\theta}_{\eta m}=\mathcal{P}_{\Xi^{N}}\left(\bm{\theta}_{\eta (m-1)}-\eta \nabla J(\bm{\theta}_{\eta m},\varsigma_{\eta m})+\sqrt{2\over \beta}\bm{\zeta}_{\eta m}  \right),
\end{align}
where $\varsigma_{\eta m}^{\otimes N}=\mathbb{P}_{\bm{\theta}_{\eta m}}$, and $\bm{\zeta}_{\eta m}=\bm{W}_{\eta m}-\bm{W}_{\eta (m-1)}$. Then, $\mathcal{T}((\bm{W}_{t})_{t\leq \eta m}))=\bar{\bm{\theta}}_{\eta m}$, and $\mathcal{T}((\tilde{\bm{W}}_{t})_{t\leq \eta m})=\tilde{\bm{\theta}}_{\eta m}$. Let $\mathbb{P}_{\tilde{\bm{\theta}}_{\eta m}}\stackrel{\rm{law}}{=} \tilde{\bm{\theta}}_{\eta m} $ and $\mathbb{P}_{\bar{\bm{\theta}}_{\eta m}}\stackrel{\rm{law}}{=} \bar{\bm{\theta}}_{\eta m}$. Furthermore, $\mathbb{P}_{\tilde{\bm{W}}}\stackrel{\rm{law}}{=} \tilde{\bm{W}}$ and $\mathbb{P}_{\bm{W}}\stackrel{\rm{law}}{=} \bm{W}$, where $\tilde{\bm{W}}=(\tilde{\bm{W}}_{\eta k})_{k\leq m}$ and $\bm{W}=(\bm{W}_{\eta k})_{0\leq k\leq m}$.

Then, for any measurable mapping $\mathcal{T}:(\Omega,\mathcal{F})\rightarrow (\Omega',\mathcal{F}')$, the following inequality holds
\begin{align}
D_{\mathrm{KL}}\left(\mathbb{P}\circ \mathcal{T}^{-1}||\mathbb{Q}\circ \mathcal{T}^{-1}\right)\leq D_{\mathrm{KL}}(\mathbb{P}||\mathbb{Q}).
\end{align}
Therefore, 
\begin{align}
\label{Eq:Plug_3}
D_{\mathrm{KL}}\Big(\mathbb{P}_{\bar{\bm{\theta}}_{\eta m}}||\mathbb{P}_{\tilde{\bm{\theta}}_{\eta m}}\Big)\leq D_{\mathrm{KL}}\Big(\mathbb{P}_{\bm{W}} ||\mathbb{P}_{\tilde{\bm{W}}}\Big).
\end{align}
Define the stopping times
\begin{align}
\tau^{k}_{n}\df \inf\left\{t\geq 0: \int_{0}^{t}G^{k}_{s}\mathrm{d}s \geq n \right\}, \quad k=1,2,\cdots,N,
\end{align}
where $\bm{G}_{s}=(G_{s}^{1},\cdots,G_{s}^{N})$ and
\begin{align}
\bm{G}_{s}\df \sqrt{\beta\over 2}\left(\nabla J(\tilde{\bm{\theta}}_{s},\rho_{s}) -\nabla J(\bm{\theta}_{s},\mu_{s})\right).
\end{align}
Furthermore, define the stopped process 
\begin{align}
\bm{U}_{t}(n)\df \bm{U}_{t\wedge \tau_{n}}\df  \int_{0}^{\eta m\wedge \tau_{n}}\bm{G}_{s}\mathrm{d}s. 
\end{align}
The processes $(\bm{U}_{t}(n))$ are adapted, locally-$\mathcal{H}_{2}$ processes, \textit{i.e.}, for every $t\geq 0$, the truncated processes $(\bm{U}_{s}(n))_{s\leq t}$ are in the class $\mathcal{H}_{2}$. Thus, the It\^{o} integrals $\int_{0}^{t}G^{k}_{s}\mathrm{d}W^{k}_{s}$ are well defined for all $k=1,2,\cdots,N$. For each coordinate $k=1,2,\cdots,N$, we define the \textit{Dol\'{e}ans-Dade exponential} as follows
\begin{align}
\label{Eq:Martingale_1}
\mathcal{E}_{t}(U^{k}(n))\df \exp\left(\int_{0}^{t\wedge \tau_{n}^{k}}G^{k}_{s}(n)\mathrm{d}W_{s}-{1\over 2}\int_{0}^{t\wedge \tau_{n}^{k}}|G^{k}_{s}(n)|^{2}\mathrm{d}s \right), \quad k=1,2,\cdots,N.
\end{align}

The following theorem provides the sufficient condition for the Dol\'{e}ans-Dade exponential to be a martingale:

\begin{theorem}\textsc{(Novikov Theorem \cite{novikov1973moment})}
	\label{Thm:Nikov}
	Consider the process $(Y_{t})_{t\geq 0}$ that is a real-valued adapted process on the probability space $(\Omega,\mathcal{F},\mathbb{P},(\mathcal{F}_{t})_{t\leq T})$ and locally $\mathcal{H}_{2}$. Furthermore, $(W_{t})_{0\leq t\leq T}$ is an adapted Wiener process. Suppose the following condition holds
	\begin{align}
	\label{Eq:Verify_Condition}
	\expect\left[\exp\left({1\over 2}\int_{0}^{t}|Y_{s}|^{2}\mathrm{d}s \right)\right]<\infty.
	\end{align}
	Then for each $t\geq 0$, the Dol\'{e}ans-Dade exponential defined as below
	\begin{align}
	\mathcal{E}_{t}(Y)=\exp\left(\int_{0}^{t}Y_{s}\mathrm{d}W_{s}-{1\over 2}\int_{0}^{t}|Y_{s}|^{2}\mathrm{d}s \right),
	\end{align} 
	is a positive martingale, under the probability measure $\mathbb{P}$. 
\end{theorem}

It is easy to see that the processes $U_{t}^{k}(n),k=1,2,\cdots,n$ defined in Eq. \eqref{Eq:Martingale_1} satisfies the condition \ref{Eq:Verify_Condition} of Theorem \ref{Thm:Nikov}, and thus the exponential process in \eqref{Eq:Martingale_1} is a martingale. This in turn allows us to employ the following change of measure argument due to Girsanov \cite{girsanov1960transforming}:

\begin{theorem}\textsc{(Girsanov's Change of Measure  \cite{girsanov1960transforming})}
	\label{Thm:Grisanov_Theorem}
	Let $(W_{t})_{t\leq T}$ be a Wiener process on the Wiener probability space $(\Omega, \mathcal{F},\mathbb{P})$. Let $(X_{t})_{t\leq T}$ be a measurable process adapted to the natural filtration of the Wiener process $\mathcal{F}_{t}=\sigma(W_{s\leq t})$ with $X_{0}=0$. If the Dol\'{e}ans-Dade exponential  $\mathcal{E}_{t}(X)$ is a strictly positive martingale, the probability measure $\mathbb{Q}$ can be defined on $(\Omega,\mathcal{F})$ via the Radon-Nikodym derivative
	\begin{align}
	\dfrac{\mathrm{d}\mathbb{Q}}{\mathrm{d}\mathbb{P}}\Big|_{\mathcal{F}_{t}}=\mathcal{E}_{t}(X).
	\end{align}
	Then for each $t\leq T$, the measure $\mathbb{Q}$ restricted to the unaugmented sigma fields $\mathcal{F}_{t}$ is equivalent to $\mathbb{P}$ restricted to $\mathcal{F}_{t}$. 
\end{theorem}

Now, consider the following \textit{change of measure}
\begin{align}
\dfrac{\mathrm{d}\mathbb{P}_{U^{k}(n)+W^{k}}}{\mathrm{d}\mathbb{P}_{W^{k}}}\Big|_{\mathcal{F}_{t}}=\mathcal{E}_{t}(G^{k}(n)), \quad k=1,2,\cdots,N,
\end{align}
which defines a new probability measure on $(\Omega,\mathcal{F})$. Furthermore, define the following Wiener process
\begin{align}
\tilde{W}_{t}^{k}(n)=U^{k}_{t}(n)+W_{t}^{k}.
\end{align}
Let $\mathbb{P}_{\bm{U}(n)+\bm{W}}\stackrel{\rm{law}}{=} \bm{U}(n)+\bm{W}$. As $n\rightarrow \infty$, $\tau^{k}_{n}\rightarrow \infty$ for all $k=1,2,\cdots,n$, and $\mathbb{P}_{\bm{U}(n)+\bm{W}}\stackrel{\mathrm{weakly}}{{\rightarrow}} \mathbb{P}_{\tilde{\bm{W}}}$. Due to the lower semi-continuity of the KL divergence, we then have that
\begin{align}
\label{Eq:Pant}
D_{\mathrm{KL}}\Big(\mathbb{P}_{\bm{W}} ||\mathbb{P}_{\tilde{\bm{W}}}\Big)&= \lim\inf_{n\rightarrow \infty}D_{\mathrm{KL}}(\mathbb{P}_{\bm{W}}||\mathbb{P}_{\bm{U}(n)+\bm{W}}).
\end{align}
Furthermore, it is easy to see that $\mathbb{P}_{\bm{U}(n)+\bm{W}}=\bigotimes_{k=1}^{N}\mathbb{P}_{U^{k}(n)+W^{k}}$ and $\mathbb{P}_{\bm{W}}=\bigotimes_{k=1}^{N}\mathbb{P}_{W^{k}}$. Due to the tensorization property of the KL divergence, we have that
\begin{align}
\label{Eq:earth}
D_{\mathrm{KL}}(\mathbb{P}_{\bm{W}}||\mathbb{P}_{\bm{U}(n)+\bm{W}})= \sum_{k=1}^{N}D_{\mathrm{KL}}(\mathbb{P}_{W^{k}}||\mathbb{P}_{U^{k}(n)+W^{k}}).
\end{align}
Therefore, we proceed from Eq. \eqref{Eq:Pant} using \eqref{Eq:earth}
\begin{align}
\nonumber
D_{\mathrm{KL}}\Big(\mathbb{P}_{\bm{W}}||\mathbb{P}_{\tilde{\bm{W}}} \Big)&=\lim\inf_{n\rightarrow \infty}\sum_{k=1}^{N}D_{\mathrm{KL}}(\mathbb{P}_{W^{k}}||\mathbb{P}_{U^{k}(n)+W^{k}})\\ \nonumber
&=-\lim\inf_{n\rightarrow \infty}\sum_{k=1}^{N}\expect\left[\log\left(\dfrac{\mathbb{P}_{U^{k}(n)+W^{k}}}{\mathbb{P}_{W^{k}}}\right)\right]\\ \nonumber
&=-\lim\inf_{n\rightarrow \infty}\sum_{k=1}^{N} \expect\Big[\log\Big(\mathcal{E}_{\eta m}(U^{k}(n))\Big)\Big]\\ \label{Eq:Martingale_1}
&=-\lim\inf_{n\rightarrow \infty}\sum_{k=1}^{N} \expect\left[\int_{0}^{\eta m\wedge \tau_{n}}G^{k}_{s}(n)\mathrm{d}W^{k}_{s}-{1\over 2}\int_{0}^{t\wedge \tau_{n}}|G^{k}_{s}(n)|^{2}\mathrm{d}s \right].
\end{align}
The integral term $\int_{0}^{t\wedge \tau_{n}}G^{k}_{s}(n)\mathrm{d}W^{k}_{s}$ is a martingale and its expectation vanishes. Therefore, Eq. \eqref{Eq:Martingale_1} reduces to 
\begin{align}
D_{\mathrm{KL}}\Big(\mathbb{P}_{\bm{W}}||\mathbb{P}_{\tilde{\bm{W}}} \Big)&=\lim\inf_{n\rightarrow \infty}\sum_{k=1}^{N} \expect\left[{1\over 2}\int_{0}^{\eta m\wedge \tau_{n}}|G^{k}_{s}(n)|^{2}\mathrm{d}s \right]\\
&\stackrel{\rm{(a)}}{=}\sum_{k=1}^{N} \expect\left[\lim\inf_{n\rightarrow \infty}{1\over 2}\int_{0}^{\eta m\wedge \tau_{n}}|G^{k}_{s}(n)|^{2}\mathrm{d}s \right]\\ \label{Eq:Plug_Inequality}
&=\dfrac{\beta}{2}\expect\left[\int_{0}^{\eta m}\big\|\nabla J(\tilde{\bm{\theta}}_{s},\nu_{s}) -\nabla J(\bm{\theta}_{s},\mu_{s})\big\|_{2}^{2}\mathrm{d}s\right],
\end{align}
where $\rm{(a)}$ follows by the monotone convergence theorem.  Plugging Inequality \eqref{Eq:Plug_Inequality} into \eqref{Eq:Plug_3} yields
\begin{align}
\nonumber
D_{\mathrm{KL}}\Big(\mathbb{P}_{\bar{\bm{\theta}}_{\eta m}}||\mathbb{P}_{\tilde{\bm{\theta}}_{\eta m}}\Big)&\leq \dfrac{\beta}{2}\expect\left[\int_{0}^{\eta m}\big\|\nabla J(\tilde{\bm{\theta}}_{s},\nu_{s}) -\nabla J(\bm{\theta}_{s},\mu_{s})\big\|_{2}^{2}\mathrm{d}s\right].
\end{align}
for all $m\in [0,T\eta^{-1}]\cap \integer$. Taking the supremum from both sides yields
\begin{align}
\nonumber
\sup_{m\in [0,{T\over \eta}]\cap \integer} D_{\mathrm{KL}}\Big(\mathbb{P}_{\bar{\bm{\theta}}_{\eta m}}||\mathbb{P}_{\tilde{\bm{\theta}}_{\eta m}}\Big)&\leq \dfrac{\beta}{2}\sup_{m\in [0,{T\eta^{-1}}]\cap \integer} \expect\left[\int_{0}^{\eta m}\big\|\nabla J(\tilde{\bm{\theta}}_{s},\nu_{s}) -\nabla J(\bm{\theta}_{s},\mu_{s})\big\|_{2}^{2}\mathrm{d}s\right]\\
&\leq \dfrac{\beta}{2} \expect\left[\sup_{m\in [0,{T\eta^{-1}}]\cap \integer}\int_{0}^{\eta m}\big\|\nabla J(\tilde{\bm{\theta}}_{s},\nu_{s}) -\nabla J(\bm{\theta}_{s},\mu_{s})\big\|_{2}^{2}\mathrm{d}s\right].
\end{align}
The discrete-time process $m\mapsto \int_{0}^{\eta m}\big\|\nabla J(\tilde{\bm{\theta}}_{s},\nu_{s}) -\nabla J(\bm{\theta}_{s},\mu_{s})\big\|_{2}^{2}\mathrm{d}s$ is a submartingale. Therefore,  invoking the discrete-time version of Doob's maximal submartingale inequality in Eq. \eqref{Eq:Doob_Boob_1} of Theorem \ref{Thm:Doob's Martingale Maximal Inequality} in Section \ref{Appendix:Proof_of_Proposition_X3} yields
\begin{align}
\nonumber
\sup_{m\in [0,{T\over \eta}]\cap \integer} D_{\mathrm{KL}}\Big(\mathbb{P}_{\bar{\bm{\theta}}_{\eta m}}||\mathbb{P}_{\tilde{\bm{\theta}}_{\eta m}}\Big)&\leq {\beta}\expect\Bigg[\int_{0}^{\eta \lfloor {T\over \eta} \rfloor }\big\|\nabla J(\tilde{\bm{\theta}}_{s},\nu_{s}) -\nabla J(\bm{\theta}_{s},\mu_{s})\big\|_{2}^{2}\mathrm{d}s \\
&\hspace{4mm}\times \log\left( \int_{0}^{\eta \lfloor {T\over \eta} \rfloor }\big\|\nabla J(\tilde{\bm{\theta}}_{s},\nu_{s}) -\nabla J(\bm{\theta}_{s},\mu_{s})\big\|_{2}^{2}\mathrm{d}s\right) \Bigg]\\ \label{Eq:we_obtain}
&\leq {\beta}\expect\left[ \left( \int_{0}^{\eta \lfloor {T\over \eta} \rfloor }\big\|\nabla J(\tilde{\bm{\theta}}_{s},\nu_{s}) -\nabla J(\bm{\theta}_{s},\mu_{s})\big\|_{2}^{2}\mathrm{d}s\right)^{2}\right],
\end{align}
where the last step is due to the basic inequality $x\log(x)\leq x^{2}-x\leq x^{2}$, and we used the fact that ${e\over e-1}\leq 2$. By the Cauchy-Schwarz inequality, we obtain
\begin{align}
\label{Eq:lethal_injection}
\left(\int_{0}^{\eta \lfloor {T\over \eta} \rfloor }\big\|\nabla J(\tilde{\bm{\theta}}_{s},\nu_{s}) -\nabla J(\bm{\theta}_{s},\mu_{s})\big\|_{2}^{2}\mathrm{d}s\right)^{2}\leq \beta T \int_{0}^{\eta \lfloor {T\over \eta} \rfloor }\big\|\nabla J(\tilde{\bm{\theta}}_{s},\nu_{s}) -\nabla J(\bm{\theta}_{s},\mu_{s})\big\|_{2}^{4}\mathrm{d}s.
\end{align} 
Combining \eqref{Eq:we_obtain} and \eqref{Eq:lethal_injection} yields
\begin{align}
\nonumber
\sup_{m\in [0,{T\over \eta}]\cap \integer} D_{\mathrm{KL}}\Big(\mathbb{P}_{\bar{\bm{\theta}}_{\eta m}}||\mathbb{P}_{\tilde{\bm{\theta}}_{\eta m}}\Big)&\leq \beta T\expect\left[  \int_{0}^{\eta \lfloor {T\over \eta} \rfloor }\big\|\nabla J(\tilde{\bm{\theta}}_{s},\nu_{s}) -\nabla J(\bm{\theta}_{s},\mu_{s})\big\|_{2}^{4}\mathrm{d}s\right]\\ \label{Eq:lets_talk_about_it}
&=\beta T\int_{0}^{\eta \lfloor {T\over \eta} \rfloor }\expect\left[\big\|\nabla J(\tilde{\bm{\theta}}_{s},\nu_{s}) -\nabla J(\bm{\theta}_{s},\mu_{s})\big\|_{2}^{4}\right]\mathrm{d}s.
\end{align}

We compute the following upper bound for the integrand using the triangle inequality
\begin{align}
\nonumber
\big\|\nabla J(\tilde{\bm{\theta}}_{s},\nu_{s}) -\nabla J(\bm{\theta}_{s},\mu_{s})\big\|_{2}&\leq \|\nabla J(\tilde{\bm{\theta}}_{s},\nu_{s})-\nabla J(\tilde{\bm{\theta}}_{s},\mu_{s})\|_{2}\\ \label{Eq:Right_hand_side_1}
&\hspace{4mm}+\|\nabla J(\tilde{\bm{\theta}}_{s},\mu_{s})-\nabla J(\bm{\theta}_{s},\mu_{s}) \|_{2}.
\end{align}
The first term on the right hand side of Eq. \eqref{Eq:Right_hand_side_1} has the following upper bound
\begin{align}
\nonumber
\|\nabla J(\tilde{\bm{\theta}}_{s},\nu_{s})-\nabla J(\tilde{\bm{\theta}}_{s},\mu_{s})\|_{2}&\leq\sum_{k=1}^{N}  \left|\nabla_{k} J(\tilde{\bm{\theta}}_{s},\nu_{s})-\nabla_{k} J(\tilde{\bm{\theta}}_{s},\mu_{s}) \right|\\ \nonumber
&= \dfrac{1}{N} \sum_{k=1}^{N}\left|\int_{\Xi}R(\tilde{\theta}_{s}^{k},\theta)\nu_{s}(\mathrm{d}\theta)-\int_{\Xi}R(\tilde{\theta}_{s}^{k},\theta)\mu_{s}(\mathrm{d}\theta)\right|\\  \label{Eq:ballsy}
&\leq \dfrac{K^{2}}{N} D_{\mathrm{BL}}(\nu_{s},\mu_{s}).
\end{align}
Now, recall that for any two probability measures $\mu,\nu\in \mathcal{M}(\mathcal{X})$, we have $D_{\mathrm{BL}}(\nu,\mu)\leq W_{1}(\nu,\mu)\leq W_{2}(\nu,\mu)$. From Eq. \eqref{Eq:ballsy} we proceed
\begin{align} 
\|\nabla J(\tilde{\bm{\theta}}_{s},\nu_{s})-\nabla J(\tilde{\bm{\theta}}_{s},\mu_{s})\|_{2}&\leq \dfrac{K^{2}}{N} W_{2}(\nu_{s},\mu_{s})\\ 
&\stackrel{\rm{(a)}}{=}\dfrac{K^{2}}{N^{2}} W_{2}(\nu_{s}^{\otimes N},\mu_{s}^{\otimes N})\\ \label{Eq:unhappy_1}
&\leq \dfrac{K^{2}}{N^{2}}\left(\expect\left[\|\tilde{\bm{\theta}}_{s}-\bm{\theta}_{s}\|^{2}_{2}\right]\right)^{1\over 2},
\end{align}
where $\mathrm{(a)}$ follows by the tensorization property of the Wasserstein distance in Theorem \ref{Thm:Tensorization}, 

Similarly, for the second term on the right hand side of Eq. \eqref{Eq:Right_hand_side_1}, the following upper bound holds
\begin{align}
\label{Eq:unhappy_2}
\|\nabla J(\tilde{\bm{\theta}}_{s},\mu_{s})- \nabla J(\bm{\theta}_{s},{\mu}_{s})\|_{2}\leq \dfrac{K^{4}(1+\gamma^{-1})}{\sqrt{N}}\|\tilde{\bm{\theta}}_{s}-\bm{\theta}_{s}\|_{2}.
\end{align}
Plugging Eqs. \eqref{Eq:unhappy_1} and \eqref{Eq:unhappy_2} into Eq. \eqref{Eq:Right_hand_side_1}, raising to the power, and taking the expectation yields
\begin{align}
\nonumber
\expect\left[\|\nabla J(\tilde{\bm{\theta}}_{s},\nu_{s})-\nabla J(\bm{\theta}_{s},\mu_{s})\|_{2}^{4}\right]\leq &\dfrac{4K^{16}}{N^{8}}\left(\expect\left[\|\tilde{\bm{\theta}}_{s}-\bm{\theta}_{s}\|^{2}_{2}\right]\right)^{2}\\ \label{Eq:PPPP}
&+\dfrac{4K^{16}(1+\gamma^{-1})^{4}}{N^{2}}\expect\left[\|\tilde{\bm{\theta}}_{s}-\bm{\theta}_{s}\|^{4}_{2}\right],
\end{align}
where we used the basic inequality $(a+b)^{n}\leq 2^{n-1}a^{n}+2^{n-1}b^{n}$. Due to Jensen's inequality $\left(\expect\left[\|\tilde{\bm{\theta}}_{s}-\bm{\theta}_{s}\|^{2}_{2}\right]\right)^{2}\leq \expect\left[\|\tilde{\bm{\theta}}_{s}-\bm{\theta}_{s}\|^{4}_{2}\right]$, and using the fact that $1/N^{8}\leq (1+\gamma^{-1})^{4}/N^{2}$ yields
\begin{align}
\label{Eq:Amin_2}
\expect\left[\|\nabla J(\tilde{\bm{\theta}}_{s},\nu_{s})-\nabla J(\bm{\theta}_{s},\mu_{s})\|_{2}^{4}\right]\leq \dfrac{8K^{16}(1+\gamma^{-1})^{4}}{N^{2}}\expect\left[\|\tilde{\bm{\theta}}_{s}-\bm{\theta}_{s}\|^{4}_{2}\right].
\end{align}

We substitute Inequality \eqref{Eq:Amin_2} into Eq. \eqref{Eq:lets_talk_about_it} 
\begin{align}
\sup_{m\in [0,{T\over \eta}]\cap \integer}D_{\mathrm{KL}}\Big(\mathbb{P}_{\tilde{\bm{\theta}}_{\eta m}}||\mathbb{P}_{\bar{\bm{\theta}}_{\eta m}}\Big)\leq \dfrac{16\beta K^{16}(1+\gamma^{-1})^{4}}{N^{2}}\int_{0}^{\eta \lfloor {T\over \eta} \rfloor}\expect\Big[\|\tilde{\bm{\theta}}_{s}-\bm{\theta}_{s}\|_{2}^{4}\Big]\mathrm{d}s.
\end{align}
From Csisz\'{a}r-Kullback-Pinsker inequality (see, \textit{e.g.}, \cite{sason2016f}), we have
\begin{align}
\left\|\mathbb{P}_{\tilde{\bm{\theta}}_{\eta m}}-\mathbb{P}_{\bar{\bm{\theta}}_{\eta m}}\right\|_{\mathrm{TV}}^{2}\leq {1\over 2}D_{\mathrm{KL}}\Big(\mathbb{P}_{\tilde{\bm{\theta}}_{\eta m}}||\mathbb{P}_{\bar{\bm{\theta}}_{\eta m}}\Big). 
\end{align}
we obtain that
\begin{align}
\label{Eq:Pinskers_Inequality}
\sup_{m\in [0,{T\over \eta}]\cap \integer}\left\|\mathbb{P}_{\tilde{\bm{\theta}}_{\eta m}}-\mathbb{P}_{\bm{\theta}_{\eta m}}\right\|_{\mathrm{TV}}^{2}\leq \dfrac{16\beta T K^{16}(1+\gamma^{-1})^{4}}{N^{2}}\int_{0}^{\eta \lfloor {T\over \eta} \rfloor}\expect\Big[\|\tilde{\bm{\theta}}_{s}-\bm{\theta}_{s}\|_{2}^{4}\Big]\mathrm{d}s.
\end{align}
To obtain an upper bound on the Wasserstein distance $W_{2}(\mathbb{P}_{\tilde{\bm{\theta}}_{\eta m}},\mathbb{P}_{\bm{\theta}_{\eta m}}\Big)$, we leverage the following result due to Villani \cite[Theorem 6.13]{villani2008optimal}:

\begin{theorem}\textsc{(An Inequality for the Wasserstein Distance, \cite[Theorem 6.13]{villani2008optimal})}
	\label{Thm:Villani}
	Let $\mu$ and $\nu$ be two probability measures on a Polish space $(\mathcal{X},d)$. Let $p\in [1,\infty)$, and $\bm{x}_{0}\in \mathcal{X}$. Then,
	\begin{align}
	W_{p}(\mu,\nu)\leq 2^{1\over q}\left(\int_{\mathcal{X}}d^{p}(\bm{x}_{0},\bm{x})\mathrm{d}|\mu-\nu|(\bm{x}) \right)^{1\over p}, \quad \dfrac{1}{p}+\dfrac{1}{q}=1.
	\end{align}
\end{theorem}
The following corollary is immediate from Theorem \ref{Thm:Villani}:
\begin{corollary}
	\label{Eq:Using_Corollary}
	Suppose the polish space $(\mathcal{X},d)$ has a finite diameter. Then, for any two probability measures $\mu$ and $\nu$ on a Polish space $(\mathcal{X},d)$ we have
	\begin{align}
	\label{Eq:corollary_W}
	W^{p}_{p}(\mu,\nu)\leq \left(\mathrm{diam}(\mathcal{X})\right)^{p} 2^{p\over q} \|\mu-\nu\|_{\mathrm{TV}},
	\end{align} 	
	where $\mathrm{diam}(\mathcal{X})= \sup_{\bm{x}_{0},\bm{x}\in \mathcal{X}}d(\bm{x}_{0},\bm{x})<\infty$.
\end{corollary}

Consider the metric space $(\Xi^{N},\|\cdot\|_{2})$. By assumption $(\mathbf{A.3})$, $\Xi=[\xi_{l},\xi_{\mathrm{u}}]$. Therefore, $\mathrm{diam}(\Xi^{N})= \mathrm{diam}(\Xi)\sqrt{N}$. Using Inequality \eqref{Eq:corollary_W} of Corollary \ref{Eq:Using_Corollary} in conjunction with Eq. \eqref{Eq:Pinskers_Inequality} yields
\begin{align}
\label{Eq:Blue_Color}
\sup_{m\in [0,{T\over \eta}]\cap \integer}W^{4}_{2}\left(\mathbb{P}_{\tilde{\bm{\theta}}_{\eta m}},\mathbb{P}_{\bar{\bm{\theta}}_{\eta m}} \right)\leq
{64\beta T (\xi_{u}-\xi_{l})^{4} K^{16}(1+\gamma^{-1})^{4}}\int_{0}^{\eta \lfloor {T\over \eta} \rfloor}\expect\Big[\|\tilde{\bm{\theta}}_{s}-\bm{\theta}_{s}\|_{2}^{4}\Big]\mathrm{d}s.
\end{align}
Recall that $\nu_{\eta m}^{\otimes N}= \mathbb{P}_{\tilde{\bm{\theta}}_{\eta m}}$ and $\rho_{\eta m}^{\otimes N}=\mathbb{P}_{\bar{\bm{\theta}}_{\eta m}}$. Using the tensorization property of the $2$-Wasserstein distance (cf.  Eq. \eqref{Eq:Suck} of Thm. \ref{Thm:Tensorization}) yields
\begin{align}
\label{Eq:Plug_1}
W^{4}_{2}\left(\mathbb{P}_{\tilde{\bm{\theta}}_{\eta m}},\mathbb{P}_{\bar{\bm{\theta}}_{\eta m}} \right)=N^{2}W_{2}^{4}(\nu_{\eta m},\rho_{\eta m}).
\end{align}
Plugging Eq. \eqref{Eq:Plug_1} into \eqref{Eq:Blue_Color} 
\begin{align}
\nonumber
\sup_{m\in [0,{T\over \eta}]\cap \integer}W_{2}^{4}(\nu_{\eta m},\rho_{\eta m})\leq \dfrac{64\beta T (\xi_{u}-\xi_{l})^{4} K^{16}(1+\gamma^{-1})^{4}}{N^2}\int_{0}^{\eta \lfloor {T\over \eta} \rfloor}\expect\Big[\|\tilde{\bm{\theta}}_{s}-\bm{\theta}_{s}\|_{2}^{4}\Big]\mathrm{d}s.
\end{align}
Taking the square root from both sides and using the fact that
\begin{align}
\sup_{m\in [0,{T\over \eta}]\cap \integer}W_{2}^{2}(\nu_{\eta m},\rho_{\eta m})\leq \left(\sup_{m\in [0,{T\over \eta}]\cap \integer}W_{2}^{4}(\nu_{\eta m},\rho_{\eta m})\right)^{1\over 2},
\end{align}
completes the proof.

\hfill $\blacksquare$

\subsection{Proof of Proposition \ref{Proposition:X3}}
\label{Appendix:Proof_of_Proposition_X3}

Before proving our results, we collect a few technical results to which we refer in the sequel, and we also give a definition. The following definition is concerned with the Skorokhod problem:
\begin{definition}\textsc{(Skorokhod problem, \cite{skorokhod1961stochastic})}
	\label{Definition:Skorokhod problem}
	Given a $d$-dimensional reflection matrix $\bm{R}$, the Skorokhod problem is the problem of constructing a map $\Psi:C((0,T],\real^{d})\rightarrow C((0,T],\real^{d})\times C((0,T],\real^{d})$ such that for every process $(\bm{Y}_{t})_{0\leq t\leq T}\in C((0,T)\times \real^{d})$, the image $(\bm{X}_{t},\bm{Z}_{t})_{0\leq t\leq T}=\Psi((\bm{Y}_{t})_{0\leq t\leq T})$ satisfies the following properties
	\begin{itemize}
		\item $\bm{X}_{t}=\bm{Y}_{t}+\bm{R}\bm{Z}_{t},\quad t\in [0,T]$.
		
		\item $\bm{Z}_{0}=\bm{0}$, and $Z_{t}^{j}$ is non-decreasing for all $j=1,2,\cdots,d$.
		
		\item $\int_{0}^{\infty}X_{t}^{j}\mathrm{d}Z_{t}^{j}=$ for all $j=1,2,\cdots,d$, where the integral is in the Stieltjes sense, which is well defined since the processes $Z_{t}^{j}$ are
		non-decreasing
	\end{itemize}
\end{definition}

We leverage the following lemma due to Tanaka \cite{tanaka2002stochastic} which establishes an upper bound on the norm of the difference between two reflected processes:

\begin{lemma}\textsc{(Tanaka \cite[Lemma 2.2.]{tanaka2002stochastic})} 
	\label{Lemma:Tanaka}	
	Let $(\bm{Y}_{t})_{0\leq t\leq T}$ and $(\tilde{\bm{Y}}_{t})_{0\leq t\leq T}$  denote two c\'{a}dl\'{a}g processes. Furthermore, let $(\bm{X}_{t},\bm{Z}_{t})_{0\leq t\leq T}$ and $(\tilde{\bm{X}}_{t},\tilde{\bm{Z}}_{t})_{0\leq t\leq T}$ denote the corresponding solutions to the Sokhrhod problem with the reflection matrix $\bm{R}=\bm{I}_{d\times d}$ in Definition \ref{Definition:Skorokhod problem}. Then, the following inequality holds
	\begin{align}
	\label{Eq:Unrival}
	\|\bm{X}_{t}-\tilde{\bm{X}}_{t}\|_{2}^{2}\leq &\|\bm{Y}_{t}-\tilde{\bm{Y}}_{t} \|_{2}^{2}+2\int_{0}^{t}\langle \bm{Y}_{t}-\tilde{\bm{Y}}_{t}-\bm{Y}_{s}+\tilde{\bm{Y}}_{s}, \mathrm{d}\bm{Z}_{s}-\mathrm{d}\tilde{\bm{Z}}_{s} \rangle,
	\end{align}	
	for all $t\in [0,T]$.
\end{lemma}

To establish our results, we also need the following maximal inequality for sub-martingales:
\begin{theorem}\textsc{(Doob's Sub-martingale Maximal Inequality, \cite[Thm. 3.4]{doob1953stochastic})} 
	\label{Thm:Doob's Martingale Maximal Inequality}	
	Consider the filtered probability space $(\Omega, \mathcal{F},(\mathcal{F}_{t})_{t\geq 0},\prob)$ and let $(M_t)_{t \ge 0}$ be a continuous $\mathcal{F}_{t}$-adapted non-negative sub-martingale.  Let $p\geq 1$ and $T>0$. If $\expect[M_{T}^{p}]<+\infty$, then we have
	\begin{subequations}
		\begin{align}
		\label{Eq:Doob_Boob}
		\expect\left[\left(\sup_{0\leq t\leq T}M_{t}\right)^{p}\right]&\leq \left({p\over p-1}\right)^{p}\expect\big[M_{T}^{p}\big],\quad p>1 \\ 	\label{Eq:Doob_Boob_1}
		\expect\left[\sup_{0\leq t\leq T}M_{t}\right]&\leq \left({e\over e-1}\right)\expect\big[M_{T}\log M_{T}\big]+\expect\big[M_{0}(1-\log M_{0})\big].
		\end{align}
	\end{subequations}
\end{theorem}

We remark that Inequality \eqref{Eq:Doob_Boob} is the classical Doob $L^{p}$-inequality, $p\in (1,\infty)$, \cite[Theorem 3.4.]{doob1953stochastic}. The second result in \eqref{Eq:Doob_Boob_1} represents the Doob $L^{1}$-inequality
in the sharp form derived by Gilat \cite{gilat1986best} from the $L\log L$ Hardy-Littlewood
inequality.

The following result is due to S{\l}omi{\'n}ski \cite{slominski2001euler}:
\begin{theorem}\textsc{(Non-Central Moments of Local Time of Semi-martingales, \cite[Thm. 2.2]{slominski2001euler})}
	\label{Thm:Higher Moments of Local Time of Semi-martingales}
	Consider the filtered probability space $(\Omega, \mathcal{F},(\mathcal{F}_{t})_{t\geq 0},\prob)$, and let $(\bm{Y}_{t})_{0\leq t\leq T}$ denotes a $\mathcal{F}_{t}$-adapted $\real^{d}$-valued semi-martingale with the following Doob-Meyer decomposition 
	\begin{align}
	\bm{Y}_{t}=\bm{Y}_{0}+\bm{M}_{t}+\bm{A}_{t},
	\end{align}
	where $(\bm{M}_{t})_{0\leq t\leq T}$ is a $\mathcal{F}_{t}$-adapted local martingale, $(\bm{A}_{t})_{0\leq t\leq T}$ is a $\mathcal{F}_{t}$-adapted process of locally bounded variation, and $\bm{Y}_{0}\in B\subset \real^{d}$ is the initial condition, confined to the convex region $B$. Let $(\bm{X}_{t},\bm{Z}_{t})_{0\leq t\leq T}$ denote the solution of the Sokhrhod problem for the process $(\bm{Y}_{t})_{0\leq t\leq T}$ with $\bm{R}=\bm{I}_{d\times d}$ in Definition \ref{Definition:Skorokhod problem}. Further, $\bm{Z}_{t}=\int_{0}^{t}\bm{n}_{t}L(\mathrm{d}s)$ is the regulator process associated with the feasible set $B$ of the process, and $\bm{n}_{t}$ and $L(s)$ are the normal vector and the local time at the boundary $\partial B$. Then, for every $p\in \integer$, every stopping time $\tau$ on $\mathcal{F}_{t}$, and any $\bm{a}\in \bar{B}\backslash \partial B$, there exists $c_{p}>0$ such that
	\begin{align}
	\expect[L^{p}(\tau)]&\leq c_{p}\left(\mathrm{dist}\Big(\bm{a},\partial B\Big)\right)^{-p} \expect\left[\sup_{0\leq t\leq \tau}\|\bm{X}_{t}-\bm{a}\|_{2}^{2p}\right]\\  \label{Eq:Inequality_Damn_Lina}
	&\leq c_{p}\left(\mathrm{dist}\Big(\bm{a},\partial B\Big)\right)^{-p}\left(\|\bm{a}- \bm{Y}_{0}\|_{2}^{2p}+\expect\left[\langle \bm{M},\bm{M} \rangle_{\tau}^{p}+ \|\bm{A}_{\tau}\|_{2}^{2p} \right]\right),
	\end{align}
	where $\langle \bm{M},\bm{M} \rangle_{\tau}\df \sum_{i=1}^{d}\langle M^{i},M^{i}\rangle_{\tau}$ is the quadratic variation, and $\mathrm{dist}(\bm{a},\partial B)\df \min_{\bm{b}\in \partial B}\|\bm{a}-\bm{b}\|_{2}$.
\end{theorem}

In the next lemma, we present a maximal inequality for the Euclidean norm of the Wiener processes:
\begin{lemma}\textsc{(A Maximal Inequality for the Wiener Processes)}
	\label{Lemma:Gaussian Concentration for the Wiener Process}
	Suppose $(\bm{W}_{t})_{0\leq t\leq T}$ is the standard Wiener process with $\bm{W}_{0}=\bm{0}$, and $\bm{W}_{t+u}-\bm{W}_{t}\sim \mathsf{N}(\bm{0},u\bm{I}_{d\times d})$. Then,
	\begin{align}
	\label{Eq:Shit_man_2}
	\expect\left[\sup_{0\leq t\leq T}\left\|\bm{W}_{t}-\bm{W}_{\eta \lfloor{t\over \eta}\rfloor}\right\|^{2p}_{2} \right]\leq  T \left(\dfrac{2p}{2p-1} \right)^{2p}2^{p-1}\eta^{2p-1}N\dfrac{\Gamma\left({N+2p\over 2}\right)}{\Gamma\left({N+2\over 2}\right)}\df A_{p}.
	\end{align}
\end{lemma}
The proof of Lemma \ref{Lemma:Gaussian Concentration for the Wiener Process} is  is due to \cite{bubeck2018sampling}  and is presented in Appendix \ref{Appendix:Gaussian Concentration for the Wiener Process} for completeness.

%

Equipped with these technical results, we are in position to prove the main result of Proposition \ref{Proposition:X3}. Consider the c\'{a}dl\'{a}g process $(\tilde{\bm{\theta}}_{t})_{0\leq t\leq T}$ in Eq. \eqref{Eq:Compared}. It can be readily verified that the process can be reformulated as follows
\begin{align}
\label{Eq:ontology}
\tilde{\bm{\theta}}_{\eta m}=\mathscr{P}_{\Xi^{N}}\left(\tilde{\bm{\theta}}_{\eta (m-1)}+\sqrt{2\over \beta}\tilde{\bm{\zeta}}_{\eta m}\right),
\end{align}
where $\tilde{\bm{\zeta}}_{\eta m}=\tilde{\bm{W}}_{\eta m}-\tilde{\bm{W}}_{\eta (m-1)}$, and 
\begin{align}
\tilde{\bm{W}}_{\eta m}=\bm{W}_{\eta m}-\sqrt{\dfrac{\beta}{2}}\int_{0}^{\eta m}\nabla  J({\bm{\theta}}_{s},\nu_{s})\mathrm{d}s.
\end{align}
Alternatively, the embedded continuous-time dynamics in Eq. \eqref{Eq:ontology} can be written as follows
\begin{align}
\label{Eq:Embedd_3}
\tilde{\bm{\theta}}_{t}&=\tilde{\bm{\theta}}_{0}+\sqrt{\dfrac{2}{\beta}}\sum_{\ell=0}^{\lfloor {t\over \eta} \rfloor}\tilde{\bm{\zeta}}_{\eta\ell}+\int_{0}^{t}\tilde{\bm{n}}_{s}\tilde{L}(\mathrm{d}s),\\
&=\tilde{\bm{\theta}}_{0}+\sqrt{\dfrac{2}{\beta}}\bm{W}_{\eta\lfloor {t\over \eta} \rfloor}- \int_{0}^{\eta\lfloor {t\over \eta} \rfloor}\nabla J(\tilde{\bm{\theta}}_{s},\nu_{s})\mathrm{d}s +\int_{0}^{t}\tilde{\bm{n}}_{s}\tilde{L}(\mathrm{d}s),
\end{align}
where $\tilde{\bm{\theta}}_{0}=\bm{\theta}_{0}$. Furthermore, the local time $\tilde{L}$ is defined as below
\begin{align}
\tilde{L}(s)\df \sum_{m=0}^{\lfloor {T\over \eta} \rfloor}\|\bm{\Delta}_{m}\|_{2}\delta_{\eta m}(s).
\end{align}
Above, $\bm{\Delta}_{m}$ is defined as follows
\begin{align}
\bm{\Delta}_{m}\df \tilde{\bm{\theta}}_{m-1}+\sqrt{2\over \beta}\tilde{\bm{\zeta}}_{m}-\mathscr{P}_{\Xi^{N}}\left(\tilde{\bm{\theta}}_{m-1}+\sqrt{2\over \beta}\tilde{\bm{\zeta}}_{m}\right).
\end{align}
Furthermore, the normal vector $(\bar{\bm{n}}_{t})_{0\leq t\leq T}$ is a piece-wise process, where $\bar{\bm{n}}_{t}=\bar{\bm{n}}_{\eta m}$ for $t\in [\eta m,\eta (m+1))$, and 
\begin{align}
\tilde{\bm{n}}_{\eta m}&\df \dfrac{\bm{\Delta}_{m}}{\|\bm{\Delta}_{m}\|_{2}}, \quad m\in [0, {T\over \eta}] \cap \integer.
\end{align}
We employ Inequality \eqref{Eq:Unrival} from Lemma \ref{Lemma:Tanaka} to compute an upper bound on the norm of the difference between the processes $(\bm{\theta}_{t})_{0\leq t\leq T}$ and $(\tilde{\bm{\theta}}_{t})_{0\leq t\leq T}$ in Eqs. \eqref{Eq:Discrete_Time_Version} and \eqref{Eq:Embedd_3}, respectively. In particular, we define the following c\'{a}dl\'{a}g processes
\begin{subequations}
	\label{Eq:Cadlag_PP}
	\begin{align}
	\bm{Y}_{t}&=\bm{\theta}_{0}+\sqrt{\dfrac{2}{\beta}}\bm{W}_{t}-\int_{0}^{t}\nabla J(\bm{\theta}_{s},\mu_{s})\mathrm{d}s,\\
	\tilde{\bm{Y}}_{t}&=\tilde{\bm{\theta}}_{0}+\sqrt{\dfrac{2}{\beta}}\bm{W}_{\eta\lfloor {t\over \eta} \rfloor}-\int_{0}^{\eta\lfloor {t\over \eta} \rfloor}\nabla J(\tilde{\bm{\theta}}_{s},\nu_{s})\mathrm{d}s.
	\end{align}
\end{subequations}
Furthermore, 
\begin{subequations}
	\begin{align}
	\bm{Z}_{t}=\int_{0}^{t}\bm{n}_{s}L(\mathrm{d}s),\quad
	\widetilde{\bm{Z}}_{t}=\int_{0}^{t}\tilde{\bm{n}}_{s}\tilde{L}(\mathrm{d}s).
	\end{align}		
\end{subequations}
Using Tanaka's Inequaltiy in Eq. \eqref{Eq:Unrival} of Lemma \ref{Lemma:Tanaka}, we derive
\begin{align}
\nonumber
\sup_{0\leq t\leq T }\left\|\tilde{\bm{\theta}}_{t}-\bm{\theta}_{t}\right\|_{2}^{2}\leq &\sup_{0\leq t\leq T }\Big\|\bm{Y}_{t}-\tilde{\bm{Y}}_{t}\Big\|_{2}^{2}\\  \label{Eq:Simplify}
&+4\sup_{0\leq t\leq T }\Big\|\bm{Y}_{t}-\tilde{\bm{Y}}_{t}\Big\|_{2}^{2} (L(T)+\tilde{L}(T)).
\end{align}
By Schwarz inequality, for any $p\in \integer$, we obtain
\begin{align}
\label{Eq:Honest_People}
\expect\left[\sup_{0\leq t\leq T }\left\|\tilde{\bm{\theta}}_{t}-\bm{\theta}_{t}\right\|_{2}^{2p} \right]\leq &2^{2p-1} \expect\left[\sup_{0\leq t\leq T }\Big\|\bm{Y}_{t}-\tilde{\bm{Y}}_{t}\Big\|_{2}^{2p}\right]\\
&+2^{2p+1}\left(\expect\left[\sup_{0\leq t\leq T }\Big\|\bm{Y}_{t}-\tilde{\bm{Y}}_{t}\Big\|_{2}^{2p}\right]\right)^{1\over 2}\left(\expect[L^{2p}(T)]+\expect[\tilde{L}^{2p}(T)] \right)^{1\over 2}.
\end{align}
We invoke Inequality \eqref{Eq:Inequality_Damn_Lina} of Theorem \ref{Thm:Higher Moments of Local Time of Semi-martingales} for the solutions of the Sokhrhod problem associated with the processes in Eq. \eqref{Eq:Cadlag_PP}, where
\begin{subequations}
	\begin{align}
	\bm{M}_{t}&=\sqrt{\dfrac{2}{\beta}}\bm{W}_{t}, \quad 	\tilde{\bm{M}}_{t}=\sqrt{\dfrac{2}{\beta}}\bm{W}_{\eta \lfloor {t\over \eta} \rfloor }, \\
	\bm{A}_{t}&=-\int_{0}^{t}\nabla J(\bm{\theta}_{s},\mu_{s})\mathrm{d}s, \quad \tilde{\bm{A}}_{t}=-\int_{0}^{\eta \lfloor{t\over \eta}\rfloor }\nabla J(\tilde{\bm{\theta}}_{s},\mu_{s})\mathrm{d}s.
	\end{align}
\end{subequations}
Then, $\langle \bm{M},\bm{M} \rangle_{t}={2\over \beta}Nt$, and $\langle \tilde{\bm{M}},\tilde{\bm{M}} \rangle_{t}={2\over \beta}\eta N\lfloor {t\over \eta} \rfloor $. Furthermore, 
\begin{align}
\nonumber
\|\bm{A}_{t}\|_{2}&\leq \left\|\int_{0}^{t}\nabla J(\bm{\theta}_{s},\mu_{s})\mathrm{d}s  \right\|_{2}\\ \nonumber
&\leq \int_{0}^{t}\|\nabla J(\bm{\theta}_{s},\mu_{s})\|_{2}\mathrm{d}s\\
&\leq \dfrac{(1+\gamma^{-1})K^{2}t}{\sqrt{N}}.
\end{align}
Similarly, we have $\|\tilde{\bm{A}}_{t}\|_{2}\leq \dfrac{(1+\gamma^{-1})K^{2}\eta \lfloor {t\over \eta}\rfloor}{\sqrt{N}}.$ For simplicity, we suppose $\bm{\theta}_{0}=\bm{\theta}_{0}\in \bar{\Xi}^{N}\backslash \partial \Xi^{N}$, and let $\bm{a}=\bm{\theta}_{0}=\tilde{\bm{\theta}}_{0}$ in Eq. \eqref{Eq:Inequality_Damn_Lina}. We then obtain that 
\begin{subequations}
	\label{Eq:dudek}
	\begin{align}
	\expect\left[L^{2p}(T)\right]\leq B_{p}, \quad
	\expect\left[\tilde{L}^{2p}(T)\right]\leq  B_{p},
	\end{align}
	where
	\begin{align}
	B_{p}\df c_{p}\left(\mathrm{dist}(\bm{\theta}_{0},\partial \Xi^{N}) \right)^{-2p}\left(\left({2\over \beta}\right)^{p}N^pT^{p} +\dfrac{K^{4p}}{N^{p}}T^{2p} \right).
	\end{align}
\end{subequations}
for some constant $c_{p}>0$ independent of $N$ and $T$, where in the last inequality we used the fact that $\eta \lfloor{T\over \eta}\rfloor\leq T$. Due to the fact that $\bm{\theta}_0=\tilde{\bm{\theta}}_{0}$, we get
\begin{align}
\nonumber
\expect\left[\sup_{0\leq t\leq T }\left\|\bm{Y}_{t}-\tilde{\bm{Y}}_{t}\right\|_{2}^{2p}\right]\leq &2^{2p-1} \expect\left[\sup_{0\leq t\leq T}\left\|\int_{0}^{t}\nabla J(\bm{\theta}_{s},\mu_{s})\mathrm{d}s-\int_{0}^{\eta \lfloor{t\over \eta}\rfloor}\nabla J(\tilde{\bm{\theta}}_{s},\nu_{s})\mathrm{d}s \right\|_{2}^{2p}\right]\\ \label{Eq:Turn_Back}
&+\dfrac{2^{4p-1}}{\beta^{2p}}\expect\left[\sup_{0\leq t\leq T}\|\bm{W}_{t}-\bm{W}_{\eta \lfloor {t\over \eta}\rfloor}\|^{2p}_{2p}\right],
\end{align}
where once again we used the geometric inequality $(a+b)^{2p}\leq 2^{2p-1}a^{2p}+2^{2p-1}b^{2p}$. For the first term in Eq. \eqref{Eq:Turn_Back} we obtain
\begin{align}
\label{Eq:Plug_2}
\left\|\int_{0}^{t}\nabla J(\bm{\theta}_{s},\mu_{s})\mathrm{d}s-\int_{0}^{\eta \lfloor{t\over \eta}\rfloor}\nabla J(\tilde{\bm{\theta}}_{s},\nu_{s})\mathrm{d}s \right\|_{2}^{2p}&\leq 2^{2p-1}\int_{\eta \lfloor {t\over \eta} \rfloor}^{t}  \left\|\nabla J(\bm{\theta}_{s},\mu_{s})\right\|_{2}^{2p}\mathrm{d}s\\ \nonumber
&+2^{2p-1}\int_{0}^{\eta \lfloor{t\over \eta} \rfloor}\left\| \nabla J(\bm{\theta}_{s},\mu_{s})-\nabla J(\tilde{\bm{\theta}}_{s},\nu_{s})\right\|_{2}^{2p}\mathrm{d}s.
\end{align}
From Eq. \eqref{Eq:Proof_System} we recall
\begin{subequations}
	\label{Eq:Shit_man}
	\begin{align}
	\left\|\nabla J(\bm{\theta}_{s},\mu_{s})\right\|_{2}^{2p}&\leq \dfrac{K^{4p}}{N^{p}}(1+\gamma^{-1})^{2p}.
	\end{align}
\end{subequations}
Furthermore, from the derivations leading to Eq. \eqref{Eq:PPPP} we obtain
\begin{subequations}	
	\begin{align}
	\nonumber
	\left\| \nabla J(\bm{\theta}_{s},\mu_{s})-\nabla J(\tilde{\bm{\theta}}_{s},\nu_{s})\right\|_{2}^{2p}&\leq 2^{2p-1}\left\|\nabla J(\bm{\theta}_{s},\mu_{s})-\nabla J(\bm{\theta}_{s},\nu_{s}) \right\|_{2}^{2p}
	\\	&\hspace{4mm}+2^{2p-1}\left\|\nabla J(\bm{\theta}_{s},\nu_{s})-\nabla J(\tilde{\bm{\theta}}_{s},\nu_{s}) \right\|_{2}^{2p}\\ \nonumber
	&\leq \dfrac{2^{2p-1}K^{4p}}{N^{4p}}\left(\expect\left[\|\tilde{\bm{\theta}}_{s}-\bm{\theta}_{s}\|^{2}_{2}\right]\right)^{p}\\  \label{Eq:Shit_man_1}
	&\hspace{4mm}+{2^{2p-1}}\left\|\nabla J(\bm{\theta}_{s},\nu_{s})-\nabla J(\tilde{\bm{\theta}}_{s},\nu_{s}) \right\|_{2}^{2p}.
	\end{align}
\end{subequations}
We plug Eqs. \eqref{Eq:Shit_man},\eqref{Eq:Shit_man_1} into Eq. \eqref{Eq:Plug_2}, take the $\sup$, and subsequently take the expectation
\begin{align}
\nonumber
&\expect\left[\sup_{0\leq t\leq T}\left\|\int_{0}^{t}\nabla J(\bm{\theta}_{s},\mu_{s})\mathrm{d}s-\int_{0}^{\eta \lfloor{t\over \eta}\rfloor}\nabla J(\tilde{\bm{\theta}}_{s},\nu_{s})\mathrm{d}s \right\|_{2}^{2p}\right] \\ \nonumber
&\leq 2^{2p-1} \dfrac{K^{4p}}{N^{p}}(1+\gamma^{-1})^{2p}+\dfrac{2^{4p-2}K^{4p}}{N^{4p}}\int_{0}^{\eta \lfloor{T\over \eta} \rfloor}\expect\left[\|\tilde{\bm{\theta}}_{s}-\bm{\theta}_{s}\|^{2p}_{2}\right]\mathrm{d}s\\
\label{Eq:SSL_0}
&\hspace{4mm}+2^{2p-1}\expect\left[\sup_{0\leq t\leq T}\int_{0}^{\eta \lfloor{t\over \eta}\rfloor}\left\|\nabla J(\bm{\theta}_{s},\nu_{s})-\nabla J(\tilde{\bm{\theta}}_{s},\nu_{s}) \right\|_{2}^{2p}\mathrm{d}s\right],
\end{align}
where in writing the upper bound, we leveraged Jensen's inequality $\left(\expect\left[\|\tilde{\bm{\theta}}_{s}-\bm{\theta}_{s}\|^{2}_{2}\right]\right)^{p}\leq \expect\left[\|\tilde{\bm{\theta}}_{s}-\bm{\theta}_{s}\|^{2p}_{2}\right]$. The process $t\mapsto \int_{0}^{\eta \lfloor{t\over \eta}\rfloor}\left\|\nabla J(\bm{\theta}_{s},\nu_{s})-\nabla J(\tilde{\bm{\theta}}_{s},\nu_{s}) \right\|_{2}^{2p}\mathrm{d}s$ is a sub-martingale. Therefore, invoking Doob's sub-martingale maximal inequality in Theorem \ref{Thm:Doob's Martingale Maximal Inequality} yields
\begin{align}
\label{Eq:my_girlfriend}
&\expect\left[\sup_{0\leq t\leq T}\int_{0}^{\eta \lfloor{t\over \eta}\rfloor}\left\|\nabla J(\bm{\theta}_{s},\nu_{s})-\nabla J(\tilde{\bm{\theta}}_{s},\nu_{s}) \right\|_{2}^{2p}\mathrm{d}s\right]\\ \nonumber
&\leq 2\expect\left[\log\left(\int_{0}^{\eta \lfloor{T\over \eta}\rfloor}\left\|\nabla J(\bm{\theta}_{s},\nu_{s})-\nabla J(\tilde{\bm{\theta}}_{s},\nu_{s}) \right\|_{2}^{2p}\mathrm{d}s\right)\int_{0}^{\eta \lfloor{T\over \eta}\rfloor}\left\|\nabla J(\bm{\theta}_{s},\nu_{s})-\nabla J(\tilde{\bm{\theta}}_{s},\nu_{s}) \right\|_{2}^{2p}\mathrm{d}s \right].
\end{align}
Let us recall the following two upper bounds 
\begin{align}
\nonumber
\left\|\nabla J(\bm{\theta}_{s},\nu_{s})-\nabla J(\tilde{\bm{\theta}}_{s},\nu_{s}) \right\|_{2}^{2p}&\leq 2^{2p-1}\|\nabla J(\bm{\theta}_{s},\nu_{s})\|_{2}^{2p}+2^{2p-1}||\nabla J(\tilde{\bm{\theta}_{s}},\nu_{s})\|_{2}^{2p}\\ \label{Eq:Plug_111}
&\leq \dfrac{2^{2p}K^{4p}(1+\gamma^{-1})^{2p}}{N^{p}}\\ \label{Eq:Plug_222}
\left\|\nabla J(\bm{\theta}_{s},\nu_{s})-\nabla J(\tilde{\bm{\theta}}_{s},\nu_{s}) \right\|_{2}^{2p}&\leq \dfrac{K^{4p}(1+\gamma^{-1})^{2p}}{N^{p}}\|\bm{\theta}_{s}- \tilde{\bm{\theta}}_{s}\|^{2p}_{2}.
\end{align}
We use the first inequality in Eq. \eqref{Eq:Plug_111} to bound the logarithm term in Eq. \eqref{Eq:my_girlfriend}, and the second inequality in Eq. \eqref{Eq:Plug_222} for the term outside of the logarithm
\begin{align}
\nonumber
&\expect\left[\sup_{0\leq t\leq T}\int_{0}^{\eta \lfloor{t\over \eta}\rfloor}\left\|\nabla J(\bm{\theta}_{s},\nu_{s})-\nabla J(\tilde{\bm{\theta}}_{s},\nu_{s}) \right\|_{2}^{2p}\mathrm{d}s\right]\\ \label{Eq:so_long}
&\leq  \dfrac{2K^{4p}(1+\gamma^{-1})^{2p}}{N^{p}}\log\left(\dfrac{2^{2p}K^{4p}T(1+\gamma^{-1})^{2p}}{N^{p}}\right)\expect\left[\int_{0}^{\eta {\lfloor{T\over \eta}\rfloor }}\|\bm{\theta}_{s}- \tilde{\bm{\theta}}_{s}\|^{2p}_{2}\mathrm{d}s \right].
\end{align}
By the Fubini-Tonelli theorem, the expectation on the right hand side of Eq. \eqref{Eq:so_long} can be moved inside the integral. Plugging the result in Eq. \eqref{Eq:SSL_0}  yields
\begin{align}
\nonumber
&\expect\left[\sup_{0\leq t\leq T}\left\|\int_{0}^{t}\nabla J(\bm{\theta}_{s},\mu_{s})\mathrm{d}s-\int_{0}^{\eta \lfloor{t\over \eta}\rfloor}\nabla J(\tilde{\bm{\theta}}_{s},\nu_{s})\mathrm{d}s \right\|_{2}^{2p}\right]\\ \label{Eq:Mama_1}
&\leq C_{p}
+D_{p}\int_{0}^{\eta \lfloor{T\over \eta} \rfloor}\expect\left[\|\tilde{\bm{\theta}}_{s}-\bm{\theta}_{s}\|^{2p}_{2}\right]\mathrm{d}s,
\end{align}
where $C_{p}$ and $D_{p}$ are defined in Eqs. \eqref{Eq:Recall_C} and \eqref{Eq:Recall_D}. We now return to Equation \eqref{Eq:Turn_Back}. By plugging Eq. \eqref{Eq:Mama_1} and \eqref{Eq:Shit_man_2} from Lemma \ref{Eq:Shit_man_2} in Eq. \eqref{Eq:Turn_Back}, we obtain
\begin{align}
\label{Eq:Have_my_back}
\expect\left[\sup_{0\leq t\leq T }\left\|\bm{Y}_{t}-\tilde{\bm{Y}}_{t}\right\|_{2}^{2p}\right]\leq & \dfrac{2^{4p-1}}{\beta^{2p}}A_{p}+2^{2p-1}C_{p}
+2^{2p-1}D_{p}\int_{0}^{\eta \lfloor{T\over \eta} \rfloor}\expect\left[\|\tilde{\bm{\theta}}_{s}-\bm{\theta}_{s}\|^{2p}_{2}\right]\mathrm{d}s.
\end{align}
Substituting Eqs. \eqref{Eq:dudek} and  \eqref{Eq:Have_my_back} into Eq. \eqref{Eq:Honest_People} yields
\begin{align}
\nonumber
\expect\left[\sup_{0\leq t\leq T }\left\|\bm{\theta}_{t}-\tilde{\bm{\theta}}_{t}\right\|_{2}^{2p}\right] \leq  &\left(\dfrac{2^{4p-1}}{\beta^{2p}}A_{p}+2^{2p-1}C_{p}\right) \sqrt{2B_{p}}
\\  \label{Eq:implies}
&+2^{2p-1}D_{p}\sqrt{2B_{p}}\int_{0}^{\eta \lfloor{T\over \eta} \rfloor}\expect\left[\|\tilde{\bm{\theta}}_{s}-\bm{\theta}_{s}\|^{2p}_{2}\right]\mathrm{d}s.
\end{align}
Now, define
\begin{align}
F(t)\df \expect\left[\sup_{0\leq s\leq t }\left\|\bm{\theta}_{s}-\tilde{\bm{\theta}}_{s}\right\|_{2}^{2p}\right].
\end{align}
Equation \eqref{Eq:implies} implies
\begin{align}
F(T)\leq \left(\dfrac{2^{4p-1}}{\beta^{2p}}A_{p}+2^{2p-1}C_{p}\right) \sqrt{2B_{p}}
+2^{2p-1}D_{p}\sqrt{2B_{p}}\int_{0}^{T}F(s)\mathrm{d}s.
\end{align}
We apply the Gronwall's inequality to obtain
\begin{align}
F(T)\leq \sqrt{2B_{p}}\left(\dfrac{2^{4p-1}}{\beta^{2p}}A_{p}+2^{2p-1}C_{p}\right) e^{2^{2p-1}TD_{p}\sqrt{2B_{p}}}.
\end{align}
This inequality completes the proof of the first part. To establish the proof of the second part, we use the fact that $\sup_{0\leq s\leq T}W_{2}^{2}(\nu_{s},\mu_{s})\leq \expect\Big[\sup_{0\leq s\leq t }\left\|\bm{\theta}_{s}-\tilde{\bm{\theta}}_{s}\right\|_{2}^{2}\Big]$. Since $\nu^{\otimes N}_{s}=\mathbb{P}_{\tilde{\bm{\theta}}_{s}}$, and $\tilde{\bm{\theta}}_{s}=\tilde{\bm{\theta}}_{\eta{\lfloor{s\over \eta}\rfloor}}$ is c\'{a}dl\'{a}g, we have $\nu_{s}=\nu_{\eta \lfloor {s\over\eta} \rfloor}$ for all $0\leq s \leq T$. Hence,
\begin{align}
\sup_{0\leq s\leq T}W_{2}^{2}\left(\nu_{\eta \lfloor {s\over\eta} \rfloor},\mu_{s}\right)&=\sup_{0\leq s\leq T}W_{2}^{2}(\nu_{s},\mu_{s}) \\
&\leq \sqrt{2B_{2}}\left(\dfrac{2^{7}}{\beta^{4}}A_{2}+2^{3}C_{2}\right) e^{2^{7}TD_{2}\sqrt{2B_{2}}}.
\end{align}

\hfill $\blacksquare$

%

\subsection{Proof of Proposition  \ref{Appendix:Proof_of_General_Kernel_LSH}}
\label{Appendix:Proof_of_General_Kernel_LSH}
Recall that a kernel function is defined via the inner product $K(\bm{x},\bm{y})=\langle\Phi(\bm{x}),\Phi(\tilde{\bm{x}})\rangle_{\mathcal{H}_{K}}$, where $\Phi(\bm{x}):\mathcal{X}\rightarrow \mathcal{H}$ is the implicit feature map which is an element of the RKHS. The random feature model of Rahimi and Recht \cite{rahimi2008random,rahimi2009weighted} relies on the following embedding of the RKHS in $L^{2}(\Omega\otimes \real,\nu\times \nu_{0})$ space
\begin{align}
K(\bm{x},\tilde{\bm{x}})=\langle\varphi(\bm{x};\bm{\omega},b),\varphi(\tilde{\bm{x}};\bm{\omega},b) \rangle_{L^{2}(\Omega\otimes \real,\nu\times \nu_{0})},
\end{align}
where $\varphi(\bm{x};\bm{\omega},b)=\cos(\langle \bm{\omega},\bm{x} \rangle+b)$. Alternatively, since $K(\bm{x},\bm{x})=K(\tilde{\bm{x}},\tilde{\bm{x}})=\phi(\bm{0})=1$, we have
\begin{align}
\nonumber
\Delta(\bm{x},\tilde{\bm{x}})&\df \|\varphi(\bm{x};\bm{\omega},b)-\varphi(\tilde{\bm{x}};\bm{\omega},b)\|_{L^{2}(\Omega\otimes \real,\nu\times \nu_{0})}\\  \label{Eq:SS}
&=\sqrt{2}(1-K(\bm{x},\tilde{\bm{x}}))^{1\over 2}.
\end{align}
Now, consider the following embedding of the kernel from ${L^{2}(\Omega\otimes \real,\nu\times \nu_{0})}$ into $\ell^{2}$ space
\begin{align}
\widehat{\Delta}_{N}(\bm{x},\tilde{\bm{x}})\df \|\bm{\varphi}_{N}(\bm{x})-\bm{\varphi}_{N}(\tilde{\bm{x}})\|_{2},
\end{align}
where we recall $\bm{\varphi}_{N}(\bm{x})\df (\cos(\langle\bm{\omega}_{k},\bm{x} \rangle+b_{k}))_{1\leq k\leq N}$, where $b_{1},\cdots,b_{N}\sim_{\text{i.i.d.}}\nu_{0}\df \mathrm{Uniform}[-\pi,\pi]$, and $\bm{\omega}_{1},\cdots,\bm{\omega}_{N}\sim_{\text{i.i.d.}}\nu$

Based on \cite[Claim 1]{rahimi2008random}, the following concentration inequality holds
\begin{align}
\nonumber
\prob\left( \sup_{\bm{x},\tilde{\bm{x}}\in \mathcal{X}}\left|K(\bm{x},\tilde{\bm{x}})-\langle\bm{\varphi}_{N}(\bm{x}),\bm{\varphi}_{N}(\tilde{\bm{x}})\rangle \right|\geq \delta \right)\leq 2^{8}\left(\dfrac{\sigma^{2}\mathrm{diam}(\mathcal{X})}{\delta} \right)^{2} \exp\left(-\dfrac{N\delta^{2}}{4(d+2)} \right),
\end{align}
where $\sigma^{2}\df \expect[\|\bm{\omega}\|_{2}^{2}]$. Applying a union bound yields the following inequality
\begin{align}
\label{Eq:Using_the_inequality}
\sup_{\bm{x},\tilde{\bm{x}}\in\mathcal{X}}\left|\Delta^{2}(\bm{x},\tilde{\bm{x}})-\widehat{\Delta}_{N}^{2}(\bm{x},\tilde{\bm{x}}) \right|\leq \sqrt{\dfrac{64(d+2)}{N}}\mathrm{ln}^{1\over 2}\left(\dfrac{2^{10}\sigma^{2}\mathrm{diam}^{2}(\mathcal{X})N}{4(d+2)}\right)\df e_{N},
\end{align}
the probability of at least $1-\rho$. Due to the stability of the Gaussian distribution and the result of \cite{datar2004locality}, we obtain
\begin{align}
\label{Eq:Decreasing_1}
\prob\left[G_{t,\bm{w}}(\bm{x})=G_{t,\bm{w}}(\tilde{\bm{x}})\right]&=\int_{0}^{W}{2\over \sqrt{2\pi \widehat{\Delta}^{2}_{N}(\bm{x},\tilde{\bm{x}})}}\exp\left(-\dfrac{s^{2}}{2\widehat{\Delta}^{2}_{N}(\bm{x},\tilde{\bm{x}})} \right)\left( 1-\dfrac{s}{W}\right)\mathrm{d}s.
\end{align}

\subsubsection{The Upper Bound} The integral on the right hand side of Eq. \eqref{Eq:Decreasing_1} is monotone decreasing in $\widehat{\Delta}^{2}_{N}(\bm{x},\tilde{\bm{x}})$. Therefore, using the inequality $\widehat{\Delta}^{2}_{N}(\bm{x},\tilde{\bm{x}})\geq \Delta^{2}(\bm{x},\tilde{\bm{x}})-e_{N}$ from Eq. \eqref{Eq:Using_the_inequality}  yields
\begin{align}
\nonumber
\prob\left[G_{t,\bm{w}}(\bm{x})=G_{t,\bm{w}}(\tilde{\bm{x}})\right]&\leq\int_{0}^{W}{2\over \sqrt{2\pi (\Delta^{2}(\bm{x},\tilde{\bm{x}})-e_{N})}}e^{-{s^{2}\over 2(\Delta^{2}(\bm{x},\tilde{\bm{x}})-e_{N})}} \left( 1-\dfrac{s}{W}\right)\mathrm{d}s\\
&\leq \int_{0}^{W}{2\over \sqrt{2\pi (\Delta^{2}(\bm{x},\tilde{\bm{x}})-e_{N})}}e^{-{s^{2}\over 2\Delta^{2}(\bm{x},\tilde{\bm{x}})}} \left( 1-\dfrac{s}{W}\right)\mathrm{d}s,
\end{align}
We proceed by using the elementary inequality $(1+x)^{-r}\leq {1\over {1+rx}}$ for $x\geq -1$ and $r\in \real\backslash (0,1)$, to obtain
\begin{align}
\nonumber
\prob\left[G_{t,\bm{w}}(\bm{x})=G_{t,\bm{w}}(\tilde{\bm{x}})\right]&\leq \dfrac{1}{1- {e_{N}\over 2\Delta^{2}(\bm{x},\tilde{\bm{x}})}}\int_{0}^{W}{2\over \sqrt{2\pi \Delta^{2}(\bm{x},\tilde{\bm{x}})}} e^{-{s^{2}\over 2\Delta^{2}(\bm{x},\tilde{\bm{x}})}} \left( 1-\dfrac{s}{W}\right)\mathrm{d}s\\ \nonumber
&= (1+\mathcal{O}(e_{N}))\int_{0}^{W}{2\over \sqrt{2\pi \Delta^{2}(\bm{x},\tilde{\bm{x}})}} e^{-{s^{2}\over 2\Delta^{2}(\bm{x},\tilde{\bm{x}})}} \left( 1-\dfrac{s}{W}\right)\mathrm{d}s,
\end{align}
provided that $e_{N}$ is sufficiently small ($N$ is sufficiently large).
\subsubsection{The Lower Bound} To derive a lower bound, we substitute $\widehat{\Delta}^{2}_{N_{0}}(\bm{x},\tilde{\bm{x}})\leq \Delta^{2}(\bm{x},\tilde{\bm{x}})+e_{N}$ which results in
\begin{align}
\nonumber
\prob\left[G_{t,\bm{w}}(\bm{x})=G_{t,\bm{w}}(\tilde{\bm{x}})\right]&\geq\int_{0}^{W}{2\over \sqrt{2\pi (\Delta^{2}(\bm{x},\tilde{\bm{x}})+e_{N})}}e^{-{s^{2}\over 2(\Delta^{2}(\bm{x},\tilde{\bm{x}})+e_{N})}} \left( 1-\dfrac{s}{W}\right)\mathrm{d}s\\ \nonumber
&\geq  \int_{0}^{W}{2\over \sqrt{2\pi (\Delta^{2}(\bm{x},\tilde{\bm{x}})+e_{N})}}e^{-{s^{2}\over 2\Delta^{2}(\bm{x},\tilde{\bm{x}})}} \left( 1-\dfrac{s}{W}\right)\mathrm{d}s,\\ \nonumber
&\geq \left( 1\over 1+{e_{N}\over \Delta^{2}(\bm{x},\tilde{\bm{x}})}\right)^{1\over 2}  \int_{0}^{W}{2\over \sqrt{2\pi \Delta^{2}(\bm{x},\tilde{\bm{x}})}}e^{-{s^{2}\over 2\Delta^{2}(\bm{x},\tilde{\bm{x}})}} \left( 1-\dfrac{s}{W}\right)\mathrm{d}s.
\end{align}
Now, using the basic inequality $\left(\dfrac{1}{1+x}\right)^{r}\leq 1-rx$ for $x<1$ and $r\in (0,1)$ yields,
\begin{align}
\nonumber
\prob\left[G_{t,\bm{w}}(\bm{x})=G_{t,\bm{w}}(\tilde{\bm{x}})\right]\geq \left( 1-\mathcal{O}(e_{N})\right)\int_{0}^{W}{2\over \sqrt{2\pi \Delta^{2}(\bm{x},\tilde{\bm{x}})}}e^{-{s^{2}\over 2\Delta^{2}(\bm{x},\tilde{\bm{x}})}} \left( 1-\dfrac{s}{W}\right)\mathrm{d}s.
\end{align}

\section{Proofs of auxiliary results}
\label{Section:Proofs of Auxiliary Results}

\subsection{Proof of Lemma \ref{Lemma:Tail Bounds for the Finite Sample Estimation Error}}
\label{Appendix:Proof_of_Thm_Tail Bounds}

Let $\bm{z}\in \mathrm{S}^{d-1}$ denote an arbitrary vector on the unit sphere. We define the following function
\begin{align}
\nonumber
Q_{\bm{z}}\big((y_{1},\bm{x}_{1}),\cdots, (y_{n},\bm{x}_{n})\big)&\df \langle \bm{z},\nabla e_{n}(\xi) \rangle\\ \nonumber
&=-\dfrac{2}{n(n-1)}\sum_{1\leq i<j\leq n}y_{i}y_{j}\|\bm{x}_{i}-\bm{x}_{j}\|^{2}_{2}e^{-\xi \|\bm{x}_{i}-\bm{x}_{j}\|_{2}^{2}}\\
&\hspace{4mm}+\expect_{P_{\bm{x},y}^{\otimes 2}}\Bigg[y\tilde{y}\|\bm{x}-\tilde{\bm{x}}\|^{2}_{2}e^{-\xi \|\bm{x}-\tilde{\bm{x}}\|_{2}^{2}}\Bigg].
\end{align} 
By definition, $\expect_{P_{\bm{x},y}}[Q_{\bm{z}}]=0$. Now, for $m\in \{1,2,\cdots,n\}$, we obtain that
\begin{align}
&\left|Q_{\bm{z}}((y_{1},\bm{x}_{1}),\cdots, (y_{m},\bm{x}_{m}),\cdots, (y_{n},\bm{x}_{n}))- Q_{\bm{z}}((y_{1},\bm{x}_{1}),\cdots, (\tilde{y}_{m},\tilde{\bm{x}}_{m}),\cdots, (y_{n},\bm{x}_{n}))\right|\\ \nonumber
&\leq \dfrac{1}{n(n-1)} \sum_{i\not = m}\left| y_{i}y_{m}\|\bm{x}_{i}-\bm{x}_{m}\|_{2}^{2}e^{-\xi \|\bm{x}_{i}-\bm{x}_{m}\|_{2}^{2}}-y_{i}\tilde{y}_{m}\|\bm{x}_{i}-\tilde{\bm{x}}_{m} \|_{2}^{2}e^{-\xi \|\bm{x}_{i}-\tilde{\bm{x}}_{m}\|_{2}^{2}}  \right| \\ \nonumber
&\leq \dfrac{1}{n(n-1)} \sum_{i\not = m}|y_{i}y_{m}| \|\bm{x}_{i}-\bm{x}_{m}\|_{2}^{2}e^{-\xi \|\bm{x}_{i}-\bm{x}_{m}\|_{2}^{2}}\\ \nonumber
&\hspace{4mm} +\dfrac{1}{n(n-1)}\sum_{i\not = m}|y_{i}\tilde{y}_{m}| \|\bm{x}_{i}-\tilde{\bm{x}}_{m}\|_{2}^{2}e^{-\xi \|\bm{x}_{i}-\tilde{\bm{x}}_{m}\|_{2}^{2}}
\\ \nonumber
&\leq \dfrac{1}{n(n-1)}\sum_{i\not = m}\|\bm{x}_{i}-\bm{x}_{m}\|_{2}^{2}+\dfrac{1}{n(n-1)}\sum_{i\not = m}\|\bm{x}_{i}-\tilde{\bm{x}}_{m}\|_{2}^{2}\\ \nonumber
&\leq \dfrac{2K}{n},
\end{align}
where the last inequality follows by Assumption \textbf{(A.1)}.
Using McDiarmid martingale's inequality \cite{mcdiarmid1989method} yields
\begin{align}
\label{Eq:Concentration_Bound_1}
\prob\left(|Q_{\bm{z}}\big((y_{1},\bm{x}_{1}),\cdots, (y_{n},\bm{x}_{n})\big)|\geq \delta \right)\leq 2\exp\left(-\dfrac{n\delta^{2}}{4K^{2}}\right),
\end{align}
for any $\delta\geq 0$. Now, for every $p\in \integer$, the $2p$-th moment of the random variable $Q_{\bm{z}}$ is given by
\begin{align}
\nonumber
\expect\Big[Q^{2p}_{\bm{z}}((y_{1},\bm{x}_{1}),\cdots,(y_{n},\bm{x}_{n}))\Big]&=\int_{\Xi}2pu^{2p-1}\prob(Q_{\bm{z}}((y_{1},\bm{x}_{1}),\cdots,(y_{n},\bm{x}_{n}))\geq u)\mathrm{d}u
\\ \nonumber &\stackrel{\rm{(a)}}{\leq}\int_{\Xi} 4pu^{2p-1}\exp\left(-\dfrac{nu^{2}}{4K^{2}}\right)\mathrm{d}u
\\ \label{Eq:Readily_Follows} &= {2\left({4K^{2}\over n}\right)^{2p}p!},
\end{align}
where $\rm{(a)}$ is due to the concentration bound in Eq. \eqref{Eq:Concentration_Bound_1}. 
Therefore, 
\begin{align}
\nonumber
\expect\big[\exp\big(Q^{2}_{\bm{z}}((y_{1},\bm{x}_{1}),\cdots,(y_{n},\bm{x}_{n}))/\sigma^{2}\big)\big]&=\sum_{p=0}^{\infty}\dfrac{1}{p!\gamma^{2p}}\expect\Big[Q^{2p}_{\bm{z}}((y_{1},\bm{x}_{1}),\cdots,(y_{n},\bm{x}_{n}))\Big]\\
\nonumber
&=1+2\sum_{p\in \integer}\left(\dfrac{4K^{2}}{n\sigma}\right)^{2p}\\ \nonumber
&=\dfrac{2}{1-(4K^{2}/n\sigma)^{2}}-1.
\end{align}
For $\sigma=4\sqrt{3}K^{2}/n$, we obtain $\expect\big[\exp\big(Q^{2}_{\bm{z}}((y_{1},\bm{x}_{1}),\cdots,(y_{n},\bm{x}_{n}))/\sigma^{2}\big)\big]\leq 2$. Therefore, $\|Q_{\bm{z}}\|_{\psi_{2}}=\|\langle \bm{z},\nabla e_{n}(\xi) \rangle \|_{\psi_{2}}\leq 4\sqrt{3}K^{2}/n$ for all $\bm{z}\in \mathrm{S}^{n-1}$ and $\xi\in \real_{+}$. Consequently, by the definition of the sub-Gaussian random vector in Eq. \eqref{Eq:Sub_Gaussian_random_vector} of Definition \ref{Definition:Sub-Gaussian Norm}, we have $\| \nabla e_{n}(\xi)\|_{\psi_{2}}\leq 4\sqrt{3}K^{2}/n$ for every $\xi\in \real_{+}$. We invoke the following lemma due to \cite[Lemma 16]{khuzani2017stochastic}:

\begin{lemma}\textsc{(The Orlicz Norm of  the Squared Vector Norms, \cite[Lemma 16]{khuzani2017stochastic})}
	\label{Lemma:The_Orlicz_Norm_of_the_Squared_Vector_Norms}
	Consider the zero-mean random vector $\bm{Z}$ satisfying $\|\bm{Z}\|_{\psi_{\nu}}\leq \beta$ for every $\nu\geq 0$. Then, $\|\|\bm{Z}\|_{2}^{2}\|_{\psi_{{\nu\over 2}}}\leq 2\cdot3^{2\over \nu}\cdot \beta^{2}$.
\end{lemma}

Using Lemma \ref{Lemma:The_Orlicz_Norm_of_the_Squared_Vector_Norms}, we now have that $\| \|\nabla e_{n}(\xi) \|_{2}^{2} \|_{\psi_{1}}\leq 64\sqrt{3}K^{2}/n^{2}$ for every $\xi\in \real_{+}$. Applying the exponential Chebyshev's inequality with $\beta=64\sqrt{3}K^{2}/n^{2}$ yields
\begin{align}
\nonumber
&\prob\Bigg(\int_{\Xi}\int_{0}^{1} \Big|\| \nabla e_{n}((1-s)\xi+s\zeta_{\ast}) \|_{2}^{2}-\expect_{\bm{x},y}[\|\nabla e_{n}((1-s)\xi+s\zeta_{\ast}) \|_{2}^{2}]\Big|\mu_{0}(\mathrm{d}\xi)\geq \delta \Bigg)  \\ \nonumber
&\leq e^{-{n^{2}\delta\over 64\sqrt{3}K^{2}}}\expect_{\bm{x},y}\left[e^{\left({n^{2}\over 64\sqrt{3} K^{2}} \int_{\real^{D}}\int_{0}^{1}\big|\| \nabla e_{n}( (1-s)\xi+s\zeta_{\ast}) \|_{2}^{2}-\expect_{\bm{x},y}[\|\nabla e_{n}((1-s)\xi+s\zeta_{\ast}) \|_{2}^{2}]\big|\mathrm{d}s\mu_{0}(\mathrm{d}\xi) \right)}\right]\\ \nonumber
&\stackrel{\rm{(a)}}{\leq} e^{-{n^{2}\delta\over 16\sqrt{3}K^{2}}} \int_{\Xi}\int_{0}^{1}\expect_{\bm{x},y}\Big[e^{{n^{2}\over 16\sqrt{3}K^{2}} (|\|\nabla e_{n}((1-s)\xi+s\zeta_{\ast})\|_{2}^{2}-\expect_{\bm{x},y}[\|\nabla e_{n}((1-s)\xi+s\zeta_{\ast})\|_{2}^{2}])|}\Big]\mathrm{d}s\mu_{0}(\mathrm{d}\xi)\\   \nonumber
&\stackrel{\rm{(b)}}{\leq} 2 e^{-{n^{2}\delta\over 16\sqrt{3}K^{2}}},
\end{align}
where $\rm{(a)}$ follows by Jensen's inequality, and $\rm{(b)}$ follows from the fact that 
\begin{align}
\expect_{\bm{x},y}\Big[e^{{n^{2}\over 16\sqrt{3}K^{2}} (|\|\nabla e_{n}((1-s)\xi+s\zeta_{\ast})\|_{2}^{2}-\expect_{\bm{x},y}[\|\nabla e_{n}((1-s)\xi+s\zeta_{\ast})\|_{2}^{2}]|)}\Big]\leq 2,
\end{align}
by definition of a sub-Gaussian random variable. Therefore,
\begin{align}
\nonumber
&\prob\left(\int_{\Xi}\int_{0}^{1} \|\nabla e_{n}((1-s)\xi+s\zeta_{\ast})\|_{2}^{2}\mathrm{d}s\mu(\mathrm{d}\xi)\geq \delta \right)\\ \label{Eq:damnation_on_lina}
&\hspace{20mm}\leq 2 \exp\left({- \dfrac{n^{2}(\delta -\int_{0}^{1}\int_{\Xi}\expect_{\bm{x},y}[\|\nabla e_{n}((1-s)\xi+s\zeta_{\ast})\|_{2}^{2}]\mathrm{d}s\mu_{0}(\mathrm{d}\xi))}{16\sqrt{3}K^{2}}}\right).
\end{align}
It now remains to compute an upper bound on the expectation $\expect_{\bm{x},y}[\|\nabla e_{n}((1-s)\xi+s\zeta_{\ast})\|_{2}^{2}]$. But this readily follows from Eq. \eqref{Eq:Readily_Follows} by letting $p=1$ and $\bm{z}={{\nabla e_{n}((1-s)\xi+s\zeta_{\ast})}\over {\| \nabla e_{n}((1-s)\xi+s\zeta_{\ast})\|_{2}}}$ as follows
\begin{align}
\nonumber
\expect_{\bm{x},y}[\|\nabla e_{n}((1-s)\xi+s\zeta_{\ast})\|_{2}^{2}]&=\expect_{\bm{x},y}\left[\left\langle{\nabla e_{n}((1-s)\xi+s\zeta_{\ast})\over \| \nabla e_{n}((1-s)\xi+s\zeta_{\ast})\|_{2}} ,\nabla e_{n}((1-s)\xi+s\zeta_{\ast}) \right\rangle^{2} \right]\\ \nonumber
&=\expect_{\bm{x},y}\Big[Q^{2}_{{\nabla e_{n}((1-s)\xi+s\zeta_{\ast})\over \| \nabla e_{n}((1-s)\xi+s\zeta_{\ast})\|_{2}}}\Big]\\
\label{Eq:Expectation_Upper_Bound}
&\leq  {32K^{4}\over n^{2}}.
\end{align}
Plugging the expectation upper bound of Eq. \eqref{Eq:Expectation_Upper_Bound} into Eq. \eqref{Eq:damnation_on_lina} completes the proof of the first part of Lemma \ref{Lemma:Tail Bounds for the Finite Sample Estimation Error}.

The second part of Lemma \ref{Lemma:Tail Bounds for the Finite Sample Estimation Error} follows by a similar approach and we thus omit the proof. $\hfill \blacksquare$

\subsection{Proof of Lemma \ref{Lemma:Gaussian Concentration for the Wiener Process}}
\label{Appendix:Gaussian Concentration for the Wiener Process}

The proof is a minor modification of the proof due to Bubeck, \textit{et al.} \cite{bubeck2018sampling}. We present the proof for completeness. By the definition of the Wiener process $\bm{W}_{t}-\bm{W}_{\eta \lfloor {t\over \eta} \rfloor}\stackrel{\mathrm{law}}{=}\bm{W}_{t-\eta\lfloor{t\over \eta} \rfloor } \sim \mathsf{N}\left(\bm{0},(t-\eta\lfloor{t\over \eta} \rfloor )\bm{I}_{N\times N}\right)$. Therefore,
\begin{align}
\expect\left[\sup_{0\leq t\leq T}\left\|\bm{W}_{t}-\bm{W}_{\eta\lfloor{t\over \eta} \rfloor }\right\|^{2p}_{2} \right]&=\expect\left[\sup_{0\leq t\leq T}\left\|\bm{W}_{t-\eta\lfloor{t\over \eta} \rfloor }\right\|^{2p}_{2} \right]\\
&=\expect\left[\max_{m\in [0,{T\over \eta}]\cap \integer}\sup_{t\in[m\eta,(m+1)\eta) }\left\|\bm{W}_{t-\eta m}\right\|^{2p}_{2} \right]\\ \label{Eq:Lab_mate}
&\leq \sum_{m=0}^{\lfloor {T\over \eta} \rfloor}\expect\left[\sup_{t\in[m\eta,(m+1)\eta) }\left\|\bm{W}_{t-\eta m}\right\|^{2p}_{2}\right].
\end{align} 
The process $t\mapsto \|\bm{W}_{t-\eta m}\|^{2p}_{2}$ is a sub-martingale on the interval $[m\eta,(m+1)\eta)$. Therefore, by Doob's sub-martingale maximal inequality Theorem \ref{Thm:Doob's Martingale Maximal Inequality}, we obtain that
\begin{align}
\label{Eq:Bala_byar}
\expect\left[\sup_{t\in[m\eta,(m+1)\eta) }\left\|\bm{W}_{t-\eta m}\right\|^{2p}_{2} \right]\leq \left(\dfrac{2p}{2p-1} \right)^{2p} \expect\left[\left\|\bm{W}_{\eta}\right\|^{2p}_{2} \right],
\end{align}
for all $p\in \integer$ and $m\in \left[0,T\eta^{-1}\right]\cap \integer$. Plugging Eq. \eqref{Eq:Bala_byar} into Eq. \eqref{Eq:Lab_mate} yields
\begin{align}
\label{Eq:PeePee}
\expect\left[\sup_{0\leq t\leq T}\left\|\bm{W}_{t}-\bm{W}_{\eta\lfloor{t\over \eta} \rfloor }\right\|^{2p}_{2} \right]\leq \left\lfloor {T\over \eta} \right\rfloor  \left(\dfrac{2p}{2p-1} \right)^{2p}\expect\left[\left\|\bm{W}_{\eta}\right\|^{2p}_{2} \right].
\end{align}
The expectation on the right hand side of Eq. \eqref{Eq:PeePee} can be evaluated using integration in the spherical coordinates 
\begin{align}
\nonumber
\expect\left[\left\|\bm{W}_{\eta}\right\|^{2p}_{2} \right]&=\dfrac{1}{(2\pi\eta)^{N}}\dfrac{N\pi^{N/2}}{\Gamma({N\over 2}+1)} \int_{0}^{\infty}e^{-{r^{2}\over 2\eta^{2}}}r^{N+2p-1}\mathrm{d}r\\ \label{Eq:PooPoo}
&=2^{p-1}\eta^{2p}N\dfrac{\Gamma\left({N+2p\over 2}\right)}{\Gamma\left({N+2\over 2}\right)},
\end{align}
where $\Gamma(\cdot)$ is the Euler's Gamma function. Combining Eq. \eqref{Eq:PooPoo} and Eq. \eqref{Eq:PeePee} gives us
\begin{align}
\expect\left[\sup_{0\leq t\leq T}\left\|\bm{W}_{t}-\bm{W}_{\eta\lfloor{t\over \eta} \rfloor }\right\|^{2p}_{2} \right]\leq  T \left(\dfrac{2p}{2p-1} \right)^{2p}2^{p-1}\eta^{2p-1}N\dfrac{\Gamma\left({N+2p\over 2}\right)}{\Gamma\left({N+2\over 2}\right)},
\end{align}
where we used the fact that $\lfloor T/\eta \rfloor\leq T/\eta$.

\hfill $\blacksquare$

\bibliographystyle{amsplain}
\bibliography{NIPS}

\end{document}